\documentclass[russian,openany,11pt%,draft
]
{book}
\usepackage[cp1251]{inputenc}
\usepackage[a5paper,left=15mm,right=15mm, top=15mm,bottom=20mm,bindingoffset=0cm]{geometry}
\usepackage{indentfirst}
\usepackage{url}

\usepackage{cmap}
\usepackage[T2A]{fontenc}

\usepackage[russian]{babel}
\usepackage{amsfonts,amssymb,amsmath,amscd,amsthm}
\usepackage[linktocpage]{hyperref}

\usepackage{remreset}
\makeatletter
\@removefromreset{section}{chapter}
\makeatother

\makeatletter
\renewcommand{\@evenhead}{\thepage \hfil Элементы математики в задачах}
\makeatother

\usepackage{etoolbox}
\patchcmd{\thebibliography}{\chapter*}{\section*}{}{}

\usepackage{makeidx}
\makeindex
\usepackage{color}
\definecolor{chocolate}{rgb}{0.48, 0.25, 0.0}
\definecolor{blue(ryb)}{rgb}{0.01, 0.28, 1.0}

\newcommand{\ull}[1]{\begingroup\color{blue(ryb)}{#1}\endgroup}
\usepackage{pgf,tikz,circuitikz}
\usetikzlibrary{arrows}

\usepackage{epsfig}
\usepackage{graphicx}
\usepackage{epic,eepic}
\graphicspath{{figures/}{fig-sturm/}{geometry/}{fig-perestanovki/}{eps-new/}{./figs/}}
\def\sseccol#1{\subsubsection*{#1}\markright{#1}}

\newif\ifasterisk

\usepackage{answers} %если хочется переместить указания в конец
\Newassociation{hintA}{hintAenv}{_hintAholder}
\Newassociation{hintB}{hintBenv}{_hintBholder}

\Newassociation{hintAA}{hintAAenv}{_hintAAholder}
\Newassociation{hintBB}{hintBBenv}{_hintBBholder}

\newtheoremstyle{mystyleA}%                % Name
  {}%                                     % Space above
  {}%                                     % Space below
  {}%                                     % Body font
  {}%                                     % Indent amount
  {\it}%                            % Theorem head font
  {.}%                                    % Punctuation after theorem head
  { }%                                    % Space after theorem head, ' ', or \newline
  {}%                                     % Theorem head spec (can be left empty, meaning `normal')
\theoremstyle{mystyleA}
\newtheorem{innerthmA}{}

\newenvironment{zadA}[1]{\renewcommand\theinnerthmA{#1}\innerthmA}{\endinnerthmA}

\newtheoremstyle{mystyleB}%                % Name
  {}%                                     % Space above
  {}%                                     % Space below
  {}%                                     % Body font
  {}%                                     % Indent amount
  {\it\bf}%                            % Theorem head font
  {.}%                                    % Punctuation after theorem head
  { }%                                    % Space after theorem head, ' ', or \newline
  {}%                                     % Theorem head spec (can be left empty, meaning `normal')
\theoremstyle{mystyleB}
\newtheorem{innerthmB}{}

\newenvironment{zadB}[1]{\renewcommand\theinnerthmB{#1}\innerthmB}{\endinnerthmB}

\def\id{\mbox{id}}
\def\st{\mbox{st}}
\def\fix{\mbox{fix}}

\def\Z{{\mathbb Z}} \def\R{{\mathbb R}} \def\N{{\mathbb N}}
\def\C{{\mathbb C}} \def\Q{{\mathbb Q}} 
 
\def\ord{\mathop{\fam0 ord}}
\def\Ree{\mbox{Re}}
\def\Imm{\mbox{Im}}

\renewcommand{\le}{\leqslant}
\renewcommand{\leq}{\leqslant}
\renewcommand{\ge}{\geqslant}
\renewcommand{\geq}{\geqslant}

\long\def\comment#1\endcomment{}
\theoremstyle{definition}

\newtheorem{pr}{}[subsection]
\renewcommand{\thepr}{\arabic{section}.\arabic{subsection}.\arabic{pr}}
\newcommand{\zvezda}{\hspace{-2mm}* }

\newcommand{\circpr}{\hspace{-2mm}$^\circ$ }

\newtheorem*{Remarks}{Замечания}
\newtheorem*{deff}{Определение}

\theoremstyle{theorem}
\newtheorem{theorem}{Теорема}[section]

\def\No{\textnumero}

\theoremstyle{thm}
\newenvironment {th*}[1]
    {\gdef\thname{#1} \begin{thn}}%
    {\end{thn}}
\newtheorem*{thn}{\thname}
\def\nopoint#1#2{}
\newcommand{\noqed}{\renewcommand{\qed}{}}
\DeclareMathOperator{\dist}{dist}

\def\Ol#1{}% для комментариев и вопросов
\usepackage{slashbox}% для диагональных линий в ячейках таблицы

\begin{document}

\def\year{2018}%
\vbox to\textheight%
{\thispagestyle{empty}

\centering{%

%{\Large Алексей Паршин}

\vskip5cm

%\vfil

{\huge %\uppercase
{Элементы математики в~задачах}\\[.5\baselineskip]
\Large через олимпиады и~кружки "--- к~профессии \par}

\vskip3\baselineskip

{\large

Под редакцией А.~А.\ Заславского,\\[.25\baselineskip] %Д.~А.\ Пермякова,\\[.25\baselineskip]
А.~Б.\ Скопенкова и М.~Б.\ Скопенкова% и~А.~В.\ Шаповалова

}

\vskip3\baselineskip

{\normalsize
Сокращенная версия, содержащая избранные материалы\\
 М.А. Берштейна, А.Д.Блинкова, В.А.Брагина, Ю.М.Бурмана, А.А.Гаврилюка, С.А.Дориченко, А.А.Заславского, А.Я.Канеля-Белова, А.А.Клячко, П.А.Кожевникова, О.А.Малиновской, Д.А.Пермякова, В.Ю.Протасова, А.Б.Скопенкова, М.Б.Скопенкова, А.В.Шаповалова, Ф.А.Шарова

\vskip2\baselineskip

Обновляемая бета-версия:
\url{http://www.mccme.ru/circles/oim/materials/sturm.pdf}
}

\par}

\vfill

\normalsize Москва\\ Издательство МЦНМО\\
\year\thispagestyle{empty}\par}

\clearpage

\thispagestyle{empty}

\vskip.15\baselineskip

\noindent\hbox to 0pt{\hskip-\leftskip\hskip4em \hss}\hskip\parindent
Элементы математики в~задачах. Через олимпиады и~кружки "--- к~профессии / Под ред. А.\,А.\ Заславского, А.\,Б.\ Скопенкова и М.\,Б.\ Скопенкова. "--- М.: МЦНМО, \year.~"--- 592\:с.

\vskip.5\baselineskip

{%\footnotesize

В данный сборник вошли материалы выездных школ по подготовке команды Москвы на Всероссийскую олимпиаду. Материалы сборника могут использоваться как школьниками для самостоятельных занятий, так и преподавателями. В большинстве материалов сборника приведены дававшиеся на занятиях задачи, а также решения или указания к ключевым задачам.

\par}

\clearpage

%\tableofcontents

{\LARGE \bf Оглавление полной версии сборника}

{
\begin{center}
\vspace{-0.8cm}
\mbox{
\normalsize \bf Материалы, находящиеся в этом файле, выделены \ull{цветом}
}
\end{center}
}

%\babel@toc {russian}{}
\contentsline {section}{\ull{\numberline {1}\IeC {\CYRO }\IeC {\cyrt } \IeC {\cyrr }\IeC {\cyre }\IeC {\cyrd }\IeC {\cyra }\IeC {\cyrk }\IeC {\cyrt }\IeC {\cyro }\IeC {\cyrr }\IeC {\cyro }\IeC {\cyrv }}}{12}{section.0.1}
\contentsline {subsection}{\ull{\numberline {1.1}\IeC {\CYRZ }\IeC {\cyra }\IeC {\cyrch }\IeC {\cyre }\IeC {\cyrm } \IeC {\cyri }~\IeC {\cyrd }\IeC {\cyrl }\IeC {\cyrya } \IeC {\cyrk }\IeC {\cyro }\IeC {\cyrg }\IeC {\cyro } \IeC {\cyrerev }\IeC {\cyrt }\IeC {\cyra } \IeC {\cyrk }\IeC {\cyrn }\IeC {\cyri }\IeC {\cyrg }\IeC {\cyra }}}{12}{subsection.0.1.1}
\contentsline {subsection}{\ull{\numberline {1.2}\IeC {\CYRI }\IeC {\cyrz }\IeC {\cyru }\IeC {\cyrch }\IeC {\cyre }\IeC {\cyrn }\IeC {\cyri }\IeC {\cyre } \IeC {\cyrp }\IeC {\cyru }\IeC {\cyrt }\IeC {\cyryo }\IeC {\cyrm } \IeC {\cyrr }\IeC {\cyre }\IeC {\cyrsh }\IeC {\cyre }\IeC {\cyrn }\IeC {\cyri }\IeC {\cyrya } \IeC {\cyri }~\IeC {\cyro }\IeC {\cyrb }\IeC {\cyrs }\IeC {\cyru }\IeC {\cyrzh }\IeC {\cyrd }\IeC {\cyre }\IeC {\cyrn }\IeC {\cyri }\IeC {\cyrya } \IeC {\cyrz }\IeC {\cyra }\IeC {\cyrd }\IeC {\cyra }\IeC {\cyrch }}}{13}{subsection.0.1.2}
\contentsline {subsection}{\ull{\numberline {1.3}\IeC {\CYRK }\IeC {\cyra }\IeC {\cyrk } \IeC {\cyru }\IeC {\cyrs }\IeC {\cyrt }\IeC {\cyrr }\IeC {\cyro }\IeC {\cyre }\IeC {\cyrn }\IeC {\cyra } \IeC {\cyrk }\IeC {\cyrn }\IeC {\cyri }\IeC {\cyrg }\IeC {\cyra }}}{14}{subsection.0.1.3}
\contentsline {subsection}{\ull{\numberline {1.4}\IeC {\CYRN }\IeC {\cyra }\IeC {\cyrp }\IeC {\cyru }\IeC {\cyrt }\IeC {\cyrs }\IeC {\cyrt }\IeC {\cyrv }\IeC {\cyri }\IeC {\cyre }. \emph {\IeC {\CYRA }.~\IeC {\CYRYA }.~\IeC {\CYRK }\IeC {\cyra }\IeC {\cyrn }\IeC {\cyre }\IeC {\cyrl }\IeC {\cyrsftsn }-\IeC {\CYRB }\IeC {\cyre }\IeC {\cyrl }\IeC {\cyro }\IeC {\cyrv }}}}{15}{subsection.0.1.4}
\contentsline {subsection}{\ull{\numberline {1.5}\IeC {\CYRO } \IeC {\cyrl }\IeC {\cyri }\IeC {\cyrt }\IeC {\cyre }\IeC {\cyrr }\IeC {\cyra }\IeC {\cyrt }\IeC {\cyru }\IeC {\cyrr }\IeC {\cyre } \IeC {\cyri }~\IeC {\cyri }\IeC {\cyrs }\IeC {\cyrt }\IeC {\cyro }\IeC {\cyrch }\IeC {\cyrn }\IeC {\cyri }\IeC {\cyrk }\IeC {\cyra }\IeC {\cyrh }}}{15}{subsection.0.1.5}
\contentsline {subsection}{\ull{\numberline {1.6}\IeC {\CYRB }\IeC {\cyrl }\IeC {\cyra }\IeC {\cyrg }\IeC {\cyro }\IeC {\cyrd }\IeC {\cyra }\IeC {\cyrr }\IeC {\cyrn }\IeC {\cyro }\IeC {\cyrs }\IeC {\cyrt }\IeC {\cyri } \IeC {\cyri }~\IeC {\cyrs }\IeC {\cyrv }\IeC {\cyre }\IeC {\cyrd }\IeC {\cyre }\IeC {\cyrn }\IeC {\cyri }\IeC {\cyrya } \IeC {\cyro }\IeC {\cyrb } \IeC {\cyra }\IeC {\cyrv }\IeC {\cyrt }\IeC {\cyro }\IeC {\cyrr }\IeC {\cyra }\IeC {\cyrh }}}{16}{subsection.0.1.6}
\contentsline {subsection}{\ull{\numberline {1.7}\IeC {\CYRV }\IeC {\cyra }\IeC {\cyrzh }\IeC {\cyrn }\IeC {\cyrery }\IeC {\cyre } \IeC {\cyrs }\IeC {\cyro }\IeC {\cyrg }\IeC {\cyrl }\IeC {\cyra }\IeC {\cyrsh }\IeC {\cyre }\IeC {\cyrn }\IeC {\cyri }\IeC {\cyrya }}}{17}{subsection.0.1.7}
\contentsline {subsection}{\ull{\numberline {1.8}\IeC {\CYRO }\IeC {\cyrs }\IeC {\cyrn }\IeC {\cyro }\IeC {\cyrv }\IeC {\cyrn }\IeC {\cyrery }\IeC {\cyre } \IeC {\cyro }\IeC {\cyrb }\IeC {\cyro }\IeC {\cyrz }\IeC {\cyrn }\IeC {\cyra }\IeC {\cyrch }\IeC {\cyre }\IeC {\cyrn }\IeC {\cyri }\IeC {\cyrya }}}{18}{subsection.0.1.8}
\contentsline {chapter}{\numberline {1}\IeC {\CYRT }\IeC {\cyre }\IeC {\cyro }\IeC {\cyrr }\IeC {\cyri }\IeC {\cyrya } \IeC {\cyrch }\IeC {\cyri }\IeC {\cyrs }\IeC {\cyre }\IeC {\cyrl }, \IeC {\cyra }\IeC {\cyrl }\IeC {\cyrg }\IeC {\cyre }\IeC {\cyrb }\IeC {\cyrr }\IeC {\cyra } \IeC {\cyri }~\IeC {\cyra }\IeC {\cyrn }\IeC {\cyra }\IeC {\cyrl }\IeC {\cyri }\IeC {\cyrz }. \emph {\IeC {\CYRA }.~\IeC {\CYRB }.~\IeC {\CYRS }\IeC {\cyrk }\IeC {\cyro }\IeC {\cyrp }\IeC {\cyre }\IeC {\cyrn }\IeC {\cyrk }\IeC {\cyro }\IeC {\cyrv }}}{21}{chapter.1}
\contentsline {section}{\numberline {2}\IeC {\CYRD }\IeC {\cyre }\IeC {\cyrl }\IeC {\cyri }\IeC {\cyrm }\IeC {\cyro }\IeC {\cyrs }\IeC {\cyrt }\IeC {\cyrsftsn } \IeC {\cyri }~\IeC {\cyrd }\IeC {\cyre }\IeC {\cyrl }\IeC {\cyre }\IeC {\cyrn }\IeC {\cyri }\IeC {\cyre } \IeC {\cyrs }~\IeC {\cyro }\IeC {\cyrs }\IeC {\cyrt }\IeC {\cyra }\IeC {\cyrt }\IeC {\cyrk }\IeC {\cyro }\IeC {\cyrm }}{21}{section.1.2}
\contentsline {subsection}{\ull{\numberline {2.1}\IeC {\CYRD }\IeC {\cyre }\IeC {\cyrl }\IeC {\cyri }\IeC {\cyrm }\IeC {\cyro }\IeC {\cyrs }\IeC {\cyrt }\IeC {\cyrsftsn } (1)}}{21}{subsection.1.2.1}
\contentsline {subsection}{\numberline {2.2}\IeC {\CYRP }\IeC {\cyrr }\IeC {\cyro }\IeC {\cyrs }\IeC {\cyrt }\IeC {\cyrery }\IeC {\cyre } \IeC {\cyrch }\IeC {\cyri }\IeC {\cyrs }\IeC {\cyrl }\IeC {\cyra } (1)}{25}{}%{subsection.1.2.2}
\contentsline {subsection}{\numberline {2.3}\IeC {\CYRN }\IeC {\CYRO }\IeC {\CYRD } \IeC {\cyri }~\IeC {\CYRN }\IeC {\CYRO }\IeC {\CYRK } (1)}{28}{}%{subsection.1.2.3}
\contentsline {subsection}{\numberline {2.4}\IeC {\CYRD }\IeC {\cyre }\IeC {\cyrl }\IeC {\cyre }\IeC {\cyrn }\IeC {\cyri }\IeC {\cyre } \IeC {\cyrs }~\IeC {\cyro }\IeC {\cyrs }\IeC {\cyrt }\IeC {\cyra }\IeC {\cyrt }\IeC {\cyrk }\IeC {\cyro }\IeC {\cyrm } \IeC {\cyri }~\IeC {\cyrs }\IeC {\cyrr }\IeC {\cyra }\IeC {\cyrv }\IeC {\cyrn }\IeC {\cyre }\IeC {\cyrn }\IeC {\cyri }\IeC {\cyrya } (1)}{30}{}%{subsection.1.2.4}
\contentsline {subsection}{\numberline {2.5}\IeC {\CYRL }\IeC {\cyri }\IeC {\cyrn }\IeC {\cyre }\IeC {\cyrishrt }\IeC {\cyrn }\IeC {\cyrery }\IeC {\cyre } \IeC {\cyrd }\IeC {\cyri }\IeC {\cyro }\IeC {\cyrf }\IeC {\cyra }\IeC {\cyrn }\IeC {\cyrt }\IeC {\cyro }\IeC {\cyrv }\IeC {\cyrery } \IeC {\cyru }\IeC {\cyrr }\IeC {\cyra }\IeC {\cyrv }\IeC {\cyrn }\IeC {\cyre }\IeC {\cyrn }\IeC {\cyri }\IeC {\cyrya } (2)}{32}{}%{subsection.1.2.5}
\contentsline {subsection}{\numberline {2.6}\IeC {\CYRK }\IeC {\cyra }\IeC {\cyrn }\IeC {\cyro }\IeC {\cyrn }\IeC {\cyri }\IeC {\cyrch }\IeC {\cyre }\IeC {\cyrs }\IeC {\cyrk }\IeC {\cyro }\IeC {\cyre } \IeC {\cyrr }\IeC {\cyra }\IeC {\cyrz }\IeC {\cyrl }\IeC {\cyro }\IeC {\cyrzh }\IeC {\cyre }\IeC {\cyrn }\IeC {\cyri }\IeC {\cyre } (2*)}{35}{}%{subsection.1.2.6}
\contentsline {subsection}{\numberline {2.7}\IeC {\CYRC }\IeC {\cyre }\IeC {\cyrl }\IeC {\cyrery }\IeC {\cyre } \IeC {\cyrt }\IeC {\cyro }\IeC {\cyrch }\IeC {\cyrk }\IeC {\cyri } \IeC {\cyrp }\IeC {\cyro }\IeC {\cyrd } \IeC {\cyrp }\IeC {\cyrr }\IeC {\cyrya }\IeC {\cyrm }\IeC {\cyro }\IeC {\cyrishrt } (2*)}{38}{}%{subsection.1.2.7}
\contentsline {section}{\numberline {3}\IeC {\CYRU }\IeC {\cyrm }\IeC {\cyrn }\IeC {\cyro }\IeC {\cyrzh }\IeC {\cyre }\IeC {\cyrn }\IeC {\cyri }\IeC {\cyre } \IeC {\cyrp }\IeC {\cyro } \IeC {\cyrp }\IeC {\cyrr }\IeC {\cyro }\IeC {\cyrs }\IeC {\cyrt }\IeC {\cyro }\IeC {\cyrm }\IeC {\cyru } \IeC {\cyrm }\IeC {\cyro }\IeC {\cyrd }\IeC {\cyru }\IeC {\cyrl }\IeC {\cyryu }}{42}{section.1.3}
\contentsline {subsection}{\ull{\numberline {3.1}\IeC {\CYRM }\IeC {\cyra }\IeC {\cyrl }\IeC {\cyra }\IeC {\cyrya } \IeC {\cyrt }\IeC {\cyre }\IeC {\cyro }\IeC {\cyrr }\IeC {\cyre }\IeC {\cyrm }\IeC {\cyra } \IeC {\CYRF }\IeC {\cyre }\IeC {\cyrr }\IeC {\cyrm }\IeC {\cyra } (2)}}{42}{subsection.1.3.1}
\contentsline {subsection}{\numberline {3.2}\IeC {\CYRP }\IeC {\cyrr }\IeC {\cyro }\IeC {\cyrv }\IeC {\cyre }\IeC {\cyrr }\IeC {\cyrk }\IeC {\cyra } \IeC {\cyrp }\IeC {\cyrr }\IeC {\cyro }\IeC {\cyrs }\IeC {\cyrt }\IeC {\cyro }\IeC {\cyrt }\IeC {\cyrery } (3*). \emph {\IeC {\CYRS }.~\IeC {\CYRV }.~\IeC {\CYRK }\IeC {\cyro }\IeC {\cyrn }\IeC {\cyrya }\IeC {\cyrg }\IeC {\cyri }\IeC {\cyrn }}}{44}{}%{subsection.1.3.2}
\contentsline {subsection}{\ull{\numberline {3.3}\IeC {\CYRK }\IeC {\cyrv }\IeC {\cyra }\IeC {\cyrd }\IeC {\cyrr }\IeC {\cyra }\IeC {\cyrt }\IeC {\cyri }\IeC {\cyrch }\IeC {\cyrn }\IeC {\cyrery }\IeC {\cyre } \IeC {\cyrv }\IeC {\cyrery }\IeC {\cyrch }\IeC {\cyre }\IeC {\cyrt }\IeC {\cyrery } (2*)}}{47}{subsection.1.3.3}
\contentsline {subsection}{\ull{\numberline {3.4}\IeC {\CYRK }\IeC {\cyrv }\IeC {\cyra }\IeC {\cyrd }\IeC {\cyrr }\IeC {\cyra }\IeC {\cyrt }\IeC {\cyri }\IeC {\cyrch }\IeC {\cyrn }\IeC {\cyrery }\IeC {\cyrishrt } \IeC {\cyrz }\IeC {\cyra }\IeC {\cyrk }\IeC {\cyro }\IeC {\cyrn } \IeC {\cyrv }\IeC {\cyrz }\IeC {\cyra }\IeC {\cyri }\IeC {\cyrm }\IeC {\cyrn }\IeC {\cyro }\IeC {\cyrs }\IeC {\cyrt }\IeC {\cyri } (3*)}}{50}{subsection.1.3.4}
\contentsline {subsection}{\ull{\numberline {3.5}\IeC {\CYRP }\IeC {\cyre }\IeC {\cyrr }\IeC {\cyrv }\IeC {\cyro }\IeC {\cyro }\IeC {\cyrb }\IeC {\cyrr }\IeC {\cyra }\IeC {\cyrz }\IeC {\cyrn }\IeC {\cyrery }\IeC {\cyre } \IeC {\cyrk }\IeC {\cyro }\IeC {\cyrr }\IeC {\cyrn }\IeC {\cyri } (3*)}}{53}{subsection.1.3.5}
\contentsline {subsection}{\ull{\numberline {3.6}\IeC {\CYRV }\IeC {\cyrery }\IeC {\cyrs }\IeC {\cyro }\IeC {\cyrk }\IeC {\cyri }\IeC {\cyre } \IeC {\cyrs }\IeC {\cyrt }\IeC {\cyre }\IeC {\cyrp }\IeC {\cyre }\IeC {\cyrn }\IeC {\cyri } (3*). \emph {\IeC {\CYRA }.~\IeC {\CYRYA }.~\IeC {\CYRK }\IeC {\cyra }\IeC {\cyrn }\IeC {\cyre }\IeC {\cyrl }\IeC {\cyrsftsn }-\IeC {\CYRB }\IeC {\cyre }\IeC {\cyrl }\IeC {\cyro }\IeC {\cyrv }}, \emph {\IeC {\CYRA }.~\IeC {\CYRB }.~\IeC {\CYRS }\IeC {\cyrk }\IeC {\cyro }\IeC {\cyrp }\IeC {\cyre }\IeC {\cyrn }\IeC {\cyrk }\IeC {\cyro }\IeC {\cyrv }}}}{55}{subsection.1.3.6}
\contentsline {section}{\numberline {4}\IeC {\CYRM }\IeC {\cyrn }\IeC {\cyro }\IeC {\cyrg }\IeC {\cyro }\IeC {\cyrch }\IeC {\cyrl }\IeC {\cyre }\IeC {\cyrn }\IeC {\cyrery } \IeC {\cyri }~\IeC {\cyrk }\IeC {\cyro }\IeC {\cyrm }\IeC {\cyrp }\IeC {\cyrl }\IeC {\cyre }\IeC {\cyrk }\IeC {\cyrs }\IeC {\cyrn }\IeC {\cyrery }\IeC {\cyre } \IeC {\cyrch }\IeC {\cyri }\IeC {\cyrs }\IeC {\cyrl }\IeC {\cyra }}{59}{section.1.4}
\contentsline {subsection}{\numberline {4.1}\IeC {\CYRR }\IeC {\cyra }\IeC {\cyrc }\IeC {\cyri }\IeC {\cyro }\IeC {\cyrn }\IeC {\cyra }\IeC {\cyrl }\IeC {\cyrsftsn }\IeC {\cyrn }\IeC {\cyrery }\IeC {\cyre } \IeC {\cyri }~\IeC {\cyri }\IeC {\cyrr }\IeC {\cyrr }\IeC {\cyra }\IeC {\cyrc }\IeC {\cyri }\IeC {\cyro }\IeC {\cyrn }\IeC {\cyra }\IeC {\cyrl }\IeC {\cyrsftsn }\IeC {\cyrn }\IeC {\cyrery }\IeC {\cyre } \IeC {\cyrch }\IeC {\cyri }\IeC {\cyrs }\IeC {\cyrl }\IeC {\cyra } (1)}{59}{}%{subsection.1.4.1}
\contentsline {subsection}{\ull{\numberline {4.2}\IeC {\CYRR }\IeC {\cyre }\IeC {\cyrsh }\IeC {\cyre }\IeC {\cyrn }\IeC {\cyri }\IeC {\cyre } \IeC {\cyru }\IeC {\cyrr }\IeC {\cyra }\IeC {\cyrv }\IeC {\cyrn }\IeC {\cyre }\IeC {\cyrn }\IeC {\cyri }\IeC {\cyrishrt } 3"~\IeC {\cyrishrt } \IeC {\cyri }~4"~\IeC {\cyrishrt } \IeC {\cyrs }\IeC {\cyrt }\IeC {\cyre }\IeC {\cyrp }\IeC {\cyre }\IeC {\cyrn }\IeC {\cyri } (2)}}{62}{subsection.1.4.2}
\contentsline {subsection}{\numberline {4.3}\IeC {\CYRT }\IeC {\cyre }\IeC {\cyro }\IeC {\cyrr }\IeC {\cyre }\IeC {\cyrm }\IeC {\cyra } \IeC {\CYRB }\IeC {\cyre }\IeC {\cyrz }\IeC {\cyru } \IeC {\cyri }~\IeC {\cyre }\IeC {\cyryo } \IeC {\cyrs }\IeC {\cyrl }\IeC {\cyre }\IeC {\cyrd }\IeC {\cyrs }\IeC {\cyrt }\IeC {\cyrv }\IeC {\cyri }\IeC {\cyrya } (2)}{68}{}%{subsection.1.4.3}
\contentsline {subsection}{\numberline {4.4}\IeC {\CYRD }\IeC {\cyre }\IeC {\cyrl }\IeC {\cyri }\IeC {\cyrm }\IeC {\cyro }\IeC {\cyrs }\IeC {\cyrt }\IeC {\cyrsftsn } \IeC {\cyrd }\IeC {\cyrl }\IeC {\cyrya } \IeC {\cyrm }\IeC {\cyrn }\IeC {\cyro }\IeC {\cyrg }\IeC {\cyro }\IeC {\cyrch }\IeC {\cyrl }\IeC {\cyre }\IeC {\cyrn }\IeC {\cyro }\IeC {\cyrv } (3*). \emph {\IeC {\CYRA }.~\IeC {\CYRYA }.~\IeC {\CYRK }\IeC {\cyra }\IeC {\cyrn }\IeC {\cyre }\IeC {\cyrl }\IeC {\cyrsftsn }-\IeC {\CYRB }\IeC {\cyre }\IeC {\cyrl }\IeC {\cyro }\IeC {\cyrv }}, \emph {\IeC {\CYRA }.~\IeC {\CYRB }.~\IeC {\CYRS }\IeC {\cyrk }\IeC {\cyro }\IeC {\cyrp }\IeC {\cyre }\IeC {\cyrn }\IeC {\cyrk }\IeC {\cyro }\IeC {\cyrv }}}{71}{}%{subsection.1.4.4}
\contentsline {subsection}{\numberline {4.5}\IeC {\CYRP }\IeC {\cyrr }\IeC {\cyri }\IeC {\cyrm }\IeC {\cyre }\IeC {\cyrn }\IeC {\cyre }\IeC {\cyrn }\IeC {\cyri }\IeC {\cyrya } \IeC {\cyrk }\IeC {\cyro }\IeC {\cyrm }\IeC {\cyrp }\IeC {\cyrl }\IeC {\cyre }\IeC {\cyrk }\IeC {\cyrs }\IeC {\cyrn }\IeC {\cyrery }\IeC {\cyrh } \IeC {\cyrch }\IeC {\cyri }\IeC {\cyrs }\IeC {\cyre }\IeC {\cyrl } (3*)}{74}{}%{subsection.1.4.5}
\contentsline {subsection}{\numberline {4.6}\IeC {\CYRT }\IeC {\cyre }\IeC {\cyro }\IeC {\cyrr }\IeC {\cyre }\IeC {\cyrm }\IeC {\cyra } \IeC {\CYRV }\IeC {\cyri }\IeC {\cyre }\IeC {\cyrt }\IeC {\cyra } \IeC {\cyri }~\IeC {\cyrs }\IeC {\cyri }\IeC {\cyrm }\IeC {\cyrm }\IeC {\cyre }\IeC {\cyrt }\IeC {\cyrr }\IeC {\cyri }\IeC {\cyrch }\IeC {\cyre }\IeC {\cyrs }\IeC {\cyrk }\IeC {\cyri }\IeC {\cyre } \IeC {\cyrm }\IeC {\cyrn }\IeC {\cyro }\IeC {\cyrg }\IeC {\cyro }\IeC {\cyrch }\IeC {\cyrl }\IeC {\cyre }\IeC {\cyrn }\IeC {\cyrery } (3*)}{77}{}%{subsection.1.4.6}
\contentsline {subsection}{\ull{\numberline {4.7}\IeC {\CYRD }\IeC {\cyri }\IeC {\cyro }\IeC {\cyrf }\IeC {\cyra }\IeC {\cyrn }\IeC {\cyrt }\IeC {\cyro }\IeC {\cyrv }\IeC {\cyrery } \IeC {\cyru }\IeC {\cyrr }\IeC {\cyra }\IeC {\cyrv }\IeC {\cyrn }\IeC {\cyre }\IeC {\cyrn }\IeC {\cyri }\IeC {\cyrya } \IeC {\cyri }~\IeC {\cyrg }\IeC {\cyra }\IeC {\cyru }\IeC {\cyrs }\IeC {\cyrs }\IeC {\cyro }\IeC {\cyrv }\IeC {\cyrery } \IeC {\cyrch }\IeC {\cyri }\IeC {\cyrs }\IeC {\cyrl }\IeC {\cyra } (4*). \emph {\IeC {\CYRA }.~\IeC {\CYRYA }.~\IeC {\CYRK }\IeC {\cyra }\IeC {\cyrn }\IeC {\cyre }\IeC {\cyrl }\IeC {\cyrsftsn }-\IeC {\CYRB }\IeC {\cyre }\IeC {\cyrl }\IeC {\cyro }\IeC {\cyrv }}}}{79}{subsection.1.4.7}
\contentsline {subsection}{\numberline {4.8}\IeC {\CYRD }\IeC {\cyri }\IeC {\cyra }\IeC {\cyrg }\IeC {\cyro }\IeC {\cyrn }\IeC {\cyra }\IeC {\cyrl }\IeC {\cyri } \IeC {\cyrp }\IeC {\cyrr }\IeC {\cyra }\IeC {\cyrv }\IeC {\cyri }\IeC {\cyrl }\IeC {\cyrsftsn }\IeC {\cyrn }\IeC {\cyrery }\IeC {\cyrh } \IeC {\cyrm }\IeC {\cyrn }\IeC {\cyro }\IeC {\cyrg }\IeC {\cyro }\IeC {\cyru }\IeC {\cyrg }\IeC {\cyro }\IeC {\cyrl }\IeC {\cyrsftsn }\IeC {\cyrn }\IeC {\cyri }\IeC {\cyrk }\IeC {\cyro }\IeC {\cyrv } (4*). \emph {\IeC {\CYRI }.~\IeC {\CYRN }.~\IeC {\CYRSH }\IeC {\cyrn }\IeC {\cyru }\IeC {\cyrr }\IeC {\cyrn }\IeC {\cyri }\IeC {\cyrk }\IeC {\cyro }\IeC {\cyrv }}}{83}{}%{subsection.1.4.8}
\contentsline {section}{\numberline {5}\IeC {\CYRR }\IeC {\cyra }\IeC {\cyrz }\IeC {\cyrr }\IeC {\cyre }\IeC {\cyrsh }\IeC {\cyri }\IeC {\cyrm }\IeC {\cyro }\IeC {\cyrs }\IeC {\cyrt }\IeC {\cyrsftsn } \IeC {\cyrv } \IeC {\cyrr }\IeC {\cyra }\IeC {\cyrd }\IeC {\cyri }\IeC {\cyrk }\IeC {\cyra }\IeC {\cyrl }\IeC {\cyra }\IeC {\cyrh }}{88}{section.1.5}
\contentsline {subsection}{\ull{\numberline {5.1}\IeC {\CYRV }\IeC {\cyrv }\IeC {\cyre }\IeC {\cyrd }\IeC {\cyre }\IeC {\cyrn }\IeC {\cyri }\IeC {\cyre }}}{88}{subsection.1.5.1}
\contentsline {subsubsection}{\ull{\numberline {5.1.1}\IeC {\CYRO } \IeC {\cyrch }\IeC {\cyryo }\IeC {\cyrm } \IeC {\cyrerev }\IeC {\cyrt }\IeC {\cyro }\IeC {\cyrt } \IeC {\cyrp }\IeC {\cyra }\IeC {\cyrr }\IeC {\cyra }\IeC {\cyrg }\IeC {\cyrr }\IeC {\cyra }\IeC {\cyrf }}}{88}{subsubsection.1.5.1.1}
\contentsline {subsubsection}{\ull{\numberline {5.1.2}\IeC {\CYRR }\IeC {\cyra }\IeC {\cyrz }\IeC {\cyrr }\IeC {\cyre }\IeC {\cyrsh }\IeC {\cyri }\IeC {\cyrm }\IeC {\cyro }\IeC {\cyrs }\IeC {\cyrt }\IeC {\cyrsftsn } \IeC {\cyrv }~\IeC {\cyrk }\IeC {\cyrv }\IeC {\cyra }\IeC {\cyrd }\IeC {\cyrr }\IeC {\cyra }\IeC {\cyrt }\IeC {\cyrn }\IeC {\cyrery }\IeC {\cyrh } \IeC {\cyrr }\IeC {\cyra }\IeC {\cyrd }\IeC {\cyri }\IeC {\cyrk }\IeC {\cyra }\IeC {\cyrl }\IeC {\cyra }\IeC {\cyrh }: \IeC {\cyrf }\IeC {\cyro }\IeC {\cyrr }\IeC {\cyrm }\IeC {\cyru }\IeC {\cyrl }\IeC {\cyri }\IeC {\cyrr }\IeC {\cyro }\IeC {\cyrv }\IeC {\cyrk }\IeC {\cyri }~(1)}}{90}{subsubsection.1.5.1.2}
\contentsline {subsubsection}{\ull{\numberline {5.1.3}\IeC {\CYRN }\IeC {\cyre }\IeC {\cyrr }\IeC {\cyra }\IeC {\cyrz }\IeC {\cyrr }\IeC {\cyre }\IeC {\cyrsh }\IeC {\cyri }\IeC {\cyrm }\IeC {\cyro }\IeC {\cyrs }\IeC {\cyrt }\IeC {\cyrsftsn } \IeC {\cyrv }~\IeC {\cyrr }\IeC {\cyra }\IeC {\cyrd }\IeC {\cyri }\IeC {\cyrk }\IeC {\cyra }\IeC {\cyrl }\IeC {\cyra }\IeC {\cyrh }: \IeC {\cyrf }\IeC {\cyro }\IeC {\cyrr }\IeC {\cyrm }\IeC {\cyru }\IeC {\cyrl }\IeC {\cyri }\IeC {\cyrr }\IeC {\cyro }\IeC {\cyrv }\IeC {\cyrk }\IeC {\cyri } (2)}}{91}{subsubsection.1.5.1.3}
\contentsline {subsubsection}{\ull{\numberline {5.1.4}\IeC {\CYRP }\IeC {\cyrl }\IeC {\cyra }\IeC {\cyrn } \IeC {\cyrp }\IeC {\cyra }\IeC {\cyrr }\IeC {\cyra }\IeC {\cyrg }\IeC {\cyrr }\IeC {\cyra }\IeC {\cyrf }\IeC {\cyra }}}{94}{subsubsection.1.5.1.4}
\contentsline {subsection}{\numberline {5.2}\IeC {\CYRV }\IeC {\cyra }\IeC {\cyrzh }\IeC {\cyrn }\IeC {\cyrery }\IeC {\cyre } \IeC {\cyro }\IeC {\cyrt }\IeC {\cyrs }\IeC {\cyrt }\IeC {\cyru }\IeC {\cyrp }\IeC {\cyrl }\IeC {\cyre }\IeC {\cyrn }\IeC {\cyri }\IeC {\cyrya }}{95}{}%{subsection.1.5.2}
\contentsline {subsubsection}{\numberline {5.2.1}\IeC {\CYRCH }\IeC {\cyre }\IeC {\cyrm } \IeC {\cyri }\IeC {\cyrn }\IeC {\cyrt }\IeC {\cyre }\IeC {\cyrr }\IeC {\cyre }\IeC {\cyrs }\IeC {\cyrn }\IeC {\cyrery } \IeC {\cyrp }\IeC {\cyrr }\IeC {\cyri }\IeC {\cyrv }\IeC {\cyro }\IeC {\cyrd }\IeC {\cyri }\IeC {\cyrm }\IeC {\cyrery }\IeC {\cyre } \IeC {\cyrd }\IeC {\cyro }\IeC {\cyrk }\IeC {\cyra }\IeC {\cyrz }\IeC {\cyra }\IeC {\cyrt }\IeC {\cyre }\IeC {\cyrl }\IeC {\cyrsftsn }\IeC {\cyrs }\IeC {\cyrt }\IeC {\cyrv }\IeC {\cyra }}{95}{}%{subsubsection.1.5.2.1}
\contentsline {subsubsection}{\numberline {5.2.2}\IeC {\CYRI }\IeC {\cyrs }\IeC {\cyrt }\IeC {\cyro }\IeC {\cyrr }\IeC {\cyri }\IeC {\cyrch }\IeC {\cyre }\IeC {\cyrs }\IeC {\cyrk }\IeC {\cyri }\IeC {\cyre } \IeC {\cyrk }\IeC {\cyro }\IeC {\cyrm }\IeC {\cyrm }\IeC {\cyre }\IeC {\cyrn }\IeC {\cyrt }\IeC {\cyra }\IeC {\cyrr }\IeC {\cyri }\IeC {\cyri }}{96}{}%{subsubsection.1.5.2.2}
\contentsline {subsubsection}{\numberline {5.2.3}\IeC {\CYRS }\IeC {\cyrv }\IeC {\cyrya }\IeC {\cyrz }\IeC {\cyrsftsn } \IeC {\cyrs }~\IeC {\cyrp }\IeC {\cyro }\IeC {\cyrs }\IeC {\cyrt }\IeC {\cyrr }\IeC {\cyro }\IeC {\cyre }\IeC {\cyrn }\IeC {\cyri }\IeC {\cyrya }\IeC {\cyrm }\IeC {\cyri } \IeC {\cyrc }\IeC {\cyri }\IeC {\cyrr }\IeC {\cyrk }\IeC {\cyru }\IeC {\cyrl }\IeC {\cyre }\IeC {\cyrm } \IeC {\cyri }~\IeC {\cyrl }\IeC {\cyri }\IeC {\cyrn }\IeC {\cyre }\IeC {\cyrishrt }\IeC {\cyrk }\IeC {\cyro }\IeC {\cyrishrt } (1)}{97}{}%{subsubsection.1.5.2.3}
\contentsline {subsection}{\numberline {5.3}\IeC {\CYRD }\IeC {\cyro }\IeC {\cyrk }\IeC {\cyra }\IeC {\cyrz }\IeC {\cyra }\IeC {\cyrt }\IeC {\cyre }\IeC {\cyrl }\IeC {\cyrsftsn }\IeC {\cyrs }\IeC {\cyrt }\IeC {\cyrv }\IeC {\cyro } \IeC {\cyrp }\IeC {\cyro }\IeC {\cyrs }\IeC {\cyrt }\IeC {\cyrr }\IeC {\cyro }\IeC {\cyri }\IeC {\cyrm }\IeC {\cyro }\IeC {\cyrs }\IeC {\cyrt }\IeC {\cyri } \IeC {\cyrv }~\IeC {\cyrt }\IeC {\cyre }\IeC {\cyro }\IeC {\cyrr }\IeC {\cyre }\IeC {\cyrm }\IeC {\cyre } \IeC {\CYRG }\IeC {\cyra }\IeC {\cyru }\IeC {\cyrs }\IeC {\cyrs }\IeC {\cyra }}{98}{}%{subsection.1.5.3}
\contentsline {subsubsection}{\ull{\numberline {5.3.1}\IeC {\CYRP }\IeC {\cyre }\IeC {\cyrr }\IeC {\cyre }\IeC {\cyrf }\IeC {\cyro }\IeC {\cyrr }\IeC {\cyrm }\IeC {\cyru }\IeC {\cyrl }\IeC {\cyri }\IeC {\cyrr }\IeC {\cyro }\IeC {\cyrv }\IeC {\cyrk }\IeC {\cyra } \IeC {\cyrp }\IeC {\cyro }\IeC {\cyrs }\IeC {\cyrt }\IeC {\cyrr }\IeC {\cyro }\IeC {\cyri }\IeC {\cyrm }\IeC {\cyro }\IeC {\cyrs }\IeC {\cyrt }\IeC {\cyri } \IeC {\cyrv }~\IeC {\cyrt }\IeC {\cyre }\IeC {\cyro }\IeC {\cyrr }\IeC {\cyre }\IeC {\cyrm }\IeC {\cyre } \IeC {\CYRG }\IeC {\cyra }\IeC {\cyru }\IeC {\cyrs }\IeC {\cyrs }\IeC {\cyra } (2)}}{98}{subsubsection.1.5.3.1}
\contentsline {subsubsection}{\ull{\numberline {5.3.2}\IeC {\CYRM }\IeC {\cyre }\IeC {\cyrt }\IeC {\cyro }\IeC {\cyrd } \IeC {\cyrr }\IeC {\cyre }\IeC {\cyrz }\IeC {\cyro }\IeC {\cyrl }\IeC {\cyrsftsn }\IeC {\cyrv }\IeC {\cyre }\IeC {\cyrn }\IeC {\cyrt } \IeC {\CYRL }\IeC {\cyra }\IeC {\cyrg }\IeC {\cyrr }\IeC {\cyra }\IeC {\cyrn }\IeC {\cyrzh }\IeC {\cyra } (2)}}{99}{subsubsection.1.5.3.2}
\contentsline {subsubsection}{\ull{\numberline {5.3.3}\IeC {\CYRD }\IeC {\cyro }\IeC {\cyrk }\IeC {\cyra }\IeC {\cyrz }\IeC {\cyra }\IeC {\cyrt }\IeC {\cyre }\IeC {\cyrl }\IeC {\cyrsftsn }\IeC {\cyrs }\IeC {\cyrt }\IeC {\cyrv }\IeC {\cyro } \IeC {\cyrp }\IeC {\cyro }\IeC {\cyrs }\IeC {\cyrt }\IeC {\cyrr }\IeC {\cyro }\IeC {\cyri }\IeC {\cyrm }\IeC {\cyro }\IeC {\cyrs }\IeC {\cyrt }\IeC {\cyri } \IeC {\cyrv }~\IeC {\cyrt }\IeC {\cyre }\IeC {\cyro }\IeC {\cyrr }\IeC {\cyre }\IeC {\cyrm }\IeC {\cyre } \IeC {\CYRG }\IeC {\cyra }\IeC {\cyru }\IeC {\cyrs }\IeC {\cyrs }\IeC {\cyra } (3)}}{105}{subsubsection.1.5.3.3}
\contentsline {subsubsection}{\numberline {5.3.4}\IeC {\CYREREV }\IeC {\cyrf }\IeC {\cyrf }\IeC {\cyre }\IeC {\cyrk }\IeC {\cyrt }\IeC {\cyri }\IeC {\cyrv }\IeC {\cyrn }\IeC {\cyrery }\IeC {\cyre } \IeC {\cyrd }\IeC {\cyro }\IeC {\cyrk }\IeC {\cyra }\IeC {\cyrz }\IeC {\cyra }\IeC {\cyrt }\IeC {\cyre }\IeC {\cyrl }\IeC {\cyrsftsn }\IeC {\cyrs }\IeC {\cyrt }\IeC {\cyrv }\IeC {\cyra } \IeC {\cyrp }\IeC {\cyro }\IeC {\cyrs }\IeC {\cyrt }\IeC {\cyrr }\IeC {\cyro }\IeC {\cyri }\IeC {\cyrm }\IeC {\cyro }\IeC {\cyrs }\IeC {\cyrt }\IeC {\cyri } (4*)}{106}{}%{subsubsection.1.5.3.4}
\contentsline {subsection}{\numberline {5.4}\IeC {\CYRZ }\IeC {\cyra }\IeC {\cyrd }\IeC {\cyra }\IeC {\cyrch }\IeC {\cyri } \IeC {\cyro }~\IeC {\cyrn }\IeC {\cyre }\IeC {\cyrr }\IeC {\cyra }\IeC {\cyrz }\IeC {\cyrr }\IeC {\cyre }\IeC {\cyrsh }\IeC {\cyri }\IeC {\cyrm }\IeC {\cyro }\IeC {\cyrs }\IeC {\cyrt }\IeC {\cyri } \IeC {\cyrv }~\IeC {\cyrr }\IeC {\cyra }\IeC {\cyrd }\IeC {\cyri }\IeC {\cyrk }\IeC {\cyra }\IeC {\cyrl }\IeC {\cyra }\IeC {\cyrh }}{114}{}%{subsection.1.5.4}
\contentsline {subsubsection}{\ull{\numberline {5.4.1}\IeC {\CYRO }\IeC {\cyrd }\IeC {\cyrn }\IeC {\cyro } \IeC {\cyri }\IeC {\cyrz }\IeC {\cyrv }\IeC {\cyrl }\IeC {\cyre }\IeC {\cyrch }\IeC {\cyre }\IeC {\cyrn }\IeC {\cyri }\IeC {\cyre } \IeC {\cyrk }\IeC {\cyrv }\IeC {\cyra }\IeC {\cyrd }\IeC {\cyrr }\IeC {\cyra }\IeC {\cyrt }\IeC {\cyrn }\IeC {\cyro }\IeC {\cyrg }\IeC {\cyro } \IeC {\cyrk }\IeC {\cyro }\IeC {\cyrr }\IeC {\cyrn }\IeC {\cyrya } (1)}}{114}{subsubsection.1.5.4.1}
\contentsline {subsubsection}{\ull{\numberline {5.4.2}\IeC {\CYRO }\IeC {\cyrd }\IeC {\cyrn }\IeC {\cyro } \IeC {\cyri }\IeC {\cyrz }\IeC {\cyrv }\IeC {\cyrl }\IeC {\cyre }\IeC {\cyrch }\IeC {\cyre }\IeC {\cyrn }\IeC {\cyri }\IeC {\cyre } \IeC {\cyrk }\IeC {\cyro }\IeC {\cyrr }\IeC {\cyrn }\IeC {\cyrya } \IeC {\cyrch }\IeC {\cyre }\IeC {\cyrt }\IeC {\cyrv }\IeC {\cyryo }\IeC {\cyrr }\IeC {\cyrt }\IeC {\cyro }\IeC {\cyrishrt } \IeC {\cyrs }\IeC {\cyrt }\IeC {\cyre }\IeC {\cyrp }\IeC {\cyre }\IeC {\cyrn }\IeC {\cyri } (2*)}}{118}{subsubsection.1.5.4.2}
\contentsline {subsubsection}{\ull{\numberline {5.4.3}\IeC {\CYRN }\IeC {\cyre }\IeC {\cyrs }\IeC {\cyrk }\IeC {\cyro }\IeC {\cyrl }\IeC {\cyrsftsn }\IeC {\cyrk }\IeC {\cyro } \IeC {\cyri }\IeC {\cyrz }\IeC {\cyrv }\IeC {\cyrl }\IeC {\cyre }\IeC {\cyrch }\IeC {\cyre }\IeC {\cyrn }\IeC {\cyri }\IeC {\cyrishrt } \IeC {\cyrk }\IeC {\cyrv }\IeC {\cyra }\IeC {\cyrd }\IeC {\cyrr }\IeC {\cyra }\IeC {\cyrt }\IeC {\cyrn }\IeC {\cyrery }\IeC {\cyrh } \IeC {\cyrk }\IeC {\cyro }\IeC {\cyrr }\IeC {\cyrn }\IeC {\cyre }\IeC {\cyrishrt } (3*)}}{121}{subsubsection.1.5.4.3}
\contentsline {subsubsection}{\numberline {5.4.4}\IeC {\CYRK } \IeC {\cyrd }\IeC {\cyro }\IeC {\cyrk }\IeC {\cyra }\IeC {\cyrz }\IeC {\cyra }\IeC {\cyrt }\IeC {\cyre }\IeC {\cyrl }\IeC {\cyrsftsn }\IeC {\cyrs }\IeC {\cyrt }\IeC {\cyrv }\IeC {\cyru } \IeC {\cyrn }\IeC {\cyre }\IeC {\cyrp }\IeC {\cyro }\IeC {\cyrs }\IeC {\cyrt }\IeC {\cyrr }\IeC {\cyro }\IeC {\cyri }\IeC {\cyrm }\IeC {\cyro }\IeC {\cyrs }\IeC {\cyrt }\IeC {\cyri } \IeC {\cyrv }~\IeC {\cyrt }\IeC {\cyre }\IeC {\cyro }\IeC {\cyrr }\IeC {\cyre }\IeC {\cyrm }\IeC {\cyre } \IeC {\CYRG }\IeC {\cyra }\IeC {\cyru }\IeC {\cyrs }\IeC {\cyrs }\IeC {\cyra } (4*)}{124}{}%{subsubsection.1.5.4.4}
\contentsline {subsubsection}{\numberline {5.4.5}\IeC {\CYRO }\IeC {\cyrd }\IeC {\cyrn }\IeC {\cyro } \IeC {\cyri }\IeC {\cyrz }\IeC {\cyrv }\IeC {\cyrl }\IeC {\cyre }\IeC {\cyrch }\IeC {\cyre }\IeC {\cyrn }\IeC {\cyri }\IeC {\cyre } \IeC {\cyrk }\IeC {\cyro }\IeC {\cyrr }\IeC {\cyrn }\IeC {\cyrya } \IeC {\cyrt }\IeC {\cyrr }\IeC {\cyre }\IeC {\cyrt }\IeC {\cyrsftsn }\IeC {\cyre }\IeC {\cyrishrt } \IeC {\cyrs }\IeC {\cyrt }\IeC {\cyre }\IeC {\cyrp }\IeC {\cyre }\IeC {\cyrn }\IeC {\cyri } (2)}{126}{}%{subsubsection.1.5.4.5}
\contentsline {subsubsection}{\numberline {5.4.6}\IeC {\CYRO }\IeC {\cyrd }\IeC {\cyrn }\IeC {\cyro } \IeC {\cyri }\IeC {\cyrz }\IeC {\cyrv }\IeC {\cyrl }\IeC {\cyre }\IeC {\cyrch }\IeC {\cyre }\IeC {\cyrn }\IeC {\cyri }\IeC {\cyre } \IeC {\cyrk }\IeC {\cyro }\IeC {\cyrr }\IeC {\cyrn }\IeC {\cyrya } \IeC {\cyrp }\IeC {\cyrr }\IeC {\cyro }\IeC {\cyrs }\IeC {\cyrt }\IeC {\cyro }\IeC {\cyrishrt } \IeC {\cyrs }\IeC {\cyrt }\IeC {\cyre }\IeC {\cyrp }\IeC {\cyre }\IeC {\cyrn }\IeC {\cyri } (3*)}{131}{}%{subsubsection.1.5.4.6}
\contentsline {subsubsection}{\numberline {5.4.7}\IeC {\CYRN }\IeC {\cyre }\IeC {\cyrs }\IeC {\cyrk }\IeC {\cyro }\IeC {\cyrl }\IeC {\cyrsftsn }\IeC {\cyrk }\IeC {\cyro } \IeC {\cyri }\IeC {\cyrz }\IeC {\cyrv }\IeC {\cyrl }\IeC {\cyre }\IeC {\cyrch }\IeC {\cyre }\IeC {\cyrn }\IeC {\cyri }\IeC {\cyrishrt } \IeC {\cyrk }\IeC {\cyro }\IeC {\cyrr }\IeC {\cyrn }\IeC {\cyre }\IeC {\cyrishrt } (4*)}{136}{}%{subsubsection.1.5.4.7}
\contentsline {subsection}{\numberline {5.5}\IeC {\CYRD }\IeC {\cyro }\IeC {\cyrk }\IeC {\cyra }\IeC {\cyrz }\IeC {\cyra }\IeC {\cyrt }\IeC {\cyre }\IeC {\cyrl }\IeC {\cyrsftsn }\IeC {\cyrs }\IeC {\cyrt }\IeC {\cyrv }\IeC {\cyra } \IeC {\cyrn }\IeC {\cyre }\IeC {\cyrr }\IeC {\cyra }\IeC {\cyrz }\IeC {\cyrr }\IeC {\cyre }\IeC {\cyrsh }\IeC {\cyri }\IeC {\cyrm }\IeC {\cyro }\IeC {\cyrs }\IeC {\cyrt }\IeC {\cyri } \IeC {\cyrv }~\IeC {\cyrr }\IeC {\cyra }\IeC {\cyrd }\IeC {\cyri }\IeC {\cyrk }\IeC {\cyra }\IeC {\cyrl }\IeC {\cyra }\IeC {\cyrh }}{137}{}%{subsection.1.5.5}
\contentsline {subsubsection}{\numberline {5.5.1}\IeC {\CYRL }\IeC {\cyre }\IeC {\cyrm }\IeC {\cyrm }\IeC {\cyra } \IeC {\cyro }~\IeC {\cyrk }\IeC {\cyra }\IeC {\cyrl }\IeC {\cyrsftsn }\IeC {\cyrk }\IeC {\cyru }\IeC {\cyrl }\IeC {\cyrya }\IeC {\cyrt }\IeC {\cyro }\IeC {\cyrr }\IeC {\cyre } \IeC {\cyri }~\IeC {\cyrp }\IeC {\cyro }\IeC {\cyrn }\IeC {\cyrya }\IeC {\cyrt }\IeC {\cyri }\IeC {\cyre } \IeC {\cyrp }\IeC {\cyro }\IeC {\cyrl }\IeC {\cyrya } (2*)}{138}{}%{subsubsection.1.5.5.1}
\contentsline {subsubsection}{\numberline {5.5.2}\IeC {\CYRD }\IeC {\cyro }\IeC {\cyrk }\IeC {\cyra }\IeC {\cyrz }\IeC {\cyra }\IeC {\cyrt }\IeC {\cyre }\IeC {\cyrl }\IeC {\cyrsftsn }\IeC {\cyrs }\IeC {\cyrt }\IeC {\cyrv }\IeC {\cyro } \IeC {\cyrn }\IeC {\cyre }\IeC {\cyrp }\IeC {\cyro }\IeC {\cyrs }\IeC {\cyrt }\IeC {\cyrr }\IeC {\cyro }\IeC {\cyri }\IeC {\cyrm }\IeC {\cyro }\IeC {\cyrs }\IeC {\cyrt }\IeC {\cyri } \IeC {\cyrv }~\IeC {\cyrt }\IeC {\cyre }\IeC {\cyro }\IeC {\cyrr }\IeC {\cyre }\IeC {\cyrm }\IeC {\cyre } \IeC {\CYRG }\IeC {\cyra }\IeC {\cyru }\IeC {\cyrs }\IeC {\cyrs }\IeC {\cyra } (3*)}{138}{}%{subsubsection.1.5.5.2}
\contentsline {subsubsection}{\numberline {5.5.3}\IeC {\CYRD }\IeC {\cyro }\IeC {\cyrk }\IeC {\cyra }\IeC {\cyrz }\IeC {\cyra }\IeC {\cyrt }\IeC {\cyre }\IeC {\cyrl }\IeC {\cyrsftsn }\IeC {\cyrs }\IeC {\cyrt }\IeC {\cyrv }\IeC {\cyro } \IeC {\cyrn }\IeC {\cyre }\IeC {\cyrr }\IeC {\cyra }\IeC {\cyrz }\IeC {\cyrr }\IeC {\cyre }\IeC {\cyrsh }\IeC {\cyri }\IeC {\cyrm }\IeC {\cyro }\IeC {\cyrs }\IeC {\cyrt }\IeC {\cyri } \IeC {\cyrv }~\IeC {\cyrv }\IeC {\cyre }\IeC {\cyrshch }\IeC {\cyre }\IeC {\cyrs }\IeC {\cyrt }\IeC {\cyrv }\IeC {\cyre }\IeC {\cyrn }\IeC {\cyrn }\IeC {\cyrery }\IeC {\cyrh } \IeC {\cyrr }\IeC {\cyra }\IeC {\cyrd }\IeC {\cyri }\IeC {\cyrk }\IeC {\cyra }\IeC {\cyrl }\IeC {\cyra }\IeC {\cyrh } (3*)}{140}{}%{subsubsection.1.5.5.3}
\contentsline {subsubsection}{\numberline {5.5.4}\IeC {\CYRD }\IeC {\cyro }\IeC {\cyrk }\IeC {\cyra }\IeC {\cyrz }\IeC {\cyra }\IeC {\cyrt }\IeC {\cyre }\IeC {\cyrl }\IeC {\cyrsftsn }\IeC {\cyrs }\IeC {\cyrt }\IeC {\cyrv }\IeC {\cyro } \IeC {\cyrn }\IeC {\cyre }\IeC {\cyrr }\IeC {\cyra }\IeC {\cyrz }\IeC {\cyrr }\IeC {\cyre }\IeC {\cyrsh }\IeC {\cyri }\IeC {\cyrm }\IeC {\cyro }\IeC {\cyrs }\IeC {\cyrt }\IeC {\cyri } \IeC {\cyrv }~\IeC {\cyrr }\IeC {\cyra }\IeC {\cyrd }\IeC {\cyri }\IeC {\cyrk }\IeC {\cyra }\IeC {\cyrl }\IeC {\cyra }\IeC {\cyrh } (4*)}{142}{}%{subsubsection.1.5.5.4}
\contentsline {subsubsection}{\numberline {5.5.5}\IeC {\CYRD }\IeC {\cyro }\IeC {\cyrk }\IeC {\cyra }\IeC {\cyrz }\IeC {\cyra }\IeC {\cyrt }\IeC {\cyre }\IeC {\cyrl }\IeC {\cyrsftsn }\IeC {\cyrs }\IeC {\cyrt }\IeC {\cyrv }\IeC {\cyro } \IeC {\cyrs }\IeC {\cyri }\IeC {\cyrl }\IeC {\cyrsftsn }\IeC {\cyrn }\IeC {\cyro }\IeC {\cyrishrt } \IeC {\cyrv }\IeC {\cyre }\IeC {\cyrshch }\IeC {\cyre }\IeC {\cyrs }\IeC {\cyrt }\IeC {\cyrv }\IeC {\cyre }\IeC {\cyrn }\IeC {\cyrn }\IeC {\cyro }\IeC {\cyrishrt } \IeC {\cyrt }\IeC {\cyre }\IeC {\cyro }\IeC {\cyrr }\IeC {\cyre }\IeC {\cyrm }\IeC {\cyrery } \IeC {\cyro }~\IeC {\cyrn }\IeC {\cyre }\-\IeC {\cyrr }\IeC {\cyra }\IeC {\cyrz }\IeC {\cyrr }\IeC {\cyre }\IeC {\cyrsh }\IeC {\cyri }\IeC {\cyrm }\IeC {\cyro }\IeC {\cyrs }\IeC {\cyrt }\IeC {\cyri } (4*)}{147}{}%{subsubsection.1.5.5.5}
\contentsline {section}{\numberline {6}\IeC {\CYRN }\IeC {\cyre }\IeC {\cyrr }\IeC {\cyra }\IeC {\cyrv }\IeC {\cyre }\IeC {\cyrn }\IeC {\cyrs }\IeC {\cyrt }\IeC {\cyrv }\IeC {\cyra }}{155}{section.1.6}
\contentsline {subsection}{\numberline {6.1}\IeC {\CYRV } \IeC {\cyrn }\IeC {\cyra }\IeC {\cyrp }\IeC {\cyrr }\IeC {\cyra }\IeC {\cyrv }\IeC {\cyrl }\IeC {\cyre }\IeC {\cyrn }\IeC {\cyri }\IeC {\cyri } \IeC {\cyrn }\IeC {\cyre }\IeC {\cyrr }\IeC {\cyra }\IeC {\cyrv }\IeC {\cyre }\IeC {\cyrn }\IeC {\cyrs }\IeC {\cyrt }\IeC {\cyrv }\IeC {\cyra } \IeC {\CYRISHRT }\IeC {\cyre }\IeC {\cyrn }\IeC {\cyrs }\IeC {\cyre }\IeC {\cyrn }\IeC {\cyra } (2)}{155}{}%{subsection.1.6.1}
\contentsline {subsection}{\numberline {6.2}\IeC {\CYRN }\IeC {\cyre }\IeC {\cyrk }\IeC {\cyro }\IeC {\cyrt }\IeC {\cyro }\IeC {\cyrr }\IeC {\cyrery }\IeC {\cyre } \IeC {\cyro }\IeC {\cyrs }\IeC {\cyrn }\IeC {\cyro }\IeC {\cyrv }\IeC {\cyrn }\IeC {\cyrery }\IeC {\cyre } \IeC {\cyrn }\IeC {\cyre }\IeC {\cyrr }\IeC {\cyra }\IeC {\cyrv }\IeC {\cyre }\IeC {\cyrn }\IeC {\cyrs }\IeC {\cyrt }\IeC {\cyrv }\IeC {\cyra } (2)}{160}{}%{subsection.1.6.2}
\contentsline {subsection}{\ull{\numberline {6.3}\IeC {\CYRP }\IeC {\cyrr }\IeC {\cyri }\IeC {\cyrm }\IeC {\cyre }\IeC {\cyrn }\IeC {\cyre }\IeC {\cyrn }\IeC {\cyri }\IeC {\cyrya } \IeC {\cyro }\IeC {\cyrs }\IeC {\cyrn }\IeC {\cyro }\IeC {\cyrv }\IeC {\cyrn }\IeC {\cyrery }\IeC {\cyrh } \IeC {\cyrn }\IeC {\cyre }\IeC {\cyrr }\IeC {\cyra }\IeC {\cyrv }\IeC {\cyre }\IeC {\cyrn }\IeC {\cyrs }\IeC {\cyrt }\IeC {\cyrv } (3*). \emph {\IeC {\CYRM }.~\IeC {\CYRA }.~\IeC {\CYRB }\IeC {\cyre }\IeC {\cyrr }\IeC {\cyrsh }\IeC {\cyrt }\IeC {\cyre }\IeC {\cyrishrt }\IeC {\cyrn }}}}{163}{subsection.1.6.3}
\contentsline {subsection}{\numberline {6.4}\IeC {\CYRG }\IeC {\cyre }\IeC {\cyro }\IeC {\cyrm }\IeC {\cyre }\IeC {\cyrt }\IeC {\cyrr }\IeC {\cyri }\IeC {\cyrch }\IeC {\cyre }\IeC {\cyrs }\IeC {\cyrk }\IeC {\cyra }\IeC {\cyrya } \IeC {\cyri }\IeC {\cyrn }\IeC {\cyrt }\IeC {\cyre }\IeC {\cyrr }\IeC {\cyrp }\IeC {\cyrr }\IeC {\cyre }\IeC {\cyrt }\IeC {\cyra }\IeC {\cyrc }\IeC {\cyri }\IeC {\cyrya } (3*)}{171}{}%{subsection.1.6.4}
\contentsline {section}{\numberline {7}\IeC {\CYRP }\IeC {\cyro }\IeC {\cyrs }\IeC {\cyrl }\IeC {\cyre }\IeC {\cyrd }\IeC {\cyro }\IeC {\cyrv }\IeC {\cyra }\IeC {\cyrt }\IeC {\cyre }\IeC {\cyrl }\IeC {\cyrsftsn }\IeC {\cyrn }\IeC {\cyro }\IeC {\cyrs }\IeC {\cyrt }\IeC {\cyri } \IeC {\cyri }~\IeC {\cyrp }\IeC {\cyrr }\IeC {\cyre }\IeC {\cyrd }\IeC {\cyre }\IeC {\cyrl }\IeC {\cyrery }}{175}{section.1.7}
\contentsline {subsection}{\ull{\numberline {7.1}\IeC {\CYRK }\IeC {\cyro }\IeC {\cyrn }\IeC {\cyre }\IeC {\cyrch }\IeC {\cyrn }\IeC {\cyrery }\IeC {\cyre } \IeC {\cyrs }\IeC {\cyru }\IeC {\cyrm }\IeC {\cyrm }\IeC {\cyrery } \IeC {\cyri }~\IeC {\cyrr }\IeC {\cyra }\IeC {\cyrz }\IeC {\cyrn }\IeC {\cyro }\IeC {\cyrs }\IeC {\cyrt }\IeC {\cyri } (3)}}{175}{subsection.1.7.1}
\contentsline {subsection}{\numberline {7.2}\IeC {\CYRL }\IeC {\cyri }\IeC {\cyrn }\IeC {\cyre }\IeC {\cyrishrt }\IeC {\cyrn }\IeC {\cyrery }\IeC {\cyre } \IeC {\cyrr }\IeC {\cyre }\IeC {\cyrk }\IeC {\cyru }\IeC {\cyrr }\IeC {\cyrr }\IeC {\cyre }\IeC {\cyrn }\IeC {\cyrt }\IeC {\cyrery } (3)}{179}{}%{subsection.1.7.2}
\contentsline {subsection}{\ull{\numberline {7.3}\IeC {\CYRK }\IeC {\cyro }\IeC {\cyrn }\IeC {\cyrk }\IeC {\cyrr }\IeC {\cyre }\IeC {\cyrt }\IeC {\cyrn }\IeC {\cyra }\IeC {\cyrya } \IeC {\cyrt }\IeC {\cyre }\IeC {\cyro }\IeC {\cyrr }\IeC {\cyri }\IeC {\cyrya } \IeC {\cyrp }\IeC {\cyrr }\IeC {\cyre }\IeC {\cyrd }\IeC {\cyre }\IeC {\cyrl }\IeC {\cyro }\IeC {\cyrv } (4*)}}{182}{subsection.1.7.3}
\contentsline {subsection}{\numberline {7.4}\IeC {\CYRK }\IeC {\cyra }\IeC {\cyrk } \IeC {\cyrk }\IeC {\cyro }\IeC {\cyrm }\IeC {\cyrp }\IeC {\cyrsftsn }\IeC {\cyryu }\IeC {\cyrt }\IeC {\cyre }\IeC {\cyrr } \IeC {\cyrv }\IeC {\cyrery }\IeC {\cyrch }\IeC {\cyri }\IeC {\cyrs }\IeC {\cyrl }\IeC {\cyrya }\IeC {\cyre }\IeC {\cyrt } \IeC {\cyrk }\IeC {\cyro }\IeC {\cyrr }\IeC {\cyre }\IeC {\cyrn }\IeC {\cyrsftsn }? (4*) \emph {\IeC {\CYRA }.~\IeC {\CYRS }.~\IeC {\CYRV }\IeC {\cyro }\IeC {\cyrr }\IeC {\cyro }\IeC {\cyrn }\IeC {\cyrc }\IeC {\cyro }\IeC {\cyrv }}, \emph {\IeC {\CYRA }.~\IeC {\CYRI }.~\IeC {\CYRS }\IeC {\cyrg }\IeC {\cyri }\IeC {\cyrb }\IeC {\cyrn }\IeC {\cyre }\IeC {\cyrv }}}{184}{}%{subsection.1.7.4}
\contentsline {subsection}{\numberline {7.5}\IeC {\CYRM }\IeC {\cyre }\IeC {\cyrt }\IeC {\cyro }\IeC {\cyrd }\IeC {\cyrery } \IeC {\cyrs }\IeC {\cyru }\IeC {\cyrm }\IeC {\cyrm }\IeC {\cyri }\IeC {\cyrr }\IeC {\cyro }\IeC {\cyrv }\IeC {\cyra }\IeC {\cyrn }\IeC {\cyri }\IeC {\cyrya } \IeC {\cyrr }\IeC {\cyrya }\IeC {\cyrd }\IeC {\cyro }\IeC {\cyrv } (4*)}{187}{}%{subsection.1.7.5}
\contentsline {subsection}{\numberline {7.6}\IeC {\CYRS }\IeC {\cyrh }\IeC {\cyro }\IeC {\cyrd }\IeC {\cyri }\IeC {\cyrm }\IeC {\cyro }\IeC {\cyrs }\IeC {\cyrt }\IeC {\cyrsftsn } \IeC {\cyrr }\IeC {\cyrya }\IeC {\cyrd }\IeC {\cyro }\IeC {\cyrv } (4*)}{192}{}%{subsection.1.7.6}
\contentsline {subsection}{\ull{\numberline {7.7}\IeC {\CYRP }\IeC {\cyrr }\IeC {\cyri }\IeC {\cyrm }\IeC {\cyre }\IeC {\cyrr }\IeC {\cyrery } \IeC {\cyrt }\IeC {\cyrr }\IeC {\cyra }\IeC {\cyrn }\IeC {\cyrs }\IeC {\cyrc }\IeC {\cyre }\IeC {\cyrn }\IeC {\cyrd }\IeC {\cyre }\IeC {\cyrn }\IeC {\cyrt }\IeC {\cyrn }\IeC {\cyrery }\IeC {\cyrh } \IeC {\cyrch }\IeC {\cyri }\IeC {\cyrs }\IeC {\cyre }\IeC {\cyrl } (3*)}}{195}{subsection.1.7.7}
\contentsline {section}{\numberline {8}\IeC {\CYRF }\IeC {\cyru }\IeC {\cyrn }\IeC {\cyrk }\IeC {\cyrc }\IeC {\cyri }\IeC {\cyri }}{199}{section.1.8}
\contentsline {subsection}{\ull{\numberline {8.1}\IeC {\CYRG }\IeC {\cyrr }\IeC {\cyra }\IeC {\cyrf }\IeC {\cyri }\IeC {\cyrk } \IeC {\cyrk }\IeC {\cyru }\IeC {\cyrb }\IeC {\cyri }\IeC {\cyrch }\IeC {\cyre }\IeC {\cyrs }\IeC {\cyrk }\IeC {\cyro }\IeC {\cyrg }\IeC {\cyro } \IeC {\cyrm }\IeC {\cyrn }\IeC {\cyro }\IeC {\cyrg }\IeC {\cyro }\IeC {\cyrch }\IeC {\cyrl }\IeC {\cyre }\IeC {\cyrn }\IeC {\cyra } (2)}}{199}{subsection.1.8.1}
\contentsline {subsection}{\ull{\numberline {8.2}\IeC {\CYREREV }\IeC {\cyrl }\IeC {\cyre }\IeC {\cyrm }\IeC {\cyre }\IeC {\cyrn }\IeC {\cyrt }\IeC {\cyrery } \IeC {\cyra }\IeC {\cyrn }\IeC {\cyra }\IeC {\cyrl }\IeC {\cyri }\IeC {\cyrz }\IeC {\cyra } \IeC {\cyrd }\IeC {\cyrl }\IeC {\cyrya } \IeC {\cyrm }\IeC {\cyrn }\IeC {\cyro }\IeC {\cyrg }\IeC {\cyro }\IeC {\cyrch }\IeC {\cyrl }\IeC {\cyre }\IeC {\cyrn }\IeC {\cyro }\IeC {\cyrv } (2)}}{203}{subsection.1.8.2}
\contentsline {subsection}{\numberline {8.3}\IeC {\CYRCH }\IeC {\cyri }\IeC {\cyrs }\IeC {\cyrl }\IeC {\cyro } \IeC {\cyrk }\IeC {\cyro }\IeC {\cyrr }\IeC {\cyrn }\IeC {\cyre }\IeC {\cyrishrt } \IeC {\cyrm }\IeC {\cyrn }\IeC {\cyro }\IeC {\cyrg }\IeC {\cyro }\IeC {\cyrch }\IeC {\cyrl }\IeC {\cyre }\IeC {\cyrn }\IeC {\cyra } (3*)}{205}{}%{subsection.1.8.3}
\contentsline {subsection}{\numberline {8.4}\IeC {\CYRO }\IeC {\cyrc }\IeC {\cyre }\IeC {\cyrn }\IeC {\cyrk }\IeC {\cyri } \IeC {\cyri }~\IeC {\cyrn }\IeC {\cyre }\IeC {\cyrr }\IeC {\cyra }\IeC {\cyrv }\IeC {\cyre }\IeC {\cyrn }\IeC {\cyrs }\IeC {\cyrt }\IeC {\cyrv }\IeC {\cyra } (4*). \emph {\IeC {\CYRV }.~\IeC {\CYRA }.~\IeC {\CYRS }\IeC {\cyre }\IeC {\cyrn }\IeC {\cyrd }\IeC {\cyre }\IeC {\cyrr }\IeC {\cyro }\IeC {\cyrv }}}{209}{}%{subsection.1.8.4}
\contentsline {subsection}{\numberline {8.5}\IeC {\CYRP }\IeC {\cyrr }\IeC {\cyri }\IeC {\cyrm }\IeC {\cyre }\IeC {\cyrn }\IeC {\cyre }\IeC {\cyrn }\IeC {\cyri }\IeC {\cyre } c\IeC {\cyru }\IeC {\cyrshch }\IeC {\cyre }\IeC {\cyrs }\IeC {\cyrt }\IeC {\cyrv }\IeC {\cyro }\IeC {\cyrv }\IeC {\cyra }\IeC {\cyrn }\IeC {\cyri }\IeC {\cyrya } \IeC {\cyrerev }\IeC {\cyrk }\IeC {\cyrs }\IeC {\cyrt }\IeC {\cyrr }\IeC {\cyre }\IeC {\cyrm }\IeC {\cyru }\IeC {\cyrm }\IeC {\cyra } (4*). \emph {\IeC {\CYRA }.~\IeC {\CYRYA }.~\IeC {\CYRK }\IeC {\cyra }\IeC {\cyrn }\IeC {\cyre }\IeC {\cyrl }\IeC {\cyrsftsn }-\IeC {\CYRB }\IeC {\cyre }\IeC {\cyrl }\IeC {\cyro }\IeC {\cyrv }}}{210}{}%{subsection.1.8.5}
\contentsline {subsection}{\ull{\numberline {8.6}\IeC {\CYRP }\IeC {\cyrr }\IeC {\cyri }\IeC {\cyrm }\IeC {\cyre }\IeC {\cyrn }\IeC {\cyre }\IeC {\cyrn }\IeC {\cyri }\IeC {\cyrya } \IeC {\cyrk }\IeC {\cyro }\IeC {\cyrm }\IeC {\cyrp }\IeC {\cyra }\IeC {\cyrk }\IeC {\cyrt }\IeC {\cyrn }\IeC {\cyro }\IeC {\cyrs }\IeC {\cyrt }\IeC {\cyri } (4*). \emph {\IeC {\CYRA }.~\IeC {\CYRYA }.~\IeC {\CYRK }\IeC {\cyra }\IeC {\cyrn }\IeC {\cyre }\IeC {\cyrl }\IeC {\cyrsftsn }-\IeC {\CYRB }\IeC {\cyre }\IeC {\cyrl }\IeC {\cyro }\IeC {\cyrv }}}}{213}{subsection.1.8.6}
\contentsline {chapter}{\numberline {2}\IeC {\CYRG }\IeC {\cyre }\IeC {\cyro }\IeC {\cyrm }\IeC {\cyre }\IeC {\cyrt }\IeC {\cyrr }\IeC {\cyri }\IeC {\cyrya }}{219}{chapter.2}
\contentsline {section}{\numberline {9}\IeC {\CYRT }\IeC {\cyrr }\IeC {\cyre }\IeC {\cyru }\IeC {\cyrg }\IeC {\cyro }\IeC {\cyrl }\IeC {\cyrsftsn }\IeC {\cyrn }\IeC {\cyri }\IeC {\cyrk }}{219}{section.2.9}
\contentsline {subsection}{\ull{\numberline {9.1}\IeC {\CYRP }\IeC {\cyrr }\IeC {\cyri }\IeC {\cyrn }\IeC {\cyrc }\IeC {\cyri }\IeC {\cyrp } \IeC {\CYRK }\IeC {\cyra }\IeC {\cyrr }\IeC {\cyrn }\IeC {\cyro } (1). \emph {\IeC {\CYRV }.~\IeC {\CYRYU }.~\IeC {\CYRP }\IeC {\cyrr }\IeC {\cyro }\IeC {\cyrt }\IeC {\cyra }\IeC {\cyrs }\IeC {\cyro }\IeC {\cyrv }}, \emph {\IeC {\CYRA }.~\IeC {\CYRA }.~\IeC {\CYRG }\IeC {\cyra }\IeC {\cyrv }\IeC {\cyrr }\IeC {\cyri }\IeC {\cyrl }\IeC {\cyryu }\IeC {\cyrk }}}}{220}{subsection.2.9.1}
\contentsline {subsection}{\numberline {9.2}\IeC {\CYRC }\IeC {\cyre }\IeC {\cyrn }\IeC {\cyrt }\IeC {\cyrr } \IeC {\cyrv }\IeC {\cyrp }\IeC {\cyri }\IeC {\cyrs }\IeC {\cyra }\IeC {\cyrn }\IeC {\cyrn }\IeC {\cyro }\IeC {\cyrishrt } \IeC {\cyro }\IeC {\cyrk }\IeC {\cyrr }\IeC {\cyru }\IeC {\cyrzh }\IeC {\cyrn }\IeC {\cyro }\IeC {\cyrs }\IeC {\cyrt }\IeC {\cyri } (2). \emph {\IeC {\CYRV }.~\IeC {\CYRYU }.~\IeC {\CYRP }\IeC {\cyrr }\IeC {\cyro }\IeC {\cyrt }\IeC {\cyra }\IeC {\cyrs }\IeC {\cyro }\IeC {\cyrv }}}{223}{}%{subsection.2.9.2}
\contentsline {subsection}{\numberline {9.3}\IeC {\CYRP }\IeC {\cyrr }\IeC {\cyrya }\IeC {\cyrm }\IeC {\cyra }\IeC {\cyrya } \IeC {\CYREREV }\IeC {\cyrishrt }\IeC {\cyrl }\IeC {\cyre }\IeC {\cyrr }\IeC {\cyra } (2). \emph {\IeC {\CYRV }.~\IeC {\CYRYU }.~\IeC {\CYRP }\IeC {\cyrr }\IeC {\cyro }\IeC {\cyrt }\IeC {\cyra }\IeC {\cyrs }\IeC {\cyro }\IeC {\cyrv }}}{226}{}%{subsection.2.9.3}
\contentsline {subsection}{\ull{\numberline {9.4}\IeC {\CYRF }\IeC {\cyro }\IeC {\cyrr }\IeC {\cyrm }\IeC {\cyru }\IeC {\cyrl }\IeC {\cyra } \IeC {\CYRK }\IeC {\cyra }\IeC {\cyrr }\IeC {\cyrn }\IeC {\cyro } (2$^*$). \emph {\IeC {\CYRA }.~\IeC {\CYRD }.~\IeC {\CYRB }\IeC {\cyrl }\IeC {\cyri }\IeC {\cyrn }\IeC {\cyrk }\IeC {\cyro }\IeC {\cyrv }}}}{227}{subsection.2.9.4}
\contentsline {subsection}{\numberline {9.5}\IeC {\CYRO }\IeC {\cyrr }\IeC {\cyrt }\IeC {\cyro }\IeC {\cyrc }\IeC {\cyre }\IeC {\cyrn }\IeC {\cyrt }\IeC {\cyrr }, \IeC {\cyro }\IeC {\cyrr }\IeC {\cyrt }\IeC {\cyro }\IeC {\cyrt }\IeC {\cyrr }\IeC {\cyre }\IeC {\cyru }\IeC {\cyrg }\IeC {\cyro }\IeC {\cyrl }\IeC {\cyrsftsn }\IeC {\cyrn }\IeC {\cyri }\IeC {\cyrk } \IeC {\cyri }~\IeC {\cyro }\IeC {\cyrk }\IeC {\cyrr }\IeC {\cyru }\IeC {\cyrzh }\IeC {\cyrn }\IeC {\cyro }\IeC {\cyrs }\IeC {\cyrt }\IeC {\cyrsftsn } \IeC {\cyrd }\IeC {\cyre }\IeC {\cyrv }\IeC {\cyrya }\IeC {\cyrt }\IeC {\cyri } \IeC {\cyrt }\IeC {\cyro }\IeC {\cyrch }\IeC {\cyre }\IeC {\cyrk }~(2). \emph {\IeC {\CYRV }.~\IeC {\CYRYU }.~\IeC {\CYRP }\IeC {\cyrr }\IeC {\cyro }\IeC {\cyrt }\IeC {\cyra }\IeC {\cyrs }\IeC {\cyro }\IeC {\cyrv }}}{232}{}%{subsection.2.9.5}
\contentsline {subsection}{\numberline {9.6}\IeC {\CYRN }\IeC {\cyre }\IeC {\cyrs }\IeC {\cyrk }\IeC {\cyro }\IeC {\cyrl }\IeC {\cyrsftsn }\IeC {\cyrk }\IeC {\cyro } \IeC {\cyrn }\IeC {\cyre }\IeC {\cyrr }\IeC {\cyra }\IeC {\cyrv }\IeC {\cyre }\IeC {\cyrn }\IeC {\cyrs }\IeC {\cyrt }\IeC {\cyrv }, \IeC {\cyrs }\IeC {\cyrv }\IeC {\cyrya }\IeC {\cyrz }\IeC {\cyra }\IeC {\cyrn }\IeC {\cyrn }\IeC {\cyrery }\IeC {\cyrh } \IeC {\cyrs }~\IeC {\cyrt }\IeC {\cyrr }\IeC {\cyre }\IeC {\cyru }\IeC {\cyrg }\IeC {\cyro }\IeC {\cyrl }\IeC {\cyrsftsn }\IeC {\cyrn }\IeC {\cyri }\IeC {\cyrk }\IeC {\cyro }\IeC {\cyrm } (3$^*$). \emph {\IeC {\CYRV }.~\IeC {\CYRYU }.~\IeC {\CYRP }\IeC {\cyrr }\IeC {\cyro }\IeC {\cyrt }\IeC {\cyra }\IeC {\cyrs }\IeC {\cyro }\IeC {\cyrv }}}{235}{}%{subsection.2.9.6}
\contentsline {subsection}{\numberline {9.7}\IeC {\CYRB }\IeC {\cyri }\IeC {\cyrs }\IeC {\cyrs }\IeC {\cyre }\IeC {\cyrk }\IeC {\cyrt }\IeC {\cyrr }\IeC {\cyri }\IeC {\cyrs }\IeC {\cyrery }, \IeC {\cyrv }\IeC {\cyrery }\IeC {\cyrs }\IeC {\cyro }\IeC {\cyrt }\IeC {\cyrery } \IeC {\cyri }~\IeC {\cyro }\IeC {\cyrp }\IeC {\cyri }\IeC {\cyrs }\IeC {\cyra }\IeC {\cyrn }\IeC {\cyrn }\IeC {\cyra }\IeC {\cyrya } \IeC {\cyro }\IeC {\cyrk }\IeC {\cyrr }\IeC {\cyru }\IeC {\cyrzh }\IeC {\cyrn }\IeC {\cyro }\IeC {\cyrs }\IeC {\cyrt }\IeC {\cyrsftsn } (2). \emph {\IeC {\CYRP }.~\IeC {\CYRA }.~\IeC {\CYRK }\IeC {\cyro }\IeC {\cyrzh }\IeC {\cyre }\IeC {\cyrv }\IeC {\cyrn }\IeC {\cyri }\IeC {\cyrk }\IeC {\cyro }\IeC {\cyrv }}}{237}{}%{subsection.2.9.7}
\contentsline {subsection}{\ull{\numberline {9.8}<<\IeC {\CYRP }\IeC {\cyro }\IeC {\cyrl }\IeC {\cyru }\IeC {\cyrv }\IeC {\cyrp }\IeC {\cyri }\IeC {\cyrs }\IeC {\cyra }\IeC {\cyrn }\IeC {\cyrn }\IeC {\cyra }\IeC {\cyrya }>> \IeC {\cyro }\IeC {\cyrk }\IeC {\cyrr }\IeC {\cyru }\IeC {\cyrzh }\IeC {\cyrn }\IeC {\cyro }\IeC {\cyrs }\IeC {\cyrt }\IeC {\cyrsftsn } (2$^*$). \emph {\IeC {\CYRP }.~\IeC {\CYRA }.~\IeC {\CYRK }\IeC {\cyro }\IeC {\cyrzh }\IeC {\cyre }\IeC {\cyrv }\IeC {\cyrn }\IeC {\cyri }\IeC {\cyrk }\IeC {\cyro }\IeC {\cyrv }}}}{242}{subsection.2.9.8}
\contentsline {subsection}{\ull{\numberline {9.9}\IeC {\CYRO }\IeC {\cyrb }\IeC {\cyro }\IeC {\cyrb }\IeC {\cyrshch }\IeC {\cyryo }\IeC {\cyrn }\IeC {\cyrn }\IeC {\cyra }\IeC {\cyrya } \IeC {\cyrt }\IeC {\cyre }\IeC {\cyro }\IeC {\cyrr }\IeC {\cyre }\IeC {\cyrm }\IeC {\cyra } \IeC {\CYRN }\IeC {\cyra }\IeC {\cyrp }\IeC {\cyro }\IeC {\cyrl }\IeC {\cyre }\IeC {\cyro }\IeC {\cyrn }\IeC {\cyra } (2$^*$). \emph {\IeC {\CYRP }.~\IeC {\CYRA }.~\IeC {\CYRK }\IeC {\cyro }\IeC {\cyrzh }\IeC {\cyre }\IeC {\cyrv }\IeC {\cyrn }\IeC {\cyri }\IeC {\cyrk }\IeC {\cyro }\IeC {\cyrv }}}}{249}{subsection.2.9.9}
\contentsline {subsection}{\numberline {9.10}\IeC {\CYRI }\IeC {\cyrz }\IeC {\cyro }\IeC {\cyrg }\IeC {\cyro }\IeC {\cyrn }\IeC {\cyra }\IeC {\cyrl }\IeC {\cyrsftsn }\IeC {\cyrn }\IeC {\cyro }\IeC {\cyre } \IeC {\cyrs }\IeC {\cyro }\IeC {\cyrp }\IeC {\cyrr }\IeC {\cyrya }\IeC {\cyrzh }\IeC {\cyre }\IeC {\cyrn }\IeC {\cyri }\IeC {\cyre } \IeC {\cyri }~\IeC {\cyrp }\IeC {\cyrr }\IeC {\cyrya }\IeC {\cyrm }\IeC {\cyra }\IeC {\cyrya } \IeC {\CYRS }\IeC {\cyri }\IeC {\cyrm }\IeC {\cyrs }\IeC {\cyro }\IeC {\cyrn }\IeC {\cyra } (3$^*$). \emph {\IeC {\CYRA }.~\IeC {\CYRV }.~\IeC {\CYRA }\IeC {\cyrk }\IeC {\cyro }\IeC {\cyrp }\IeC {\cyrya }\IeC {\cyrn }}}{257}{}%{subsection.2.9.10}
\contentsline {section}{\numberline {10}\IeC {\CYRO }\IeC {\cyrk }\IeC {\cyrr }\IeC {\cyru }\IeC {\cyrzh }\IeC {\cyrn }\IeC {\cyro }\IeC {\cyrs }\IeC {\cyrt }\IeC {\cyrsftsn }}{269}{section.2.10}
\contentsline {subsection}{\numberline {10.1}\IeC {\CYRP }\IeC {\cyrr }\IeC {\cyro }\IeC {\cyrs }\IeC {\cyrt }\IeC {\cyre }\IeC {\cyrishrt }\IeC {\cyrsh }\IeC {\cyri }\IeC {\cyre } \IeC {\cyrs }\IeC {\cyrv }\IeC {\cyro }\IeC {\cyrishrt }\IeC {\cyrs }\IeC {\cyrt }\IeC {\cyrv }\IeC {\cyra } \IeC {\cyro }\IeC {\cyrk }\IeC {\cyrr }\IeC {\cyru }\IeC {\cyrzh }\IeC {\cyrn }\IeC {\cyro }\IeC {\cyrs }\IeC {\cyrt }\IeC {\cyri } (1). \emph {\IeC {\CYRA }.~\IeC {\CYRD }.~\IeC {\CYRB }\IeC {\cyrl }\IeC {\cyri }\IeC {\cyrn }\IeC {\cyrk }\IeC {\cyro }\IeC {\cyrv }}}{269}{}%{subsection.2.10.1}
\contentsline {subsection}{\numberline {10.2}\IeC {\CYRV }\IeC {\cyrp }\IeC {\cyri }\IeC {\cyrs }\IeC {\cyra }\IeC {\cyrn }\IeC {\cyrn }\IeC {\cyrery }\IeC {\cyrishrt } \IeC {\cyru }\IeC {\cyrg }\IeC {\cyro }\IeC {\cyrl } (1). \emph {\IeC {\CYRA }.~\IeC {\CYRD }.~\IeC {\CYRB }\IeC {\cyrl }\IeC {\cyri }\IeC {\cyrn }\IeC {\cyrk }\IeC {\cyro }\IeC {\cyrv }}, \emph {\IeC {\CYRD }.~\IeC {\CYRA }.~\IeC {\CYRP }\IeC {\cyre }\IeC {\cyrr }\IeC {\cyrm }\IeC {\cyrya }\IeC {\cyrk }\IeC {\cyro }\IeC {\cyrv }}}{274}{}%{subsection.2.10.2}
\contentsline {subsection}{\numberline {10.3}\IeC {\CYRV }\IeC {\cyrp }\IeC {\cyri }\IeC {\cyrs }\IeC {\cyra }\IeC {\cyrn }\IeC {\cyrn }\IeC {\cyrery }\IeC {\cyre } \IeC {\cyri }~\IeC {\cyro }\IeC {\cyrp }\IeC {\cyri }\IeC {\cyrs }\IeC {\cyra }\IeC {\cyrn }\IeC {\cyrn }\IeC {\cyrery }\IeC {\cyre } \IeC {\cyro }\IeC {\cyrk }\IeC {\cyrr }\IeC {\cyru }\IeC {\cyrzh }\IeC {\cyrn }\IeC {\cyro }\IeC {\cyrs }\IeC {\cyrt }\IeC {\cyri } (2). \emph {\IeC {\CYRA }.~\IeC {\CYRA }.~\IeC {\CYRG }\IeC {\cyra }\IeC {\cyrv }\IeC {\cyrr }\IeC {\cyri }\IeC {\cyrl }\IeC {\cyryu }\IeC {\cyrk }}}{279}{}%{subsection.2.10.3}
\contentsline {subsection}{\numberline {10.4}\IeC {\CYRR }\IeC {\cyra }\IeC {\cyrd }\IeC {\cyri }\IeC {\cyrk }\IeC {\cyra }\IeC {\cyrl }\IeC {\cyrsftsn }\IeC {\cyrn }\IeC {\cyra }\IeC {\cyrya } \IeC {\cyro }\IeC {\cyrs }\IeC {\cyrsftsn } (2). \emph {\IeC {\CYRI }.~\IeC {\CYRN }.~\IeC {\CYRSH }\IeC {\cyrn }\IeC {\cyru }\IeC {\cyrr }\IeC {\cyrn }\IeC {\cyri }\IeC {\cyrk }\IeC {\cyro }\IeC {\cyrv }}, \emph {\IeC {\CYRA }.~\IeC {\CYRI }.~\IeC {\CYRZ }\IeC {\cyra }\IeC {\cyrs }\IeC {\cyro }\IeC {\cyrr }\IeC {\cyri }\IeC {\cyrn }}}{281}{}%{subsection.2.10.4}
\contentsline {subsection}{\numberline {10.5}\IeC {\CYRK }\IeC {\cyra }\IeC {\cyrs }\IeC {\cyra }\IeC {\cyrn }\IeC {\cyri }\IeC {\cyre } (2). \emph {\IeC {\CYRI }.~\IeC {\CYRN }.~\IeC {\CYRSH }\IeC {\cyrn }\IeC {\cyru }\IeC {\cyrr }\IeC {\cyrn }\IeC {\cyri }\IeC {\cyrk }\IeC {\cyro }\IeC {\cyrv }}, \emph {\IeC {\CYRA }.~\IeC {\CYRZ }\IeC {\cyra }\IeC {\cyrs }\IeC {\cyro }\IeC {\cyrr }\IeC {\cyri }\IeC {\cyrn }}}{282}{}%{subsection.2.10.5}
\contentsline {subsection}{\ull{\numberline {10.6}\IeC {\CYRT }\IeC {\cyre }\IeC {\cyro }\IeC {\cyrr }\IeC {\cyre }\IeC {\cyrm }\IeC {\cyrery } \IeC {\CYRP }\IeC {\cyrt }\IeC {\cyro }\IeC {\cyrl }\IeC {\cyre }\IeC {\cyrm }\IeC {\cyre }\IeC {\cyrya } \IeC {\cyri }~\IeC {\CYRK }\IeC {\cyre }\IeC {\cyrz }\IeC {\cyri } (3$^*$). \emph {\IeC {\CYRA }.~\IeC {\CYRD }.~\IeC {\CYRB }\IeC {\cyrl }\IeC {\cyri }\IeC {\cyrn }\IeC {\cyrk }\IeC {\cyro }\IeC {\cyrv }}, \emph {\IeC {\CYRA }.~\IeC {\CYRA }.~\IeC {\CYRZ }\IeC {\cyra }\IeC {\cyrs }\IeC {\cyrl }\IeC {\cyra }\IeC {\cyrv }\IeC {\cyrs }\IeC {\cyrk }\IeC {\cyri }\IeC {\cyrishrt }}}}{284}{subsection.2.10.6}
\contentsline {subsubsection}{\ull{\numberline {10.6.1}\IeC {\CYRT }\IeC {\cyre }\IeC {\cyro }\IeC {\cyrr }\IeC {\cyre }\IeC {\cyrm }\IeC {\cyra } \IeC {\CYRP }\IeC {\cyrt }\IeC {\cyro }\IeC {\cyrl }\IeC {\cyre }\IeC {\cyrm }\IeC {\cyre }\IeC {\cyrya }}}{284}{subsubsection.2.10.6.1}
\contentsline {subsubsection}{\ull{\numberline {10.6.2}\IeC {\CYRT }\IeC {\cyre }\IeC {\cyro }\IeC {\cyrr }\IeC {\cyre }\IeC {\cyrm }\IeC {\cyra } \IeC {\CYRK }\IeC {\cyre }\IeC {\cyrz }\IeC {\cyri }}}{285}{subsubsection.2.10.6.2}
\contentsline {section}{\numberline {11}\IeC {\CYRG }\IeC {\cyre }\IeC {\cyro }\IeC {\cyrm }\IeC {\cyre }\IeC {\cyrt }\IeC {\cyrr }\IeC {\cyri }\IeC {\cyrch }\IeC {\cyre }\IeC {\cyrs }\IeC {\cyrk }\IeC {\cyri }\IeC {\cyre } \IeC {\cyrp }\IeC {\cyrr }\IeC {\cyre }\IeC {\cyro }\IeC {\cyrb }\IeC {\cyrr }\IeC {\cyra }\IeC {\cyrz }\IeC {\cyro }\IeC {\cyrv }\IeC {\cyra }\IeC {\cyrn }\IeC {\cyri }\IeC {\cyrya }}{292}{section.2.11}
\contentsline {subsection}{\ull{\numberline {11.1}\IeC {\CYRP }\IeC {\cyrr }\IeC {\cyri }\IeC {\cyrm }\IeC {\cyre }\IeC {\cyrn }\IeC {\cyre }\IeC {\cyrn }\IeC {\cyri }\IeC {\cyrya } \IeC {\cyrd }\IeC {\cyrv }\IeC {\cyri }\IeC {\cyrzh }\IeC {\cyre }\IeC {\cyrn }\IeC {\cyri }\IeC {\cyrishrt }. (1) \emph {\IeC {\CYRA }.~\IeC {\CYRD }.~\IeC {\CYRB }\IeC {\cyrl }\IeC {\cyri }\IeC {\cyrn }\IeC {\cyrk }\IeC {\cyro }\IeC {\cyrv }}}}{292}{subsection.2.11.1}
\contentsline {subsection}{\numberline {11.2}\IeC {\CYRK }\IeC {\cyrl }\IeC {\cyra }\IeC {\cyrs }\IeC {\cyrs }\IeC {\cyri }\IeC {\cyrf }\IeC {\cyri }\IeC {\cyrk }\IeC {\cyra }\IeC {\cyrc }\IeC {\cyri }\IeC {\cyrya } \IeC {\cyrd }\IeC {\cyrv }\IeC {\cyri }\IeC {\cyrzh }\IeC {\cyre }\IeC {\cyrn }\IeC {\cyri }\IeC {\cyrishrt } \IeC {\cyrp }\IeC {\cyrl }\IeC {\cyro }\IeC {\cyrs }\IeC {\cyrk }\IeC {\cyro }\IeC {\cyrs }\IeC {\cyrt }\IeC {\cyri } (2). \emph {\IeC {\CYRA }.~\IeC {\CYRB }.~\IeC {\CYRS }\IeC {\cyrk }\IeC {\cyro }\IeC {\cyrp }\IeC {\cyre }\IeC {\cyrn }\IeC {\cyrk }\IeC {\cyro }\IeC {\cyrv }}}{300}{}%{subsection.2.11.2}
\contentsline {subsection}{\numberline {11.3}\IeC {\CYRK }\IeC {\cyrl }\IeC {\cyra }\IeC {\cyrs }\IeC {\cyrs }\IeC {\cyri }\IeC {\cyrf }\IeC {\cyri }\IeC {\cyrk }\IeC {\cyra }\IeC {\cyrc }\IeC {\cyri }\IeC {\cyrya } \IeC {\cyrd }\IeC {\cyrv }\IeC {\cyri }\IeC {\cyrzh }\IeC {\cyre }\IeC {\cyrn }\IeC {\cyri }\IeC {\cyrishrt } \IeC {\cyrp }\IeC {\cyrr }\IeC {\cyro }\IeC {\cyrs }\IeC {\cyrt }\IeC {\cyrr }\IeC {\cyra }\IeC {\cyrn }\IeC {\cyrs }\IeC {\cyrt }\IeC {\cyrv }\IeC {\cyra } (3*). \emph {\IeC {\CYRA }.~\IeC {\CYRB }.~\IeC {\CYRS }\IeC {\cyrk }\IeC {\cyro }\IeC {\cyrp }\IeC {\cyre }\IeC {\cyrn }\IeC {\cyrk }\IeC {\cyro }\IeC {\cyrv }}}{302}{}%{subsection.2.11.3}
\contentsline {subsection}{\numberline {11.4}\IeC {\CYRP }\IeC {\cyrr }\IeC {\cyri }\IeC {\cyrm }\IeC {\cyre }\IeC {\cyrn }\IeC {\cyre }\IeC {\cyrn }\IeC {\cyri }\IeC {\cyre } \IeC {\cyrp }\IeC {\cyro }\IeC {\cyrd }\IeC {\cyro }\IeC {\cyrb }\IeC {\cyri }\IeC {\cyrya } \IeC {\cyri }~\IeC {\cyrg }\IeC {\cyro }\IeC {\cyrm }\IeC {\cyro }\IeC {\cyrt }\IeC {\cyre }\IeC {\cyrt }\IeC {\cyri }\IeC {\cyri } (1). \emph {\IeC {\CYRA }.~\IeC {\CYRD }.~\IeC {\CYRB }\IeC {\cyrl }\IeC {\cyri }\IeC {\cyrn }\IeC {\cyrk }\IeC {\cyro }\IeC {\cyrv }}}{304}{}%{subsection.2.11.4}
\contentsline {subsection}{\ull{\numberline {11.5}\IeC {\CYRP }\IeC {\cyro }\IeC {\cyrv }\IeC {\cyro }\IeC {\cyrr }\IeC {\cyro }\IeC {\cyrt }\IeC {\cyrn }\IeC {\cyra }\IeC {\cyrya } \IeC {\cyrg }\IeC {\cyro }\IeC {\cyrm }\IeC {\cyro }\IeC {\cyrt }\IeC {\cyre }\IeC {\cyrt }\IeC {\cyri }\IeC {\cyrya } (2). \emph {\IeC {\CYRP }.~\IeC {\CYRA }.~\IeC {\CYRK }\IeC {\cyro }\IeC {\cyrzh }\IeC {\cyre }\IeC {\cyrv }\IeC {\cyrn }\IeC {\cyri }\IeC {\cyrk }\IeC {\cyro }\IeC {\cyrv }}}}{312}{subsection.2.11.5}
\contentsline {subsubsection}{\ull{\numberline {11.5.1}\IeC {\CYRV }\IeC {\cyrv }\IeC {\cyro }\IeC {\cyrd }\IeC {\cyrn }\IeC {\cyrery }\IeC {\cyre } \IeC {\cyrz }\IeC {\cyra }\IeC {\cyrd }\IeC {\cyra }\IeC {\cyrch }\IeC {\cyri }: \IeC {\cyrn }\IeC {\cyre }\IeC {\cyrm }\IeC {\cyrn }\IeC {\cyro }\IeC {\cyrg }\IeC {\cyro } \IeC {\cyro }~\IeC {\cyrv }\IeC {\cyre }\IeC {\cyrl }\IeC {\cyro }\IeC {\cyrs }\IeC {\cyri }\IeC {\cyrp }\IeC {\cyre }\IeC {\cyrd }\IeC {\cyri }\IeC {\cyrs }\IeC {\cyrt }\IeC {\cyra }\IeC {\cyrh }}}{312}{subsubsection.2.11.5.1}
\contentsline {subsubsection}{\ull{\numberline {11.5.2}\IeC {\CYRO }\IeC {\cyrs }\IeC {\cyrn }\IeC {\cyro }\IeC {\cyrv }\IeC {\cyrn }\IeC {\cyrery }\IeC {\cyre } \IeC {\cyrz }\IeC {\cyra }\IeC {\cyrd }\IeC {\cyra }\IeC {\cyrch }\IeC {\cyri }}}{313}{subsubsection.2.11.5.2}
\contentsline {subsubsection}{\ull{\numberline {11.5.3}\IeC {\CYRD }\IeC {\cyro }\IeC {\cyrp }\IeC {\cyro }\IeC {\cyrl }\IeC {\cyrn }\IeC {\cyri }\IeC {\cyrt }\IeC {\cyre }\IeC {\cyrl }\IeC {\cyrsftsn }\IeC {\cyrn }\IeC {\cyrery }\IeC {\cyre } \IeC {\cyrz }\IeC {\cyra }\IeC {\cyrd }\IeC {\cyra }\IeC {\cyrch }\IeC {\cyri }}}{314}{subsubsection.2.11.5.3}
\contentsline {subsection}{\numberline {11.6}\IeC {\CYRP }\IeC {\cyro }\IeC {\cyrd }\IeC {\cyro }\IeC {\cyrb }\IeC {\cyri }\IeC {\cyre } (1). \emph {\IeC {\CYRA }.~\IeC {\CYRB }.~\IeC {\CYRS }\IeC {\cyrk }\IeC {\cyro }\IeC {\cyrp }\IeC {\cyre }\IeC {\cyrn }\IeC {\cyrk }\IeC {\cyro }\IeC {\cyrv }}}{319}{}%{subsection.2.11.6}
\contentsline {subsection}{\numberline {11.7}\IeC {\CYRS }\IeC {\cyrzh }\IeC {\cyra }\IeC {\cyrt }\IeC {\cyri }\IeC {\cyre } \IeC {\cyrk }~\IeC {\cyrp }\IeC {\cyrr }\IeC {\cyrya }\IeC {\cyrm }\IeC {\cyro }\IeC {\cyrishrt } (2). \emph {\IeC {\CYRA }.~\IeC {\CYRYA }.~\IeC {\CYRK }\IeC {\cyra }\IeC {\cyrn }\IeC {\cyre }\IeC {\cyrl }\IeC {\cyrsftsn }-\IeC {\CYRB }\IeC {\cyre }\IeC {\cyrl }\IeC {\cyro }\IeC {\cyrv }}}{320}{}%{subsection.2.11.7}
\contentsline {subsection}{\numberline {11.8}\IeC {\CYRP }\IeC {\cyra }\IeC {\cyrr }\IeC {\cyra }\IeC {\cyrl }\IeC {\cyrl }\IeC {\cyre }\IeC {\cyrl }\IeC {\cyrsftsn }\IeC {\cyrn }\IeC {\cyra }\IeC {\cyrya } \IeC {\cyrp }\IeC {\cyrr }\IeC {\cyro }\IeC {\cyre }\IeC {\cyrk }\IeC {\cyrc }\IeC {\cyri }\IeC {\cyrya } \IeC {\cyri }~\IeC {\cyra }\IeC {\cyrf }\IeC {\cyrf }\IeC {\cyri }\IeC {\cyrn }\IeC {\cyrn }\IeC {\cyrery }\IeC {\cyre } \IeC {\cyrp }\IeC {\cyrr }\IeC {\cyre }\IeC {\cyro }\IeC {\cyrb }\IeC {\cyrr }\IeC {\cyra }\IeC {\cyrz }\IeC {\cyro }\IeC {\cyrv }\IeC {\cyra }\IeC {\cyrn }\IeC {\cyri }\IeC {\cyrya } (2). \emph {\IeC {\CYRA }.~\IeC {\CYRB }.~\IeC {\CYRS }\IeC {\cyrk }\IeC {\cyro }\IeC {\cyrp }\IeC {\cyre }\IeC {\cyrn }\IeC {\cyrk }\IeC {\cyro }\IeC {\cyrv }}}{321}{}%{subsection.2.11.8}
\contentsline {subsection}{\numberline {11.9}\IeC {\CYRC }\IeC {\cyre }\IeC {\cyrn }\IeC {\cyrt }\IeC {\cyrr }\IeC {\cyra }\IeC {\cyrl }\IeC {\cyrsftsn }\IeC {\cyrn }\IeC {\cyra }\IeC {\cyrya } \IeC {\cyrp }\IeC {\cyrr }\IeC {\cyro }\IeC {\cyre }\IeC {\cyrk }\IeC {\cyrc }\IeC {\cyri }\IeC {\cyrya } \IeC {\cyri }~\IeC {\cyrp }\IeC {\cyrr }\IeC {\cyro }\IeC {\cyre }\IeC {\cyrk }\IeC {\cyrt }\IeC {\cyri }\IeC {\cyrv }\IeC {\cyrn }\IeC {\cyrery }\IeC {\cyre } \IeC {\cyrp }\IeC {\cyrr }\IeC {\cyre }\IeC {\cyro }\IeC {\cyrb }\IeC {\cyrr }\IeC {\cyra }\IeC {\cyrz }\IeC {\cyro }\IeC {\cyrv }\IeC {\cyra }\IeC {\cyrn }\IeC {\cyri }\IeC {\cyrya } (3). \emph {\IeC {\CYRA }.~\IeC {\CYRB }.~\IeC {\CYRS }\IeC {\cyrk }\IeC {\cyro }\IeC {\cyrp }\IeC {\cyre }\IeC {\cyrn }\IeC {\cyrk }\IeC {\cyro }\IeC {\cyrv }}}{325}{}%{subsection.2.11.9}
\contentsline {subsection}{\numberline {11.10}\IeC {\CYRI }\IeC {\cyrn }\IeC {\cyrv }\IeC {\cyre }\IeC {\cyrr }\IeC {\cyrs }\IeC {\cyri }\IeC {\cyrya } (2). \emph {\IeC {\CYRA }.~\IeC {\CYRB }.~\IeC {\CYRS }\IeC {\cyrk }\IeC {\cyro }\IeC {\cyrp }\IeC {\cyre }\IeC {\cyrn }\IeC {\cyrk }\IeC {\cyro }\IeC {\cyrv }}}{328}{}%{subsection.2.11.10}
\contentsline {section}{\numberline {12}\IeC {\CYRA }\IeC {\cyrf }\IeC {\cyrf }\IeC {\cyri }\IeC {\cyrn }\IeC {\cyrn }\IeC {\cyra }\IeC {\cyrya } \IeC {\cyri }~\IeC {\cyrp }\IeC {\cyrr }\IeC {\cyro }\IeC {\cyre }\IeC {\cyrk }\IeC {\cyrt }\IeC {\cyri }\IeC {\cyrv }\IeC {\cyrn }\IeC {\cyra }\IeC {\cyrya } \IeC {\cyrg }\IeC {\cyre }\IeC {\cyro }\IeC {\cyrm }\IeC {\cyre }\IeC {\cyrt }\IeC {\cyrr }\IeC {\cyri }\IeC {\cyrya }}{334}{section.2.12}
\contentsline {subsection}{\numberline {12.1}\IeC {\CYRB }\IeC {\cyru }\IeC {\cyrr }\IeC {\cyrya } \IeC {\cyrn }\IeC {\cyra } \IeC {\CYRM }\IeC {\cyra }\IeC {\cyrs }\IeC {\cyrs }\IeC {\cyro }\IeC {\cyrv }\IeC {\cyro }\IeC {\cyrm } \IeC {\cyrp }\IeC {\cyro }\IeC {\cyrl }\IeC {\cyre } (2). \emph {\IeC {\CYRA }.~\IeC {\CYRA }.~\IeC {\CYRG }\IeC {\cyra }\IeC {\cyrv }\IeC {\cyrr }\IeC {\cyri }\IeC {\cyrl }\IeC {\cyryu }\IeC {\cyrk }}}{335}{}%{subsection.2.12.1}
\contentsline {subsection}{\numberline {12.2}\IeC {\CYRD }\IeC {\cyrv }\IeC {\cyro }\IeC {\cyrishrt }\IeC {\cyrn }\IeC {\cyrery }\IeC {\cyre } \IeC {\cyro }\IeC {\cyrt }\IeC {\cyrn }\IeC {\cyro }\IeC {\cyrsh }\IeC {\cyre }\IeC {\cyrn }\IeC {\cyri }\IeC {\cyrya } (2). \emph {\IeC {\CYRA }.~\IeC {\CYRA }.~\IeC {\CYRG }\IeC {\cyra }\IeC {\cyrv }\IeC {\cyrr }\IeC {\cyri }\IeC {\cyrl }\IeC {\cyryu }\IeC {\cyrk }}}{338}{}%{subsection.2.12.2}
\contentsline {subsection}{\ull{\numberline {12.3}\IeC {\CYRP }\IeC {\cyro }\IeC {\cyrl }\IeC {\cyrya }\IeC {\cyrr }\IeC {\cyrn }\IeC {\cyro }\IeC {\cyre } \IeC {\cyrs }\IeC {\cyro }\IeC {\cyro }\IeC {\cyrt }\IeC {\cyrv }\IeC {\cyre }\IeC {\cyrt }\IeC {\cyrs }\IeC {\cyrt }\IeC {\cyrv }\IeC {\cyri }\IeC {\cyre } (2). \emph {\IeC {\CYRA }.~\IeC {\CYRA }.~\IeC {\CYRG }\IeC {\cyra }\IeC {\cyrv }\IeC {\cyrr }\IeC {\cyri }\IeC {\cyrl }\IeC {\cyryu }\IeC {\cyrk }}, \emph {\IeC {\CYRP }.~\IeC {\CYRA }.~\IeC {\CYRK }\IeC {\cyro }\IeC {\cyrzh }\IeC {\cyre }\IeC {\cyrv }\IeC {\cyrn }\IeC {\cyri }\IeC {\cyrk }\IeC {\cyro }\IeC {\cyrv }}}}{343}{}%{subsection.2.12.3}
\contentsline {section}{\ull{\numberline {13}\IeC {\CYRK }\IeC {\cyro }\IeC {\cyrm }\IeC {\cyrp }\IeC {\cyrl }\IeC {\cyre }\IeC {\cyrk }\IeC {\cyrs }\IeC {\cyrn }\IeC {\cyrery }\IeC {\cyre } \IeC {\cyrch }\IeC {\cyri }\IeC {\cyrs }\IeC {\cyrl }\IeC {\cyra } \IeC {\cyri }~\IeC {\cyrg }\IeC {\cyre }\IeC {\cyro }\IeC {\cyrm }\IeC {\cyre }\IeC {\cyrt }\IeC {\cyrr }\IeC {\cyri }\IeC {\cyrya } (3). \emph {\IeC {\CYRA }.~\IeC {\CYRA }.~\IeC {\CYRZ }\IeC {\cyra }\IeC {\cyrs }\IeC {\cyrl }\IeC {\cyra }\IeC {\cyrv }\IeC {\cyrs }\IeC {\cyrk }\IeC {\cyri }\IeC {\cyrishrt }}}}{351}{section.2.13}
\contentsline {subsection}{\ull{\numberline {13.1}\IeC {\CYRK }\IeC {\cyro }\IeC {\cyrm }\IeC {\cyrp }\IeC {\cyrl }\IeC {\cyre }\IeC {\cyrk }\IeC {\cyrs }\IeC {\cyrn }\IeC {\cyrery }\IeC {\cyre } \IeC {\cyrch }\IeC {\cyri }\IeC {\cyrs }\IeC {\cyrl }\IeC {\cyra } \IeC {\cyri }~\IeC {\cyrerev }\IeC {\cyrl }\IeC {\cyre }\IeC {\cyrm }\IeC {\cyre }\IeC {\cyrn }\IeC {\cyrt }\IeC {\cyra }\IeC {\cyrr }\IeC {\cyrn }\IeC {\cyra }\IeC {\cyrya } \IeC {\cyrg }\IeC {\cyre }\IeC {\cyro }\IeC {\cyrm }\IeC {\cyre }\IeC {\cyrt }\IeC {\cyrr }\IeC {\cyri }\IeC {\cyrya }.}}{352}{}%{subsection.2.13.1}
\contentsline {subsection}{\ull{\numberline {13.2}\IeC {\CYRK }\IeC {\cyro }\IeC {\cyrm }\IeC {\cyrp }\IeC {\cyrl }\IeC {\cyre }\IeC {\cyrk }\IeC {\cyrs }\IeC {\cyrn }\IeC {\cyrery }\IeC {\cyre } \IeC {\cyrch }\IeC {\cyri }\IeC {\cyrs }\IeC {\cyrl }\IeC {\cyra } \IeC {\cyri }~\IeC {\cyrk }\IeC {\cyrr }\IeC {\cyru }\IeC {\cyrg }\IeC {\cyro }\IeC {\cyrv }\IeC {\cyrery }\IeC {\cyre } \IeC {\cyrp }\IeC {\cyrr }\IeC {\cyre }\IeC {\cyro }\IeC {\cyrb }\IeC {\cyrr }\IeC {\cyra }\IeC {\cyrz }\IeC {\cyro }\IeC {\cyrv }\IeC {\cyra }\IeC {\cyrn }\IeC {\cyri }\IeC {\cyrya }.}}{355}{subsection.2.13.2}
\contentsline {section}{\numberline {14}\IeC {\CYRP }\IeC {\cyro }\IeC {\cyrs }\IeC {\cyrt }\IeC {\cyrr }\IeC {\cyro }\IeC {\cyre }\IeC {\cyrn }\IeC {\cyri }\IeC {\cyrya } \IeC {\cyri }~\IeC {\cyrg }\IeC {\cyre }\IeC {\cyro }\IeC {\cyrm }\IeC {\cyre }\IeC {\cyrt }\IeC {\cyrr }\IeC {\cyri }\IeC {\cyrch }\IeC {\cyre }\IeC {\cyrs }\IeC {\cyrk }\IeC {\cyri }\IeC {\cyre }\\ \IeC {\cyrm }\IeC {\cyre }\IeC {\cyrs }\IeC {\cyrt }\IeC {\cyra } \IeC {\cyrt }\IeC {\cyro }\IeC {\cyrch }\IeC {\cyre }\IeC {\cyrk }}{359}{section.2.14}
\contentsline {subsection}{\numberline {14.1}\IeC {\CYRG }\IeC {\cyre }\IeC {\cyro }\IeC {\cyrm }\IeC {\cyre }\IeC {\cyrt }\IeC {\cyrr }\IeC {\cyri }\IeC {\cyrch }\IeC {\cyre }\IeC {\cyrs }\IeC {\cyrk }\IeC {\cyri }\IeC {\cyre } \IeC {\cyrm }\IeC {\cyre }\IeC {\cyrs }\IeC {\cyrt }\IeC {\cyra } \IeC {\cyrt }\IeC {\cyro }\IeC {\cyrch }\IeC {\cyre }\IeC {\cyrk } (1). \emph {\IeC {\CYRA }.~\IeC {\CYRD }.~\IeC {\CYRB }\IeC {\cyrl }\IeC {\cyri }\IeC {\cyrn }\IeC {\cyrk }\IeC {\cyro }\IeC {\cyrv }}}{359}{}%{subsection.2.14.1}
\contentsline {subsection}{\numberline {14.2}\IeC {\CYRZ }\IeC {\cyra }\IeC {\cyrd }\IeC {\cyra }\IeC {\cyrch }\IeC {\cyri } \IeC {\cyrn }\IeC {\cyra } \IeC {\cyrp }\IeC {\cyro }\IeC {\cyrs }\IeC {\cyrt }\IeC {\cyrr }\IeC {\cyro }\IeC {\cyre }\IeC {\cyrn }\IeC {\cyri }\IeC {\cyre } \IeC {\cyri }~\IeC {\CYRG }\IeC {\CYRM }\IeC {\CYRT }, \IeC {\cyrs }\IeC {\cyrv }\IeC {\cyrya }\IeC {\cyrz }\IeC {\cyra }\IeC {\cyrn }\IeC {\cyrn }\IeC {\cyrery }\IeC {\cyre } \IeC {\cyrs }~\IeC {\cyrp }\IeC {\cyrl }\IeC {\cyro }\IeC {\cyrshch }\IeC {\cyra }\IeC {\cyrd }\IeC {\cyrya }\IeC {\cyrm }\IeC {\cyri }~(1). \emph {\IeC {\CYRA }.~\IeC {\CYRD }.~\IeC {\CYRB }\IeC {\cyrl }\IeC {\cyri }\IeC {\cyrn }\IeC {\cyrk }\IeC {\cyro }\IeC {\cyrv }}}{367}{}%{subsection.2.14.2}
\contentsline {subsection}{\ull{\numberline {14.3}\IeC {\CYRP }\IeC {\cyro }\IeC {\cyrs }\IeC {\cyrt }\IeC {\cyrr }\IeC {\cyro }\IeC {\cyre }\IeC {\cyrn }\IeC {\cyri }\IeC {\cyrya }. \IeC {\CYRYA }\IeC {\cyrshch }\IeC {\cyri }\IeC {\cyrk } \IeC {\cyri }\IeC {\cyrn }\IeC {\cyrs }\IeC {\cyrt }\IeC {\cyrr }\IeC {\cyru }\IeC {\cyrm }\IeC {\cyre }\IeC {\cyrn }\IeC {\cyrt }\IeC {\cyro }\IeC {\cyrv } (2). \emph {\IeC {\CYRA }.~\IeC {\CYRA }.~\IeC {\CYRG }\IeC {\cyra }\IeC {\cyrv }\IeC {\cyrr }\IeC {\cyri }\IeC {\cyrl }\IeC {\cyryu }\IeC {\cyrk }}}}{373}{subsection.2.14.3}
\contentsline {subsection}{\numberline {14.4}\IeC {\CYRD }\IeC {\cyro }\IeC {\cyrp }\IeC {\cyro }\IeC {\cyrl }\IeC {\cyrn }\IeC {\cyri }\IeC {\cyrt }\IeC {\cyre }\IeC {\cyrl }\IeC {\cyrsftsn }\IeC {\cyrn }\IeC {\cyrery }\IeC {\cyre } \IeC {\cyrp }\IeC {\cyro }\IeC {\cyrs }\IeC {\cyrt }\IeC {\cyrr }\IeC {\cyro }\IeC {\cyre }\IeC {\cyrn }\IeC {\cyri }\IeC {\cyrya } (2$^*$). \emph {\IeC {\CYRI }.~\IeC {\CYRN }.~\IeC {\CYRSH }\IeC {\cyrn }\IeC {\cyru }\IeC {\cyrr }\IeC {\cyrn }\IeC {\cyri }\IeC {\cyrk }\IeC {\cyro }\IeC {\cyrv }}}{377}{}%{subsection.2.14.4}
\contentsline {section}{\numberline {15}\IeC {\CYRS }\IeC {\cyrt }\IeC {\cyre }\IeC {\cyrr }\IeC {\cyre }\IeC {\cyro }\IeC {\cyrm }\IeC {\cyre }\IeC {\cyrt }\IeC {\cyrr }\IeC {\cyri }\IeC {\cyrya }}{385}{section.2.15}
\contentsline {subsection}{\ull{\numberline {15.1}\IeC {\CYRR }\IeC {\cyri }\IeC {\cyrs }\IeC {\cyro }\IeC {\cyrv }\IeC {\cyra }\IeC {\cyrn }\IeC {\cyri }\IeC {\cyre } (2). \emph {\IeC {\CYRA }.~\IeC {\CYRB }.~\IeC {\CYRS }\IeC {\cyrk }\IeC {\cyro }\IeC {\cyrp }\IeC {\cyre }\IeC {\cyrn }\IeC {\cyrk }\IeC {\cyro }\IeC {\cyrv }}}}{385}{subsection.2.15.1}
\contentsline {subsection}{\numberline {15.2}\IeC {\CYRP }\IeC {\cyrr }\IeC {\cyra }\IeC {\cyrv }\IeC {\cyri }\IeC {\cyrl }\IeC {\cyrsftsn }\IeC {\cyrn }\IeC {\cyrery }\IeC {\cyre } \IeC {\cyrm }\IeC {\cyrn }\IeC {\cyro }\IeC {\cyrg }\IeC {\cyro }\IeC {\cyrg }\IeC {\cyrr }\IeC {\cyra }\IeC {\cyrn }\IeC {\cyrn }\IeC {\cyri }\IeC {\cyrk }\IeC {\cyri } (3)}{387}{}%{subsection.2.15.2}
\contentsline {subsubsection}{\numberline {15.2.1}\IeC {\CYRV }\IeC {\cyrp }\IeC {\cyri }\IeC {\cyrs }\IeC {\cyra }\IeC {\cyrn }\IeC {\cyrn }\IeC {\cyrery }\IeC {\cyre } \IeC {\cyri }~\IeC {\cyro }\IeC {\cyrp }\IeC {\cyri }\IeC {\cyrs }\IeC {\cyra }\IeC {\cyrn }\IeC {\cyrn }\IeC {\cyrery }\IeC {\cyre }. \emph {\IeC {\CYRA }.~\IeC {\CYRYA }.~\IeC {\CYRK }\IeC {\cyra }\IeC {\cyrn }\IeC {\cyre }\IeC {\cyrl }\IeC {\cyrsftsn }-\IeC {\CYRB }\IeC {\cyre }\IeC {\cyrl }\IeC {\cyro }\IeC {\cyrv }}}{387}{}%{subsubsection.2.15.2.1}
\contentsline {subsubsection}{\numberline {15.2.2}\IeC {\CYRS }\IeC {\cyra }\IeC {\cyrm }\IeC {\cyro }\IeC {\cyrs }\IeC {\cyro }\IeC {\cyrv }\IeC {\cyrm }\IeC {\cyre }\IeC {\cyrshch }\IeC {\cyre }\IeC {\cyrn }\IeC {\cyri }\IeC {\cyrya }. \emph {\IeC {\CYRA }.~\IeC {\CYRB }.~\IeC {\CYRS }\IeC {\cyrk }\IeC {\cyro }\IeC {\cyrp }\IeC {\cyre }\IeC {\cyrn }\IeC {\cyrk }\IeC {\cyro }\IeC {\cyrv }}}{390}{}%{subsubsection.2.15.2.2}
\contentsline {subsection}{\numberline {15.3}\IeC {\CYRM }\IeC {\cyrn }\IeC {\cyro }\IeC {\cyrg }\IeC {\cyro }\IeC {\cyrm }\IeC {\cyre }\IeC {\cyrr }\IeC {\cyrsftsn }\IeC {\cyre } (4$^*$). \emph {\IeC {\CYRA }.~\IeC {\CYRYA }.~\IeC {\CYRK }\IeC {\cyra }\IeC {\cyrn }\IeC {\cyre }\IeC {\cyrl }\IeC {\cyrsftsn }-\IeC {\CYRB }\IeC {\cyre }\IeC {\cyrl }\IeC {\cyro }\IeC {\cyrv }}}{392}{}%{subsection.2.15.3}
\contentsline {subsubsection}{\ull{\numberline {15.3.1}\IeC {\CYRP }\IeC {\cyrr }\IeC {\cyro }\IeC {\cyrs }\IeC {\cyrt }\IeC {\cyre }\IeC {\cyrishrt }\IeC {\cyrsh }\IeC {\cyri }\IeC {\cyre } \IeC {\cyrm }\IeC {\cyrn }\IeC {\cyro }\IeC {\cyrg }\IeC {\cyro }\IeC {\cyrg }\IeC {\cyrr }\IeC {\cyra }\IeC {\cyrn }\IeC {\cyrn }\IeC {\cyri }\IeC {\cyrk }\IeC {\cyri } \IeC {\cyrv }~\IeC {\cyrm }\IeC {\cyrn }\IeC {\cyro }\IeC {\cyrg }\IeC {\cyro }\IeC {\cyrm }\IeC {\cyre }\IeC {\cyrr }\IeC {\cyrn }\IeC {\cyro }\IeC {\cyrm } \IeC {\cyrp }\IeC {\cyrr }\IeC {\cyro }\IeC {\cyrs }\IeC {\cyrt }\IeC {\cyrr }\IeC {\cyra }\IeC {\cyrn }\IeC {\cyrs }\IeC {\cyrt }\IeC {\cyrv }\IeC {\cyre }. \emph {\IeC {\CYRYU }.~\IeC {\CYRM }.~\IeC {\CYRB }\IeC {\cyru }\IeC {\cyrr }\IeC {\cyrm }\IeC {\cyra }\IeC {\cyrn }}, \emph {\IeC {\CYRA }.~\IeC {\CYRYA }.~\IeC {\CYRK }\IeC {\cyra }\IeC {\cyrn }\IeC {\cyre }\IeC {\cyrl }\IeC {\cyrsftsn }-\IeC {\CYRB }\IeC {\cyre }\IeC {\cyrl }\IeC {\cyro }\IeC {\cyrv }}}}{392}{subsubsection.2.15.3.1}
\contentsline {subsubsection}{\ull{\numberline {15.3.2}\IeC {\CYRM }\IeC {\cyrn }\IeC {\cyro }\IeC {\cyrg }\IeC {\cyro }\IeC {\cyrm }\IeC {\cyre }\IeC {\cyrr }\IeC {\cyrn }\IeC {\cyrery }\IeC {\cyre } \IeC {\cyro }\IeC {\cyrb }\IeC {\cyrhrdsn }\IeC {\cyryo }\IeC {\cyrm }\IeC {\cyrery }}}{397}{subsubsection.2.15.3.2}
\contentsline {subsubsection}{\numberline {15.3.3}\IeC {\CYRO }\IeC {\cyrb }\IeC {\cyrhrdsn }\IeC {\cyryo }\IeC {\cyrm }\IeC {\cyrery } \IeC {\cyri }~\IeC {\cyrs }\IeC {\cyre }\IeC {\cyrch }\IeC {\cyre }\IeC {\cyrn }\IeC {\cyri }\IeC {\cyrya }}{399}{}%{subsubsection.2.15.3.3}
\contentsline {subsubsection}{\numberline {15.3.4}\IeC {\CYRD }\IeC {\cyrv }\IeC {\cyre } \IeC {\cyrz }\IeC {\cyra }\IeC {\cyrd }\IeC {\cyra }\IeC {\cyrch }\IeC {\cyri } \IeC {\cyrd }\IeC {\cyrl }\IeC {\cyrya } \IeC {\cyri }\IeC {\cyrs }\IeC {\cyrs }\IeC {\cyrl }\IeC {\cyre }\IeC {\cyrd }\IeC {\cyro }\IeC {\cyrv }\IeC {\cyra }\IeC {\cyrn }\IeC {\cyri }\IeC {\cyrya }}{400}{}%{subsubsection.2.15.3.4}
\contentsline {subsubsection}{\numberline {15.3.5}\IeC {\CYRR }\IeC {\cyra }\IeC {\cyrz }\IeC {\cyrb }\IeC {\cyri }\IeC {\cyre }\IeC {\cyrn }\IeC {\cyri }\IeC {\cyre } \IeC {\cyrn }\IeC {\cyra } \IeC {\cyrch }\IeC {\cyra }\IeC {\cyrs }\IeC {\cyrt }\IeC {\cyri } \IeC {\cyrm }\IeC {\cyre }\IeC {\cyrn }\IeC {\cyrsftsn }\IeC {\cyrsh }\IeC {\cyre }\IeC {\cyrg }\IeC {\cyro } \IeC {\cyrd }\IeC {\cyri }\IeC {\cyra }\IeC {\cyrm }\IeC {\cyre }\IeC {\cyrt }\IeC {\cyrr }\IeC {\cyra }. \emph {\IeC {\CYRA }.~\IeC {\CYRM }.~\IeC {\CYRR }\IeC {\cyra }\IeC {\cyrishrt }\IeC {\cyrg }\IeC {\cyro }\IeC {\cyrr }\IeC {\cyro }\IeC {\cyrd }\IeC {\cyrs }\IeC {\cyrk }\IeC {\cyri }\IeC {\cyrishrt }}}{401}{}%{subsubsection.2.15.3.5}
\contentsline {section}{\numberline {16}\IeC {\CYRR }\IeC {\cyra }\IeC {\cyrz }\IeC {\cyrn }\IeC {\cyrery }\IeC {\cyre } \IeC {\cyrz }\IeC {\cyra }\IeC {\cyrd }\IeC {\cyra }\IeC {\cyrch }\IeC {\cyri } \IeC {\cyrp }\IeC {\cyro } \IeC {\cyrg }\IeC {\cyre }\IeC {\cyro }\IeC {\cyrm }\IeC {\cyre }\IeC {\cyrt }\IeC {\cyrr }\IeC {\cyri }\IeC {\cyri }}{407}{section.2.16}
\contentsline {subsection}{\numberline {16.1}\IeC {\CYRG }\IeC {\cyre }\IeC {\cyro }\IeC {\cyrm }\IeC {\cyre }\IeC {\cyrt }\IeC {\cyrr }\IeC {\cyri }\IeC {\cyrch }\IeC {\cyre }\IeC {\cyrs }\IeC {\cyrk }\IeC {\cyri }\IeC {\cyre } \IeC {\cyrz }\IeC {\cyra }\IeC {\cyrd }\IeC {\cyra }\IeC {\cyrch }\IeC {\cyri } \IeC {\cyrn }\IeC {\cyra } \IeC {\cyrerev }\IeC {\cyrk }\IeC {\cyrs }\IeC {\cyrt }\IeC {\cyrr }\IeC {\cyre }\IeC {\cyrm }\IeC {\cyra }\IeC {\cyrl }\IeC {\cyrsftsn }\IeC {\cyrn }\IeC {\cyrery }\IeC {\cyre } \IeC {\cyrz }\IeC {\cyrn }\IeC {\cyra }\IeC {\cyrch }\IeC {\cyre }\IeC {\cyrn }\IeC {\cyri }\IeC {\cyrya } (2). \emph {\IeC {\CYRA }.~\IeC {\CYRD }.~\IeC {\CYRB }\IeC {\cyrl }\IeC {\cyri }\IeC {\cyrn }\IeC {\cyrk }\IeC {\cyro }\IeC {\cyrv }}}{407}{}%{subsection.2.16.1}
\contentsline {subsection}{\numberline {16.2}\IeC {\CYRP }\IeC {\cyrl }\IeC {\cyro }\IeC {\cyrshch }\IeC {\cyra }\IeC {\cyrd }\IeC {\cyri } (2). \emph {\IeC {\CYRA }.~\IeC {\CYRD }.~\IeC {\CYRB }\IeC {\cyrl }\IeC {\cyri }\IeC {\cyrn }\IeC {\cyrk }\IeC {\cyro }\IeC {\cyrv }}}{413}{}%{subsection.2.16.2}
\contentsline {subsection}{\numberline {16.3}\IeC {\CYRK }\IeC {\cyro }\IeC {\cyrn }\IeC {\cyri }\IeC {\cyrch }\IeC {\cyre }\IeC {\cyrs }\IeC {\cyrk }\IeC {\cyri }\IeC {\cyre } \IeC {\cyrs }\IeC {\cyre }\IeC {\cyrch }\IeC {\cyre }\IeC {\cyrn }\IeC {\cyri }\IeC {\cyrya } (3$^*$). \emph {\IeC {\CYRA }.~\IeC {\CYRV }.~\IeC {\CYRA }\IeC {\cyrk }\IeC {\cyro }\IeC {\cyrp }\IeC {\cyrya }\IeC {\cyrn }}}{423}{}%{subsection.2.16.3}
\contentsline {subsection}{\numberline {16.4}\IeC {\CYRK }\IeC {\cyrr }\IeC {\cyri }\IeC {\cyrv }\IeC {\cyro }\IeC {\cyrl }\IeC {\cyri }\IeC {\cyrn }\IeC {\cyre }\IeC {\cyrishrt }\IeC {\cyrn }\IeC {\cyrery }\IeC {\cyre } \IeC {\cyrt }\IeC {\cyrr }\IeC {\cyre }\IeC {\cyru }\IeC {\cyrg }\IeC {\cyro }\IeC {\cyrl }\IeC {\cyrsftsn }\IeC {\cyrn }\IeC {\cyri }\IeC {\cyrk }\IeC {\cyri } \IeC {\cyri }~\IeC {\cyrn }\IeC {\cyre }\IeC {\cyre }\IeC {\cyrv }\IeC {\cyrk }\IeC {\cyrl }\IeC {\cyri }\IeC {\cyrd }\IeC {\cyro }\IeC {\cyrv }\IeC {\cyra } \IeC {\cyrg }\IeC {\cyre }\IeC {\cyro }\IeC {\cyrm }\IeC {\cyre }\IeC {\cyrt }\IeC {\cyrr }\IeC {\cyri }\IeC {\cyrya } (3$^*$). \emph {\IeC {\CYRM }.~\IeC {\CYRB }.~\IeC {\CYRS }\IeC {\cyrk }\IeC {\cyro }\IeC {\cyrp }\IeC {\cyre }\IeC {\cyrn }\IeC {\cyrk }\IeC {\cyro }\IeC {\cyrv }}}{434}{}%{subsection.2.16.4}
\contentsline {chapter}{\numberline {3}\IeC {\CYRK }\IeC {\cyro }\IeC {\cyrm }\IeC {\cyrb }\IeC {\cyri }\IeC {\cyrn }\IeC {\cyra }\IeC {\cyrt }\IeC {\cyro }\IeC {\cyrr }\IeC {\cyri }\IeC {\cyrk }\IeC {\cyra }}{441}{chapter.3}
\contentsline {section}{\numberline {17}\IeC {\CYRP }\IeC {\cyro }\IeC {\cyrd }\IeC {\cyrs }\IeC {\cyrch }\IeC {\cyre }\IeC {\cyrt }\IeC {\cyrery } \IeC {\cyrv }~\IeC {\cyrk }\IeC {\cyro }\IeC {\cyrm }\IeC {\cyrb }\IeC {\cyri }\IeC {\cyrn }\IeC {\cyra }\IeC {\cyrt }\IeC {\cyro }\IeC {\cyrr }\IeC {\cyri }\IeC {\cyrk }\IeC {\cyre }}{441}{section.3.17}
\contentsline {subsection}{\ull{\numberline {17.1}\IeC {\CYRP }\IeC {\cyro }\IeC {\cyrd }\IeC {\cyrs }\IeC {\cyrch }\IeC {\cyre }\IeC {\cyrt }\IeC {\cyrery } \IeC {\cyrch }\IeC {\cyri }\IeC {\cyrs }\IeC {\cyrl }\IeC {\cyra } \IeC {\cyrs }\IeC {\cyrp }\IeC {\cyro }\IeC {\cyrs }\IeC {\cyro }\IeC {\cyrb }\IeC {\cyro }\IeC {\cyrv } (1). \emph {\IeC {\CYRA }.~\IeC {\CYRA }.~\IeC {\CYRG }\IeC {\cyra }\IeC {\cyrv }\IeC {\cyrr }\IeC {\cyri }\IeC {\cyrl }\IeC {\cyryu }\IeC {\cyrk }}, \emph {\IeC {\CYRD }.~\IeC {\CYRA }.~\IeC {\CYRP }\IeC {\cyre }\IeC {\cyrr }\IeC {\cyrm }\IeC {\cyrya }\IeC {\cyrk }\IeC {\cyro }\IeC {\cyrv }}}}{441}{subsection.3.17.1}
\contentsline {subsection}{\numberline {17.2}\IeC {\CYRN }\IeC {\cyra }\IeC {\cyrb }\IeC {\cyro }\IeC {\cyrr }\IeC {\cyrery } \IeC {\cyrp }\IeC {\cyro }\IeC {\cyrd }\IeC {\cyrm }\IeC {\cyrn }\IeC {\cyro }\IeC {\cyrzh }\IeC {\cyre }\IeC {\cyrs }\IeC {\cyrt }\IeC {\cyrv } (2). \emph {\IeC {\CYRD }.~\IeC {\CYRA }.~\IeC {\CYRP }\IeC {\cyre }\IeC {\cyrr }\IeC {\cyrm }\IeC {\cyrya }\IeC {\cyrk }\IeC {\cyro }\IeC {\cyrv }}}{445}{}%{subsection.3.17.2}
\contentsline {subsection}{\ull{\numberline {17.3}\IeC {\CYRF }\IeC {\cyro }\IeC {\cyrr }\IeC {\cyrm }\IeC {\cyru }\IeC {\cyrl }\IeC {\cyra } \IeC {\cyrv }\IeC {\cyrk }\IeC {\cyrl }\IeC {\cyryu }\IeC {\cyrch }\IeC {\cyre }\IeC {\cyrn }\IeC {\cyri }\IeC {\cyrishrt } \IeC {\cyri }~\IeC {\cyri }\IeC {\cyrs }\IeC {\cyrk }\IeC {\cyrl }\IeC {\cyryu }\IeC {\cyrch }\IeC {\cyre }\IeC {\cyrn }\IeC {\cyri }\IeC {\cyrishrt } (2). \emph {\IeC {\CYRD }.~\IeC {\CYRA }.~\IeC {\CYRP }\IeC {\cyre }\IeC {\cyrr }\IeC {\cyrm }\IeC {\cyrya }\IeC {\cyrk }\IeC {\cyro }\IeC {\cyrv }}}}{448}{subsection.3.17.3}
\contentsline {subsection}{\numberline {17.4}\IeC {\CYRN }\IeC {\cyre }\IeC {\cyrs }\IeC {\cyrk }\IeC {\cyro }\IeC {\cyrl }\IeC {\cyrsftsn }\IeC {\cyrk }\IeC {\cyro } \IeC {\cyrv }\IeC {\cyrz }\IeC {\cyrg }\IeC {\cyrl }\IeC {\cyrya }\IeC {\cyrd }\IeC {\cyro }\IeC {\cyrv } \IeC {\cyrn }\IeC {\cyra } \IeC {\cyrch }\IeC {\cyri }\IeC {\cyrs }\IeC {\cyrl }\IeC {\cyra } \IeC {\CYRK }\IeC {\cyra }\IeC {\cyrt }\IeC {\cyra }\IeC {\cyrl }\IeC {\cyra }\IeC {\cyrn }\IeC {\cyra }. \emph {\IeC {\CYRG }.~\IeC {\CYRB }.~\IeC {\CYRSH }\IeC {\cyra }\IeC {\cyrb }\IeC {\cyra }\IeC {\cyrt }}{{}{}} }{455}{}%{subsection.3.17.4}
\contentsline {section}{\numberline {18}\IeC {\CYRP }\IeC {\cyrr }\IeC {\cyri }\IeC {\cyrn }\IeC {\cyrc }\IeC {\cyri }\IeC {\cyrp } \IeC {\CYRD }\IeC {\cyri }\IeC {\cyrr }\IeC {\cyri }\IeC {\cyrh }\IeC {\cyrl }\IeC {\cyre } \IeC {\cyri }~\IeC {\cyri }\IeC {\cyrn }\IeC {\cyrd }\IeC {\cyru }\IeC {\cyrk }\IeC {\cyrc }\IeC {\cyri }\IeC {\cyrya }}{462}{section.3.18}
\contentsline {subsection}{\numberline {18.1}\IeC {\CYRP }\IeC {\cyrr }\IeC {\cyri }\IeC {\cyrn }\IeC {\cyrc }\IeC {\cyri }\IeC {\cyrp } \IeC {\CYRD }\IeC {\cyri }\IeC {\cyrr }\IeC {\cyri }\IeC {\cyrh }\IeC {\cyrl }\IeC {\cyre } (1). \emph {\IeC {\CYRA }.~\IeC {\CYRYA }.~\IeC {\CYRK }\IeC {\cyra }\IeC {\cyrn }\IeC {\cyre }\IeC {\cyrl }\IeC {\cyrsftsn }-\IeC {\CYRB }\IeC {\cyre }\IeC {\cyrl }\IeC {\cyro }\IeC {\cyrv }}}{462}{}%{subsection.3.18.1}
\contentsline {subsection}{\numberline {18.2}\IeC {\CYRP }\IeC {\cyrr }\IeC {\cyra }\IeC {\cyrv }\IeC {\cyri }\IeC {\cyrl }\IeC {\cyro } \IeC {\cyrk }\IeC {\cyrr }\IeC {\cyra }\IeC {\cyrishrt }\IeC {\cyrn }\IeC {\cyre }\IeC {\cyrg }\IeC {\cyro } (2). \emph {\IeC {\CYRA }.~\IeC {\CYRYA }.~\IeC {\CYRK }\IeC {\cyra }\IeC {\cyrn }\IeC {\cyre }\IeC {\cyrl }\IeC {\cyrsftsn }-\IeC {\CYRB }\IeC {\cyre }\IeC {\cyrl }\IeC {\cyro }\IeC {\cyrv }}}{466}{}%{subsection.3.18.2}
\contentsline {subsection}{\numberline {18.3}\IeC {\CYRC }\IeC {\cyri }\IeC {\cyrk }\IeC {\cyrl }\IeC {\cyri }\IeC {\cyrch }\IeC {\cyrn }\IeC {\cyro }\IeC {\cyrs }\IeC {\cyrt }\IeC {\cyrsftsn } I (2){{}{}} . \emph {\IeC {\CYRA }.~\IeC {\CYRYA }.~\IeC {\CYRK }\IeC {\cyra }\IeC {\cyrn }\IeC {\cyre }\IeC {\cyrl }\IeC {\cyrsftsn }-\IeC {\CYRB }\IeC {\cyre }\IeC {\cyrl }\IeC {\cyro }\IeC {\cyrv }}}{468}{}%{subsection.3.18.3}
\contentsline {subsection}{\numberline {18.4}\IeC {\CYRC }\IeC {\cyri }\IeC {\cyrk }\IeC {\cyrl }\IeC {\cyri }\IeC {\cyrch }\IeC {\cyrn }\IeC {\cyro }\IeC {\cyrs }\IeC {\cyrt }\IeC {\cyrsftsn } II (2). \emph {\IeC {\CYRP }.~\IeC {\CYRA }.~\IeC {\CYRK }\IeC {\cyro }\IeC {\cyrzh }\IeC {\cyre }\IeC {\cyrv }\IeC {\cyrn }\IeC {\cyri }\IeC {\cyrk }\IeC {\cyro }\IeC {\cyrv }}}{471}{}%{subsection.3.18.4}
\contentsline {subsection}{\numberline {18.5}\IeC {\CYRK }\IeC {\cyro }\IeC {\cyrn }\IeC {\cyre }\IeC {\cyrch }\IeC {\cyrn }\IeC {\cyro }\IeC {\cyre } \IeC {\cyri }~\IeC {\cyrs }\IeC {\cyrch }\IeC {\cyryo }\IeC {\cyrt }\IeC {\cyrn }\IeC {\cyro }\IeC {\cyre } (2). \emph {\IeC {\CYRP }.~\IeC {\CYRA }.~\IeC {\CYRK }\IeC {\cyro }\IeC {\cyrzh }\IeC {\cyre }\IeC {\cyrv }\IeC {\cyrn }\IeC {\cyri }\IeC {\cyrk }\IeC {\cyro }\IeC {\cyrv }}}{475}{}%{subsection.3.18.5}
\contentsline {subsection}{\numberline {18.6}\IeC {\CYRN }\IeC {\cyre }\IeC {\cyrm }\IeC {\cyrn }\IeC {\cyro }\IeC {\cyrg }\IeC {\cyro } \IeC {\cyri }\IeC {\cyrn }\IeC {\cyrd }\IeC {\cyru }\IeC {\cyrk }\IeC {\cyrc }\IeC {\cyri }\IeC {\cyri } \IeC {\cyri }~\IeC {\cyrp }\IeC {\cyre }\IeC {\cyrr }\IeC {\cyre }\IeC {\cyrb }\IeC {\cyro }\IeC {\cyrr }\IeC {\cyra } (3). \emph {\IeC {\CYRI }.~\IeC {\CYRN }.~\IeC {\CYRSH }\IeC {\cyrn }\IeC {\cyru }\IeC {\cyrr }\IeC {\cyrn }\IeC {\cyri }\IeC {\cyrk }\IeC {\cyro }\IeC {\cyrv }}}{481}{}%{subsection.3.18.6}
\contentsline {section}{\numberline {19}\IeC {\CYRG }\IeC {\cyrr }\IeC {\cyra }\IeC {\cyrf }\IeC {\cyrery }. \emph {\IeC {\CYRD }.~\IeC {\CYRA }.~\IeC {\CYRP }\IeC {\cyre }\IeC {\cyrr }\IeC {\cyrm }\IeC {\cyrya }\IeC {\cyrk }\IeC {\cyro }\IeC {\cyrv }}, \emph {\IeC {\CYRA }.~\IeC {\CYRB }.~\IeC {\CYRS }\IeC {\cyrk }\IeC {\cyro }\IeC {\cyrp }\IeC {\cyre }\IeC {\cyrn }\IeC {\cyrk }\IeC {\cyro }\IeC {\cyrv }}}{485}{section.3.19}
\contentsline {subsection}{\numberline {19.1}\IeC {\CYRG }\IeC {\cyrr }\IeC {\cyra }\IeC {\cyrf }\IeC {\cyrery } \IeC {\cyrp }\IeC {\cyro }\IeC {\cyrd } \IeC {\cyrsh }\IeC {\cyru }\IeC {\cyrb }\IeC {\cyro }\IeC {\cyrishrt } (1)}{485}{}%{subsection.3.19.1}
\contentsline {subsection}{\numberline {19.2}\IeC {\CYRP }\IeC {\cyro }\IeC {\cyrd }\IeC {\cyrs }\IeC {\cyrch }\IeC {\cyryo }\IeC {\cyrt }\IeC {\cyrery } \IeC {\cyrv }~\IeC {\cyrg }\IeC {\cyrr }\IeC {\cyra }\IeC {\cyrf }\IeC {\cyra }\IeC {\cyrh } (1)}{488}{}%{subsection.3.19.2}
\contentsline {subsection}{\numberline {19.3}\IeC {\CYRP }\IeC {\cyru }\IeC {\cyrt }\IeC {\cyri } \IeC {\cyrv }~\IeC {\cyrg }\IeC {\cyrr }\IeC {\cyra }\IeC {\cyrf }\IeC {\cyra }\IeC {\cyrh } (2)}{491}{}%{subsection.3.19.3}
\contentsline {section}{\numberline {20}\IeC {\CYRK }\IeC {\cyro }\IeC {\cyrn }\IeC {\cyrs }\IeC {\cyrt }\IeC {\cyrr }\IeC {\cyru }\IeC {\cyrk }\IeC {\cyrc }\IeC {\cyri }\IeC {\cyri } \IeC {\cyri }~\IeC {\cyri }\IeC {\cyrn }\IeC {\cyrv }\IeC {\cyra }\IeC {\cyrr }\IeC {\cyri }\IeC {\cyra }\IeC {\cyrn }\IeC {\cyrt }\IeC {\cyrery }}{494}{section.3.20}
\contentsline {subsection}{\ull{\numberline {20.1}\IeC {\CYRK }\IeC {\cyro }\IeC {\cyrn }\IeC {\cyrs }\IeC {\cyrt }\IeC {\cyrr }\IeC {\cyru }\IeC {\cyrk }\IeC {\cyrc }\IeC {\cyri }\IeC {\cyri }{{}{}} \ (1). \emph {\IeC {\CYRA }.~\IeC {\CYRV }.~\IeC {\CYRSH }\IeC {\cyra }\IeC {\cyrp }\IeC {\cyro }\IeC {\cyrv }\IeC {\cyra }\IeC {\cyrl }\IeC {\cyro }\IeC {\cyrv }{{}{}} }}}{494}{subsection.3.20.1}
\contentsline {subsection}{\numberline {20.2}\IeC {\CYRI }\IeC {\cyrn }\IeC {\cyrv }\IeC {\cyra }\IeC {\cyrr }\IeC {\cyri }\IeC {\cyra }\IeC {\cyrn }\IeC {\cyrt }\IeC {\cyrery } I (1). \emph {\IeC {\CYRA }.~\IeC {\CYRYA }.~\IeC {\CYRK }\IeC {\cyra }\IeC {\cyrn }\IeC {\cyre }\IeC {\cyrl }\IeC {\cyrsftsn }-\IeC {\CYRB }\IeC {\cyre }\IeC {\cyrl }\IeC {\cyro }\IeC {\cyrv }}}{509}{}%{subsection.3.20.2}
\contentsline {subsection}{\numberline {20.3}\IeC {\CYRI }\IeC {\cyrn }\IeC {\cyrv }\IeC {\cyra }\IeC {\cyrr }\IeC {\cyri }\IeC {\cyra }\IeC {\cyrn }\IeC {\cyrt }\IeC {\cyrery } II (1){{}{}} . \emph {\IeC {\CYRA }.~\IeC {\CYRV }.~\IeC {\CYRSH }\IeC {\cyra }\IeC {\cyrp }\IeC {\cyro }\IeC {\cyrv }\IeC {\cyra }\IeC {\cyrl }\IeC {\cyro }\IeC {\cyrv }}}{512}{}%{subsection.3.20.3}
\contentsline {subsection}{\numberline {20.4}\IeC {\CYRR }\IeC {\cyra }\IeC {\cyrs }\IeC {\cyrk }\IeC {\cyrr }\IeC {\cyra }\IeC {\cyrs }\IeC {\cyrk }\IeC {\cyri }}{522}{}%{subsection.3.20.4}
\contentsline {subsubsection}{\numberline {20.4.1}\IeC {\CYRZ }\IeC {\cyra }\IeC {\cyrm }\IeC {\cyro }\IeC {\cyrshch }\IeC {\cyre }\IeC {\cyrn }\IeC {\cyri }\IeC {\cyrya } (1). \emph {\IeC {\CYRA }.~\IeC {\CYRYA }.~\IeC {\CYRK }\IeC {\cyra }\IeC {\cyrn }\IeC {\cyre }\IeC {\cyrl }\IeC {\cyrsftsn }-\IeC {\CYRB }\IeC {\cyre }\IeC {\cyrl }\IeC {\cyro }\IeC {\cyrv }}}{522}{}%{subsubsection.3.20.4.1}
\contentsline {subsubsection}{\numberline {20.4.2}\IeC {\CYRT }\IeC {\cyra }\IeC {\cyrb }\IeC {\cyrl }\IeC {\cyri }\IeC {\cyrc }\IeC {\cyrery } (2){{}{}} . \emph {\IeC {\CYRD }.~\IeC {\CYRA }.~\IeC {\CYRP }\IeC {\cyre }\IeC {\cyrr }\IeC {\cyrm }\IeC {\cyrya }\IeC {\cyrk }\IeC {\cyro }\IeC {\cyrv }}}{523}{}%{subsubsection.3.20.4.2}
\contentsline {subsection}{\ull{\numberline {20.5}\IeC {\CYRP }\IeC {\cyro }\IeC {\cyrl }\IeC {\cyru }\IeC {\cyri }\IeC {\cyrn }\IeC {\cyrv }\IeC {\cyra }\IeC {\cyrr }\IeC {\cyri }\IeC {\cyra }\IeC {\cyrn }\IeC {\cyrt }\IeC {\cyrery }{{}{}} (1). \emph {\IeC {\CYRA }.~\IeC {\CYRV }.~\IeC {\CYRSH }\IeC {\cyra }\IeC {\cyrp }\IeC {\cyro }\IeC {\cyrv }\IeC {\cyra }\IeC {\cyrl }\IeC {\cyro }\IeC {\cyrv }}}}{524}{subsection.3.20.5}
\contentsline {section}{\numberline {21}\IeC {\CYRA }\IeC {\cyrl }\IeC {\cyrg }\IeC {\cyro }\IeC {\cyrr }\IeC {\cyri }\IeC {\cyrt }\IeC {\cyrm }\IeC {\cyrery }}{536}{section.3.21}
\contentsline {subsection}{\ull{\numberline {21.1}\IeC {\CYRI }\IeC {\cyrg }\IeC {\cyrr }\IeC {\cyrery } (1){{}{}} . \emph {\IeC {\CYRD }.~\IeC {\CYRA }.~\IeC {\CYRP }\IeC {\cyre }\IeC {\cyrr }\IeC {\cyrm }\IeC {\cyrya }\IeC {\cyrk }\IeC {\cyro }\IeC {\cyrv }}, \emph {\IeC {\CYRM }.~\IeC {\CYRB }.~\IeC {\CYRS }\IeC {\cyrk }\IeC {\cyro }\IeC {\cyrp }\IeC {\cyre }\IeC {\cyrn }\IeC {\cyrk }\IeC {\cyro }\IeC {\cyrv }}, \emph {\IeC {\CYRA }.~\IeC {\CYRV }.~\IeC {\CYRSH }\IeC {\cyra }\IeC {\cyrp }\IeC {\cyro }\IeC {\cyrv }\IeC {\cyra }\IeC {\cyrl }\IeC {\cyro }\IeC {\cyrv }}}}{536}{subsection.3.21.1}
\contentsline {subsection}{\numberline {21.2}\IeC {\CYRI }\IeC {\cyrn }\IeC {\cyrf }\IeC {\cyro }\IeC {\cyrr }\IeC {\cyrm }\IeC {\cyra }\IeC {\cyrc }\IeC {\cyri }\IeC {\cyro }\IeC {\cyrn }\IeC {\cyrn }\IeC {\cyrery }\IeC {\cyre } \IeC {\cyrz }\IeC {\cyra }\IeC {\cyrd }\IeC {\cyra }\IeC {\cyrch }\IeC {\cyri } (2). \emph {\IeC {\CYRA }.~\IeC {\CYRYA }.~\IeC {\CYRK }\IeC {\cyra }\IeC {\cyrn }\IeC {\cyre }\IeC {\cyrl }\IeC {\cyrsftsn }-\IeC {\CYRB }\IeC {\cyre }\IeC {\cyrl }\IeC {\cyro }\IeC {\cyrv }}}{549}{}%{subsection.3.21.2}
\contentsline {subsection}{\numberline {21.3}\IeC {\CYRK }\IeC {\cyro }\IeC {\cyrd }\IeC {\cyrery }, \IeC {\cyri }\IeC {\cyrs }\IeC {\cyrp }\IeC {\cyrr }\IeC {\cyra }\IeC {\cyrv }\IeC {\cyrl }\IeC {\cyrya }\IeC {\cyryu }\IeC {\cyrshch }\IeC {\cyri }\IeC {\cyre } \IeC {\cyro }\IeC {\cyrsh }\IeC {\cyri }\IeC {\cyrb }\IeC {\cyrk }\IeC {\cyri } (2). \emph {\IeC {\CYRM }.~\IeC {\CYRB }.~\IeC {\CYRS }\IeC {\cyrk }\IeC {\cyro }\IeC {\cyrp }\IeC {\cyre }\IeC {\cyrn }\IeC {\cyrk }\IeC {\cyro }\IeC {\cyrv }}}{553}{}%{subsection.3.21.3}
\contentsline {subsection}{\numberline {21.4}\IeC {\CYRB }\IeC {\cyru }\IeC {\cyrl }\IeC {\cyre }\IeC {\cyrv } \IeC {\cyrk }\IeC {\cyru }\IeC {\cyrb } (2). \emph {\IeC {\CYRA }.~\IeC {\CYRB }.~\IeC {\CYRS }\IeC {\cyrk }\IeC {\cyro }\IeC {\cyrp }\IeC {\cyre }\IeC {\cyrn }\IeC {\cyrk }\IeC {\cyro }\IeC {\cyrv }}}{556}{}%{subsection.3.21.4}
\contentsline {subsection}{\numberline {21.5}\IeC {\CYRV }\IeC {\cyrery }\IeC {\cyrr }\IeC {\cyra }\IeC {\cyrz }\IeC {\cyri }\IeC {\cyrm }\IeC {\cyro }\IeC {\cyrs }\IeC {\cyrt }\IeC {\cyrsftsn } \IeC {\cyrd }\IeC {\cyrl }\IeC {\cyrya } \IeC {\cyrf }\IeC {\cyru }\IeC {\cyrn }\IeC {\cyrk }\IeC {\cyrc }\IeC {\cyri }\IeC {\cyrishrt } \IeC {\cyra }\IeC {\cyrl }\IeC {\cyrg }\IeC {\cyre }\IeC {\cyrb }\IeC {\cyrr }\IeC {\cyrery } \IeC {\cyrl }\IeC {\cyro }\IeC {\cyrg }\IeC {\cyri }\IeC {\cyrk }\IeC {\cyri }. \emph {\IeC {\CYRA }.~\IeC {\CYRB }.~\IeC {\CYRS }\IeC {\cyrk }\IeC {\cyro }\IeC {\cyrp }\IeC {\cyre }\IeC {\cyrn }\IeC {\cyrk }\IeC {\cyro }\IeC {\cyrv }}}{561}{}%{subsection.3.21.5}
\contentsline {subsubsection}{\numberline {21.5.1}\IeC {\CYRP }\IeC {\cyrr }\IeC {\cyri }\IeC {\cyrm }\IeC {\cyre }\IeC {\cyrr }\IeC {\cyrery } \IeC {\cyri }~\IeC {\cyro }\IeC {\cyrp }\IeC {\cyrr }\IeC {\cyre }\IeC {\cyrd }\IeC {\cyre }\IeC {\cyrl }\IeC {\cyre }\IeC {\cyrn }\IeC {\cyri }\IeC {\cyrya }}{561}{}%{subsubsection.3.21.5.1}
\contentsline {subsubsection}{\numberline {21.5.2}\IeC {\CYRT }\IeC {\cyre }\IeC {\cyro }\IeC {\cyrr }\IeC {\cyre }\IeC {\cyrm }\IeC {\cyra } \IeC {\CYRP }\IeC {\cyro }\IeC {\cyrs }\IeC {\cyrt }\IeC {\cyra } (2*)}{563}{}%{subsubsection.3.21.5.2}
\contentsline {subsection}{\numberline {21.6}\IeC {\CYRS }\IeC {\cyrl }\IeC {\cyro }\IeC {\cyrzh }\IeC {\cyrn }\IeC {\cyro }\IeC {\cyrs }\IeC {\cyrt }\IeC {\cyrsftsn } \IeC {\cyrs }\IeC {\cyru }\IeC {\cyrm }\IeC {\cyrm }\IeC {\cyri }\IeC {\cyrr }\IeC {\cyro }\IeC {\cyrv }\IeC {\cyra }\IeC {\cyrn }\IeC {\cyri }\IeC {\cyrya }{{}{}} {}. \emph {\IeC {\CYRYU }.~\IeC {\CYRG }.~\IeC {\CYRK }\IeC {\cyru }\IeC {\cyrd }\IeC {\cyrr }\IeC {\cyrya }\IeC {\cyrsh }\IeC {\cyro }\IeC {\cyrv }}, \emph {\IeC {\CYRA }.~\IeC {\CYRB }.~\IeC {\CYRS }\IeC {\cyrk }\IeC {\cyro }\IeC {\cyrp }\IeC {\cyre }\IeC {\cyrn }\IeC {\cyrk }\IeC {\cyro }\IeC {\cyrv }}}{567}{}%{subsection.3.21.6}
\contentsline {subsubsection}{\numberline {21.6.1}\IeC {\CYRV }\IeC {\cyrv }\IeC {\cyro }\IeC {\cyrd }\IeC {\cyrn }\IeC {\cyrery }\IeC {\cyre } \IeC {\cyrz }\IeC {\cyra }\IeC {\cyrd }\IeC {\cyra }\IeC {\cyrch }\IeC {\cyri } (2)}{567}{}%{subsubsection.3.21.6.1}
\contentsline {subsubsection}{\numberline {21.6.2}\IeC {\CYRO }\IeC {\cyrp }\IeC {\cyrr }\IeC {\cyre }\IeC {\cyrd }\IeC {\cyre }\IeC {\cyrl }\IeC {\cyre }\IeC {\cyrn }\IeC {\cyri }\IeC {\cyrya } \IeC {\cyri }~\IeC {\cyrp }\IeC {\cyrr }\IeC {\cyri }\IeC {\cyrm }\IeC {\cyre }\IeC {\cyrr }\IeC {\cyrery } (3*)}{568}{}%{subsubsection.3.21.6.2}
\contentsline {subsubsection}{\numberline {21.6.3}\IeC {\CYRA }\IeC {\cyrs }\IeC {\cyri }\IeC {\cyrm }\IeC {\cyrp }\IeC {\cyrt }\IeC {\cyro }\IeC {\cyrt }\IeC {\cyri }\IeC {\cyrch }\IeC {\cyre }\IeC {\cyrs }\IeC {\cyrk }\IeC {\cyri }\IeC {\cyre } \IeC {\cyro }\IeC {\cyrc }\IeC {\cyre }\IeC {\cyrn }\IeC {\cyrk }\IeC {\cyri } (4*)}{570}{}%{subsubsection.3.21.6.3}
\contentsline {section}{\ull{\numberline {22}\IeC {\CYRV }\IeC {\cyre }\IeC {\cyrr }\IeC {\cyro }\IeC {\cyrya }\IeC {\cyrt }\IeC {\cyrn }\IeC {\cyro }\IeC {\cyrs }\IeC {\cyrt }\IeC {\cyrsftsn }{{}{}} . \emph {\IeC {\CYRA }.~\IeC {\CYRA }.~\IeC {\CYRZ }\IeC {\cyra }\IeC {\cyrs }\IeC {\cyrl }\IeC {\cyra }\IeC {\cyrv }\IeC {\cyrs }\IeC {\cyrk }\IeC {\cyri }\IeC {\cyrishrt }}}}{579}{section.3.22}
\contentsline {subsection}{\ull{\numberline {22.1}\IeC {\CYRK }\IeC {\cyrl }\IeC {\cyra }\IeC {\cyrs }\IeC {\cyrs }\IeC {\cyri }\IeC {\cyrch }\IeC {\cyre }\IeC {\cyrs }\IeC {\cyrk }\IeC {\cyro }\IeC {\cyre } \IeC {\cyro }\IeC {\cyrp }\IeC {\cyrr }\IeC {\cyre }\IeC {\cyrd }\IeC {\cyre }\IeC {\cyrl }\IeC {\cyre }\IeC {\cyrn }\IeC {\cyri }\IeC {\cyre } \IeC {\cyrv }\IeC {\cyre }\IeC {\cyrr }\IeC {\cyro }\IeC {\cyrya }\IeC {\cyrt }\IeC {\cyrn }\IeC {\cyro }\IeC {\cyrs }\IeC {\cyrt }\IeC {\cyri } (1).}}{579}{subsection.3.22.1}
\contentsline {subsection}{\ull{\numberline {22.2}\IeC {\CYRB }\IeC {\cyro }\IeC {\cyrl }\IeC {\cyre }\IeC {\cyre } \IeC {\cyro }\IeC {\cyrb }\IeC {\cyrshch }\IeC {\cyre }\IeC {\cyre } \IeC {\cyro }\IeC {\cyrp }\IeC {\cyrr }\IeC {\cyre }\IeC {\cyrd }\IeC {\cyre }\IeC {\cyrl }\IeC {\cyre }\IeC {\cyrn }\IeC {\cyri }\IeC {\cyre } \IeC {\cyrv }\IeC {\cyre }\IeC {\cyrr }\IeC {\cyro }\IeC {\cyrya }\IeC {\cyrt }\IeC {\cyrn }\IeC {\cyro }\IeC {\cyrs }\IeC {\cyrt }\IeC {\cyri } (1)}}{583}{subsection.3.22.2}
\contentsline {subsection}{\ull{\numberline {22.3}\IeC {\CYRN }\IeC {\cyre }\IeC {\cyrz }\IeC {\cyra }\IeC {\cyrv }\IeC {\cyri }\IeC {\cyrs }\IeC {\cyri }\IeC {\cyrm }\IeC {\cyro }\IeC {\cyrs }\IeC {\cyrt }\IeC {\cyrsftsn } \IeC {\cyri }~\IeC {\cyru }\IeC {\cyrs }\IeC {\cyrl }\IeC {\cyro }\IeC {\cyrv }\IeC {\cyrn }\IeC {\cyra }\IeC {\cyrya } \IeC {\cyrv }\IeC {\cyre }\IeC {\cyrr }\IeC {\cyro }\IeC {\cyrya }\IeC {\cyrt }\IeC {\cyrn }\IeC {\cyro }\IeC {\cyrs }\IeC {\cyrt }\IeC {\cyrsftsn } (1)}}{585}{subsection.3.22.3}
\contentsline {subsection}{\ull{\numberline {22.4}\IeC {\CYRS }\IeC {\cyrl }\IeC {\cyru }\IeC {\cyrch }\IeC {\cyra }\IeC {\cyrishrt }\IeC {\cyrn }\IeC {\cyrery }\IeC {\cyre } \IeC {\cyrv }\IeC {\cyre }\IeC {\cyrl }\IeC {\cyri }\IeC {\cyrch }\IeC {\cyri }\IeC {\cyrn }\IeC {\cyrery } (3)}}{591}{subsection.3.22.4}
\contentsline {subsection}{\ull{\numberline {22.5}\IeC {\CYRI }\IeC {\cyrs }\IeC {\cyrp }\IeC {\cyrery }\IeC {\cyrt }\IeC {\cyra }\IeC {\cyrn }\IeC {\cyri }\IeC {\cyrya } \IeC {\CYRB }\IeC {\cyre }\IeC {\cyrr }\IeC {\cyrn }\IeC {\cyru }\IeC {\cyrl }\IeC {\cyrl }\IeC {\cyri } (3)}}{595}{subsection.3.22.5}
\contentsline {subsection}{\numberline {22.6}\IeC {\CYRS }\IeC {\cyrl }\IeC {\cyru }\IeC {\cyrch }\IeC {\cyra }\IeC {\cyrishrt }\IeC {\cyrn }\IeC {\cyrery }\IeC {\cyre } \IeC {\cyrb }\IeC {\cyrl }\IeC {\cyru }\IeC {\cyrzh }\IeC {\cyrd }\IeC {\cyra }\IeC {\cyrn }\IeC {\cyri }\IeC {\cyrya } \IeC {\cyri }~\IeC {\cyrerev }\IeC {\cyrl }\IeC {\cyre }\IeC {\cyrk }\IeC {\cyrt }\IeC {\cyrr }\IeC {\cyri }\IeC {\cyrch }\IeC {\cyre }\IeC {\cyrs }\IeC {\cyrk }\IeC {\cyri }\IeC {\cyre } \IeC {\cyrc }\IeC {\cyre }\IeC {\cyrp }\IeC {\cyri }{{}{}} \ (3). \emph {\IeC {\CYRA }.~\IeC {\CYRA }.~\IeC {\CYRZ }\IeC {\cyra }\IeC {\cyrs }\IeC {\cyrl }\IeC {\cyra }\IeC {\cyrv }\IeC {\cyrs }\IeC {\cyrk }\IeC {\cyri }\IeC {\cyrishrt }}, \emph {\IeC {\CYRM }.~\IeC {\CYRB }.~\IeC {\CYRS }\IeC {\cyrk }\IeC {\cyro }\IeC {\cyrp }\IeC {\cyre }\IeC {\cyrn }\IeC {\cyrk }\IeC {\cyro }\IeC {\cyrv }}, \emph {\IeC {\CYRA }.~\IeC {\CYRV }.~\IeC {\CYRU }\IeC {\cyrs }\IeC {\cyrt }\IeC {\cyri }\IeC {\cyrn }\IeC {\cyro }\IeC {\cyrv }}}{597}{}%{subsection.3.22.6}
\contentsline {subsection}{\numberline {22.7}\IeC {\CYRT }\IeC {\cyre }\IeC {\cyro }\IeC {\cyrr }\IeC {\cyri }\IeC {\cyrya } \IeC {\cyrv }\IeC {\cyre }\IeC {\cyrr }\IeC {\cyro }\IeC {\cyrya }\IeC {\cyrt }\IeC {\cyrn }\IeC {\cyro }\IeC {\cyrs }\IeC {\cyrt }\IeC {\cyre }\IeC {\cyrishrt } \IeC {\cyri }~\IeC {\cyrk }\IeC {\cyro }\IeC {\cyrm }\IeC {\cyrb }\IeC {\cyri }\IeC {\cyrn }\IeC {\cyra }\IeC {\cyrt }\IeC {\cyro }\IeC {\cyrr }\IeC {\cyrn }\IeC {\cyra }\IeC {\cyrya } \IeC {\cyrg }\IeC {\cyre }\IeC {\cyro }\IeC {\cyrm }\IeC {\cyre }\IeC {\cyrt }\IeC {\cyrr }\IeC {\cyri }\IeC {\cyrya } (4*). \emph {\IeC {\CYRA }.~\IeC {\CYRM }.~\IeC {\CYRR }\IeC {\cyra }\IeC {\cyrishrt }\IeC {\cyrg }\IeC {\cyro }\IeC {\cyrr }\IeC {\cyro }\IeC {\cyrd }\IeC {\cyrs }\IeC {\cyrk }\IeC {\cyri }\IeC {\cyrishrt }}}{616}{}%{subsection.3.22.7}
\contentsline {section}{\numberline {23}\IeC {\CYRP }\IeC {\cyre }\IeC {\cyrr }\IeC {\cyre }\IeC {\cyrs }\IeC {\cyrt }\IeC {\cyra }\IeC {\cyrn }\IeC {\cyro }\IeC {\cyrv }\IeC {\cyrk }\IeC {\cyri }. \emph {\IeC {\CYRA }.~\IeC {\CYRB }.~\IeC {\CYRS }\IeC {\cyrk }\IeC {\cyro }\IeC {\cyrp }\IeC {\cyre }\IeC {\cyrn }\IeC {\cyrk }\IeC {\cyro }\IeC {\cyrv }}}{620}{section.3.23}
\contentsline {subsection}{\ull{\numberline {23.1}\IeC {\CYRP }\IeC {\cyro }\IeC {\cyrr }\IeC {\cyrya }\IeC {\cyrd }\IeC {\cyro }\IeC {\cyrk }, \IeC {\cyrt }\IeC {\cyri }\IeC {\cyrp }, \IeC {\cyrs }\IeC {\cyro }\IeC {\cyrp }\IeC {\cyrr }\IeC {\cyrya }\IeC {\cyrzh }\IeC {\cyryo }\IeC {\cyrn }\IeC {\cyrn }\IeC {\cyro }\IeC {\cyrs }\IeC {\cyrt }\IeC {\cyrsftsn } (1)}}{621}{subsection.3.23.1}
\contentsline {subsection}{\numberline {23.2}\IeC {\CYRCH }\IeC {\cyryo }\IeC {\cyrt }\IeC {\cyrn }\IeC {\cyro }\IeC {\cyrs }\IeC {\cyrt }\IeC {\cyrsftsn } \IeC {\cyrp }\IeC {\cyre }\IeC {\cyrr }\IeC {\cyre }\IeC {\cyrs }\IeC {\cyrt }\IeC {\cyra }\IeC {\cyrn }\IeC {\cyro }\IeC {\cyrv }\IeC {\cyrk }\IeC {\cyri } (1)}{624}{}%{subsection.3.23.2}
\contentsline {subsection}{\ull{\numberline {23.3}\IeC {\CYRK }\IeC {\cyro }\IeC {\cyrm }\IeC {\cyrb }\IeC {\cyri }\IeC {\cyrn }\IeC {\cyra }\IeC {\cyrt }\IeC {\cyro }\IeC {\cyrr }\IeC {\cyri }\IeC {\cyrk }\IeC {\cyra } \IeC {\cyrk }\IeC {\cyrl }\IeC {\cyra }\IeC {\cyrs }\IeC {\cyrs }\IeC {\cyro }\IeC {\cyrv } \IeC {\cyrerev }\IeC {\cyrk }\IeC {\cyrv }\IeC {\cyri }\IeC {\cyrv }\IeC {\cyra }\IeC {\cyrl }\IeC {\cyre }\IeC {\cyrn }\IeC {\cyrt }\IeC {\cyrn }\IeC {\cyro }\IeC {\cyrs }\IeC {\cyrt }\IeC {\cyri } (2)}}{626}{subsection.3.23.3}
\contentsline {section}{\numberline {24}\IeC {\CYRG }\IeC {\cyrr }\IeC {\cyru }\IeC {\cyrp }\IeC {\cyrp }\IeC {\cyrery }. \emph {\IeC {\CYRV }.~\IeC {\CYRB }\IeC {\cyrr }\IeC {\cyra }\IeC {\cyrg }\IeC {\cyri }\IeC {\cyrn }}, \emph {\IeC {\CYRA }.~\IeC {\CYRK }\IeC {\cyrl }\IeC {\cyrya }\IeC {\cyrch }\IeC {\cyrk }\IeC {\cyro }}, \emph {\IeC {\CYRA }.~\IeC {\CYRS }\IeC {\cyrk }\IeC {\cyro }\IeC {\cyrp }\IeC {\cyre }\IeC {\cyrn }\IeC {\cyrk }\IeC {\cyro }\IeC {\cyrv }}}{632}{section.3.24}
\contentsline {subsection}{\ull{\numberline {24.1}\IeC {\CYRZ }\IeC {\cyra }\IeC {\cyrch }\IeC {\cyre }\IeC {\cyrm }, \IeC {\cyrd }\IeC {\cyrl }\IeC {\cyrya } \IeC {\cyrk }\IeC {\cyro }\IeC {\cyrg }\IeC {\cyro } \IeC {\cyri }~\IeC {\cyrk }\IeC {\cyra }\IeC {\cyrk } \IeC {\cyru }\IeC {\cyrs }\IeC {\cyrt }\IeC {\cyrr }\IeC {\cyro }\IeC {\cyre }\IeC {\cyrn } \IeC {\cyrerev }\IeC {\cyrt }\IeC {\cyro }\IeC {\cyrt } \IeC {\cyrp }\IeC {\cyra }\IeC {\cyrr }\IeC {\cyra }\IeC {\cyrg }\IeC {\cyrr }\IeC {\cyra }\IeC {\cyrf }}}{632}{subsection.3.24.1}
\contentsline {subsection}{\numberline {24.2}\IeC {\CYRK }\IeC {\cyra }\IeC {\cyrk } \IeC {\cyrp }\IeC {\cyrr }\IeC {\cyri }\IeC {\cyrd }\IeC {\cyru }\IeC {\cyrm }\IeC {\cyra }\IeC {\cyrt }\IeC {\cyrsftsn }}{635}{}%{subsection.3.24.2}
\contentsline {subsubsection}{\ull{\numberline {24.2.1}\IeC {\CYRP }\IeC {\cyro }\IeC {\cyrs }\IeC {\cyrt }\IeC {\cyra }\IeC {\cyrn }\IeC {\cyro }\IeC {\cyrv }\IeC {\cyrk }\IeC {\cyra } \IeC {\cyrz }\IeC {\cyra }\IeC {\cyrd }\IeC {\cyra }\IeC {\cyrch }\IeC {\cyri } (2)}}{635}{subsubsection.3.24.2.1}
\contentsline {subsubsection}{\numberline {24.2.2}\IeC {\CYRP }\IeC {\cyrr }\IeC {\cyri }\IeC {\cyrm }\IeC {\cyre }\IeC {\cyrr }\IeC {\cyrery } \IeC {\cyrg }\IeC {\cyrr }\IeC {\cyru }\IeC {\cyrp }\IeC {\cyrp } (2)}{636}{}%{subsubsection.3.24.2.2}
\contentsline {subsubsection}{\numberline {24.2.3}\IeC {\CYRD }\IeC {\cyro }\IeC {\cyrk }\IeC {\cyra }\IeC {\cyrzh }\IeC {\cyre }\IeC {\cyrm } \IeC {\cyri }~\IeC {\cyrp }\IeC {\cyrr }\IeC {\cyri }\IeC {\cyrm }\IeC {\cyre }\IeC {\cyrn }\IeC {\cyri }\IeC {\cyrm } \IeC {\cyrt }\IeC {\cyre }\IeC {\cyro }\IeC {\cyrr }\IeC {\cyre }\IeC {\cyrm }\IeC {\cyru } \IeC {\CYRL }\IeC {\cyra }\IeC {\cyrg }\IeC {\cyrr }\IeC {\cyra }\IeC {\cyrn }\IeC {\cyrzh }\IeC {\cyra } (2)}{638}{}%{subsubsection.3.24.2.3}
\contentsline {subsubsection}{\numberline {24.2.4}\IeC {\CYRP }\IeC {\cyrr }\IeC {\cyri }\IeC {\cyrm }\IeC {\cyre }\IeC {\cyrn }\IeC {\cyri }\IeC {\cyrm } \IeC {\cyrs }\IeC {\cyro }\IeC {\cyrp }\IeC {\cyrr }\IeC {\cyrya }\IeC {\cyrzh }\IeC {\cyre }\IeC {\cyrn }\IeC {\cyri }\IeC {\cyre } (3)}{640}{}%{subsubsection.3.24.2.4}
\contentsline {subsubsection}{\numberline {24.2.5}\IeC {\CYRM }\IeC {\cyra }\IeC {\cyrk }\IeC {\cyrs }\IeC {\cyri }\IeC {\cyrm }\IeC {\cyra }\IeC {\cyrl }\IeC {\cyrsftsn }\IeC {\cyrn }\IeC {\cyrery }\IeC {\cyre } \IeC {\cyrp }\IeC {\cyro }\IeC {\cyrd }\IeC {\cyrg }\IeC {\cyrr }\IeC {\cyru }\IeC {\cyrp }\IeC {\cyrp }\IeC {\cyrery } \IeC {\cyri }~\IeC {\cyrc }\IeC {\cyre }\IeC {\cyrn }\IeC {\cyrt }\IeC {\cyrr } (4*)}{641}{}%{subsubsection.3.24.2.5}
\contentsline {subsection}{\numberline {24.3}\IeC {\CYRI }\IeC {\cyrt }\IeC {\cyro }\IeC {\cyrg }: \IeC {\cyrf }\IeC {\cyro }\IeC {\cyrr }\IeC {\cyrm }\IeC {\cyru }\IeC {\cyrl }\IeC {\cyri }\IeC {\cyrr }\IeC {\cyro }\IeC {\cyrv }\IeC {\cyrk }\IeC {\cyra } \IeC {\cyri }~\IeC {\cyrd }\IeC {\cyro }\IeC {\cyrk }\IeC {\cyra }\IeC {\cyrz }\IeC {\cyra }\IeC {\cyrt }\IeC {\cyre }\IeC {\cyrl }\IeC {\cyrsftsn }\IeC {\cyrs }\IeC {\cyrt }\IeC {\cyrv }\IeC {\cyro }}{647}{}%{subsection.3.24.3}
\contentsline {subsubsection}{\numberline {24.3.1}\IeC {\CYRF }\IeC {\cyro }\IeC {\cyrr }\IeC {\cyrm }\IeC {\cyru }\IeC {\cyrl }\IeC {\cyri }\IeC {\cyrr }\IeC {\cyro }\IeC {\cyrv }\IeC {\cyrk }\IeC {\cyra } \IeC {\cyro }\IeC {\cyrs }\IeC {\cyrn }\IeC {\cyro }\IeC {\cyrv }\IeC {\cyrn }\IeC {\cyro }\IeC {\cyrg }\IeC {\cyro } \IeC {\cyrr }\IeC {\cyre }\IeC {\cyrz }\IeC {\cyru }\IeC {\cyrl }\IeC {\cyrsftsn }\IeC {\cyrt }\IeC {\cyra }\IeC {\cyrt }\IeC {\cyra } (2)}{647}{}%{subsubsection.3.24.3.1}
\contentsline {subsubsection}{\numberline {24.3.2}\IeC {\CYRD }\IeC {\cyro }\IeC {\cyrk }\IeC {\cyra }\IeC {\cyrz }\IeC {\cyra }\IeC {\cyrt }\IeC {\cyre }\IeC {\cyrl }\IeC {\cyrsftsn }\IeC {\cyrs }\IeC {\cyrt }\IeC {\cyrv }\IeC {\cyro } \IeC {\cyrch }\IeC {\cyra }\IeC {\cyrs }\IeC {\cyrt }\IeC {\cyri } <<\IeC {\cyrt }\IeC {\cyro }\IeC {\cyrl }\IeC {\cyrsftsn }\IeC {\cyrk }\IeC {\cyro } \IeC {\cyrt }\IeC {\cyro }\IeC {\cyrg }\IeC {\cyrd }\IeC {\cyra }>> (3*)}{648}{}%{subsubsection.3.24.3.2}
\contentsline {subsubsection}{\numberline {24.3.3}\IeC {\CYRD }\IeC {\cyro }\IeC {\cyrk }\IeC {\cyra }\IeC {\cyrz }\IeC {\cyra }\IeC {\cyrt }\IeC {\cyre }\IeC {\cyrl }\IeC {\cyrsftsn }\IeC {\cyrs }\IeC {\cyrt }\IeC {\cyrv }\IeC {\cyro } \IeC {\cyrch }\IeC {\cyra }\IeC {\cyrs }\IeC {\cyrt }\IeC {\cyri } <<\IeC {\cyrt }\IeC {\cyro }\IeC {\cyrg }\IeC {\cyrd }\IeC {\cyra }>> (4*)}{649}{}%{subsubsection.3.24.3.3}
\contentsline {section}{\numberline {25}\IeC {\CYRK }\IeC {\cyro }\IeC {\cyrm }\IeC {\cyrb }\IeC {\cyri }\IeC {\cyrn }\IeC {\cyra }\IeC {\cyrt }\IeC {\cyro }\IeC {\cyrr }\IeC {\cyrn }\IeC {\cyra }\IeC {\cyrya } \IeC {\cyrg }\IeC {\cyre }\IeC {\cyro }\IeC {\cyrm }\IeC {\cyre }\IeC {\cyrt }\IeC {\cyrr }\IeC {\cyri }\IeC {\cyrya }}{654}{section.3.25}
\contentsline {subsection}{\numberline {25.1}\IeC {\CYRO } \IeC {\cyrk }\IeC {\cyro }\IeC {\cyrv }\IeC {\cyrr }\IeC {\cyro }\IeC {\cyrv }\IeC {\cyrery }\IeC {\cyrh } \IeC {\cyrd }\IeC {\cyro }\IeC {\cyrr }\IeC {\cyro }\IeC {\cyrzh }\IeC {\cyrk }\IeC {\cyra }\IeC {\cyrh } \IeC {\cyri }~\IeC {\cyrs }\IeC {\cyra }\IeC {\cyrl }\IeC {\cyrf }\IeC {\cyre }\IeC {\cyrt }\IeC {\cyrk }\IeC {\cyra }\IeC {\cyrh } (2). \emph {\IeC {\CYRP }.~\IeC {\CYRA }.~\IeC {\CYRK }\IeC {\cyro }\IeC {\cyrzh }\IeC {\cyre }\IeC {\cyrv }\IeC {\cyrn }\IeC {\cyri }\IeC {\cyrk }\IeC {\cyro }\IeC {\cyrv }}}{654}{}%{subsection.3.25.1}
\contentsline {subsection}{\numberline {25.2}\IeC {\CYRT }\IeC {\cyre }\IeC {\cyro }\IeC {\cyrr }\IeC {\cyre }\IeC {\cyrm }\IeC {\cyra } \IeC {\CYRH }\IeC {\cyre }\IeC {\cyrl }\IeC {\cyrl }\IeC {\cyri } (2). \emph {\IeC {\CYRA }.~\IeC {\CYRV }.~\IeC {\CYRA }\IeC {\cyrk }\IeC {\cyro }\IeC {\cyrp }\IeC {\cyrya }\IeC {\cyrn }}}{662}{}%{subsection.3.25.2}
\contentsline {subsection}{\numberline {25.3}\IeC {\CYRM }\IeC {\cyrn }\IeC {\cyro }\IeC {\cyrg }\IeC {\cyro }\IeC {\cyru }\IeC {\cyrg }\IeC {\cyro }\IeC {\cyrl }\IeC {\cyrsftsn }\IeC {\cyrn }\IeC {\cyri }\IeC {\cyrk }\IeC {\cyri } \IeC {\cyrn }\IeC {\cyra } \IeC {\cyrk }\IeC {\cyrl }\IeC {\cyre }\IeC {\cyrt }\IeC {\cyrch }\IeC {\cyra }\IeC {\cyrt }\IeC {\cyro }\IeC {\cyrishrt } \IeC {\cyrb }\IeC {\cyru }\IeC {\cyrm }\IeC {\cyra }\IeC {\cyrg }\IeC {\cyre } (2). \emph {\IeC {\CYRV }.~\IeC {\CYRV }.~\IeC {\CYRP }\IeC {\cyrr }\IeC {\cyra }\IeC {\cyrs }\IeC {\cyro }\IeC {\cyrl }\IeC {\cyro }\IeC {\cyrv }}, \emph {\IeC {\CYRM }.~\IeC {\CYRB }.~\IeC {\CYRS }\IeC {\cyrk }\IeC {\cyro }\IeC {\cyrp }\IeC {\cyre }\IeC {\cyrn }\IeC {\cyrk }\IeC {\cyro }\IeC {\cyrv }}}{665}{}%{subsection.3.25.3}
\contentsline {subsection}{\numberline {25.4}\IeC {\CYRP }\IeC {\cyrr }\IeC {\cyri }\IeC {\cyrn }\IeC {\cyrc }\IeC {\cyri }\IeC {\cyrp } \IeC {\CYRD }\IeC {\cyri }\IeC {\cyrr }\IeC {\cyri }\IeC {\cyrh }\IeC {\cyrl }\IeC {\cyre } \IeC {\cyrn }\IeC {\cyra } \IeC {\cyrp }\IeC {\cyrr }\IeC {\cyrya }\IeC {\cyrm }\IeC {\cyro }\IeC {\cyrishrt } (3). \emph {\IeC {\CYRA }.~\IeC {\CYRYA }.~\IeC {\CYRK }\IeC {\cyra }\IeC {\cyrn }\IeC {\cyre }\IeC {\cyrl }\IeC {\cyrsftsn }-\IeC {\CYRB }\IeC {\cyre }\IeC {\cyrl }\IeC {\cyro }\IeC {\cyrv }}}{681}{}%{subsection.3.25.4}
\contentsline {subsection}{\numberline {25.5}\IeC {\CYRP }\IeC {\cyrr }\IeC {\cyri }\IeC {\cyrn }\IeC {\cyrc }\IeC {\cyri }\IeC {\cyrp } \IeC {\CYRD }\IeC {\cyri }\IeC {\cyrr }\IeC {\cyri }\IeC {\cyrh }\IeC {\cyrl }\IeC {\cyre } \IeC {\cyri }~\IeC {\cyre }\IeC {\cyrg }\IeC {\cyro } \IeC {\cyrp }\IeC {\cyrr }\IeC {\cyri }\IeC {\cyrm }\IeC {\cyre }\IeC {\cyrn }\IeC {\cyre }\IeC {\cyrn }\IeC {\cyri }\IeC {\cyrya } \IeC {\cyrv }~\IeC {\cyrg }\IeC {\cyre }\IeC {\cyro }\IeC {\cyrm }\IeC {\cyre }\IeC {\cyrt }\IeC {\cyrr }\IeC {\cyri }\IeC {\cyri }{{}{}} {} (3). \emph {\IeC {\CYRI }.~\IeC {\CYRV }.~\IeC {\CYRA }\IeC {\cyrr }\IeC {\cyrzh }\IeC {\cyra }\IeC {\cyrn }\IeC {\cyrc }\IeC {\cyre }\IeC {\cyrv }}}{682}{}%{subsection.3.25.5}
\contentsline {subsection}{\numberline {25.6}\IeC {\CYRF }\IeC {\cyra }\IeC {\cyrz }\IeC {\cyro }\IeC {\cyrv }\IeC {\cyrery }\IeC {\cyre } \IeC {\cyrp }\IeC {\cyrr }\IeC {\cyro }\IeC {\cyrs }\IeC {\cyrt }\IeC {\cyrr }\IeC {\cyra }\IeC {\cyrn }\IeC {\cyrs }\IeC {\cyrt }\IeC {\cyrv }\IeC {\cyra } (3). \emph {\IeC {\CYRA }.~\IeC {\CYRYA }.~\IeC {\CYRK }\IeC {\cyra }\IeC {\cyrn }\IeC {\cyre }\IeC {\cyrl }\IeC {\cyrsftsn }-\IeC {\CYRB }\IeC {\cyre }\IeC {\cyrl }\IeC {\cyro }\IeC {\cyrv }}}{690}{}%{subsection.3.25.6}
\contentsline {subsection}{\numberline {25.7}\IeC {\CYRL }\IeC {\cyri }\IeC {\cyrn }\IeC {\cyre }\IeC {\cyrishrt }\IeC {\cyrn }\IeC {\cyro }\IeC {\cyre } \IeC {\cyrv }\IeC {\cyra }\IeC {\cyrr }\IeC {\cyrsftsn }\IeC {\cyri }\IeC {\cyrr }\IeC {\cyro }\IeC {\cyrv }\IeC {\cyra }\IeC {\cyrn }\IeC {\cyri }\IeC {\cyre } (3). \emph {\IeC {\CYRA }.~\IeC {\CYRYA }.~\IeC {\CYRK }\IeC {\cyra }\IeC {\cyrn }\IeC {\cyre }\IeC {\cyrl }\IeC {\cyrsftsn }-\IeC {\CYRB }\IeC {\cyre }\IeC {\cyrl }\IeC {\cyro }\IeC {\cyrv }}}{692}{}%{subsection.3.25.7}
\contentsline {subsection}{\ull{\numberline {25.8}\IeC {\CYRS }\IeC {\cyro }\IeC {\cyrb }\IeC {\cyre }\IeC {\cyrr }\IeC {\cyri } \IeC {\cyrk }\IeC {\cyrv }\IeC {\cyra }\IeC {\cyrd }\IeC {\cyrr }\IeC {\cyra }\IeC {\cyrt } (3*). \emph {\IeC {\CYRM }.~\IeC {\CYRB }.~\IeC {\CYRS }\IeC {\cyrk }\IeC {\cyro }\IeC {\cyrp }\IeC {\cyre }\IeC {\cyrn }\IeC {\cyrk }\IeC {\cyro }\IeC {\cyrv }}, \emph {\IeC {\CYRO }.~\IeC {\CYRA }.~\IeC {\CYRM }\IeC {\cyra }\IeC {\cyrl }\IeC {\cyri }\IeC {\cyrn }\IeC {\cyro }\IeC {\cyrv }\IeC {\cyrs }\IeC {\cyrk }\IeC {\cyra }\IeC {\cyrya }}, \emph {\IeC {\CYRS }.~\IeC {\CYRA }.~\IeC {\CYRD }\IeC {\cyro }\IeC {\cyrr }\IeC {\cyri }\IeC {\cyrch }\IeC {\cyre }\IeC {\cyrn }\IeC {\cyrk }\IeC {\cyro }}, \emph {\IeC {\CYRF }.~\IeC {\CYRA }.~\IeC {\CYRSH }\IeC {\cyra }\IeC {\cyrr }\IeC {\cyro }\IeC {\cyrv }}}}{694}{subsection.3.25.8}
\contentsline {subsection}{\numberline {25.9}\IeC {\CYRM }\IeC {\cyro }\IeC {\cyrzh }\IeC {\cyrn }\IeC {\cyro } \IeC {\cyrl }\IeC {\cyri } \IeC {\cyri }\IeC {\cyrz } \IeC {\cyrt }\IeC {\cyre }\IeC {\cyrt }\IeC {\cyrr }\IeC {\cyra }\IeC {\cyrerev }\IeC {\cyrd }\IeC {\cyrr }\IeC {\cyra } \IeC {\cyrs }\IeC {\cyrd }\IeC {\cyre }\IeC {\cyrl }\IeC {\cyra }\IeC {\cyrt }\IeC {\cyrsftsn } \IeC {\cyrk }\IeC {\cyru }\IeC {\cyrb }?{{}{}} {} (3). \emph {\IeC {\CYRM }.~\IeC {\CYRV }.~\IeC {\CYRP }\IeC {\cyrr }\IeC {\cyra }\IeC {\cyrs }\IeC {\cyro }\IeC {\cyrl }\IeC {\cyro }\IeC {\cyrv }}, \emph {\IeC {\CYRM }.~\IeC {\CYRB }.~\IeC {\CYRS }\IeC {\cyrk }\IeC {\cyro }\IeC {\cyrp }\IeC {\cyre }\IeC {\cyrn }\IeC {\cyrk }\IeC {\cyro }\IeC {\cyrv }}}{709}{}%{subsection.3.25.9}
\contentsline {chapter}{\numberline {4}\IeC {\CYRO } \IeC {\cyrp }\IeC {\cyrr }\IeC {\cyre }\IeC {\cyrp }\IeC {\cyro }\IeC {\cyrd }\IeC {\cyra }\IeC {\cyrv }\IeC {\cyra }\IeC {\cyrn }\IeC {\cyri }\IeC {\cyri }. \emph {\IeC {\CYRA }.~\IeC {\CYRB }.~\IeC {\CYRS }\IeC {\cyrk }\IeC {\cyro }\IeC {\cyrp }\IeC {\cyre }\IeC {\cyrn }\IeC {\cyrk }\IeC {\cyro }\IeC {\cyrv }}}{724}{chapter.4}
\contentsline {section}{\ull{\numberline {26}\IeC {\CYRO }\IeC {\cyrl }\IeC {\cyri }\IeC {\cyrm }\IeC {\cyrp }\IeC {\cyri }\IeC {\cyra }\IeC {\cyrd }\IeC {\cyrery } \IeC {\cyri }~\IeC {\cyrm }\IeC {\cyra }\IeC {\cyrt }\IeC {\cyre }\IeC {\cyrm }\IeC {\cyra }\IeC {\cyrt }\IeC {\cyri }\IeC {\cyrk }\IeC {\cyra }}}{724}{section.4.26}
\contentsline {section}{\ull{\numberline {27}\IeC {\CYRN }\IeC {\cyra }\IeC {\cyrch }\IeC {\cyri }\IeC {\cyrn }\IeC {\cyra }\IeC {\cyrt }\IeC {\cyrsftsn } \IeC {\cyrs }~\IeC {\cyrya }\IeC {\cyrz }\IeC {\cyrery }\IeC {\cyrk }\IeC {\cyra } \IeC {\cyri }\IeC {\cyrl }\IeC {\cyri } \IeC {\cyrs }\IeC {\cyro }\IeC {\cyrd }\IeC {\cyre }\IeC {\cyrr }\IeC {\cyrzh }\IeC {\cyra }\IeC {\cyrn }\IeC {\cyri }\IeC {\cyrya }?}}{726}{section.4.27}
\contentsline {section}{\numberline {28}\IeC {\CYRO } \IeC {\cyrn }\IeC {\cyre }\IeC {\cyro }\IeC {\cyrb }\IeC {\cyrh }\IeC {\cyro }\IeC {\cyrd }\IeC {\cyri }\IeC {\cyrm }\IeC {\cyro }\IeC {\cyrs }\IeC {\cyrt }\IeC {\cyri } \IeC {\cyrm }\IeC {\cyro }\IeC {\cyrt }\IeC {\cyri }\IeC {\cyrv }\IeC {\cyri }\IeC {\cyrr }\IeC {\cyro }\IeC {\cyrv }\IeC {\cyro }\IeC {\cyrk }}{730}{}%{section.4.28}
\contentsline {subsection}{\numberline {28.1}<<\IeC {\CYRZ }\IeC {\cyra }>> \IeC {\cyri }~<<\IeC {\cyrp }\IeC {\cyrr }\IeC {\cyro }\IeC {\cyrt }\IeC {\cyri }\IeC {\cyrv }>> \IeC {\cyrm }\IeC {\cyro }\IeC {\cyrt }\IeC {\cyri }\IeC {\cyrv }\IeC {\cyri }\IeC {\cyrr }\IeC {\cyro }\IeC {\cyrv }\IeC {\cyro }\IeC {\cyrk }}{731}{}%{subsection.4.28.1}
\contentsline {subsection}{\numberline {28.2}\IeC {\CYRO } \IeC {\cyrm }\IeC {\cyro }\IeC {\cyrt }\IeC {\cyri }\IeC {\cyrv }\IeC {\cyri }\IeC {\cyrr }\IeC {\cyro }\IeC {\cyrv }\IeC {\cyrk }\IeC {\cyra }\IeC {\cyrh } \IeC {\cyrt }\IeC {\cyre }\IeC {\cyro }\IeC {\cyrr }\IeC {\cyri }\IeC {\cyri } \IeC {\CYRG }\IeC {\cyra }\IeC {\cyrl }\IeC {\cyru }\IeC {\cyra }}{733}{}%{subsection.4.28.2}
\contentsline {subsection}{\numberline {28.3}\IeC {\CYRP }\IeC {\cyro }\IeC {\cyrch }\IeC {\cyre }\IeC {\cyrm }\IeC {\cyru } \IeC {\cyrn }\IeC {\cyre } \IeC {\cyrp }\IeC {\cyrr }\IeC {\cyri }\IeC {\cyrn }\IeC {\cyri }\IeC {\cyrm }\IeC {\cyra }\IeC {\cyre }\IeC {\cyrt }\IeC {\cyrs }\IeC {\cyrya } \IeC {\cyrm }\IeC {\cyro }\IeC {\cyrt }\IeC {\cyri }\IeC {\cyrv }\IeC {\cyri }\IeC {\cyrr }\IeC {\cyro }\IeC {\cyrv }\IeC {\cyra }\IeC {\cyrn }\IeC {\cyrn }\IeC {\cyro }\IeC {\cyre } \IeC {\cyri }\IeC {\cyrz }\IeC {\cyrl }\IeC {\cyro }\IeC {\cyrzh }\IeC {\cyre }\IeC {\cyrn }\IeC {\cyri }\IeC {\cyre }?}{734}{}%{subsection.4.28.3}
\contentsline {subsubsection}{\numberline {28.3.1}\IeC {\CYRO }\IeC {\cyrt }\IeC {\cyrz }\IeC {\cyrery }\IeC {\cyrv }}{734}{}%{subsubsection.4.28.3.1}
\contentsline {subsubsection}{\numberline {28.3.2}\IeC {\CYRK }\IeC {\cyro }\IeC {\cyrm }\IeC {\cyrm }\IeC {\cyre }\IeC {\cyrn }\IeC {\cyrt }\IeC {\cyra }\IeC {\cyrr }\IeC {\cyri }\IeC {\cyri } \IeC {\cyrk }~\IeC {\cyro }\IeC {\cyrt }\IeC {\cyrz }\IeC {\cyrery }\IeC {\cyrv }\IeC {\cyru }}{735}{}%{subsubsection.4.28.3.2}
\contentsline {subsubsection}{\numberline {28.3.3}\IeC {\CYRD }\IeC {\cyrr }\IeC {\cyru }\IeC {\cyrg }\IeC {\cyri }\IeC {\cyre } \IeC {\cyrv }\IeC {\cyrery }\IeC {\cyrs }\IeC {\cyrk }\IeC {\cyra }\IeC {\cyrz }\IeC {\cyrery }\IeC {\cyrv }\IeC {\cyra }\IeC {\cyrn }\IeC {\cyri }\IeC {\cyrya }}{738}{}%{subsubsection.4.28.3.3}
\contentsline {section}{\numberline {29}\IeC {\CYRK }\IeC {\cyrr }\IeC {\cyru }\IeC {\cyrzh }\IeC {\cyrk }\IeC {\cyri } \IeC {\cyri }~\IeC {\cyro }\IeC {\cyrl }\IeC {\cyri }\IeC {\cyrm }\IeC {\cyrp }\IeC {\cyri }\IeC {\cyra }\IeC {\cyrd }\IeC {\cyrery } \IeC {\cyrk }\IeC {\cyra }\IeC {\cyrk } \IeC {\cyrp }\IeC {\cyru }\IeC {\cyrt }\IeC {\cyrsftsn } \IeC {\cyrv }~\IeC {\cyrm }\IeC {\cyra }\IeC {\cyrt }\IeC {\cyre }\IeC {\cyrm }\IeC {\cyra }\IeC {\cyrt }\IeC {\cyri }\IeC {\cyrk }\IeC {\cyru } \IeC {\cyri }~\IeC {\cyrk }\IeC {\cyra }\IeC {\cyrk } \IeC {\cyrs }\IeC {\cyrp }\IeC {\cyro }\IeC {\cyrr }\IeC {\cyrt }. \emph {\IeC {\CYRA }.~\IeC {\CYRYA }.~\IeC {\CYRK }\IeC {\cyra }\IeC {\cyrn }\IeC {\cyre }\IeC {\cyrl }\IeC {\cyrsftsn }-\IeC {\CYRB }\IeC {\cyre }\IeC {\cyrl }\IeC {\cyro }\IeC {\cyrv }}, \emph {\IeC {\CYRA }.\tmspace +\thinmuskip {.1667em}\IeC {\CYRI }.\tmspace +\thinmuskip {.1667em}\IeC {\CYRB }\IeC {\cyru }\IeC {\cyrf }\IeC {\cyre }\IeC {\cyrt }\IeC {\cyro }\IeC {\cyrv }}}{745}{}%{section.4.29}
\contentsline {subsection}{\numberline {29.1}\IeC {\CYRV }\IeC {\cyrv }\IeC {\cyre }\IeC {\cyrd }\IeC {\cyre }\IeC {\cyrn }\IeC {\cyri }\IeC {\cyre }}{745}{}%{subsection.4.29.1}
\contentsline {subsection}{\numberline {29.2}\IeC {\CYRS }\IeC {\cyrp }\IeC {\cyro }\IeC {\cyrr }\IeC {\cyrt }\IeC {\cyri }\IeC {\cyrv }\IeC {\cyrn }\IeC {\cyrery }\IeC {\cyrishrt } \IeC {\cyrp }\IeC {\cyro }\IeC {\cyrd }\IeC {\cyrh }\IeC {\cyro }\IeC {\cyrd }}{745}{}%{subsection.4.29.2}
\contentsline {subsection}{\numberline {29.3}\IeC {\CYRO }\IeC {\cyrl }\IeC {\cyri }\IeC {\cyrm }\IeC {\cyrp }\IeC {\cyri }\IeC {\cyra }\IeC {\cyrd }\IeC {\cyra } \IeC {\cyrk }\IeC {\cyra }\IeC {\cyrk } \IeC {\cyrp }\IeC {\cyru }\IeC {\cyrt }\IeC {\cyrsftsn } \IeC {\cyrv }~\IeC {\cyrm }\IeC {\cyra }\IeC {\cyrt }\IeC {\cyre }\IeC {\cyrm }\IeC {\cyra }\IeC {\cyrt }\IeC {\cyri }\IeC {\cyrk }\IeC {\cyru }}{747}{}%{subsection.4.29.3}

\newpage

%%%!!!
%\chapter*{Введение}
%\addcontentsline{toc}{chapt}{Введение}
%\markboth{Введение}{Введение}
%\clearpage\thispagestyle{empty}\vspace*{10pt}

%\Ol{1. Фамилии с ё или без? Федоров, Золотарев, Соловьев, Чеботарёв, Терешин, Личев, Шевелев, Кузьмин, Толмачев}

\section{От редакторов}\label{0int}
%\addcontentsline{toc}{subsection}{От редакторов}

\medskip

\subsection{Зачем и~для кого эта книга}

Глубокое понимание математики полезно и~математику, и~профессионалу в~наукоёмкой отрасли. В~частности, <<профессия>> в~названии этой книги не обязательно означает профессию математика.

Эта книга предназначена для старшеклассников и~младшекурсников (в~частности, ориентированных на олимпиады). См.~подробнее \S\,\ref{oim} <<Олимпиады и~математика>>. Книгу можно использовать как для самостоятельных занятий, так и~для преподавания.

Книга содержит наиболее стандартный <<базовый>> материал (впрочем, частично, скорее, для повторения, чем для первоначального изучения). Основное содержание книги составляет более сложный материал. Некоторые темы малоизвестны в~традиции математических кружков, но полезны как для математического образования, так и~для подготовки к~олимпиадам.

Книга основана на занятиях, проведённых авторами в~разное время в~школе им.~А.~Н.~Колмогорова (СУНЦ МГУ), школе \No~1543 г.~Москвы, летней школе <<Современная математика>>, Кировской и~Костромской летних математических школах, Московской выездной олимпиадной школе, в~кружках <<Математический семинар>> и~<<Олимпиады и~математика>>, на летней конференции Турнира городов, при подготовке команды России к~международной математической олимпиаде,
в~системе дистанционного обучения математике МИОО, а~также в~Независимом московском университете и~на математическом факультете Высшей школы экономики.

Книга доступна уже старшеклассникам, интересующимся математикой\footnote{Часть материала в~некоторых кружках и~летних школах изучается теми, кто только знакомится с~математикой (например, 6-классниками). Однако приводимое изложение рассчитано на читателя, уже имеющего хотя бы минимальную математическую культуру. Заниматься с~6-классниками нужно по-другому, см., например, \cite{GIF}.}. Приводятся почти все определения, не входящие в~школьную программу. Если где-то нужны дополнительные сведения, то приводятся ссылки.

При этом многие темы трудны, если изучать их <<с нуля>>. Однако \emph{последовательность изложения} помогает преодолевать трудности. В~то же время многие темы \emph{независимы} друг от друга. См.~подробнее п.\;\ref{s:intstr} <<Как устроена книга>>.

\subsection{Изучение путём решения и~обсуждения задач}\label{s:intpro}

Мы следуем традиции изучения материала в~виде решения и~обсуждения задач. Эти задачи подобраны так, что в~процессе их решения читатель (точнее, решатель) освоит основы важных теорий "--- как классических, так и~современных. Основные идеи демонстрируются по одной и~на <<олимпиадных>> примерах, т.\,е. на простейших частных случаях, свободных от технических деталей. Этим мы показываем, \emph{как можно придумать} эти теории. См.~подробнее \S\,\ref{oim} <<Олимпиады и~математика>>.

Обучение путём решения задач не только характерно для серьёзного изучения математики, но и~продолжает древнюю культурную традицию. Например, послушники дзенских монастырей обучаются, размышляя над загадками, данными им наставниками. Впрочем, эти загадки являются скорее парадоксами, а~не задачами.
См.~подробнее \cite{Su}; ср. \cite[с.\,26--33]{Pl}. А вот некоторые <<математические>> примеры: \cite{Ar, BS, GDI, KK, Pr07-1, PS, SC, Sk09, Va, Zv}; кое-где не только приведены задачи, но и~изложены \emph{принципы отбора} удачных задач.

Учиться, решая задачи, трудно. В частности, потому, что такое обучение обычно
не создаёт \emph{иллюзию} понимания. Однако усилия сполна вознаграждаются глубоким пониманием материала "--- в~первую очередь, умением проводить аналогичные (и~даже не очень аналогичные) рассуждения. Кое-где вслед за великими математиками в~процессе изучения интересных задач читатель увидит, как естественно возникают важные понятия и~теории. Надеемся, это поможет ему совершить собственные настолько же полезные открытия (не обязательно в~математике)!

Для решения задач достаточно понимания их условий. Другие знания и~теории не нужны. (Впрочем, такие знания и~теории как раз появляются при решении подобранных задач.) Но может потребоваться владение другими частями книги, что
отражено в~подсказках и~указаниях.

К важнейшим задачам приводятся подсказки, указания, решения и~ответы. Они расположены в~конце каждого пункта. Однако к~ним стоит обращаться после прорешивания каждой задачи.

Если задача выделена словом <<теорема>> (<<лемма>>, <<следствие>> и~т.\,д.) и~жирным шрифтом, то её утверждение важное.

Как правило, мы приводим \emph{формулировку} красивого или важного утверждения (в~виде задачи) перед его \emph{доказательством}. В~таких случаях для доказательства утверждения могут потребоваться следующие задачи. Это всегда явно оговаривается в~подсказках, а~иногда и~прямо в~тексте. Поэтому если некоторая задача не получается, то читайте дальше. (На занятии задача-подсказка выдаётся только тогда, когда ученик немного подумал над самой задачей.) Такой процесс обучения полезен, поскольку моделирует реальную исследовательскую ситуацию. См.~подробнее \S\,28%\ref{s:motiv1} 
<<О необходимости мотивировок>>.

Всё это "--- попытка продемонстрировать занятие в~виде \emph{диалога}, основанного на решении и~обсуждении задач. Подробнее см.~\cite{KK15}.

\subsection{Как устроена книга}\label{s:intstr}

Книгу не обязательно изучать подряд. Читатель может выбрать удобную ему последовательность изучения (или вовсе опустить некоторые пункты) на основании приводимого плана. Для занятия кружка можно использовать любой пункт (или подпункт) книги.

Книга разбита на главы, параграфы и~пункты (некоторые пункты разбиты на подпункты). Структура параграфов приблизительно описана в~их начале. Если в~задаче используется материал другого пункта, то можно либо игнорировать эту задачу, либо посмотреть то место, на которое приводится ссылка. Это даёт большую свободу читателю при изучении книги, но одновременно может требовать его внимательности.

Пункты внутри каждого параграфа расположены примерно в~порядке возрастания сложности материала. Цифры в~скобках после названия пункта означают его <<относительный уровень>>: 1 "--- самый простой, 4 "--- самый сложный. Первые пункты (не отмеченные звёздочкой) являются базовыми; если не указано противное, с~них можно начать изучение главы. А~к~остальным пунктам (отмеченным звёздочкой) можно возвращаться потом; если не указано противное, то они независимы друг от друга. При изучении полезно \emph{возвращаться} к~пройденному материалу, но на новом уровне. Поэтому разные пункты одного параграфа можно изучать \emph{не подряд}, а~с~перерывами на другие темы.

Обозначения, используемые в~разных главах книги, приведены в~конце введения.
Понятия и~обозначения, используемые в~некоторой главе, вводятся в~начале главы.

Последняя глава составлена из заметок об общих принципах преподавания, адресованных прежде всего учителям. Возможно, заметки окажутся полезными и~ученикам.

В конце книги есть предметный указатель. Жирным шрифтом выделены номера страниц, на которых приводятся \emph{формальные определения} понятий.

Обновляемая электронная версия части книги, выложенная с~разрешения издательства:\\
\url{http://www.mccme.ru/circles/oim/materials/sturm.pdf}.

%%%!!!
\subsection{Напутствие. \emph{А.~Я.~Канель-Белов}}
%\addcontentsline{toc}{section}{Напутствие. А.~Я.~Канель-Белов}

Для успешного решения задач математических олимпиад высшего уровня необходимы в~первую очередь общеукрепляющие средства: хорошая проработка алгебры (культура алгебраических преобразований), проработка школьной геометрии. Задачи этих олимпиад (кроме первых задач) практически всегда предполагают смешанный сценарий решения; редки задачи на применение некоторого метода или идеи в~чистом виде. Решению таких <<смешанных>> задач должна предшествовать работа с~ключевыми задачами, в~которых идеи работают в~чистом виде. См., например, литературу к~п.\;\ref{s:intpro} или настоящий сборник.

\subsection{О литературе и~источниках}

В конце каждого параграфа приводится литература, относящаяся ко всему параграфу, и~отдельно литература по каждому пункту. Ссылка на книгу~\cite{GKP}, относящуюся и~к~комбинаторике, и~к~алгебре, приведена в~этом параграфе. Мы старались указать не только литературу, использованную при подготовке конкретного материала, но также и~жемчужины научно-популярного жанра на изучаемые темы. Мы надеемся, что наш список литературы, хотя бы в первом приближении, сможет стать путеводителем в море научно-популярной литературы по математике. Однако в~него наверняка не вошли многие замечательные материалы, ввиду необъятности их количества. Важно, что обращение к~литературе не нужно для решения задач, если явно не указано обратное.

В книге использованы материалы из сборника~\cite{ZPS}, см.~также~\cite{U}. При этом большинство задач существенно переработано и~добавлены новые. Многие материалы по комбинаторике, в~том числе <<базовые>>, перемещены в~\cite{GDI}.
Удалены <<разные задачи>> и~материалы, перепечатанные из других источников.
По последним приводятся ссылки.

Многие задачи не оригинальны, но первоисточник (даже если его можно установить) обычно не указывается.
Ссылка после условия задачи указывает источник, из которого взяты задача.
Эти ссылки приведены для того, чтобы читатель смог сравнить своё решение с~приведённым там.
Если мы знали, что пересечение какого-то пункта с~каким-то источником велико, то упоминали об этом.

%обозначают номер приведённой здесь задачи в~книге \cite{prasolov}
%(кроме материалов А.~Д.~Блинкова; разумеется, остальные задачи также не являются оригинальными).

Мы не даём ссылок на интернет-версии статей в~журналах <<Квант>> и~<<Математическое Просвещение>>, их можно найти на сайтах

\url{http://kvant.ras.ru}, \url{http://kvant.mccme.ru},

\url{http://www.mccme.ru/free-books/matpros.html}.

\subsection{Благодарности и~сведения об авторах}

Мы благодарим за серьёзную работу авторов материалов.
Благодарим за полезные замечания рецензентов книги Е.~А.~Авксентьева%(\S\,\ref{preobr},\ref{raznye})
, А.~В.~Антропова% (\S\,\ref{s:numb}, \ref{s:feresi})
, Е.~В.~Бакаева% (\S\,\ref{gmt})
, В.~Н.~Дубровского% (\S\,\ref{preobr})
, К.~А.~Кнопа% (\S\,\ref{s:alg})
, Д.~В.~Мусатова% (п.\,\ref{s:logic})
, Л.~Э.~Медникова% (глава~3)
, А.~А.~Полянского% (\S\,\ref{affin},\ref{stereo})
, А.~И.~Сгибнева% (\S\,\ref{s:compol}, \ref{s:fun}, \ref{oim}, \ref{phil-met}, \ref{s:motiv1})
, С.~Л.~Табачникова% (\S\,\ref{s:rad})
, А.~И.~Храброва% (\S\,\ref{s:ineq},\ref{s:seq})
, Г.~И.~Шарыгина% (\S\,\ref{treug},\ref{okruzhn})
и~Д.~Э.~Шноля% (\S\,\ref{compgeom})
, а~также анонимных рецензентов отдельных материалов. Благодарим А.~Я.~Канеля-Белова и~А.~В.~Шаповалова, авторов большого количества материалов, высказавших также ряд полезных идей и~замечаний. Благодарим Д.~А.~Пермякова, редактора книги~\cite{ZPS}. Благодарим учеников за каверзные вопросы и~указания на неточности.
Благодарим Е.~С.~Горскую и~П.~А.~Широкова за подготовку многих рисунков.
Благодарности по отдельным материалам приводятся прямо в~них.

Мы приносим извинения за допущенные неточности и~будем благодарны читателям за указания на них.

Главы 1, 2 и~3 редактировали А.~Б.~Скопенков, А.~А.~Заславский и~М.~Б.~Скопенков соответственно. Мы организовали рецензирование материалов главы~4, но  редактировали их сами авторы.

%%%!!!
%Какие пункты отмечать звёздочкой, решал редактор соответствующей главы\-.

М.~Б.~Скопенков и~А.~Б.~Скопенков частично поддержаны грантами фонда Саймонса и~фонда <<Династия>>. М.~Б.~ Скопенков частично поддержан грантом Президента РФ МК-6137.2016.1.

\textbf{Места работы и~интернет-страницы.}

А.~А.~Заславский: ЦЭМИ РАН, школа 1543.

А.~Б.~Скопенков: Московский физико-технический институт (ГУ) и~Независимый московский университет, \url{www.mccme.ru/~skopenko}.

М.~Б.~Скопенков: Национальный исследовательский университет Высшая школа экономики (факультет математики) и~Институт проблем передачи информации РАН, \url{http://skopenkov.ru}.

\subsection{Важные соглашения}\label{s:intrem}

Пункты внутри каждого параграфа расположены примерно в~порядке возрастания сложности материала. Цифры в~скобках после названия пункта означают его <<относительный уровень>>: 1 "--- самый простой, 4 "--- самый сложный. Первые пункты (не отмеченные звёздочкой) являются базовыми;
если не указано противное, с~них можно начать изучение главы.
А к~остальным пунктам (отмеченным звёздочкой) можно возвращаться потом; если не указано противное, то они независимы друг от друга.

Номера задач обозначаются жирным шрифтом.
Если условие задачи является формулировкой утверждения, то в~задаче требуется это утверждение доказать.
\emph{Загадкой} называется не сформулированный чётко вопрос; здесь нужно придумать и~чёткую формулировку, и~доказательство, ср.~\cite{VIN}.
В задачах, отмеченных кружочком $^\circ$, требуется привести только ответ без доказательства.
Наиболее трудные задачи отмечены звёздочкой *.
Если в~условии задачи написано <<найдите>>, то нужно дать ответ без знака суммы и~многоточия.
\emph{Указание и~решение} к~задаче может опираться на \emph{подсказку} к~ней.

Если некоторая задача не получается, то читайте дальше "--- следующие задачи могут оказаться подсказками.

\subsection{Основные обозначения}\label{s1-8}
%\addcontentsline{toc}{subsection}{Основные обозначения}

%\begin{bul}
$\bullet$ %\item
$[x]$ "--- (нижняя) целая часть числа $x$.

$\bullet$ %\item
$\{x\}$ "--- дробная часть числа $x$.

$\bullet$ %\item
$d \mid n$, или $n \,\text{\raisebox{-1pt}{$\vdots$}}\, d$ "--- число $n$ \emph{делится} на число $d$, т.\,е. $d\ne0$ и~существует такое целое $k$, что $n=kd$ (число $d$ называется \emph{делителем} числа $n$).

%%%!!!
%\blue{$\bullet$ %\item
%${\mathcal R}_n$ "--- множество $\{1,2,\ldots,n\}$. \textbf{"--- УБРАТЬ ЭТУ СТРОКУ!}}

$\bullet$ %\item
$\R$, $\Q$, $\Z$ "--- множества всех действительных, рациональных и~целых чисел соответственно.

$\bullet$ %\item
$\Z_2$ "--- множество $\{0, 1\}$ остатков от деления на 2 с~операциями сложения и~умножения по модулю 2.

$\bullet$ %\item
$\Z_m$ "--- множество $\{0, 1, \ldots, m-1\}$ остатков от деления на $m$
с~операциями сложения и~умножения по модулю $m$. (Заметим, что специалисты по алгебре чаще обозначают это множество $\Z/m\Z$, а~через $\Z_m$ обозначают
множество \emph{целых $m$"~адических чисел} для простого $m$.)

$\bullet$ %\item
$\binom{n}{k}$ "--- количество $k$"~элементных подмножеств $n$"~элементного  множества (другое обозначение: $C_n^k$).

$\bullet$ %\item
$|X|$ "--- число элементов во множестве $X$.

$\bullet$ %\item
$A - B = \{ x \mid x \in A \text{ и~} x \notin B \}$ "--- разность множеств $A$ и~$B$.

$\bullet$ %\item
$A \sqcup B $ "--- дизъюнктное объединение множеств $A$ и~$B$,
т.\,е. объединение $A \cup B$ непересекающихся множеств $A$ и~$B$.

$\bullet$ %\item
$A \subset B $ "--- <<множество $A$ содержится в~множестве $B$>>.
(В некоторых других книгах это обозначают $A \subseteq B $, а~$A \subset B $ означает <<множество $A$ содержится в~множестве $B$ и~не равно $B$>>.)

$\bullet$ %\item
Фраза <<обозначим $x=a$>> сокращается до $x:=a$.

$\bullet$ %\item
$\id$ "--- отображение множества в~себя, переводящее каждый элемент в~себя (тождественное).
%\end{bul}

%\chapter[Теория чисел, алгебра и~анализ]{Теория чисел, алгебра и~анализ \emph{А.~Б.~Скопенков}}
\begingroup

\chapter[Теория чисел, алгебра и~анализ. \emph{А.~Б.~Скопенков}]{Теория чисел, алгебра и~анализ}

\vspace*{-3.5\baselineskip}

\begin{center}
{\large\itshape\bfseries А.~Б.~Скопенков}
\end{center}

\vspace*{3\baselineskip}

{\it Умение преобразовывать алгебраические выражения "--- одно из базовых.
Его недостает <<олимпиадникам>>, из-за его отсутствия часто возникают нелепые и~обидные ошибки.
Поэтому для успешного решения задач алгебраического и~теоретико-числового типа рекомендуем нарабатывать культуру
арифметических выкладок.}

\section{Делимость и~деление с~остатком}\label{s:numb}

Из этого параграфа далее используются в~основном алгоритм Евклида и~его применения (задачи 2.5.7 и 2.5.9%\ref{numbdio-lem} и~\ref{numbdio-euc}
), язык сравнений (п.\;2.4%\ref{s:numbrem} 
 <<Деление с~остатком и~сравнения>>) и~простые факты (типа задач \ref{numbdiv-true} и~2.3.2%\ref{s:numbgcd}.2
 ).

В этом параграфе латинскими буквами обозначаются \emph{целые} числа.
Многие решения написаны с~использованием текстов М.~А.~Пра\-со\-лова.

\subsection{Делимость (1)}\label{s:numbdiv}

\begin{pr}\label{numbdiv-priz}

(a)
Сформулируйте и~докажите признаки делимости на 2, 4, 5, 10, 3, 9, 11.

(b)
Делится ли число $11 \ldots 1$ из 1993 единиц на 111111?

(c)
Число $1\ldots1$ из 2001 единиц делится на 37.
\end{pr}

\begin{pr}\label{numbdiv-even}
Если $a$ делится на 2 и~не делится на 4, то количество чётных делителей числа $a$ равно количеству его нечётных делителей.
\end{pr}

\begin{pr}\label{numbdiv-true}
Какие из следующих утверждений верны для любых $a,b$:

(a) $2|(a^2-a)$;\quad (b) $4|(a^4-a)$;\quad (c) $6|(a^3-a)$;\quad (d) $30|(a^5-a)$;

(e) если $c|a$ и~$c|b$, то $c|(a+b)$;

(f) если $b|a$, то $bc|ac$ для любого $c\ne0$;

(g) если $bc|ac$ для некоторого $c$, то $b|a$?
\end{pr}

При решении задачи \ref{numbdiv-true}\,(c) вы использовали следующий факт \ref{numbdiv-prod}\,(a). Докажите его по определению делимости, не используя единственности разложения на простые множители (задача 2.2.7%\ref{numbpri-dec}
\,(c))! Использование единственности может привести к~порочному кругу, ведь обычно при доказательстве единственности используется факт, близкий к~утверждению~\ref{numbdiv-prod}\,(b).

\begin{pr}\label{numbdiv-prod}
(a) Если число $a$ делится на 2 и~на 3, то $a$ делится на 6.

(b) Если число $a$ делится на 2, на 3 и~на 5, то $a$ делится на 30.

(c) Если число $a$ делится на 17 и~на 19, то $a$ делится на 323.
\end{pr}

\begin{pr}\label{numbdiv-ferm}
(a) Если $k$ не кратно ни 2, ни 3, ни 5, то $k^4-1$ кратно 240.

(b) Если $a+b+c$ делится на 6, то и~$a^3+b^3+c^3$ делится на 6.

(c) Если $a+b+c$ делится на 30, то и~$a^5+b^5+c^5$ делится на 30.

(d) Если $n\ge0$, то $20^{2n}+16^{2n}-3^{2n}-1$ делится на 323.
\end{pr}

\subsubsection*{Подсказки}

\paragraph*{\ref{numbdiv-prod}.} (a) Имеем $3a-2a=a$, поэтому $a$ делится на~6.

\sseccol{Указания, ответы и~решения}

\paragraph*{\ref{numbdiv-priz}.} (a) Для доказательства нижеприведённых признаков обозначим упоминаемое в~них число через
 $$
 n=\pm(10^ma_m+10^{m-1}a_{m-1}+\ldots+10a_1+a_0)
 $$
для некоторых $a_i$, $0\le a_i\le9$.

\section{Умножение по простому модулю}\label{s:feresi}

Из этого параграфа далее используются в~основном теорема Ферма"--~Эйлера
(задачи \ref{numbfer-fer} и~\ref{numbfer-eul}) и~теорема о~первообразном корне (задача \ref{s3.5.6}\,(b)). Впрочем, при применении теоремы о~первообразном корне понимать её доказательство не обязательно.

В этом параграфе латинскими буквами обозначаются \emph{целые числа} или \emph{вычеты} по простому модулю $p$ (что именно "--- видно из контекста).

\subsection{Малая теорема Ферма (2)}\label{s:numbfer}

\begin{pr}\label{numbfer-fer}
(a) Обозначим~$\Z_{97}=\{0,1,\ldots,96\}$. Определим отображение~$f\colon \Z_{97}\to \Z_{97}$ так: $f(a)$ равно остатку от~деления числа~$14a$ на~$97$. Тогда~$f$ "--- взаимно однозначное соответствие.

\emph{Обсуждение.} Достаточно доказать либо сюръективность, либо инъективность. Обычно доказывают \emph{инъективность}. Но необходимая для этого основная лемма арифметики обычно доказывается через разрешимость уравнения $97x+14y=1$, из которой сразу вытекает \emph{сюръективность}.

(b) Справедливо соотношение $(14\cdot1)\cdot(14\cdot2)\cdot\ldots\cdot(14\cdot96)\equiv96!\pmod{97}$.

(c) Справедливо соотношение $14^{96}\equiv1\pmod{97}$.

(d) {\bf Малая теорема Ферма.} Если $p$ простое, то $n^p-n$ делится на $p$ для любого целого $n$.\index{Теорема!малая теорема Ферма}

\emph{Alio modo.} Если $p$ простое и~$n$ не делится на $p$, то $n^{p-1}-1$ делится на $p$.

(f) Для простого~$p$ число~$\binom{p}{k}$ делится на~$p$ для любого $k=1,2,\ldots,p-1$. (Из~этого получается иное "--- по индукции "--- доказательство малой теоремы Ферма.)
\end{pr}

\begin{pr}\label{numbfer-rem}
Найдите остаток от деления
\quad

(a) $2^{100}$ на 101;  \quad
(b) $3^{102}$ на 101; \quad
(c) $8^{900}$ на 29;

%\parbox{2.5cm}{(a) $2^{100}$ на 101;} \quad
%\parbox{2.5cm}{(b) $3^{102}$ на 101; AS: убрал parbox, а то херня получалась} \quad
%(c) $8^{900}$ на 29;

\parbox{2.5cm}{(d) $3^{2000}$ на 43;} \quad
\parbox{2.5cm}{(e) $7^{60}$ на 143;} \quad
(f) $2^{60}+6^{50}$ на 143.
\end{pr}

\begin{pr}\label{numbfer-app}
(a) Если $p$ простое и~$p>2$, то $7^p-5^p-2$ делится на $6p$.

(b) Число $111\ldots11$ из 2002 единиц делится на 2003.

(c) Если $p$ и~$q$ "--- различные простые числа, то $p^q+q^p-p-q$ делится на $pq$.

(d) Число $30^{239}+239^{30}$ составное.

(e) Если $p$ простое, то длина периода десятичной дроби $1/p$ делит $p-1$.
\end{pr}

\begin{pr}\label{numbfer-ord}
Для простого $p$ и~$a$, не делящегося на $p$, назовём \emph{порядком} $\ord\,a=\ord_pa$ числа (или вычета) $a$ по модулю $p$ наименьшее $k>0$, для которого $a^k\equiv1\pmod p$:\index{Порядок|вычета по простому модулю|textbf}
 $$
 \ord\,a=\ord\phantom{}_pa:=\min\{k\ge1 \mid a^k\equiv1\!\!\pmod p\}.
 $$

(a) Множество $\{m\ge0\colon a^m\equiv1\pmod p\}$ состоит из целых неотрицательных чисел, кратных $\ord a$.

(b) Если $a^m\equiv a^n\pmod p$, то $m-n$ делится на $\ord a$.

(c) {\bf Лемма.} Число $p-1$ делится на $\ord a$. \quad

(d) Если $\ord x$ и~$\ord y$ взаимно просты, то $\ord(xy)=\ord x\cdot\ord y$.

(e) Для любых ли $a,x,p$ верно, что $a\ord_p x^a=\ord_p x$?
\end{pr}

Заметим, что по простому модулю можно определить деление и~отрицательные степени. Аналоги утверждений \ref{numbfer-ord}\,(a, b) справедливы для отрицательных степеней.

\begin{pr}\label{numbfer-eul}
В этой задаче буквами $p,q,p_1,\ldots,p_k$ обозначаются различные простые числа.

(a) Если $p\ne q$ и~$n$ не делится ни на $p$, ни на $q$, то $n^{(p-1)(q-1)}-1$ делится на $pq$.

(b) Если $n$ не делится на $p$, то $n^{p^\alpha(p-1)}-1$ делится на $p^{\alpha+1}$.

(c) \textbf{Теорема Эйлера.}
Если $n$ взаимно просто с~$m=p_1^{\alpha_1}\cdot\ldots\cdot p_k^{\alpha_k}$ и~$\varphi(m):=(p_1-1)p_1^{\alpha_1-1}\cdot\ldots\cdot (p_k-1)p_k^{\alpha_k-1}$,
то $n^{\varphi(m)}-1$ делится на $m$.\index{Теорема!Эйлера}

(d) Число $\varphi(m)$ равно количеству чисел от 1 до $m$, взаимно простых с~$m$.
\end{pr}

\begin{pr}\label{numbfer-cod}
(Загадка.) Известно, что $n$ "--- нечётное число от 3 до 47, не делящееся на~5. Как быстро вычислять неизвестное $n$ по известному $n^7 \bmod50$?
\end{pr}

Решение этой загадки показывает, почему для шифрования так важно быстро находить разложение числа на простые множители, в~частности быстро распознавать простоту числа.

\newpage
\setcounter{page}{47}
\refstepcounter{subsection}
\subsection{Квадратичные вычеты (2*)}\label{s:numquad}

Цель этого цикла задач "--- мотивировать и~обсудить проблему разрешимости сравнения $x^2\equiv a\pmod p$ для простого $p$.
В этом пункте через $p$ обозначается нечётное простое число.

\begin{pr}\label{s3.3.1}
(a) Какие остатки могут давать квадраты целых чисел при делении на 3, 4, 5, 6, 7, 8, 9, 10?

(b) Если $a^2+b^2$ делится на 3 (на 7), то $a$ и~$b$ делятся на 3 (на 7).

(c) Число вида $4k+3$ не представимо в~виде суммы двух квадратов.

(d) Существует бесконечно много чисел, не пред\-ставимых в~виде суммы трёх квадратов.
\end{pr}

\begin{pr}\label{s3.3.2}
Решите уравнения в~целых числах:

(a) $x_1^2+x_2^2+x_3^2+x_4^2+x_5^2=y^2$ (в~нечётных числах);

(b) $3x=5y^2+4y-1$;
\quad
(c) $x^2+y^2=3z^2$;
\quad
(d) $2^x+1=3y^2$;
\quad

(e) $x^2=2003y-1$;
\quad
(f) $x^2+1=py$, где $p=4k+3$;
\end{pr}

\begin{pr}\label{s3.3.3}
(a) Если $p=4k+3$ делит $a^2+b^2$, то $p|a$ и~$p|b$.

(b) Число, в~каноническое разложение которого некоторый простой делитель вида $4k+3$ входит в~нечётной степени, не представимо в~виде суммы двух квадратов (целых чисел).

(c)$^*$ Уравнение $x^2+1=py$ разрешимо в~целых числах при $p=4k+1$
(и~неразрешимо при $p=4k+3$).

(d)$^*$ Любое простое число вида $4k+1$ представимо в~виде суммы двух квадратов.

(e)$^*$ Число, в~каноническое разложение которого любой простой делитель вида $4k+3$ входит в~чётной степени,
представимо в~виде суммы двух квадратов.

(f) Простых чисел вида $4k+1$ бесконечно много.
\end{pr}

Доказательство Дон Загира утверждения (d) можно найти в~книге~\cite{Pr07-1}.

\begin{pr}[]\label{s3.3.4}
(Загадка.) <<Сведите>> уравнение $py=at^2+bt+c$, $a\ne0$, к~сравнению $x^2\equiv k\pmod p$.
\end{pr}

Остаток $a\ne0$ называется \emph{квадратичным вычетом \textup{(}квадратичным невычетом\textup{)} по модулю $p$}, если сравнение $x^2\equiv a(p)$ разрешимо (неразрешимо). Слова <<по модулю $p$>> далее опускаются.
\index{Квадратичный|вычет|textbf}\index{Квадратичный|невычет|textbf}

\begin{pr}\label{s3.3.5}
(a) Приведите пример таких $a$ и~$p$, что оба числа $a$ и~$-a$ являются квадратичными вычетами.

(b) Если $a$ не делится на $p$, то сравнение $x^2\equiv a^2(p)$ имеет ровно два решения.

(c) {\bf Лемма.} Число квадратичных вычетов равно числу квадратичных невычетов и~равно $\frac{p-1}2$.
\end{pr}

\begin{pr}\label{s3.3.6}
(a) {\bf Лемма.} Для любого $a\ne0$ существует и~единственно такое~$b$, что $ab\equiv1(p)$.

Обозначение: $b=a^{-1}$.

(b) Решите сравнение $x\equiv x^{-1}(p)$.

(c) \textbf{Теорема Вильсона.} Число $(p-1)!+1$ делится на $p$.
\end{pr}

\begin{pr}\label{s3.3.7}
(a) Если $a\ne0$ "--- квадратичный вычет, то $a^{-1}$ тоже квадратичный вычет.

(b) Число квадратичных вычетов чётно тогда и~только тогда, когда $-1$ является квадратичным вычетом.
\end{pr}

\begin{pr}\label{s3.3.8}
\textbf{Лемма.} (a) Произведение двух квадратичных вычетов является квадратичным вычетом.

(b) Произведение квадратичного вычета и~квадратичного невычета является
квадратичным невычетом.

(c) Произведение двух квадратичных невычетов является квадратичным вычетом.
\end{pr}

\newpage
\setcounter{page}{50}
\subsection{Квадратичный закон взаимности (3*)}\label{s:numbrec}

Здесь строится алгоритм выяснения разрешимости сравнения $x^2\equiv a\pmod p$ для простого~$p$. Используется п.\;\ref{s:numquad} <<Квадратичные вычеты>>.

\begin{pr}\label{s3.4.1}
Если число $p=8k+5$ простое, то

(a) $2^{4k+2}\equiv-1\pmod p$;

(b) уравнение $x^2-2=py$ неразрешимо в~целых числах.
\end{pr}

\begin{pr}\label{s3.4.2}
Если число $p=8k+1$ простое, то

(a) $2^{4k}\equiv1\pmod p$;

(b) уравнение $x^2-2=py$ разрешимо в~целых числах.
\end{pr}

\begin{pr}\label{s3.4.3}
(a) Если число $p=8k\pm1$ простое, то $2^{(p-1)/2}\equiv1\pmod p$.

(b) Если число $p=8k\pm3$ простое, то $2^{(p-1)/2}\equiv-1\pmod p$.

(c) Для каких простых $p$ разрешимо в~целых числах уравнение $x^2-2=py$?
\end{pr}

\begin{pr}\label{s3.4.4}
(а) Если число $p=12k\pm1$ простое, то $3^{(p-1)/2}\equiv1\pmod p$.

(b) Если число $p=12k\pm5$ простое, то $3^{(p-1)/2}\equiv-1\pmod p$.

(c) Для каких простых $p$ разрешимо в~целых числах уравнение $x^2-3=py$?
\end{pr}

\begin{pr}\label{s3.4.5}
Для нечётного простого числа $p$ рассмотрим \emph{символ Лежандра}
 $$
 \Big(\frac ap\Big):=\begin{cases}
 +1, & \text{$a$ "--- квадратичный вычет по модулю $p$};\\[-1mm]
 -1, & \text{$a$ "--- квадратичный невычет по модулю $p$.}
 \end{cases}
 $$
Например, $\Big(\frac2p\Big)=(-1)^{(p^2-1)/8}$ по задаче \ref{s3.4.3} и~$\Big(\frac{ab}p\Big)=\Big(\frac ap\Big)\Big(\frac bp\Big)$ по задаче \ref{s3.3.8}.
%\end{pr}

(a) {\bf Критерий Эйлера.} Справедливо соотношение
$$\Big(\frac ap\Big)\equiv a^{\frac{p-1}2}\pmod p.$$

(b) {\bf Лемма Гаусса.} Справедливо соотношение
$$\Big(\frac ap\Big) = (-1)^{\sum\limits_{x=1}^{(p-1)/2}[\frac{2ax}p]}.$$

\end{pr}

\newpage
\setcounter{page}{53}

\subsection{Первообразные корни (3*)}\label{s:numbroo}

\begin{pr}\label{s3.5.1}
(2--7) Сформулируйте и~обоснуйте алгоритм решения сравнения $a^x\equiv b\ (m)$ для заданных $a,b$, взаимно простых с~заданным $m\in\{2,3,4,5,6,7\}$.

(Решение такого сравнения "--- одна из основных мотивировок этого занятия.)
\end{pr}

\begin{pr}\label{s3.5.2}
(a) Если $(a,35)=1$, то $a^{12}\equiv1\pmod{35}$.

(b) Если $m$ делится на два различных простых нечётных числа и~$(a,m)=1$, то $a^{\frac{\varphi(m)}2}\equiv1\pmod m$.
\end{pr}

Пусть $(g,m)=1$.
Вычет $g$ называется \emph{первообразным корнем} по модулю $m$,
если остатки от деления на $m$ чисел $g^1,g^2,\ldots,g^{\varphi(m)}\equiv1$ различны.
\index{Первообразный корень по модулю|textbf}
Например,

%\begin{bul}
$\bullet$ %\item
число 2 является первообразным корнем по модулю 5, а~число 4 "--- нет;

$\bullet$ %\item
по задаче \ref{s3.5.2}\,(b) если $m$ делится на два различных простых нечётных числа, то не существует первообразного корня по модулю~$m$.
%\end{bul}

\begin{pr}\label{s3.5.3}
Докажите существование первообразного корня по простому модулю следующего вида:\nopagebreak\smallskip

(a) 257; \quad (b) $2^l+1$; \quad (c) $2^k\cdot 3^l+1$; \quad (d) 151; \quad (e) $2^k\cdot 3^l\cdot5^m+1$.
\end{pr}

Простой метод решения пунктов (a), (b), (c), не проходит для (d),~(e). Продемонстрируем метод решения пунктов (d), (e) на примерах.

\begin{pr}\label{s3.5.4}
(a) Вычет $g$ "--- первообразный корень по модулю 97 тогда и~только тогда, когда ни $g^3$, ни $g^{32}$ не сравнимы с~1 по модулю~97.

(b) Сравнение $x^3\equiv1\pmod{97}$ имеет ровно 3 решения.

(c) Сравнение $x^{32}\equiv1\pmod{97}$ имеет ровно 32 решения.

(d) Существует первообразный корень по модулю 97.

(e) Количество первообразных корней по модулю 97 равно 63.
\end{pr}

\begin{pr}\label{s3.5.5}
(a) Вычет $g$ "--- первообразный корень по модулю 151 тогда и~только тогда, когда ни $g^2$, ни $g^3$, ни $g^{25}$ не сравнимы с~1 по модулю~151.

(b) Сравнение $x^k\equiv1\pmod{151}$ имеет ровно $k$ решений для $k\in\{30,50,75\}$.
%%%!!!
%\smallskip

(c) Имеет место равносильность
 $$\begin{cases}x^{30}\equiv1\pmod{151},\\
 x^{50}\equiv1\pmod{151}
 \end{cases}\Leftrightarrow x^{10}\equiv1\pmod{151}.
 $$

(d) Существует первообразный корень по модулю 151.

(e) Количество первообразных корней по модулю 151 равно 40.
\end{pr}

\begin{pr}\label{s3.5.6}
(a)
Если $p$ простое и~$p-1$ делится на $d$, то сравнение $x^d\equiv1\pmod p$ имеет ровно $d$ решений.

(b) \textbf{Теорема о~первообразном корне.} Для любого простого $p$ существует число
$g$, для которого остатки от деления на $p$ чисел $g^1,g^2,g^3,\ldots,g^{p-1}=1$ различны.
\index{Теорема!о первообразном корне}

(c) Сколько существует первообразных корней по простому модулю $p$?
\end{pr}

\sseccol{Указания, ответы и~решения}

\paragraph*{\ref{s3.5.3}.} (b) Если первообразного корня нет, то сравнение $x^{2^{l-1}}\!\equiv1\!\pmod p$ имеет $p-1=2^l>2^{l-1}$ решений.

\paragraph*{\ref{s3.5.6}.}
(a) Заметьте, что многочлен $x^{p-1}-1$ над $\Z_p$ имеет ровно $p-1$ корень и~делится на $x^d-1$. Докажите, что если многочлен степени $a$ имеет ровно $a$ корней и~делится на многочлен степени $b$, то этот многочлен степени $b$ имеет ровно $b$ корней.

Другое решение можно получить, заметив, что если $p=kd$, то для любого $a$
сравнение $y^k\equiv a\pmod p$ имеет не более $k$ решений.

(c) \emph{Ответ}: $\varphi(p-1)$.

%\newpage
\subsection{Высокие степени (3*). \emph{А.~Я.~Канель-Белов}, \emph{А.~Б.~Скопенков}}\label{s:numbhi}

\begin{pr}\label{s3.6.0}
(a) Для любых $n$ и нечетного $k$ число $k^{2^n}-1$ делится на $2^{n+2}$.

(b) Для любого $n$ число $2^{3\cdot7^n}-1$ делится на $7^{n+1}$.
\end{pr}

\begin{pr}\label{s3.6.1}
При каких $a$ \quad

(a) $2^a -1$ делится на $3^{100}$; \quad
(b) $2^a +1$ делится на $3^{100}$;

(c) $5^a-1$ делится на $2^{100}$; \quad
(d) $2^a -1$ делится на $5^{100}$?
\end{pr}

Утверждение \ref{s3.6.0}\,(a) означает, что ни при каком $n\ge3$ не существует первообразного корня по модулю $2^n$ (см. определение в п. \ref{s:numbroo}).
Ответы к задачам \ref{s3.6.1}.(a),(d),(c) и утверждение \ref{s3.6.0}.(b) означают, что для любого $n$ число 2 является первообразным корнем
по модулю $3^n$ и по модулю $5^n$, а числа 5 и 2 не являются первообразными корнями по модулю $2^n$ и по модулю $7^n$.

%далее переставлено

\begin{pr}\label{s3.6.2}
(a) Найдите первообразный корень по модулю $7^{100}$.

(b) \textbf{Теорема.} Первообразные корни существуют только по модулям $2,4,p^n,2p^n$.
\end{pr}

\begin{pr}\label{s3.6.4} Пусть $p>2$ простое, $g$ "--- первообразный корень по модулю~$p$ и $g^{p-1}-1$ не делится на $p^2$.
Тогда  $g$ "--- первообразный корень по модулю

(a) $p^2$; \quad (b) $p^3$; \quad (c) $p^n$ для любого $n$.
\end{pr}

\begin{pr}\label{s3.6.4}
Пусть $p>2$ простое.

(a) Если $g$ "--- первообразный корень по модулю $p$, то одно из чисел $g^{p-1}-1$ и $(g+p)^{p-1}-1$ не делится на $p^2$.

%(a) существуют такие $t$ и~$u$, что $(g+pt)^{p-1}=1+pu$ и~$u$ не делится на $p$.

(b) Если $g$ "--- первообразный корень по модулю $p^2$, то $g$ "--- первообразный корень по модулю $p^n$ для любого~$n$.

(c) Для любого целого положительного $n$ существует первообразный корень по модулю $p^n$.

(d) Tо же по модулю $2p^n$.
\end{pr}

\begin{pr}\label{s3.6.6}
{\bf Лемма об уточнении показателя.}
Пусть $p$ "--- простое число, $p>2$ или $n>1$, $q$ не делится на~$p$ и~$x-1$ делится на $p^n$, но не на $p^{n+1}$.

(a) Число $x^q-1$ делится на $p^n$, но не на $p^{n+1}$.

(b) Число $x^p-1$ делится на $p^{n+1}$, но не на $p^{n+2}$.

(c) Число $x^{p^kq}-1$ делится на $p^{n+k}$, но не на $p^{n+k+1}$.
\index{Лемма!об уточнении показателя}

(Близкое утверждение называется \emph{леммой Гензеля}.)
\end{pr}

\begin{pr}\label{s3.6.7}
Найдите длину периода дроби\ \ (a) $1/3^{100}$;\quad (b) $1/7^{100}$.
\end{pr}

\newpage
\setcounter{page}{62}
\refstepcounter{section}
\refstepcounter{subsection}
\sectionmark{Многочлены и комплексные числа}
\subsection{Решение уравнений 3"~й и~4"~й степени (2)}\label{mot34}

Благодарю О.~Е.~Орёл за полезные обсуждения.

%AS переставил

Приведённый здесь материал важен и~широко известен, но не входит в~школьную или университетскую программу.
Отличие приводимого рассуждения от встречающихся в~других источниках в~том, что (вместо немотивированных замен)
уравнения естественно сводятся к~таким, которые ясно, как решать.

%не заменять на: мы можем легко решить.

Например, уравнение $x^2+4x-1=0$ сводится к уравнению $y^2-5=0$ заменой переменной $y=x+2$.

\begin{pr}\label{2zamena}
(a) Уравнение $x^3+3x^2+5x+7=0$ <<сводится>> заменой переменных к~уравнению $y^3+py+q=0$ с~некоторыми числами $p,q$.

(b) Уравнение $ax^3+bx^2+cx+d=0$ при $a\ne0$ <<сводится>> заменой переменных к~уравнению $y^3+py+q=0$ с~некоторыми числами $p,q$.

(c) Уравнение $ax^4+bx^3+cx^2+dx+e=0$ при $a\ne0$ <<сводится>> заменой переменных к~уравнению $y^4+py^2+qy+r=0$ с некоторыми числами $p,q,r$.
\end{pr}

\begin{pr}\label{2cubic}
(a) Докажите, что $\sqrt[3]{2+\sqrt5}-\sqrt[3]{\sqrt5-2}=1$.

(b) Найдите хотя бы одно решение уравнения $x^3-3\sqrt[3]2x+3=0$.

\emph{Указание}. \textbf{Метод дель Ферро.} Так как
 $$
 (u+v)^3=u^3+v^3+3uv(u+v),
 $$
то число $u+v$ является корнем уравнения $x^3-3uvx-(u^3+v^3)=0$.
\index{Метод!дель Ферро}

(c) Решите уравнение $x^3-3\sqrt[3]2x+3=0$.
\end{pr}

\begin{pr}\label{2decom}
(a) Разложите на множители выражение $a^3+b^3+c^3-3abc$.

(b) Справедливо неравенство $a^2+b^2+c^2\ge ab+bc+ca$. Когда достигается равенство?

(c) Справедливо неравенство $a^3+b^3+c^3\ge3abc$ при $a,b,c>0$.

(d) Разложите выражение $a^3+b^3+c^3-3abc$ на линейные множители с~комплексными коэффициентами.
\end{pr}

Задачи этого пункта о~комплексных числах можно пропустить. Но для их решения необходимы лишь
минимальные сведения о~комплексных числах, достаточно уметь решать задачи 4.5.1 и 4.5.2.

\begin{pr}\label{2dferro}
(a) Сформулируйте и~докажите теоремы, описывающие все вещественные (все комплексные) решения уравнения $x^2+px+q=0$.

(b) Сформулируйте и~докажите теоремы, описывающие все вещественные (все комплексные) решения уравнения $x^3+px+q=0$ в~том случае, когда работает метод дель Ферро (см.~задачу \ref{2cubic}).
А~при каком условии на $p,q$ применим этот метод для вещественных решений,
если квадратные корни разрешается извлекать только из положительных чисел?

%(b) То же для комплексных решений.

(c) Составьте алгоритм (точного, или символьного) нахождения всех вещественных корней уравнения $ax^3+bx^2+cx+d=0$, где $a\ne 0$.
\end{pr}

При решении некоторых кубических уравнений методом дель Ферро в~формулах неожиданным образом возникают комплексные числа "--- как раз тогда, когда все корни исходного уравнения вещественны.
Такие уравнения можно также решать следующим <<чисто вещественным>> методом.
Он также интересен тем, что подводит к~\emph{трансцендентным методам} решения уравнений \cite{PrSo}.

%%%!!!

\begin{pr}\label{2vieta} (a) Решите уравнение $4x^3-3x=\frac12$.
\index{Метод!Виета}

(b) Решите уравнение $x^3-3x-1=0$.

(c) Используя функции $\cos$ и~$\arccos$, напишите общую формулу для
решения уравнения $x^3+px+q=0$ методом, намеченным в~этой задаче.
При каком условии уравнение $x^3+px+q=0$ решается этим методом?
\end{pr}

\begin{pr}\label{2fourth}
Решите уравнение

\parbox{4.5cm}{(a) $(x^2+2)^2=9(x-1)^2$;}
(b) $x^4+4x-1=0$;

\parbox{4.5cm}{(c) $x^4+2x^2-8x-4=0$;}
(d) $x^4-12x^2-24x-14=0$.
\end{pr}

\emph{Указание к~задаче} \ref{2fourth}\,(b). \textbf{Метод Феррари.} Подберите такие $\alpha$, $b$, $c$, что
 $$
 x^4+4x-1=(x^2+\alpha)^2-(bx+c)^2.
 $$
Для этого найдите хотя бы одно $\alpha$, для которого квадратный трёхчлен $(x^2+\alpha)^2-(x^4+{4x-1})$ от $x$ является полным квадратом. Для этого найдите дискриминант этого квадратного трёхчлена. Он является кубическим многочленом от $\alpha$ и~называется \emph{кубической резольвентой} многочлена $x^4+4x-1$.
\index{Метод!Феррари}

\newpage
\setcounter{page}{79}
\refstepcounter{subsection}
\refstepcounter{subsection}
\refstepcounter{subsection}
\refstepcounter{subsection}
\subsection{Диофантовы уравнения и~гауссовы числа (4*). \emph{А.~Я.~Канель-Белов}}\label{s:polgau}

Всем хорошо знаком алгоритм Евклида. Даны два числа $a$, $b$. Из них выбирается большее, из большего вычитается меньшее, большее заменяется на разность, и~с~новой парой чисел производится та же процедура. См.~задачу 2.5.9%\ref{numbdio-euc}
\,(b). С~помощью алгоритма Евклида доказываются арифметические свойства чисел и~это Вы изучали раньше (см. п.\;2.5%\ref{s:numbdio} 
<<Линейные диофантовы уравнения>> и~п.\;4.4%\ref{s:poleup} 
<<Делимость для многочленов>>). Приведём принципиально новые (для большинства школьников) его применения.

\begin{pr}\label{polgau-equ}
Решите уравнения в~целых числах:

(a) $x^2+4=y^3$; \quad (b) $x^2+2=y^n$; \quad (c)$^*$ $x^3+y^3=z^3$.
\end{pr}

Попробуйте порешать их, не читая дальнейшего! Впрочем, у~Вас вряд ли получится. Возвращайтесь к~этой задаче по мере чтения дальнейшего материала.

При решении уравнения $x^2+4=y^3$ в~целых числах хочется действовать так:
$x^2+4=(x+2i)(x-2i)$. При нечётном $x$ оба эти множители взаимно просты, и~потому оба являются кубами. Из этого получается решение. (Случай чётного $x$ хитрее: обе скобки могут делиться на $(1+i)^3$.) Попробуйте довести решение до конца, а~затем сравнить с~приведённым в~конце темы.

Одним словом, хочется наслаждаться дополнительными возможностями при разложении на множители за счёт использования \emph{гауссовых чисел}, т.\,е. чисел вида $a+bi$ с~целыми $a$ и~$b$. Однако не всё коту масленица "--- так получается не всегда (см.~задачи 2.2.7%\ref{numbpri-dec}
\,(b) и~\ref{polgau-exa}\,(b)), но иногда получается. Чтобы применять разложение на множители для решения уравнений, нужна \emph{однозначность разложения на простые множители}. Если она имеет место, то мы имеем всё те же арифметические удовольствия, что и~для целых чисел. Следующая задача показывает удивительный факт: для \emph{арифметических} удовольствий достаточно доказать \emph{геометрический} факт о~возможности деления с~остатком.

\begin{pr}\label{polgau-pro}
Гауссово число называется \emph{простым}, если оно не разлагается на два множителя,
каждый из которых отличен от $\pm1$ и~$\pm i$.
В этой задаче латинские буквы обозначают гауссовы числа.

(a) Однозначность разложения на простые множители вытекает из
следующего свойства (аналога леммы Евклида 2.5.7 %\ref{numbdio-lem}
\,(с)).

\textbf{Факториальность.} Для любых $a,b$ если простое число $p$ делит $ab$, то $p$ делит $a$ или $p$ делит $b$.

(b) Факториальность вытекает из следующего свойства (аналога леммы о~представлении НОД 2.5.7% \ref{numbdio-lem}
\,(a)).

\textbf{Главноидеальность.} Для любых $a,b$ существуют такие $x,y$, что $xa+yb=\gcd(a,b)$. (Дайте определение наибольшего общего делителя $gcd(a,b)$ чисел $a,b$ самостоятельно!)

(c) Главноидеальность обеспечивается следующим свойством (аналогом теоремы о~делении с~остатком 2.4.1 %\ref{numbdiv-rema}
\,(b)).
\index{Теорема!о делении с~остатком!гауссовых чисел}
\index{Алгоритм Евклида!для гауссовых чисел}

\textbf{Евклидовость.} Для любых $b\ne0$ и~$a$ существует такое $k$, что $|a-kb|<|b|$.

\end{pr}

\begin{pr}\label{polgau-exa}
Верна ли евклидовость (и, значит, факториальность!) для множества $\Z[\xi]$ чисел вида $a+b\xi$ с~целыми $a,b$, если $\xi$ есть

(a) $\sqrt{-2}$; \quad (b) $\sqrt{-3}$; \quad (c) $(1-\sqrt{-3})/2$; \quad (d) $(1-\sqrt{-5})/2$; \quad

(e) $(1-\sqrt{-7})/2$?
\end{pr}

\begin{pr}\label{p4-7-4}
(a) Никакое простое число вида $4k-1$ не разлагается в~сумму двух квадратов.

(b) Любое простое число вида $4k+1$ разлагается в~сумму двух квадратов, причём ровно одним способом.

(b) Существует целое число, ровно 1024 способами разлагающееся в~сумму двух квадратов.

Эту задачу проще решать без гауссовых чисел (см.~п.\;\ref{s:numquad}), однако полезно потренироваться в~их применении!
\end{pr}

%Подробнее см.~\cite[\S\,4]{Po}. См.~также задачу \ref{poleupol-per}.

\newpage
\setcounter{page}{88}
\section{Разрешимость в радикалах}\label{s:rad}

%%%!!!
%удалил comment

Основное содержание этого параграфа "--- простые элементарные доказательства знаменитых теорем Гаусса, Абеля, Галуа и~Кронекера о~построимости правильных многоугольников и~неразрешимости уравнений в~радикалах. На примере этих доказательств иллюстрируются некоторые основные идеи алгебры.
Определения построимости и~разрешимости в~радикалах, а~также формулировки указанных теорем приводятся.
Этот параграф адресован всем любителям изложения глубоких идей на примерах красивых результатов и~доказательств: старшеклассникам, студентам, учителям и~профессиональным математикам. Хороший опыт в~работе с~комплексными числами и~многочленами получит и~тот, кто не дойдёт до полного доказательства основных результатов.

%%%!!!Для понимания доказательств достаточно знакомства с~многочленами \blue{(п. \ref{s:ratir}, \ref{s:polbez}, %\ref{s:poleup})} и~умения извлекать корни из комплексных чисел \blue{(задача \ref{compl-tri})}.

\subsection{Введение}\label{s:intr}

\subsubsection{О чём этот параграф}

Основное содержание этого параграфа "--- простые элементарные доказательства

%\begin{bul}
$\bullet$ %\item
теоремы Гаусса о~построимости правильных многоугольников
(и~даже более сильного результата "--- теоремы Гаусса о~понижении);
\index{Теорема!Гаусса о~построимости правильных многоугольников}

$\bullet$ %\item
существования уравнения 3"~й степени, неразрешимого в~\emph{вещественных} радикалах (и~даже более сильного результата "--- сильной вещественной теоремы о~неразрешимости);
\index{Теорема!о неразрешимости в~вещественных радикалах}

$\bullet$ %\item
теоремы Галуа о~существовании уравнения 5"~й степени, неразрешимого в~\emph{комплексных} радикалах (и~даже более сильного результата "--- теоремы Кронекера).
\index{Теорема!Галуа о~неразрешимости в~радикалах}
\index{Теорема!Кронекера о~неразрешимости в~радикалах}
%\end{bul}

Определения построимости и~разрешимости в~радикалах, а~также формулировки указанных теорем приведены в~п.\;\ref{intsqu} и~\ref{intrad}. Я~не привожу историю этих знаменитых теорем, отсылая заинтересованного читателя к~текстам~\cite{Gi,Gi1,Ma}.

Приводимые доказательства интересны тем, что для их понимания достаточно уметь доказывать иррациональность (п.\;4.1%\ref{s:ratir}
), делить многочлены с~остатком (п.\;4.3%\ref{s:polbez}
 и~задачи 4.4.3 и 4.4.4%\ref{poleup-rem}, \ref{poleupol-mi}
 ), извлекать корни из комплексных чисел (задача 4.5.4 % \ref{compl-tri}) 
 и~решать системы линейных уравнений. Эти прямые доказательства интересны тем, что на них ясно видны базовые идеи важной \emph{теории Галуа}
(см.~п.\;5.2.1%\ref{motint} 
и~\S\,\ref{phil-met}).

Приводимые доказательства не претендуют на новизну (хотя в~этом тексте имеется много методических находок, см.~п.\;5.2.1 и 5.2.2%\ref{motint} и~\ref{phil-hist}
).
Однако, к~сожалению, они малоизвестны. Как следствие, малоизвестно, что не только решать квадратные и~кубические уравнения, но и~доказывать указанные
теоремы экономнее не строя и~затем применяя теорию Галуа (как, например, в~\cite{Kh13, Ki05}), а~напрямую "--- но при этом, конечно, открывая и~используя базовые идеи этой теории.

Параграф адресован тем, кому интересен хотя бы один из этих результатов.
Разбор доказательств (или их начала) полезен для закрепления тем <<иррациональность>>, <<многочлены>>, <<комплексные числа>> и~<<основы линейной алгебры>>. Старшеклассник найдёт здесь задачи для исследования, не претендующие на научную новизну, подробнее см.~п.\;5.2.1%\ref{motint}
. Пункты \ref{gaucon} и~\ref{pro} могут быть интересны профессиональному математику.

Как устроен параграф, написано в~п.\;\ref{intplan}. Впрочем, начать изучать параграф можно не с~п.\;\ref{s:intr}, а~с~решения задач в~п.\;5.4%\ref{non}
, поскольку большинство из них использует предыдущий материал только в~качестве мотивировки.

Благодарю А.~Я.~Белова-Канеля, И.~И.~Богданова, Э.~Б.~Винберга, В.~В.~Волкова, М.~Н.~Вялого, А.~С.~Голованова, П.~А.~Дергача, Д.~Зунга, А.~А.~Казначеева, А.~Л.~Канунникова, Г.~А.~Мерзона, А.~А.~Пахарева, В.~В.~Прасолова, А.~Д.~Руховича, Л.~М.~Самойлова, М.~Б.~Скопенкова, Г.~Р.~Челнокова, Л.~А.~Шабанова, В.~В.~Шувалова и~особенно В.~А.~Клепцына за полезные замечания и~предложения.

\subsubsection{Разрешимость в~квадратных радикалах: формулировки~(1)}\label{intsqu}

Известно, что
%$$
\begin{equation*}
\begin{gathered}
 \cos\frac{2\pi}3=-\frac12,\quad \cos\frac{2\pi}4=0, \quad \cos\frac{2\pi}5=\frac{\sqrt5-1}4,\\[1mm]
 \cos\frac{2\pi}6=\frac12, \quad \cos\frac{2\pi}8=\frac1{\sqrt2}.
\end{gathered}\tag*{($*$)}\label{zv1}
\end{equation*}
%$$
Как обобщить эти формулы (используя только четыре арифметических действия и~извлечения корней)? Для формализации этого вопроса введём следующие определения.

%AS вернул авторскую версию, ибо было неправильно
%\begin{deff}

\smallskip
{\bf Определение вещественного калькулятора.} Рассмотрим калькулятор с~кнопками
 $$
 1,\quad +,\quad -,\quad \times,\quad :\quad \text{и}\quad \sqrt[n]{\phantom{n}}\quad
 \text{для любого $n$}.
 $$
Калькулятор вычисляет числа с~абсолютной точностью и~имеет неограниченную память. При делении на 0 он выдаёт ошибку.
\index{Калькулятор!вещественный|textbf}

\emph{Вещественный} калькулятор оперирует с~вещественными числами
и~при извлечении корня чётной степени из отрицательного числа выдаёт ошибку.
%\end{deff}

%\begin{deff}
\smallskip
{\bf Определение вещественной построимости.}
Вещественное чис\-ло называется \emph{вещественно построимым}, если его можно получить на вещественном калькуляторе так, чтобы при этом извлекались корни только второй степени (т.\,е. получить из 1 при помощи сложений, вычитаний, умножений, делений и~извлечений квадратного корня из положительных чисел).
\index{Число!вещественно построимое|textbf}
%\end{deff}

\smallskip
Например, вещественно построимы числа
 $$
 \sqrt[4]2=\sqrt{\sqrt2},\ \ \sqrt{2\sqrt3},\ \ \sqrt2+\sqrt3,\ \ \sqrt{1+\sqrt2},\ \ 1+\sqrt{3-2\sqrt2},\ \ \frac1{1+\sqrt2},
 $$
числа из формулы~$(*)$ и~даже число $\cos\frac{\pi}{60}=\cos 3^\circ.$
Про последнее число это не совсем очевидно (но мы это увидим в~п.\;\ref{gauref}).

Вопрос об обобщении формул $(*)$ формализуется так: для каких $n$ число $\cos(2\pi/n)$ вещественно построимо? Ответ даётся следующей теоремой.

\begin{theorem}[Гаусса] Число $\cos(2\pi/n)$ вещественно построимо тогда и~только тогда\textup{,} когда $n=2^\alpha p_1\ldots p_l,$ где $p_1,\ldots,p_l$ "--- различные простые числа вида $2^{2^s}+1$.
\index{Теорема!Гаусса о~построимости правильных многоугольников}
\end{theorem}

Построимость в~теореме Гаусса доказана в~п.\;\ref{gauref} и~\ref{gaucon}, а~непостроимость "--- в~п.\;\ref{progau}.

Вещественная построимость числа равносильна его \emph{построимости циркулем и~линейкой}. Поэтому теорема Гаусса равносильна критерию построимости циркулем и~линейкой правильных многоугольников. Мы обсудим эту равносильность в~п.\;5.2.3%\ref{motcon}
; впрочем, она не будет использоваться в~остальном тексте.

Из вещественной непостроимости числа $\cos(2\pi/9)$ вытекает следующий результат.

\smallskip
{\bf Следствие.}  {\it Трисекция угла невозможна на вещественном калькуляторе\textup{,} если можно извлекать корни только второй степени\textup{,} или\textup{,} формально\textup{,} на нём невозможно получить число $\cos(\alpha/3)$\textup{,} имея число $\cos\alpha$ \textup{(}например\textup{,} для $\alpha=2\pi/3)$.}

\begin{Remarks}
(a) Строго говоря, теорема Гаусса не даёт настоящего решения проблемы построимости, поскольку неизвестно, какие числа вида $2^{2^s}+1$ являются простыми. Однако теорема Гаусса даёт, например, быстрый
алгоритм выяснения построимости числа $\cos(2\pi/n)$.

(b) Для практики приближённые методы вычисления тригонометрических функций и~решения уравнений более полезны, чем радикальные формулы. Кроме того, уравнения степени выше 4 разрешимы при помощи трансцедентных функций (см.~метод Виета в~п.\;\ref{mot34} и~\cite{PrSo}; о~развитии этих идей рассказывается, например, в~\cite{Sk10}). Однако проблема разрешимости в~радикалах интересна как пробная задача современных теорий символьных вычислений и~сложности вычислений.
\end{Remarks}

\subsubsection{Неразрешимость в~радикалах: формулировки (2)}\label{intrad}

Следующее утверждение даёт достаточное условие разрешимости уравнений третьей степени <<в вещественных радикалах>>.

\smallskip
{\bf Утверждение о~разрешимости в~вещественных радикалах.}
{\it Если многочлен третьей степени с~рациональными коэффициентами имеет ровно один вещественный корень\textup{,} то этот корень можно получить на вещественном калькуляторе\textup{\footnote{Стандартная терминология: уравнение \emph{разрешимо в~вещественных радикалах}.}}. Более того\textup{,} это можно сделать так\textup{,} чтобы извлечение корня происходило только два раза\textup{,} один раз второй и~один раз третьей степени.}

\smallskip
Это утверждение доказывается \emph{методом дель Ферро} (оно вытекает из теоремы, приведённой в~подсказке к~задаче \ref{2dferro}, и~результата задачи \ref{2numbroot}\,(d)).

\begin{theorem}[о~неразрешимости в~вещественных радикалах]
Существует многочлен $3$"~й степени с~рациональными коэффициентами \textup{(}например\textup{,} $x^3-3x+1),$ ни один из корней которого невозможно получить на вещественном калькуляторе.
\index{Tеорема!о неразрешимости в~вещественных радикалах}
\end{theorem}

Эта теорема доказана в~п.\;\ref{prorea}.

\smallskip
{\bf Следствие.} {\it Трисекция угла невозможна на вещественном калькуляторе\textup{,} или\textup{,} формально\textup{,} на нём невозможно получить число $\cos(\alpha/3),$ имея число $\cos\alpha$ \textup{(}например\textup{,} для $\alpha=2\pi/3)$.}

\begin{proof}
По формуле 4.1.5%\ref{3alpha}
\,(a) косинуса тройного угла каждое из чисел $\cos(2\pi/9)$, $\cos(8\pi/9)$, $\cos(14\pi/9)$ удовлетворяет уравнению $8y^3-6y+1=0$. Замена $x=2y$ превращает его в~уравнение $x^3-3x+1=0$. Значит, по теореме ни одно из них невозможно получить на вещественном калькуляторе.
\end{proof}

Перейдём теперь к~формулам, которые могут содержать комплексные числа.

%\begin{deff}
\smallskip
{\bf Определение комплексного калькулятора.}
\emph{Комплексный} калькулятор имеет те же кнопки, что и~вещественный, но оперирует с~комплексными числами и~при нажатии кнопки $\sqrt[n]{\phantom{n}}$ выдаёт все значения корня. На комплексном калькуляторе \emph{можно получить число}, если на нём можно получить множество чисел, содержащих заданное число.
\index{Калькулятор!комплексный|textbf}
%\end{deff}

\smallskip
Оказывается, уравнение третьей степени (например, $x^3-3x+1$), неразрешимое на вещественном калькуляторе, разрешимо на комплексном.

\smallskip
{\bf Утверждение о~разрешимости в~комплексных радикалах.}
{\it Все корни любого многочлена третьей или четвёртой степени с~рациональными коэффициентами можно получить на комплексном калькуляторе\textup{\footnote{Стандартная терминология: уравнение \emph{разрешимо в~радикалах}.}}. Более того\textup{,} это можно сделать так\textup{,} чтобы извлечение корня происходило только

%\begin{bul}
$\bullet$ %\item
два раза\textup{,} причём один раз третьей степени и~один раз второй "--- для многочлена третьей степени\textup{;}

$\bullet$ %\item
четыре раза\textup{,} причём один раз третьей степени и~три раза второй "--- для многочлена четвёртой степени.
%\end{bul}
}

\smallskip
Это утверждение доказывается \emph{методами дель Ферро и~Феррари} (оно вытекает из теорем, приведённых в~указании к~задачам \ref{2dferro} и~4.2.7%\ref{2ferrari}
).

Однако аналог этого утверждения для более высоких степеней неверен.

\begin{theorem}[Галуа]
Существует многочлен $5$"~й степени с~рациональными коэффициентами \textup{(}например\textup{,} $x^5-4x+2$\textup{),} ни один из корней которого невозможно получить на комплексном калькуляторе\textup{\footnote{Немного ранее была доказана более слабая теорема П.~Руффини"--~Н.~Х.~Абеля. Она сложнее формулируется~\cite{Al, FT, Sk11}, но именно она решила знаменитую проблему о~разрешимости уравнений в~радикалах.}}.
\index{Теорема!Галуа о~неразрешимости в~радикалах}
\end{theorem}

Из приведённых теорем о~неразрешимости тривиально следует, что \emph{для любого $n\ge3$ $(n\ge5)$ существует многочлен $n$"~й степени\textup{,} один из корней которого невозможно получить на вещественном \textup{(}комплексном\textup{)} калькуляторе.} Более сложно доказывается аналог этого утверждения с~заменой слов <<один из корней>> на <<ни один из корней>>. При этом корни \emph{некоторых} уравнений высоких степеней вполне могут получаться на калькуляторе, см.~например, п.\;\ref{gaucon}.

Теорема Галуа вытекает из следующего результата. Он интересен и~нетривиален даже для многочленов пятой степени.

\begin{theorem}[Кронекера] Если многочлен простой степени неприводим над $\Q,$ имеет более одного вещественного корня и~хотя бы один невещественный\textup{,} то ни один из его корней невозможно получить на комплексном калькуляторе.
\index{Теорема!Кронекера о~неразрешимости в~радикалах}
\end{theorem}

Эта теорема доказана в~п.\;5.5.4%\ref{progal}
. Для её доказательства необходимо следующее обобщение теоремы Гаусса (которое доказывается аналогично, см.~задачу~\ref{lower}). Обозначим
 $$
 \varepsilon_n:=\cos(2\pi/n)+i\sin(2\pi/n).
 $$

\begin{theorem}[Гаусса о~понижении]
\textup{(a)} Если $n$ простое\textup{,} то на комплексном калькуляторе можно получить $\varepsilon_n$ так\textup{,} чтобы корни извлекались только $(n-1)$"~й степени.

\textup{(b)} Для любого $n$ на комплексном калькуляторе можно получить~$\varepsilon_n$ так\textup{,} чтобы корни извлекались только степеней\textup{,} строго меньших $n$.
\end{theorem}

Вещественный аналог теоремы Кронекера следующий.

\begin{theorem}[Сильная вещественная теорема о~неразрешимости]
Если многочлен простой нечётной степени неприводим над $\Q$ и~имеет более одного вещественного корня\textup{,} то ни один из его корней невозможно получить на вещественном калькуляторе.
\end{theorem}

Доказательство этой <<вещественной>> теоремы приведено в~п.\;5.5.5%\ref{prostr}
.
Оно \emph{сложнее} доказательства <<комплексной>> теоремы Кронекера.

\subsubsection{План параграфа}\label{intplan}

Этот параграф не обязательно изучать подряд. Читатель может выбрать удобную ему последовательность изучения (или вовсе опустить некоторые пункты) на основании приводимого плана. К~плану разумно вернуться, если читатель потеряет нить изложения.

Пункт 5.2%\ref{mot} 
 независим от остального текста (т.\,е. он не используется в~остальном тексте и~для его изучения достаточно прочитать п.\;\ref{s:intr}).
В п.\;5.2.3%\ref{motcon}
 приводится переформулировка теоремы Гаусса (упомянутая в~п.\;\ref{intsqu}).

\newpage
\setcounter{page}{98}
\refstepcounter{subsection}
\subsection{Доказательство построимости в~теореме Гаусса}\label{gau}

В п.\;\ref{gauref} и~п.\;\ref{gaucon} доказана построимость в~теореме Гаусса.\break
В~п.\;\ref{gaures} идеи из п.\;\ref{gaucon} иллюстрируются на примерах и~задачах. Материал п.\;5.3.4%\ref{gaueff}
 не используется далее.

\subsubsection{Переформулировка построимости в~теореме Гаусса (2)}\label{gauref}

Начнём с~простых задач, подводящих к~основному результату этого пункта "--- лемме о~комплексификации.

\begin{pr}
Число $\cos(2\pi/n)$ вещественно построимо для $n=$3, 4, 5, 6, 8, 10, 15.
\end{pr}

\begin{pr}
{\bf Лемма об умножении (вещественная версия).}

(a) Если число $\cos(2\pi/n)$ вещественно построимо, то число\break $\cos(\pi/n)$ вещественно построимо.

(b) Если числа $\cos(2\pi/n)$ и~$\cos(2\pi/m)$ вещественно построимы и~$m,n$ взаимно просты, то число $\cos(2\pi/mn)$ вещественно построимо.
\end{pr}

Из этой леммы вытекает, что вещественная построимость в~теореме Гаусса следует из вещественной построимости чисел $\cos(2\pi/n)$ для простых $n$ вида $2^{2^s}+1$.

%\begin{deff}
{\bf Определение построимости.}
Комплексное число называется \emph{построимым}, если его можно получить на комплексном калькуляторе так,
чтобы при этом извлекались корни только второй степени.
\index{Число!построимое|textbf}
%\end{deff}

\begin{pr}\label{g-eps}
Число $\cos(2\pi/n)$ построимо тогда и~только тогда, когда число $\varepsilon_n:=\cos(2\pi/n)+i\sin(2\pi/n)$ построимо.
\end{pr}

\begin{pr}\label{g-com} (a)
{\bf Лемма о~комплексификации.}
Комплексное число построимо тогда и~только тогда, когда его вещественная и~мнимая части вещественно построимы.
\index{Число!вещественно построимое}

(b)$^*$ Можно ли получить число $e$ на калькуляторе, если уже есть число $e+\pi i$? (Используйте без доказательства тот факт, что числа $e$ и $\pi$ невозможно получить на калькуляторе.)
\end{pr}

Из этой леммы вытекает, что вещественное число построимо тогда и~только тогда, когда оно вещественно построимо\footnote{Заметим, что на комплексном калькуляторе нет кнопок $\Ree$ и~$\Imm$. Однако их можно <<реализовать>>, доказав, что если можно получить число $z$, то можно получить и~$\overline z$.
Но так будет доказана \emph{построимость} вещественной и~мнимой частей, а~не их \emph{вещественная построимость}. Для доказательства вещественной построимости нужно научиться извлекать корень из комплексного числа при помощи вещественного калькулятора. Это возможно только для корней второй степени. Если в~определении построимости и~вещественной построимости допускать извлечения корней третьей степени, то аналог леммы о~комплексификации будет неверен (ибо $\varepsilon_9\in\sqrt[3]{\sqrt[3]{1}}$ получается на комплексном калькуляторе с~извлечением корней только третьей степени, а~$\cos(2\pi/9)$ не получается на вещественном калькуляторе, см.~п.\;\ref{intrad}).}.
Значит, построимость в~теореме Гаусса достаточно доказать с~заменой <<вещественной построимости>> на <<построимость>>.

\subsubsection{Метод резольвент Лагранжа (2)}\label{gaures}

\index{Метод!резольвент Лагранжа}

\begin{pr} \label{postr7}
(a) Число $\varepsilon_5$ построимо.

(b) На комплексном калькуляторе можно получить число $\varepsilon_7$ так, что при этом извлекаются только корни второй и~третьей степени.

(b$')^*$ Можно ли получить на комплексном калькуляторе число $\varepsilon_7$ так, что при этом извлекаются только корни второй и~третьей степени, причём только по одному разу?

(c) На комплексном калькуляторе можно получить число $\varepsilon_{11}$ так, что при этом извлекаются только корни второй и~пятой степени.

(d) Число $\varepsilon_{17}$ построимо.
\end{pr}

Пункты (a), (b) и~(b$'$) можно решить непосредственно. Для решения пунктов (c), (d) уже нужна новая идея, изложенная далее. Оказывается, в~задачах \ref{postr7}\,(c, d) (и~во многих других ситуациях!) вместо работы с~набором корней удобнее работать с~некоторыми выражениями от корней "---
\emph{резольвентами Лагранжа}, которые мы скоро определим.

\begin{pr}\label{ident} Решите системы уравнений ($x,y,z,t$ "--- неизвестные, $a,b$, $c,d$ известны, $\varepsilon_3=\cfrac{-1+i\sqrt3}2$\,):
%AS убрал знаки препинания в конце
\[
\begin{gathered}
 \text{(a)}\ \ \left\{\arraycolsep=0pt\begin{array}{l} x+y+z+t=a,\\
x+y-z-t=b,\\
x-y+z-t=c,\\
x-y-z+t=d;
\end{array}\right.
\qquad
\text{(b)}\ \ \left\{\arraycolsep=0pt\begin{array}{l} x+y+z+t=a,\\
x+iy-z-it=b,\\
x-y+z-t=c,\\
x-iy-z+it=d;
\end{array}\right.\\
 \text{(c)} \ \ \left\{\arraycolsep=0pt\begin{array}{l}
 x+y+z=a,\\
 x+\varepsilon_3 y+\varepsilon_3^2z=b,\\
 x+\varepsilon_3^2y+\varepsilon_3 z=c.
 \end{array}\right.
\end{gathered}
\]
\end{pr}

Выражения из задачи \ref{ident} и~называются \emph{резольвентами Лагранжа}. Они <<лучше>> корней, поскольку <<симметричнее>> в~следующем смысле.

%AS вернул
\smallskip
\emph{Решение кубического уравнения при помощи резольвент Лагранжа.}
Для нахождения корней $x,y,z$ кубического уравнения достаточно найти выражения $a,b,c$ из задачи \ref{ident}\,(с). (Заметим, что метод дель Ферро из п.\;\ref{mot34} фактически приводит к~тому же результату.) По теореме Виета $a=a(x,y,z)$ "--- коэффициент уравнения. При замене $x\leftrightarrow y$ многочлен $b=b(x,y,z)$ переходит в~$\varepsilon_3 c$, а~$c=c(x,y,z)$ в~$\varepsilon_3^2 b$ (проверьте!). Значит, многочлены $bc$ и~$b^3+c^3$ не меняются при этой замене.
Аналогично они не меняются при замене $z\leftrightarrow y$. Поэтому многочлены $bc$ и~$b^3+c^3$ \emph{симметрические}, т.\,е. не меняется при любой перестановке переменных. Тогда из теоремы Виета и~теоремы о~представимости симметрического многочлена в~виде многочлена от элементарных симметрических многочленов (утверждение 4.6.3%\ref{repsymth}
\,(c)) следует, что эти многочлены от $x,y,z$ представляются в~виде многочленов от коэффициентов уравнения. Теперь, решая квадратное уравнение, можно получить $b^3$ и~$c^3$. Далее легко получить сами $b$ и~$c$.

%AS вернул
\smallskip
\emph{Решение уравнения $4$"~й степени при помощи резольвент Лагранжа.} Для нахождения корней $x,y,z,t$ уравнения 4"~й степени достаточно найти выражения $a,b,c,d$ от корней из задачи \ref{ident}\,(a). По теореме Виета $a$ "--- коэффициент уравнения. При замене $x\leftrightarrow y$ многочлены $c^2$ и~$d^2$ меняются местами, а~многочлен $b^2$ переходит в~себя. При циклической замене $x\rightarrow y\rightarrow z\rightarrow t\rightarrow x$ многочлены $b^2$ и~$d^2$ меняются местами, а~многочлен $c^2$ переходит в~себя. Значит, многочлены $b^2,c^2,d^2$ переставляются при любой перестановке переменных. Поэтому виетовские многочлены от них, т.\,е.
 $$
 b^2+c^2+d^2,\quad b^2c^2+b^2d^2+c^2d^2,\quad b^2c^2d^2,
 $$
симметрические. Тогда эти многочлены от $x,y,z$ представляются в~виде многочленов от коэффициентов уравнения. Теперь, решая кубическое уравнение, можно получить сами $b^2,c^2,d^2$. Далее легко получить $b,c,d$.

%AS вернул
\smallskip
Сообразите, почему же этот метод не работает для общего уравнения 5"~й степени!

\begin{pr}\label{5solh} (a) Если $x_1,\ldots,x_5$ --- корни многочлена $f\in\Q[x]$ 5-й степени, то
 $$
 T(y):=\prod\limits_{\tau\in S_5}
 (y-x_{\tau(1)}-\varepsilon_5x_{\tau(2)}-\varepsilon_5^2x_{\tau(3)}-\varepsilon_5^3x_{\tau(4)}-\varepsilon_5^4 x_{\tau(5)}) \in \Q[\varepsilon_5][y].
 $$

(b) Для некоторого $G\in\Q[\varepsilon_5][y]$ выполнено равенство $T(y)=G(y^5)$. Такой многочлен $G$ называется {\it разрешающим многочленом} для~$f$.

(c)* На комплексном калькуляторе можно получить все корни разрешающего многочлена для $f(x)=x^5+15x+11$ (а~значит, и~самого многочлена~$f$).
\end{pr}

\begin{pr}\zvezda\label{5sol}
На комплексном калькуляторе можно получить хотя бы один корень уравнения $x^5+ax+b=0$ для

(a) $(a,b)=(15,11)$; \quad (b) $(a,b)=(-5,52)$; \quad
(c) $(a,b)=(35,36)$; \quad

(d) $(a,b)=\frac{(15\pm20c,44\mp8c)}{c^2+1}$, $c\in\Q$, $c\ge0$. \quad

(Для других $(a,b)$ этого сделать нельзя \cite{PrSo}.)
\end{pr}

\begin{proof}[Идея доказательства построимости числа $\varepsilon:=\varepsilon_5$]
Во-первых,
 $$
 T_0:=\varepsilon+\varepsilon^2+\varepsilon^4+\varepsilon^8=-1.
 $$
Сначала докажем построимость числа
 $$
 T_2:=\varepsilon-\varepsilon^2+\varepsilon^4-\varepsilon^8.
 $$
При замене $\varepsilon$ на $\varepsilon^2$ число $T_2$ переходит в~$-T_2$. Значит, $T_2^2$ не меняется при этой замене. Поэтому $T_2^2$ не меняется при двукратной и трехкратной таких заменах, т.\,е. при заменах $\varepsilon$ на $\varepsilon^4$ и $\varepsilon$ на $\varepsilon^8=\varepsilon^3$. Итак, для любого $k$ число $T_2^2$ не меняется при замене $\varepsilon$ на $\varepsilon^k$.

Раскроем скобки в произведении $T_2^2$ и заменим $\varepsilon^5$ на 1.
Получим равенство
$$
T_2^2=a_0+a_1\varepsilon+a_2\varepsilon^2+a_3\varepsilon^3+a_4\varepsilon^4\quad\text{для некоторых}\quad a_k\in\Z.
$$
Так как для любого $k$ число $T_2^4$ не меняется при замене $\varepsilon$ на $\varepsilon^k$, то $a_1=a_2=a_3=a_4$. Поэтому $T_2^2=a_0-a_1\in \Z$. Значит, $T_2$ построимо.

Обозначим
$$
T_1:=\varepsilon+i\varepsilon^2-\varepsilon^4-i\varepsilon^8\qquad\text{и}\qquad
T_3:=\varepsilon-i\varepsilon^2-\varepsilon^4+i\varepsilon^8.
$$
Тогда $T_0+T_1+T_2+T_3=4\varepsilon$.
Поэтому достаточно доказать построимость чисел $T_1$ и $T_3$.
Сделаем это для $T_1$; доказательство для $T_3$ аналогично.

При замене $\varepsilon$ на $\varepsilon^2$ число $T_1$ переходит в~$-iT_1$. Значит, $T_1^4$ при этой замене не меняется. Поэтому $T_1^4$ не меняется при двукратной и~трёхкратной замене такого вида, т.\,е. при замене $\varepsilon$ на $\varepsilon^4$ и~$\varepsilon$ на $\varepsilon^8=\varepsilon^3$. Итак, для любого $k$ число $T_1^4$ не меняется при замене $\varepsilon$ на $\varepsilon^k$.

Раскроем скобки в~произведении $T_1^4$ и~заменим $\varepsilon^5$ на~1. Получим равенство
 $$
 T_1^4=a_0+a_1\varepsilon+a_2\varepsilon^2+a_3\varepsilon^3+a_4\varepsilon^4\quad\text{для некоторых}\quad a_k\in\Z+i\Z.
 $$
Так как для любого $k$ число $T_1^4$ не меняется при замене $\varepsilon$ на $\varepsilon^k$, то $a_1=a_2=a_3=a_4$. Поэтому $T_1^4=a_0-a_1\in \Z+i\Z$. Значит, $T_1$ построимо.
\noqed\end{proof}

В приведённом рассуждении нужно обосновать вывод равенства $a_1=a_2=a_3=a_4$ и~строго определить, что такое <<замена $\varepsilon$ на $\varepsilon^2$>>. Обоснование для общего случая трудное; читатель может найти пример такого рассуждения в~\cite[\S\,24]{E1}. Поэтому вместо того, чтобы его приводить, мы немного изменим доказательство; именно этим изменением приводимое доказательство отличается от данного в~\cite{E1}, \cite[\S\,6.4]{PrSo}. Для этого вместо того, чтобы работать с~\emph{числами}, мы будем работать с~\emph{многочленами} и~подставлять в~них $\varepsilon$ в~качестве аргумента.

\begin{pr}\label{gaures-i}
Обозначим $T_1(x):=x+ix^2-x^4-ix^8$. Тогда\footnote{Два многочлена называются \emph{сравнимыми по модулю многочлена $x^5-1$}, если их разность делится на $x^5-1$.\index{Сравнимость!многочленов|textbf}}

(a) $iT_1(x^2)\equiv T_1(x)\bmod(x^5-1)$;

(b) $T_1^4(x^2)\equiv T_1^4(x)\bmod(x^5-1)$;

(c) $T_1^4(x^k)\equiv T_1^4(x)\bmod(x^5-1)$ для любого $k$.
\end{pr}

\begin{proof}[Доказательство построимости числа $\varepsilon:=\varepsilon_5$]
Определим многочлен $T_1(x):=x+ix^2-x^4-ix^8$.
Определим многочлены $T_0(x)$, $T_2(x)$ и~$T_3(x)$ формулами, аналогичными вышенаписанным.
Как и~выше, $(T_0+T_1+T_2+T_3)(\varepsilon)=4\varepsilon$.
Поэтому достаточно доказать построимость каждого из чисел $T_r(\varepsilon)$, $r=1,2,3$.
Имеем
 $$
\begin{gathered}
 iT_1(x^2)\underset{x^5-1}\equiv T_1(x)\quad\Longrightarrow\quad T_1^4(x^2)\underset{x^5-1}\equiv T_1^4(x)\quad\Longrightarrow\\
\Longrightarrow\quad T_1^4(x^k)\underset{x^5-1}\equiv T_1^4(x)\quad\text{для любого }k.
\end{gathered}
$$
Возьмём многочлен $a_0+a_1x+a_2x^2+a_3x^3+a_4x^4$ с~коэффициентами в~$\Z+i\Z$,
сравнимый с~$T_1^4(x)$ по модулю $x^5-1$.

Тогда $a_1=a_2=a_3=a_4$.
Поэтому\footnote{Другой способ, предложенный М.~Ягудиным:
\begin{multline*}
 T_1^4(\varepsilon)=a_0+a_1\varepsilon+a_2\varepsilon^2+a_3\varepsilon^3+a_4\varepsilon^4=
 a_0+a_1\varepsilon^2+a_2\varepsilon^4+a_3\varepsilon+a_4\varepsilon^3={}\\
 =a_0+a_1\varepsilon^3+a_2\varepsilon+a_3\varepsilon^4+a_4\varepsilon^2=
 a_0+a_1\varepsilon^4+a_2\varepsilon^3+a_3\varepsilon^2+a_4\varepsilon.
\end{multline*}
Суммируя эти выражения, получим $4T_1^4(\varepsilon)=a_0-a_1-a_2-a_3-a_4\in \Z+i\Z$.}
$T_1^4(\varepsilon)=a_0-a_1\in \Z+i\Z$.

Значит, $T_1(\varepsilon)$ построимо. Аналогично $T_2(\varepsilon)$ и~$T_3(\varepsilon)$ построимы.~\end{proof}

\newpage
\subsubsection{Доказательство построимости в~теореме Гаусса (3)}\label{gaucon}

Читатель, изучивший (точнее, прорешавший) два предыдущих пункта, подготовлен к~доказательству. Напомним, что формально оно независимо от двух предыдущих пунктов.

\smallskip
{\bf Лемма об умножении.} {\it \textup{(a)} Если $\varepsilon_n$ построимо\textup{,} то $\varepsilon_{2n}$ построимо.

\textup{(b)} Если $\varepsilon_n$ и~$\varepsilon_m$ построимы и~$m,n$ взаимно просты\textup{,} то $\varepsilon_{mn}$ построимо.}

\begin{proof}[Доказательство\nopoint] получается из формул $\varepsilon_{2n}\in\sqrt{\varepsilon_n}$ и~$\varepsilon_{mn}=\varepsilon_m^x\varepsilon_n^y$, где $x$ и~$y$ "--- целые числа, для которых $nx+my=1$.
\end{proof}

При решении задач \ref{postr7}\,(a) мы использовали различность остатков от деления чисел $2,2^2,2^3,2^4$ на~5. При решении задач \ref{postr7}\,(c, d) и~5.3.10%\ref{7}
\,(a) мы использовали аналогичное свойство чисел 2 и~11, 6 и~17, 3 и~7. Для общего случая необходимо следующее обобщение.

\begin{theorem}[о~первообразном корне] Для любого простого $p$ существует число $g,$ для которого остатки от деления на $p$ чисел $g^1,g^2,g^3,\ldots,g^{p-1}$ различны.
\end{theorem}

\emph{Указание к~доказательству для $p=2^m+1$ \textup{(}только этот случай нужен для теоремы Гаусса\textup{)}.} Если первообразного корня нет, то сравнение $x^{2^{m-1}}\equiv1\bmod p$ имеет $p-1=2^m>2^{m-1}$ решений.
Это противоречит теореме Безу.

Заинтересованный читатель может получить и~полное доказательство, см.~п.\;\ref{s:numbroo}.

\begin{proof}[Доказательство построимости в~теореме Гаусса]
\index{Теорема!Гаусса о~построимости правильных многоугольников}
По лемме \ref{g-com} о~комплексификации и~по лемме об умножении достаточно доказать, что $\varepsilon_n$ построимо для любого простого $n=2^{2^s}+1$.
Так как $n-1=2^m$, то по лемме об умножении $\beta:=\varepsilon_{n-1}$ построимо.
Обозначим
 $$
 \Z[\beta]:=\{a_0+a_1\beta+a_2\beta^2+\ldots+a_{n-2}\beta^{n-2}\ |\ a_0,\ldots,a_{n-2}\in\Z\}.
 $$
Обозначим через $g$ первообразный корень по модулю $n$. Для $r=0,1,2,\ldots,n-2$, обозначим
 $$
 T_r(x):=x+\beta^rx^g+\beta^{2r}x^{g^2}+\ldots+\beta^{(n-2)r}x^{g^{n-2}}\in\Z[\beta][x].
 $$
Тогда $(T_0+T_1+\ldots+T_{n-2})(\varepsilon)=(n-1)\varepsilon.$ Кроме того, $T_0(\varepsilon)=-1$. Поэтому достаточно доказать построимость каждого из чисел $T_r(\varepsilon)$, $r=1,2,\ldots,n-2$. Имеем
 $$
\begin{gathered}
 \beta^rT_r(x^g)\underset{x^n-1}\equiv T_r(x)\quad\Longrightarrow\quad T_r^{n-1}(x^g)\underset{x^n-1}\equiv T_r^{n-1}(x) \quad\Longrightarrow\\
 \Longrightarrow\quad T_r^{n-1}(x^k)\underset{x^n-1}\equiv T_r^{n-1}(x)\quad\text{для любого }k.
\end{gathered}
$$
Возьмём многочлен $a_0+a_1x+a_2x^2+\ldots+a_{n-1}x^{n-1}$ с~коэффициентами в~$\Z[\beta]$, сравнимый с~$T_r^{n-1}(x)$ по модулю $x^n-1$. Тогда
$a_1=a_2=\ldots=a_{n-1}$. Поэтому $T_r^{n-1}(\varepsilon)=a_0-a_1\in \Z[\beta]$.
Значит, $T_r(\varepsilon)$ построимо.~\end{proof}

\begin{pr}\zvezda\label{lower} Докажите теорему Гаусса о~понижении.
\end{pr}

\subsubsection*{Подсказки}

\paragraph*{\ref{lower}.}
(a) Аналогично доказательству построимости в~теореме Гаусса.

(b) Докажем теорему при помощи индукции по $n$.

Если $n=ab$ для некоторых целых $a,b$, $0<a,b<n$, то шаг индукции следует из равенства $\varepsilon_n=\sqrt[a]\varepsilon_b$.

Если же $n$ простое, то шаг индукции следует из п.\,(a).

\newpage

\setcounter{page}{114}
\refstepcounter{subsection}
\subsubsection{Одно извлечение квадратного корня (1)}\label{non12}

%\Opensolutionfile{_hintAholder}
%\Opensolutionfile{_hintBholder}

Перед решением задач этого пункта полезно прорешать п.\;4.1%\ref{s:ratir}
.

\begin{pr}\label{number2} Представимо ли следующее число в~виде $a+\sqrt b$, где $a,b\in\Q$:

\parbox{2.7cm}{(a) $\sqrt{3+2\sqrt2}$;}
(b) $\frac1{7+5\sqrt2}$;

\parbox{2.7cm}{(c) $\sqrt[3]{7+5\sqrt2}$;}
(d) $\cos(2\pi/5)$; \quad
(e) $\sqrt[3]2$; \quad
(f) $\sqrt2+\sqrt[3]2$; \quad

\parbox{2.7cm}{(g) $\cos(2\pi/9)$;}
(h)$^*$\ $\sqrt{2+\sqrt2}$; \quad
(i)$^*$\ $\cos(2\pi/7)$.

\end{pr}

Задачи \ref{number2} и~\ref{utv2} интересны в~связи с~неразрешимостью в~радикалах, поскольку нам нужно придумать многочлен, корни которого невозможно получить на калькуляторе, а~числа из задачи \ref{number2} являются корнями многочленов (подумайте, каких).

\begin{pr}\label{sopr2}
Пусть $r\in\R-\Q$ и~$r^2\in\Q$.

(a) \textbf{Лемма о~неприводимости.} Многочлен $x^2-r^2$ неприводим над~$\Q$.

(b) \textbf{Лемма о~линейной независимости.} Если $a,b\in\Q$ и~$a+br=0$, то $a=b=0$.

(c) Если многочлен имеет корень $r$, то этот многочлен делится на $x^2-r^2$.

(d) \textbf{Теорема о~сопряжении.}
Если многочлен имеет корень $r$, то корнем этого многочлена является также число $-r$.
\index{Теорема!о сопряжении}

(e) \textbf{Следствие.}
Если $a,b\in\Q$ и~многочлен имеет корень $a+br$, то корнем этого многочлена является также число $a-br$.

(f) \textbf{Следствие.} Если $a,b\in\Q$ и~кубический многочлен имеет корень $a+br$, то он имеет рациональный корень.

\end{pr}

\begin{pr}\label{utv2}
\textbf{Утверждение.} Если многочлен степени выше второй неприводим над $\Q$, то ни один из его корней не представим в~виде $a\pm\sqrt b$, где $a,b\in\Q$.
\end{pr}

\begin{pr}\label{calc2}
\textbf{Лемма о~калькуляторе.} Пусть $F\in\{\R,\C\}$.
Число, которое можно получить на $F$"~калькуляторе так, чтобы извлечение корня происходило только один раз, причём второй степени, имеет вид $a\pm\sqrt b$, где $a,b\in\Q$ и~если $F=\R$, то $b>0$.
\end{pr}

\newpage
\setcounter{page}{137}
\subsection{Доказательства неразрешимости в~радикалах}\label{pro}

Читатель, изучивший (точнее, прорешавший) предыдущий параграф, подготовлен к~доказательствам.
Напомним, что формально они независимы от предыдущего параграфа.
Для понимания этого параграфа достаточно прочитать п.\;\ref{intsqu} и~\ref{intrad}
(кроме того, в~п.\;\ref{progau} понадобится простая лемма \ref{g-com} о~комплексификации).
%(она доказана в~\ref{gausol}),
%а в~\S\,\ref{progal} "--- теорема Гаусса \ref{lower}\,(a) о~понижении
%(она доказывается аналогично доказательству построимости в~теореме Гаусса, приведённому в~\ref{gaucon}).
%??? док-во теоремы Гаусса о~понижении.
%Пункты \ref{progau}, \ref{prorea} и~\ref{progal} формально независимы друг от друга.
%Неформально же, для каждого следующего пункта полезно прочитать предыдущий.
\emph{План доказательства} в~каждом пункте получается из текста пункта пропуском доказательств лемм.
%В п.\;\ref{prostr} приводится дополнительный материал.
%важнейшие идеи доказательства неразрешимости в~радикалах продемонстрированы на вещественном случае.

\subsubsection{Лемма о~калькуляторе и~понятие поля (2*)}\label{profie}

Если $F\subset\C$, $r\in\C$ и~$r^q\in F$ для некоторого целого положительного $q$, то обозначим
 $$
 F[r]:=\{a_0+a_1r+a_2r^2+\ldots+a_{q-1}r^{q-1}\ |\ a_0,\ldots,a_{q-1}\in F\}.
 $$

\smallskip
{\bf Лемма о~калькуляторе.} {\it Пусть $F\in\{\R,\C\}$. Число $x\in F$ можно получить на $F$"~калькуляторе тогда и~только тогда\textup{,} когда существуют такие $r_1,\ldots r_{s-1}\in F$ и~такие простые $q_1,\ldots q_{s-1},$ что
 $$
 \Q=F_1\subset F_2\subset F_3\subset \ldots\subset F_{s-1}\subset F_s\ni x,
 $$
где $r_k^{q_k}\in F_k,$ $r_k\not\in F_k$ и~$F_{k+1}=F_k[r_k]$ для любого $k=1,\ldots,s-1$.}
\index{Калькулятор!вещественный}
\index{Калькулятор!комплексный}

\smallskip
Такая последовательность называется \emph{башней \textup{(радикальных)} расширений}.

Эта лемма доказывается несложно (аналогично леммам о~калькуляторе из п.\;5.4%\ref{non}
).

В этом пункте \emph{полем} называется подмножество множества $\C$, замкнутое относительно операций сложения, умножения, вычитания и~деления на ненулевое число. Общепринятое название: числовое поле (а~\emph{полем} в~математике называется немного другой объект). Это понятие полезно для нас тем, что теорема делении с~остатком верна для многочленов с~коэффициентами в~поле.

Если $F$ "--- поле, $q$ простое, $r\not\in F$ и~либо $F=\Q$, либо $\varepsilon_q\in F$, то многочлен $t^q-r^q$ неприводим над $F$ (это фактически доказано в~леммах о~линейной независимости далее в~п.\;\ref{prorea} и~5.5.4%\ref{progal}
 аналогично лемме о~неприводимости 5.4.22%\ref{idea}
 \,(a) и~задаче 5.4.27%\ref{mailemtow}
 ). Тогда $F[r]$ "--- поле.

Напомним, что $\varepsilon_n:=\cos\frac{2\pi}n+i\sin\frac{2\pi}n.$

\subsubsection{Доказательство непостроимости в~теореме Гаусса (3*)}\label{progau}
\index{Теорема!Гаусса о~построимости правильных многоугольников}

Так как $\varepsilon_n=\varepsilon_{nk}^k$, то из построимости числа $\varepsilon_{nk}$ вытекает
построимость числа $\varepsilon_n$. Поэтому и~по лемме \ref{g-com} о~комплексификации для доказательства вещественной непостроимости в~теореме Гаусса достаточно показать, что $\varepsilon_n$ непостроимо для

(A) простого числа $n$, не представимого в~виде $2^m+1$;

(B) квадрата простого числа, т.\,е. $n=p^2$.

\smallskip
{\bf Лемма о~степенях двойки.} {\it Если неприводимый над $\Q$ многочлен $P$ с~рациональными коэффициентами имеет построимый корень\textup{,} то $\deg P$ есть степень двойки.}

\smallskip
Эта лемма доказана далее. Для её применения нужны следующие результаты.

\smallskip
{\bf Признак Эйзенштейна.} {\it Пусть $p$ простое. Если для многочлена с~целыми коэффициентами старший коэффициент не делится на $p,$ остальные делятся на $p,$ а~свободный член не делится на $p^2,$ то этот многочлен неприводим над $\Z$.}

\smallskip
{\bf Лемма Гаусса.} {\it Если многочлен с~целыми коэффициентами неприводим над $\Z,$ то он неприводим и~над~$\Q$.}

\smallskip
И признак Эйзенштейна, и~лемма Гаусса легко доказываются переходом к~многочленам с~коэффициентами $\Z_p$ (для леммы Гаусса рассмотрим разложение $P=P_1P_2$ данного полинома $P$ над $\Q$, возьмём такие целые $n_1$ и~$n_2$, что и~$n_1P_1$, и~$n_2P_2$ имеют целые коэффициенты, и~возьмём простой делитель $p$ числа $n_1n_2$).

\begin{proof}[Доказательство непостроимости числа $\varepsilon_n$]
Непостроимость чис\-ла $\varepsilon_n$ следует из леммы о~калькуляторе и~леммы о~степенях двойки для корня $\varepsilon_n$ многочлена

%\begin{bul}

$\bullet$ %\item
$P(x):=x^{n-1}+x^{n-2}+\ldots+x+1$ в~случае (A) и~

$\bullet$ %\item
$P(x):=x^{p(p-1)}+x^{p(p-2)}+\ldots+x^p+1$ в~случае (B).

%\end{bul}

Неприводимость этих многочленов $P(x)$ над $\Q$ вытекает из их неприводимости над $\Z$ и~леммы Гаусса.
Неприводимость этих многочленов $P(x)$ над $\Z$ вытекает из неприводимости многочленов $P(x+1)$ над $\Z$.
Последняя неприводимость доказывается применением признака Эйзенштейна.
Выполнение предположения <<если>> для многочленов $P(x+1)$ легко проверяется с~помощью сравнения
$(a+b)^p\equiv a^p+b^p\bmod p$.
\end{proof}

Лемма о~степенях двойки является случаем $k=1$ следующего утверждения.

\smallskip
{\bf Обобщённая лемма о~степенях двойки.} {\it Если
 $$
 \Q=F_1\subset F_2\subset F_3\subset \ldots\subset F_{s-1}\subset F_s\ni \alpha,
 $$
где $r_k^2\in F_k,$ $r_k\not\in F_k$ и~$F_{k+1}=F_k[r_k]$ для любого $k=1,\ldots,s-1,$ то для каждого $k=1,2,\ldots,s$ степень любого неприводимого над $F_k$ многочлена с~коэффициентами из $F_k$ и~корнем $\alpha$ есть степень двойки.}

\begin{proof}
Индукция по $k$ вниз. База $k=s$ очевидна. Докажем шаг. Для $j\in\{k,k+1\}$ обозначим через $P_j$ любой неприводимый над $F_k$ многочлен с~коэффициентами из $F_k$ и~корнем $\alpha$. Предположение индукции заключается в~том, что $\deg P_{k+1}$ есть степень двойки. Будем рассматривать делимость и~НОД в~$F_{k+1}$.
Так как $P_k(\alpha)=0$, то многочлен $P_k$ делится на $P_{k+1}$.

Будем использовать следующий простой результат.

\smallskip
{\bf Лемма о~сопряжении.} {\it Пусть $F\subset\C$ "--- поле\textup{,} $r\in\C-F$ и~$r^2\in F$. Определим отображение сопряжения $\overline\cdot\,\colon F[r]\to F[r]$ формулой\break $\overline{x+yr}{:=}x-yr.$ Это отображение корректно определено\textup{,} $\overline{z+w}=\overline z+\overline w$ и~$\overline{zw}=\overline z\cdot\overline w$.}

\smallskip
Применим эту лемму о~сопряжении к~$F=F_k$ и~$F[r]=F_{k+1}$. Получим, что $P_k=\overline{P_k}$ делится на $\overline{P_{k+1}}$. Обозначим $D:=\gcd(P_{k+1},\overline{P_{k+1}})$. Так как многочлен $P_{k+1}$ неприводим над $F_{k+1}$ и~делится на $D$, то либо $D=1$, либо $P_{k+1}=D$.

Во втором случае из того, что $\overline D=D$, получаем, что $P_{k+1}=D\in F_k[x]$. Значит, $P_k=P_{k+1}$ и~шаг индукции доказан.

В первом случае $P_k$ делится на $M:=P_{k+1}\overline{P_{k+1}}$. Из того, что $\overline M =M$, следует, что $M\in F_k[x]$. Так как многочлен $P_k$ неприводим над $F_k$, то $P_k=M$. Значит, $\deg P_k=2\deg P_{k+1}$ есть степень двойки.
\end{proof}

 %\newpage
\subsubsection{Доказательство неразрешимости в~вещественных радикалах (3*)}\label{prorea}
\index{Теорема!о неразрешимости в~вещественных радикалах}

\smallskip
{\bf Основная лемма (вещественный случай).}
{\it Пусть $q$ простое, $F\subset\R$ "--- поле, $r\in\R-F$ и~$r^q\in F$.

\textup{(a) (линейная независимость).}
Если $P(r)=0$ для некоторого многочлена $P\in F[\varepsilon_q][t]$ степени меньше $q,$ то $P=0$.

\textup{(b) (сопряжение).} Если $P\in F[\varepsilon_q][t]$ и~$P(r)=0,$ то $P(r\varepsilon_q^k)=0$ для любого $k=0,1,\ldots,q-1$.}
\index{Теорема!о сопряжении}

\begin{proof}[Доказательство части \textup{(a)}]
Оба многочлена $P$ и~$t^q-r^q$ с~коэффициентами из $F[\varepsilon_q]$ имеют корень~$r$. Значит, их НОД имеет корень $r$ и~степень $k$, $0<k\le\deg P<q$.
Все корни многочлена $t^q-r^q$ есть $r,r\varepsilon_q,r\varepsilon_q^2,\ldots,r\varepsilon_q^{q-1}$. Свободный член наименьшего общего делителя равен произведению некоторых $k$ из этих корней. Тогда $r^k\in F[\varepsilon_q]$. Так как $q$ простое, то $kx+qy=1$ для некоторых целых $x,y$.
Тогда $r = (r^k)^x (r^q)^y\in F[\varepsilon_q]$.

Поэтому\footnote{Другая запись этого абзаца с~использованием понятия размерности: тогда $\dim_F F[r]\le \dim_F F[\varepsilon_q]\le q-1$.}
 $r^2,r^3,\ldots,r^{q-1}\in F[\varepsilon_q]$. Составим таблицу $a_{kl}\in F$ размера $q\times (q-1)$ из разложений чисел $r^k$ по степеням числа $\varepsilon_q$:
 $$
 r^k=\sum\limits_{l=0}^{q-2}a_{kl}\varepsilon_q^l,\quad 0\le k\le q-1.
 $$
При помощи нескольких операций прибавления к~одной строке другой, умноженной на число из $F$, можно получить таблицу с~нулевой строкой.

Значит, имеется ненулевой многочлен $P_1$ степени меньше $q$ с~коэффициентами из $F$ и~корнем~$r$. Дальнейшие рассуждения аналогичны первому абзацу этого доказательства. Нужно только заменить $P$ на $P_1$ и~$F[\varepsilon_q]$ на $F$.
Получаем, что $r\in F$, "--- противоречие.

\begin{proof}[Доказательство части \textup{(b)}]
Так как $P(r)=0$, то остаток от деления многочлена $P(t)$ на $t^q-r^q$ принимает значение 0 в~точке~$r$. Значит, по части (a) этот остаток равен нулю. Отсюда вытекает заключение части~(b).
\noqed\end{proof}

\emph{Доказательство теоремы о~неразрешимости в~вещественных радикалах.}
Предположим, напротив, что некоторый корень $x_0$ уравнения $x^3-3x+1=0$ можно получить на вещественном калькуляторе. Тогда по лемме о~калькуляторе для $F=\R$ существует наименьшее~$s$, для которого найдётся башня расширений, последнее поле $F_s\subset\R$ которой содержит некоторый корень $x_1$ уравнения $x^3-3x+1=0$ (возможно, $x_1\ne x_0$). Обозначим $F:=F_{s-1}$, $q:=q_{s-1}$ и~$r:=r_{s-1}$. Тогда $x_1=h(r)$ для некоторого многочлена $h$ с~коэффициентами в~$F$ степени больше 0 и~меньше~$q$.

Применим часть (b) основной леммы (вещественный случай) к~многочлену $P(t):=h(t)^3-3h(t)+1$. Так как $h(r)^3-3h(r)+1=0$, то $h(r\varepsilon_q^k)$ является корнем уравнения $x^3-3x+1=0$ для любого $k=0,1,\ldots,q-1$. Если $h(r\varepsilon_q^k)=h(r\varepsilon_q^l)$ для некоторых $k,l$, $0\le k<l\le q-1$, то по части (a) основной леммы (вещественный случай) получим, что $\deg h=0$ "--- противоречие. Итак, числа $h(r\varepsilon_q^k)$, $0\le k\le q-1$, "--- попарно различные корни уравнения $x^3-3x+1=0$. Значит, $q=2$ или $q=3$.

Если $q=2$, то по теореме Виета третий корень уравнения $x^3-3x+1=0$ равен $-2h(0)\in F$ "--- противоречие с~минимальностью числа~$s$.

Если $q=3$, то возьмём $h_0+h_1t+h_2\in F$, для которых $h_0+h_1t+h_2t^2=h(t)$. Так как
 $$
 h(r\varepsilon_3)\in\{2\cos(2\pi/9),2\cos(8\pi/9),2\cos(14\pi/9)\}\subset\R,
 $$
то $0=\Imm h(r\varepsilon_3)=\frac{\sqrt3}2(h_1r-h_2r^2)$. Так как $r\not\in F$, то $h_1=h_2=0$. Противоречие с~неравенством $\deg h>0$\footnote{\emph{Другое завершение доказательства для $q=3$.}
Если $q=3$, то из равенство $\overline{\varepsilon_3}=\varepsilon_3^2$ вытекает, что $\overline{h(r\varepsilon_3)}=h(r\varepsilon_3^2)$. Это противоречит вещественности и~различности последних двух чисел.}.
\end{proof}

\newpage
\setcounter{page}{163}
\refstepcounter{section}
\refstepcounter{subsection}
\refstepcounter{subsection}
\subsection{Применения основных неравенств (3*). \emph{М.~А.~Берштейн}}\label{s:ineqsim}
\sectionmark{Неравенства}

Автор благодарен А.~Берштейну, А.~Дудко, В.~Карайко, К.~Кнопу и~В.~Франку, которые научили его почти всему, что здесь написано.

\begin{pr}\label{ineqsim-pow}
Для целых $a, b, c>0$ выполнено неравенство
 $$
 \Big(\frac{a^2+b^2+c^2}{a+b+c}\Big)^{a+b+c} \ge a^ab^bc^c \ge\Big(\frac{a+b+c}{3}\Big)^{a+b+c}.
 $$
\end{pr}

\begin{pr}\label{ineqsim-cau}
(a) \textbf{Весовое неравенство Коши.} Если $a_1>0,\ldots,a_n>0$ и~$a_1+\ldots+a_n=1$, то $a_1x_1+\ldots+a_nx_n\ge x_1^{a_1} \cdot \ldots \cdot x_n^{a_n}$.\index{Неравенство!Коши!весовое}
(Это обобщение неравенства Юнга 6.2.5%\ref{ineqsim-jng}
\,(b).)

(b)$^*$ Определим \emph{среднее степенное порядка $m$ чисел $x_1,\ldots,x_n$ с~весами $a_1,\ldots,a_n>0$}, где $a_1+\ldots+a_n=1$, как
\[
\begin{gathered}
 S_m:=\root m\of{a_1x_1^m+\ldots+a_nx_n^m}\quad\text{при}\ \ m\ne0,\quad
 S_0:=x_1^{a_1}\cdot\ldots\cdot x_n^{a_n},\\
 S_{-\infty}:=\min\{x_1,\ldots,x_n\}\quad \text{и}\quad S_{+\infty}:=\max\{x_1,\ldots,x_n\}.
\end{gathered}
\]
Докажите, что $S_a\le S_b$ при $a\le b$ для любых $a,b\in\R\cup\{-\infty,+\infty\}$.

(c) Остаётся ли справедливым неравенство $S_a\le S_b$ при $a\le b$ и~любых положительных значениях переменных $x_1,\ldots,x_n$, если одно из чисел $a_i$ меньше нуля?
\end{pr}

\begin{pr}\label{ineqsim-cyc} Докажите неравенства

(a) $\frac{a_1^2}{a_2}+\frac{a_2^2}{a_3}+\ldots+\frac{a_n^2}{a_1} \ge a_1+a_2+\ldots+a_n$;

(b) $\frac{a_1^2}{a_1+a_2}+\frac{a_2^2}{a_2+a_3}+\ldots+
\frac{a_n^2}{a_n+a_1} \ge \frac{1}{2}(a_1+a_2+\ldots+a_n).$
\end{pr}

\begin{pr}\label{ineqsim-2} Докажите неравенство
 $$
 \frac{a^2}{b(a+c)}+\frac{b^2}{c(b+d)}+\frac{c^2}{d(c+a)}+\frac{d^2}{a(d+b)} \ge 2.
 $$
\end{pr}

\begin{pr}\label{ineqsim-3mur} Докажите неравенства

(a) $a^3b+b^3c+c^3a \ge abc(a+b+c);$

(b) $a^3b^2+b^3c^2+c^3a^2 \ge abc(ab+bc_ca)$.
\end{pr}

\comment
\section{Последовательности и~пределы}\label{s:seq}

Этот параграф почти независим от остальной части книги. В~других местах из него используются лишь простые факты.
\endcomment

\newpage
\setcounter{page}{175}
\refstepcounter{section}
\subsection{Конечные суммы и~разности (3)}\label{s:seqfin}

\emph{Последовательностью сумм} последовательности $\{a_n\}_{n=1}^\infty$
называется последовательность $b_n=\Sigma a_n\!:=\!a_1\!+\!\ldots\!+\!a_n$, а~\emph{последовательностью разностей} "--- последовательность $c_n\!=\!\Delta a_n\!:=\!a_{n+1}\!-\!a_n$.
\index{Последовательность!сумм|textbf}%
\index{Последовательность!разностей|textbf}%

\vspace{-0.3cm}
Например, $\Delta 2^n=2^n$ и~$\Sigma2^n=2^{n+1}-2$.

(Сумма и~разность "--- аналоги \emph{интеграла} и~\emph{производной}.)

В этом пункте $n$ обозначает номер члена последовательности, <<по которому>> берётся сумма и~разность. Так, например, $\Delta 2^k=0$.

\begin{pr}\label{seqfin-del}
Найдите

(a) $\Delta n^k$ для каждого целого $k\ge-1$; \quad
(b) $\Delta\cos n$; \quad
(c) $\Delta(n\cdot2^n)$.
\end{pr}

\begin{pr}\label{seqfin-sum}
Найдите

(a) $\Sigma\sin n$; \quad
(b) $\Sigma\frac1{n(n+1)\ldots(n+k)}$ для каждого целого $k>0$.
\end{pr}

\begin{pr}\label{seqfin-que}
Какие из указанных равенств выполняются для некоторой непостоянной последовательности $a_n$:

\parbox{2.5cm}{(a) $\Delta a_n=0$;}
\parbox{2.75cm}{(b) $\Delta a_n=1$;}
(c) $\Delta a_n=a_n$; \quad

\parbox{2.5cm}{(d) $\Sigma a_n=a_n$;}
\parbox{2.75cm}{(e) $\Sigma\Delta a_n=a_n$;}
(f) $\Delta\Sigma a_n=a_n$?
\end{pr}

\begin{pr}\label{seqfin-lei}
(a) Найдите $\sum\limits_{k=0}^n(-1)^kk^2\binom{n}{k}$.

(b) {\bf Лемма.} $k$-я разность многочлена $k$"~й степени есть постоянная, а~$(k+1)$-я равна 0.

(c) (Загадка.) Выразите $\Delta^ka_n$ через $a_n,a_{n+1},\ldots,a_{n+k}$.

(d) {\bf Лемма.} Равенство $\Delta^ka_n=0$ имеет место тогда и~только тогда, когда $a_n$ "--- многочлен от $n$ степени не выше $k-1$.

(e) Для некоторого многочлена $P_\lambda(n)$,
имеющего степень $l$ при $\lambda\ne1$ и~степень $l-1$ при $\lambda\ne1$, выполняется равенство  $\Delta(n^l\lambda^n)=P_\lambda(n)\lambda^n$.

(f) {\bf Формула Лейбница.} Справедливо равенство
 $$
 \Delta(a_nb_n)=a_{n+1}\Delta b_n+b_n\Delta a_n.
 $$
\index{Формула!Лейбница}%

\vspace{-0.6cm}
(g)$^*$ Сформулируйте и~докажите аналогичную формулу для $\Delta^l(a_nb_n)$.
\end{pr}

\begin{pr}\label{seqfin-pow}
(a) Найдите $\Sigma n^k$ для $k=0,1,2,3,4$.

(b) {\bf Лемма.} Последовательность сумм многочлена степени $k\ge0$ есть многочлен степени $k+1$.
\end{pr}

\begin{pr}\label{seqfin-qua}
(a) Найдите $\Sigma(n\cdot2^n)$.

\end{pr}

\newpage
\setcounter{page}{182}
\refstepcounter{subsection}
\subsection{Конкретная теория пределов (4*)}\label{s:analim}

Задачи этого пункта интересны не только как простейший способ разобраться в~теории пределов. Похожие задачи о~конкретных, хотя и~грубых оценках часто возникают и~на олимпиадах, и~в~прикладной математике, и~в~теоретической математике.

В решении этих задач нельзя пользоваться функциями $\sqrt[n] x$, $a^x$, $\log_ax$, $\arcsin x$ и~т.\,п. без определения этих функций (поскольку для их определения "--- например, для доказательства существования такого~$x$, что
$x^2=2$, "--- фактически нужно эти задачи решить). Исключение: если некоторая функция используется в~условии, то её можно использовать и~в~решении.
Можно пользоваться без доказательства свойствами неравенств.

\begin{pr}\label{analim} Найдите хотя бы одно такое $N$, чтобы для любого $n>N$ выполнялось неравенство $a_n>10^9$, если

(a) $a_n=\sqrt n$; \quad
(b) $a_n=n^2-3n+5$; \quad
(c) $a_n=1{,}02^n$;

(d) $a_n=1+\frac12+\frac13+\frac14+\ldots+\frac1n$.
\end{pr}

%%%!!!
\begin{pr}\label{analim-bern}
{\bf Неравенство Бернулли.} Докажите, что $(1+x)^a\ge 1+ax$ для любых $x\ge-1$ и

(a) целого $a\ge1$; \quad
(b) рационального $a\ge1$;
\quad

(c) действительного $a\ge1$.
\index{Неравенство!Бернулли}
\end{pr}

\begin{pr}\label{analim-an}
Найдите хотя бы одну пару таких $a$ и~$N$, чтобы для любого $n>N$ выполнялось неравенство $|a_n-a|<10^{-8}$, если

(a) $a_n=\frac{n^2-n+28}{n-2n^2}$; \quad
(b) $\sqrt{5+\dfrac{2}{n}}$; \quad
(c) $a_n=n\bigg(\!\sqrt{1+\frac1n}-1\!\bigg)$;

\end{pr}

\newpage
\setcounter{page}{195}
\refstepcounter{subsection}
\refstepcounter{subsection}
\refstepcounter{subsection}
\subsection{Примеры трансцендентных чисел (3*)}\label{s:trans}

\begin{pr}\label{trans-irr}
Следующие числа иррациональны:

(a) $e:=\sum\limits_{n=0}^{\infty}\frac1{n!}$;
\quad
(b) $\lambda:=\sum\limits_{n=0}^{\infty}2^{-n!}$;
\quad
(c) $\mu:=\sum\limits_{n=0}^{\infty}2^{-2^n}$.

(Используемые здесь бесконечные суммы определены в~п.\;7.5%\ref{s:anasum}
.)
\end{pr}

\begin{pr}\label{trans-qua}
(e), ($\lambda$), ($\mu$)
Ни одно из чисел $e$, $\lambda$, $\mu$ не является корнем квадратного уравнения с~целыми коэффициентами.
\end{pr}

\begin{pr}\label{p7-7-3} Для любого рационального числа $p/q$, не являющегося корнем многочлена $f$ степени $t$ с~целыми коэффициентами, выполнено неравенство $|f(p/q)|\ge q^{-t}$.
\end{pr}

Число $x$ называется \emph{трансцендентным}, если оно не является корнем уравнения $a_tx^t+a_{t-1}x^{t-1}+\ldots+a_1x+a_0=0$ с~целыми коэффициентами $a_t\ne0$, $a_{t-1},\ldots,a_0$.
\index{Число!трансцендентное|textbf}

\begin{pr}\label{lambda}
(a) \textbf{Теорема Лиувилля.} Число $\lambda$ трансцендентно.

(b) {\bf Общая теорема Лиувилля.} Для любых многочлена степени $t$ с~рациональными коэффициентами и~его иррационального корня $\alpha$ существует такое $C>0$, что для любых целых $p,q$ выполнено неравенство $\Big|\alpha-\frac pq\Big|>Cq^{-t}$.
\index{Теорема!Лиувилля}
\end{pr}

\begin{pr}\label{nu}
(a) Число $\mu$ не является корнем кубического уравнения с~целыми коэффициентами.

(b) Справедливо равенство $\mu^q=\sum\limits_{n=0}^\infty d_n(q)2^{-n}$,
где $d_n(q)$ есть количество упорядоченных представлений числа $n$ в~виде
суммы $q$ степеней двойки (не обязательно различных степеней):
 $$
 d_n(q)=\#\{(w_1,\ldots,w_q)\in\Z^q\mid n=2^{w_1}+\ldots+2^{w_q}\text{ и~}w_1,\ldots,w_q>0\}.
 $$
Например, $d_3(2)=2$, поскольку $3=2^0+2^1=2^1+2^0$. По определению полагаем $d_0(0)=1$.

(c) Справедливо неравенство $d_n(q)\le (q!)^2$.

(d) Число $\mu$ трансцендентно.
\end{pr}

Следующая задача "--- удачная тема для исследовательских работ старшеклассников. Описание удачных примеров этой деятельности читатель может найти в~материалах Московской математической конференции школьников~[M]%\cite{M}
. Пункты (a), (b), (c) заведомо не претендуют на научную новизну. Решение остальных пунктов мне неизвестно, но наверняка доступно сильному старшекласснику.

\begin{pr}\label{unso} Является ли число $\sum\limits_{n=0}^{\infty}a_n$ трансцендентным, если

(a) $a_n=2^{-3^n}$;

(b) $a_n=d_n2^{-2^n}$ для некоторой ограниченной последовательности $d_n>0$ целых чисел;

(c) $a_n=2^{-f_n}$, где $f_{n+2}=f_{n+1}+f_n$, $f_0=f_1=1$ "--- последовательность Фибоначчи;

(d)$^*$ $a_n=2^{-[1,1^n]}$;
\quad
(e)$^*$ $a_n=n2^{-2^n}$;
\quad
(f)$^*$ $a_n=2^{n-2^n}$;
\quad

(g)$^*$ $a_n=(-1)^n2^{-2^n}$.
\end{pr}

\subsubsection*{Подсказки}

\paragraph*{\ref{trans-irr}.}
(a) Предположим, напротив, что существует линейный многочлен
$f(x)=bx+c$ с~целыми коэффициентами $b\ne0$ и~$c$, для которого $f(e)=0$. Обозначим $e_s=\sum\limits_{n=0}^s\frac1{n!}$. Заметим, что уравнение $bx+c=0$ имеет только один корень, значит, $f(e_s)\ne0$. Мы получим противоречие из следующих неравенств для $s=2|b|$:
 $$
 \frac1{s!} \le |f(e_s)| = |f(e)-f(e_s)| =  |b|\cdot(e-e_s)< \frac{2|b|}{(s+1)!}.
 $$

(b, c) Обобщение вышеприведённого доказательства (см.~указание к~п.\,(a)) работает.

\comment
\section{Функции}\label{s:fun}

Этот параграф почти независим с~остальной частью книги.
В~других местах из него используются лишь простые факты.

В этом параграфе, если не оговорено противное, под \emph{многочленом} понимается многочлен с~\emph{вещественными} коэффициентами и~латинские буквы обозначают \emph{вещественные} числа.
\endcomment

\newpage
\refstepcounter{section}
\setcounter{page}{199}
\subsection{График кубического многочлена (2)}
\sectionmark{Функции}

Известно, что график любого квадратного трёхчлена имеет ось симметрии.

\begin{pr}\label{2center}
(a) График любого кубического многочлена имеет центр симметрии.

(b) Найдите координаты центра симметрии графика функции $y=-2x^3-6x^2+4$.

(c) Верно ли, что график любого многочлена 4"~й степени имеет ось симметрии?
\end{pr}

Известно, что квадратное уравнение $ax^2+bx+c=0$ имеет два решения при $D>0$, имеет одно решение при $D=0$ и~не имеет решений при $D<0$.
Здесь $D=b^2-4ac$. Способ нахождения \emph{количества} решений кубического уравнения \emph{без решения самого уравнения} легко вывести напрямую (задача \ref{2numbroot} ниже). В~частности, для решения следующих задач не требуется знать формулы для корней кубического уравнения; более того, решения, не использующие этих формул, \emph{проще} вывода указанных формул. Ср. с~задачами 8.1.7%\ref{numro4}
, \ref{desca} ниже.

\begin{pr}\label{2numroo}
Сколько (вещественных) решений имеет уравнение

(a) $x^3+2x+7=0$; \quad (b) $x^3-4x-1=0$?
\end{pr}

\begin{pr}\zvezda
\textbf{Теорема о~промежуточном значении.}
Для многочлена $f$ и~чисел $a<b$ если $f(a)>0>f(b)$, то существует такое $c\in(a,b)$, что $f(c) =0$.
\end{pr}
\index{Теорема!о промежуточном значении}

Этой теоремой можно пользоваться в~дальнейшем без доказательства.

\begin{pr}\label{numeasy}
(a) Уравнение $x^3+x+q=0$ имеет ровно одно решение при любом $q$.

(b) При каком условии на $p,q$ уравнение $x^3+px+q=0$ имеет ровно два решения?

(c) Выразите эти два решения через $p,q$.
\end{pr}

\begin{pr}\label{2numbroot}
(a) Для функции $f(x)=x^3-6x+2$ найдите промежутки возрастания и~убывания.

%(b) Для той же функции найдите наибольшее и~наименьшее значения на отрезке $[0,3]$.

%(c) При каких $q$ уравнение $x^3-x+q=0$ имеет ровно одно решение?
\end{pr}

\newpage
\setcounter{page}{203}
\subsection{Элементы анализа для многочленов (2)}\label{s:funpol}

\begin{pr}\label{desca}
(a) {\bf Правило знаков Декарта.} Число \emph{положительных} решений уравнения $p_nx^n + \ldots + p_1x + p_0 = 0$ не превосходит числа перемен знака
в последовательности $p_0,\ldots, p_n$.
%, из которой вычеркнуты нули.
\index{Правило!знаков Декарта}

(b) Как аналогично правилу знаков Декарта оценить количество \emph{отрицательных} корней данного многочлена?

(c)$^*$ Как аналогично правилу знаков Декарта оценить количество корней данного многочлена на данном промежутке $[a,b]$?

(d) {\bf Неравенства Маклорена.} Для $x_1,\ldots,x_n>0$ обозначим
 $$
 M_k=\sqrt[k]{\frac{\sum\limits_{i_1<\ldots<i_k}x_{i_1}\cdot\ldots\cdot x_{i_k}}{\binom{n}{k}}}.
 $$
(Заметьте, что $M_1$ "--- это среднее арифметическое и~$M_n$ "--- среднее геометрическое.) Тогда $M_1\ge \ldots \ge M_n$.
\end{pr}

Для решения этих и~многих других задач необходимо следующее понятие.

\emph{Производной} $f'$ многочлена $f$ называется многочлен, полученный подстановкой $y = x$ в~многочлен $\frac{f(y) - f(x)}{y - x}$ от двух переменных $x$, $y$. (Сообразите, почему это многочлен.)
\index{Производная!многочлена|textbf}

Геометрический смысл: уравнение \emph{касательной} к~графику многочлена $f$ в~точке $a$ есть $y = f'(a)(x-a)+f(a)$. (Формально, это можно воспринимать как определение касательной.)

\begin{pr}\label{difrul} Приведите чёткие формулировки утверждений, скрывающиеся за следующими формулами, и~докажите их:

(a) $(f + g)' = f' + g'$; \quad
(b) $(af)'=af'$; \quad
(c) $(x^n)' = nx^{n-1}$;

(d) $(p_nx^n + \ldots + p_1x + p_0)' = np_nx^{n-1}+(n-1)p_{n-1}x^{n-2} + \ldots + p_1$ (при $n = 0$ это выражение равно~$0$).

(e) {\bf Формула Лейбница.} Справедливо равенство $(fg)' = f'g + fg'$.
\index{Формула!Лейбница}
\end{pr}

\newpage
\setcounter{page}{213}
\refstepcounter{subsection}
\refstepcounter{subsection}
\refstepcounter{subsection}
\subsection{Применения компактности (4*). \emph{А.~Я.~Канель-Белов}}\label{s:funcomp}

В этом пункте задачи посложнее и~подсказок поменьше. Однако он будет интересен читателю, так как, насколько нам известно, такая подборка интересных задач
по этой важной теме впервые публикуется в~неспециальной литературе.

\begin{pr}\label{funcomp-cha}
\textbf{Близкая идея в~конечном случае.} Запись числа состоит из нулей и~единиц. Любой фрагмент <<10>> числа заменяют на <<0001>>. Докажите, что рано или поздно заменять будет нечего.
\end{pr}

\begin{pr}\label{funcomp-ide}
\textbf{Идея компактности.} (a) Известно, что человечество живёт вечно, а~число людей в~каждом поколении конечно. Докажите, что найдётся бесконечная мужская цепочка наследников.

(b) В~бесконечном парламенте у~каждого парламентария не более трёх врагов. Докажите, что парламент можно разбить на две палаты так, что у~каждого
парламентария будет не более одного врага в~своей палате. (Для конечного парламента эта задача разбирается в~задаче \ref{pri1-7-7} пункта <<Полуинварианты>>.)

(c) Известно, что любую \emph{конечную} карту на плоскости можно правильно раскрасить в~$4$ цвета. Докажите, что тогда \emph{произвольную} карту на плоскости также можно правильно раскрасить в~$4$ цвета. (Страны можно считать многоугольниками. Раскраска называется \emph{правильной}, если
любые две страны с~общим участком границы раскрашены в~разные цвета.)

(d) (Загадка.) Прочитайте п.\;\ref{ss:polu} <<Полуинварианты>>. Какие утверждения верны для бесконечных множеств, а~какие нет?
\end{pr}

\begin{pr}\label{p8-6-3}
Для любых $M$ и~$k$ найдётся достаточно большое $v$ с~таким свойством: если все рёбра полного графа с~$v$ вершинами покрашены в~$M$ цветов, то найдётся полный подграф с~$k$ вершинами, все рёбра которого покрашены в~один цвет.
\end{pr}

\begin{pr}\label{p8-6-4}
Из любой бесконечной последовательности целых чисел можно выбрать подпоследовательность либо так, чтобы каждый её член делился на предыдущий, либо так, чтобы ни один член не делился на другой.
\end{pr}

\begin{pr}\label{p8-6-5}
На плоскости отмечено бесконечное множество точек, никакие три из которых не лежат на одной прямой. Тогда найдётся выпуклая фигура, граница которой проходит через его бесконечное подмножество.
\end{pr}

Идеи, с~помощью которых доказываются остановки процессов, зачастую
работают вместе с~идеей компактности, с~которой они являются своего рода родственниками.

\begin{pr}\label{p8-6-6}
Существует ли такое $n$, что любое рациональное число между 0 и~1 представимо в~виде $\sum\limits_{i=1}^n \frac{1}{a_i}$, где $0<a_i \in\Z$?
\end{pr}

\endgroup

\begingroup

\chapter{Геометрия}
%\addcontentsline{toc}{section}{Геометрия}
\setcounter{page}{219}

Как правило, параграфы и~пункты этой главы можно изучать независимо друг от друга и~от остальных частей книги. В~тех случаях, когда для решения задач какого-нибудь пункта желательно знакомство с~другими материалами, это указывается в~начале пункта. Если задачу можно решать разными методами, она приводится в~пункте, посвящённом одному из них, а~о возможности других решений говорится в~комментарии. Помимо обозначений, принятых во всей книге, в~данной главе везде, где не оговорено обратное, используются принятые в~геометрии обозначения элементов треугольника, описанные в~начале параграфа <<Треугольник>>.

\section{Треугольник}\label{treug}

\emph{Всюду в~данной главе\textup{,} кроме специально оговорённых случаев\emph{,} используются следующие обозначения}: $ABC$ "--- данный треугольник, $A_i$, $B_i$, $C_i$, $i=1$, $2$, $\ldots,$ "--- точки на сторонах $BC$,$CA$ и~$AB$ соответственно (или на продолжениях этих сторон, если
это оговорено в~условии задачи); $\omega$ "--- вписанная
окружность, $I$ "--- её центр, $r$ "--- её радиус; $\Omega$ "---
описанная окружность, $O$ "--- её центр, $R$ "--- её радиус; $G$ "--- точка пересечения медиан (центр тяжести, центроид), $H$ "---
точка пересечения высот (ортоцентр). Проведём биссектрисы $AI$,
$BI$, $CI$ до пересечения с~$\Omega$ в~точках $A'$, $B'$, $C'$
соответственно. Таким образом, $A'$, $B'$, $C'$ "--- середины дуг
$AB$, $BC$, $CA$. \emph{Ортотреугольник} "--- треугольник
с~вершинами в~основаниях высот, \emph{серединный треугольник} "---
треугольник с~вершинами в~серединах сторон данного треугольника.
Перпендикуляр, опущенный из точки~$A$ на $BC$, обозначается
$h(A,BC)$.\index{Окружность!вписанная}\index{Окружность!описанная}%
\index{Центр!тяжести}\index{Центроид}\index{Ортоцентр}%
\index{Ортотреугольник}\index{Треугольник!серединный}%

{\it Окружностью $ABC$} называется окружность, описанная вокруг треугольника $ABC$.

\subsection{Принцип Карно (1). \emph{В.~Ю.~Протасов}, \emph{А.~А.~Гаврилюк}}\label{s9.1}

%\Opensolutionfile{_hintBholder}

\begin{pr}\label{pr1-1-1} \textbf{Теорема Карно.}\index{Теорема!Карно} В~точках $A_1$, $B_1$, $C_1$, лежащих на сторонах треугольника $ABC$ или на их
продолжениях, восставлены перпендикуляры к~этим сторонам. Докажите,
что они пересекаются в~одной точке тогда и~только тогда, когда
\[
C_1A^2 - C_1B^2 + A_1B^2 - A_1C^2 + B_1C^2 - B_1A^2 = 0.
\]
\end{pr}

\begin{pr}\label{pr1-1-2}
Сформулируйте и~докажите обобщённую теорему Карно для произвольных точек плоскости $A_1$, $B_1$, $C_1$, не обязательно лежащих на прямых, содержащих стороны треугольника $ABC$.
\end{pr}

\begin{pr}\circpr\label{pr1-1-11}
В~каком из следующих случаев перпендикуляры, восставленные
к~сторонам треугольника в~указанных точках, могут не пересекаться
в~одной точке:

%\begin{arab}
%\item
1) $A_1$, $B_1$, $C_1$ "--- точки касания сторон с~вписанной
окружностью;

%\item
2) $A_2$, $B_2$, $C_2$ "--- точки касания сторон с~соответствующими
вневписанными окружностями;

%\item
3) $A_3$, $B_3$, $C_3$ "--- основания биссектрис треугольника?
%\end{arab}
\end{pr}

\begin{pr}\label{pr1-1-3}
Пусть вневписанная окружность треугольника касается его стороны $AB$ в~точке $C_1$ и~касается продолжений двух других сторон. Аналогично определяются точки $A_1$ и~$B_1$. Докажите, что перпендикуляры, восставленные к~сторонам треугольника в~точках $A_1$, $B_1$, $C_1$ пересекаются в~одной точке.
\end{pr}

\begin{pr}\label{pr1-1-4}
На плоскости даны три пересекающиеся окружности. Докажите, что три
их общие хорды пересекаются в~одной точке.

\emph{Примечание}. Это утверждение обычно доказывают, используя
понятие степени точки (см.~п.\;10.4%\ref{s10.4}
). Однако его
легко вывести и~из обобщённой теоремы Карно.
\end{pr}

\newpage
\setcounter{page}{227}
\refstepcounter{subsection}
\refstepcounter{subsection}
\subsection{Формула Карно (2$^*$). \emph{А.~Д.~Блинков}}\label{s9.4}

%\Opensolutionfile{_hintAholder}
%\Opensolutionfile{_hintBholder}

\textbf{Формула Карно} (по имени французского математика, физика и~политического деятеля Лазаря Карно, 1753--1823) утверждает, что в~остроугольном треугольнике сумма расстояний от центра описанной
окружности до сторон треугольника равняется сумме радиусов описанной
и вписанной окружностей, т.\,е. $OM_1+OM_2+OM_3=R+r$, где $M_1$,
$M_2$, $M_3$ "--- середины $BC$, $CA$, $AB$ соответственно. Её
доказательство с~помощью теоремы Птолемея приводится в~п.\;\ref{s10.6}~<<Теоремы Птолемея и~Кези>>. Здесь мы рассмотрим её применения и~ещё один
способ её доказательства, в~процессе которого будут получены другие
важные факты.

\begin{pr}\label{pr1-4-1}
Пусть биссектриса угла $A$ пересекает окружность, описанную около
треугольника $ABC$, в~точке $W$, а~точка $D$ диаметрально
противоположна точке $W$. Докажите, что

(a) $M_1W=(r_a-r)/2$;

(b) $M_1D=(r_b+r_c)/2$, где $r$, $r_a$, $r_b$, $r_c$ "--- радиусы
вписанной и~вневписанных окружностей.
\end{pr}

\begin{pr}\label{pr1-4-2}
Докажите формулу Карно.

Рассмотрим теперь несколько задач на применение формулы Карно. Если
явно не оговорено обратное, то треугольник, заданный в~условии, остроугольный.
\end{pr}

\begin{pr}\label{pr1-4-3}
Докажите, что сумма расстояний от вершин треугольника до ортоцентра
равна сумме диаметров его вписанной и~описанной окружностей.
\end{pr}

\begin{pr}\label{pr1-4-4}
Докажите, что в~треугольнике $ABC$ выполняются неравенства

(a) $AH+BH+CH\leq 3R$;

(b) $3OH\geq R-2r$.
\end{pr}

\begin{pr}\label{pr1-4-5}
(a) Докажите, что $m_a+m_b+m_c\leq\frac92R$, где $m_a$, $m_b$ и~$m_c$ "--- длины медиан треугольника.

(b) Пусть в~треугольнике $ABC$ биссектрисы углов $A$, $B$ и~$C$
пересекают описанную окружность в~точках $W_1$, $W_2$ и~$W_3$
соответственно. Докажите, что $AW_1+BW_2+CW_3\leq 6{,}5R-r$.
\end{pr}

\begin{pr}\label{pr1-4-6}
(a) Докажите, что для углов треугольника выполняется неравенство
$$
\frac{3r}{R}\leq\cos A+\cos B+\cos C\leq\frac32.
$$

(b) Пусть $AH_1$, $BH_2$ и~$CH_3$ "--- высоты треугольника $ABC$.
Выразите сумму диаметров окружностей, описанных около треугольников
$AH_2H_3$, $BH_1H_3$ и~$CH_1H_2$, через $R$ и~$r$.
\end{pr}

\begin{pr}\label{pr1-4-7}
В окружность радиуса $R$ вписан треугольник, а~в каждый сегмент,
ограниченный стороной треугольника и~меньшей из дуг окружности,
вписана окружность наибольшего возможного радиуса. Найдите сумму
диаметров трёх получившихся окружностей и~радиуса окружности,
вписанной в~треугольник.
\end{pr}

\begin{pr}\label{pr1-4-8}
(a) Докажите, что в~треугольнике $ABC$ выполняется равенство
$$
a(OM_2+OM_3)+b(OM_1+OM_3)+c(OM_1+OM_2)=2pR.
$$

(b) \textbf{Неравенство Эрдёша.} Пусть $h_a$ "--- наибольшая высота
треугольника $ABC$. Докажите, что $h_a\geq R+r$.\index{Неравенство!Эрдёша}
\end{pr}

\begin{pr}\label{pr1-4-9}
(a) Выведите аналоги формулы Карно для прямоугольного и~тупоугольного
треугольников.

(b) Четырёхугольник $ABCD$ вписанный. Пусть $r_1$ и~$r_2$ "---
радиусы окружностей, вписанных в~треугольники $ABC$ и~$ADC$, а~$r_3$
и $r_4$ "--- радиусы окружностей, вписанных в~треугольники $ABD$ и~$CBD$. Докажите, что $r_1+r_2=r_3+r_4$.
\end{pr}

\begin{pr}\label{pr1-4-10}
Пусть $d$, $d_1$, $d_2$ и~$d_3$ "--- расстояния от центра $O$
окружности, описанной около треугольника, до центров его вписанной и~вневписанных окружностей. Докажите, что
$$
R^2=\frac{d^2+d_1^2+d_2^2+d_3^2}{12}.
$$
\end{pr}

\begin{pr}\label{pr1-4-11}
(a) Докажите, что если точка принадлежит отрезку, соединяющему
основания двух биссектрис треугольника, то сумма расстояний от этой
точки до двух сторон треугольника равна расстоянию от неё до третьей
стороны.

\comment
(b) Пусть центр окружности, описанной около треугольника, лежит на
отрезке, соединяющем основания двух биссектрис. Докажите, что
расстояние от ортоцентра треугольника до одной из его вершин равно
$R+r$.
\endcomment
\end{pr}

\newpage
\setcounter{page}{242}
\refstepcounter{subsection}
\refstepcounter{subsection}
\refstepcounter{subsection}
\subsection{<<Полувписанная>> окружность (2$^*$). \emph{П.~А.~Кожевников}}\label{s9.8}

%\Opensolutionfile{_hintAholder}
%\Opensolutionfile{_hintBholder}

Пусть $A'$ и~$A''$ "--- середины дуг $BC$ описанной окружности
$\Omega$, соответственно не содержащей и~содержащей точку $A$; $B'$
и~$B''$, $C'$ и~$C''$ определяются аналогично.

Рассмотрим окружность $S_A$ (назовём её \emph{полувписанной}),
касающуюся сторон $AB$, $AC$ и~окружности $\Omega$ (внутренним
образом). Основными в~этой серии являются следующие факты:
\index{Окружность!полувписанная}

"--* прямая, проходящая через точки касания полувписанной окружности со сторонами, содержит точку~$I$;

"--* точка касания полувписанной окружности с~окружностью $\Omega$
лежит на прямой~$A''I$.

\subsection*{Основная серия-1}

Докажите следующие утверждения.

\begin{pr}\label{pr1-8-1}
Пусть перпендикуляр к~биссектрисе $AI$, проведённый через точку~$I$,
пересекает $AB$ и~$AC$ в~точках $K$ и~$L$ соответственно. Тогда
окружности $BKI$, $CLI$ и~$\Omega$ пересекаются в~одной точке~$T$.

%\Ol{fig1.eps заменить шрифт!!!!!!!!!!!!!!!!!!!}
\end{pr}

\begin{pr}\label{pr1-8-2}
Точки $T$, $I$, $A''$ лежат на одной прямой.
\end{pr}

\begin{pr}\label{pr1-8-3}
Точки $T$, $K$, $C'$ лежат на одной прямой.
\end{pr}

\begin{pr}\label{pr1-8-4}
Точки $K$, $L$ и~$T$ являются точками касания окружности $S_A$ с~прямыми $AB$, $AC$ и~окружностью $\Omega $.
\end{pr}

\begin{pr}\label{pr1-8-5}
(a) Прямая $CC'$ касается окружности $TBKI$.

(b) Точка $T$ "--- центр поворотной гомотетии, переводящей треугольник
$BKI$ в~треугольник $ILC$. \index{Гомотетия!поворотная}
\end{pr}

\subsection*{Основная серия-2}

\begin{pr}\label{pr1-8-6}
Прямая $AT$ проходит через центр гомотетии с~положительным коэффициентом,
переводящей окружность $\omega $ в~$\Omega$.
\end{pr}

\begin{pr}\label{pr1-8-7}
Пусть $A_1$ и~$A_2$ "--- точки касания вписанной и~вневписанной
окружностей со стороной $BC$ соответственно. Тогда

(a) $AA'$ "--- биссектриса угла $TAA_2$;

(b) $\angle BTA_1 = \angle ABC$. %(\emph{Задача \textup{4.7.7} из \textup{\cite{sb1}}}.)
\end{pr}

\begin{pr}\label{pr1-8-8}
Пусть $AT$ пересекает $KL$ в~точке $Z$. Тогда $\angle BZK = \angle
CZL$.%(\emph{Задача \textup{4.7.5} из \textup{\cite{sb1}}}.)
\end{pr}

\begin{pr}\label{pr1-8-9}
Прямые $KL$, $TA'$ и~$BC$ пересекаются в~одной точке или параллельны. (\emph{И.~Шарыгин.})
\end{pr}

\begin{pr}\label{pr1-8-10}
Точка пересечения $Y_A$ из предыдущей задачи и~точки $Y_B$,~$Y_C$,
определённые аналогичным образом, лежат на одной прямой.
\end{pr}

\subsection*{Дополнительные задачи-1}

\begin{pr}\label{pr1-8-11}
Пусть $P$ "--- произвольная точка на дуге $BA'C$.

(a) Пусть $P_b = BB'\cap PC'$, $P_c = CC'\cap PB'$. Тогда окружность
$PP_bP_c$ проходит через $T$. %(\emph{\textup{\cite{olshar2}}\textup{,} задача \textup{8.8, 2013}~г.})

(b) Пусть $J_b$ и~$J_c$ "--- центры вписанных окружностей
треугольников $PAB$ и~$PAC$. Тогда окружность $PJ_bJ_c$ проходит
через $T$. %(\emph{Задача \textup{4.7.9} из \textup{\cite{sb1}}}.)

(c) Пусть касательные к~$\omega$ из точки $P$ пересекают $BC$ в~точках $U_1$ и~$U_2$. Тогда окружность $PU_1U_2$ проходит через~$T$.
%(\emph{Задача \textup{4.7.10} из \textup{\cite{sb1}}}.)

(d) Пусть прямые, проходящие через $I$ параллельно биссектрисам углов
между прямыми $AP$ и~$BC$ пересекают прямую $BC$ в~точках $V_1$ и~$V_2$ Тогда окружность $PV_1V_2$ проходит через $T$. %(\emph{См.~частный случай задачи \textup{4.7.18} из \textup{\cite{sb1}}.})
\end{pr}

\subsection*{Дополнительные задачи-2}

Следующие задачи "--- про <<обобщённые полувписанные>> окружности, т.\,е. окружности, касающиеся двух прямых и~окружности.

\begin{pr}\label{pr1-8-12}
Пусть $D$ "--- точка на стороне $AC$ треугольника $ABC$, и~пусть $S_1$ "--- окружность, касающаяся окружности $\Omega$ внутренним образом в~точке $R$, а~также отрезков $BD$ и~$AD$ в~точках $M$ и~$N$ соответственно.

(a) Докажите, что точки $B$, $M$, $I$, $R$ лежат на одной окружности.

(b) \textbf{Лемма Саваямы.} Прямая $MN$ проходит через центр $I$
вписанной окружности $\omega$ треугольника~$ABC$.
\index{Лемма!Саваямы}
\end{pr}

\newpage
\setcounter{page}{249}
\subsection{Обобщённая теорема Наполеона (2$^*$). \emph{П.~А.~Кожевников}}\label{s9.9}

%\Opensolutionfile{_hintAholder}
%\Opensolutionfile{_hintBholder}

Классическая теорема Наполеона гласит, что центры правильных
треугольников, построенных на сторонах произвольного треугольника
вне его, являются вершинами равностороннего треугольника.
\index{Теорема!Наполеона}\index{Теорема!Наполеона!обобщённая}

Теорема Наполеона является частным случаем утверждения задачи 13.1.7% \ref{pr2-20-7}
, п.~<<Комплексные числа и~геометрия>>. В~этом пункте мы докажем
другое обобщение теоремы Наполеона.

Предлагаем для решения серию задач, внешне не имеющих никакой связи
с~теоремой Наполеона. Можно решать задачи любыми методами, а~затем
познакомиться с~обобщением теоремы Наполеона и~получить решения
задач как следствия этого сильного факта.

\subsection*{Вводные задачи}

\begin{pr}\label{pr1-9-1}
Докажите, что центры квадратов, построенных на сторонах параллелограмма вне его, являются вершинами квадрата.
\end{pr}

\begin{pr}\label{pr1-9-2}
На боковых сторонах трапеции $ABCD$ построены треугольники $ABE$
и~$CDF$ так, что $AE\parallel CF$ и~$BE\parallel DF$. Докажите, что
если $E$ лежит на стороне $CD$, то $F$ лежит на стороне $AB$.
\end{pr}

\begin{pr}\label{pr1-9-3}
(a) Две окружности пересекаются в~точках $A$ и~$B$. Через точку $A$
проведена прямая, вторично пересекающая первую окружность в~точке
$C$, а~вторую "--- в~точке $D$ (можно считать, что точки $C$ и~$D$
лежат по разные стороны от точки $A$). Пусть $M$ и~$N$ "--- середины
дуг $BC$ и~$BD$, не содержащих точку $A$, а~$K$ "--- середина
отрезка $CD$. Докажите, что угол $MKN$ прямой. (\emph{Д.~Терёшин. Всероссийская математическая олимпиада \textup{1997}~г.})

(b) Круг поделили хордой $AB$ на два круговых сегмента и~один из них
повернули вокруг точки $A$ на некоторый угол. Пусть при этом
повороте точка $B$ перешла в~точку $D$. Докажите, что отрезки,
соединяющие середины дуг сегментов с~серединой отрезка $BD$,
перпендикулярны друг другу. %(\emph{З.~Насыров}, \cite{Kvant92-2}.)
\end{pr}

\begin{pr}\label{pr1-9-4}
Через вершину~$A$ треугольника~$ABC$ проведены прямые~$l_1$ и~$l_2$,
симметричные относительно биссектрисы угла~$A$. Докажите, что
проекции точек~$B$ и~$C$ на~$l_1$ и~$l_2$ соответственно, середина
стороны~$BC$ и~основание высоты, опущенной из~вершины~$A$, лежат
на~одной окружности.
\end{pr}

\begin{pr}\label{pr1-9-5}
Во вписанном четырёхугольнике $ABCD$ диагонали пересекаются в~точке
$E$, точки $K$ и~$M$ "--- середины сторон $AB$ и~$CD$, точки $L$
и~$N$ "--- проекции точки $E$ на $BC$ и~$AD$. Докажите, что $KM\perp LN$.
\end{pr}

\begin{pr}\label{pr1-9-6}
По двум окружностям, пересекающимся в~точках $P$ и~$Q$, одновременно
начали движение с~равными угловыми скоростями из точки $P$ два
велосипедиста $A$ и~$B$: один едет по часовой стрелке, другой "---
против часовой стрелки. Докажите, что $A$ и~$B$ всё время
равноудалены от фиксированной точки. %(\emph{Задача о~велосипедистах, случай движения в~разные стороны. Н.~Васильев}, \emph{И.~Шарыгин},  \cite{Kvant79-12}.)
\end{pr}

\begin{pr}\label{pr1-9-7}
В остроугольном треугольнике $ABC$ точка $K$ "--- середина $AC$. На
сторонах $AB$ и~$BC$ как на основаниях внутрь треугольника построены
равнобедренные треугольники $ABM$ и~$BCN$ так, что $AM=BM$, $\angle
AMB=\angle AKB$ и~$BN=CN$, $\angle BNC=\angle BKC$. Докажите, что
окружность, описанная около треугольника $MNK$, касается стороны
$AC$. %(\emph{А.~Антропов}, \emph{М.~Урьев}, \cite{Kvant15-2}.)
\end{pr}

\begin{pr}\label{pr1-9-8}
На описанной окружности треугольника $ABC$ взяты точки $A_1$, $B_1$,
$C_1$ так, что $AA_1$, $BB_1$ и~$CC_1$ пересекаются в~одной точке.
При отражении точек $A_1$, $B_1$, $C_1$ относительно сторон $BC$, $CA$,
$AB$ соответственно получаются точки $A_2$, $B_2$, $C_2$. Докажите,
что треугольники $A_1B_1C_1$ и~$A_2B_2C_2$ подобны. %(\emph{А.~Заславский}, \cite{olshar1}.)
\end{pr}

\subsection*{Формулировка и~доказательство обобщённой теоремы
Наполеона}

Через $\angle(\vec{a}, \vec{b})$ будем обозначать угол поворота от вектора $\vec{a}\not=\overrightarrow{0}$ до вектора $\vec{b}\not=\overrightarrow{0}$,
отсчитываемый против часовой стрелки. Этот угол определён с~точностью до прибавления $2\pi k$, $k\in \Z$. Например, равенство $\angle(\vec{a}, \vec{b})\equiv 0 \pmod{2\pi}$ означает, что $\angle(\vec{a}, \vec{b})=2k\pi$ для некоторого $k\in\Z$.

\begin{pr}\zvezda\label{pr1-9-9}
\textbf{Обобщённая теорема Наполеона.} Пусть на сторонах
треугольника $ABC$ построены такие треугольники (возможно,
вырожденные) $BCA_1,$ $CAB_1,$ $ABC_1,$ что выполнены следующие
условия: \index{Теорема!Наполеона!обобщённая}

%\begin{arab}
%\item
1)
$\angle (\overrightarrow{A_1B}, \overrightarrow{A_1C})+
\angle (\overrightarrow{B_1C}, \overrightarrow{B_1A}) + \angle
(\overrightarrow{C_1A}, \overrightarrow{C_1B})\equiv 0 \pmod{2\pi};$

%\item
2)
 $AB_1\cdot BC_1\cdot CA_1 = BA_1\cdot CB_1\cdot AC_1.$
%\end{arab}

Тогда углы треугольника $A_1B_1C_1$ находятся из равенств
\begin{align*}
&\angle (\overrightarrow{A_1C_1}, \overrightarrow{A_1B_1})\equiv
\angle (\overrightarrow{BC_1}, \overrightarrow{BA}) +
\angle (\overrightarrow{CA}, \overrightarrow{CB_1}) \pmod{2\pi};\\
&\angle (\overrightarrow{B_1A_1}, \overrightarrow{B_1C_1}) \equiv
\angle (\overrightarrow{CA_1}, \overrightarrow{CB}) +
\angle (\overrightarrow{AB}, \overrightarrow{AC_1}) \pmod{2\pi};\\
&\angle (\overrightarrow{C_1B_1}, \overrightarrow{C_1A_1}) \equiv
\angle (\overrightarrow{AB_1}, \overrightarrow{AC}) + \angle
(\overrightarrow{BC}, \overrightarrow{BA_1}) \pmod{2\pi}.
\end{align*}
\end{pr}

\emph{Примечание.} В~теореме предполагается, что точка $A_1$ отлична от
$B$, $C$, $B_1$, $C_1$ и~т.\,д. Однако допускается, что вершины
каких-то из треугольников $BCA_1$, $CAB_1$, $ABC_1$ и~$A_1B_1C_1$
лежат на одной прямой. В~этом случае говорят, что соответствующий
треугольник является вырожденным, а~его углы считают равными (с~точностью до $2\pi$) 0, 0 и~$\pi$.

Таким образом, в~теореме утверждается, что при выполнении условий 1 и~2 углы треугольника $A_1B_1C_1$ зависят лишь от углов треугольников, построенных на сторонах треугольника $ABC$, и~не зависят от углов самого треугольника $ABC$. Условие теоремы может быть описано также таким изящным образом (см.~[Bel]%\cite{sb8}
): пусть даны точки $M, N, P, T$ и~на сторонах треугольника $ABC$ строятся треугольники $ABC_1$, $BCA_1$, $CAB_1$, подобные с~сохранением ориентации соответственно треугольникам $MNT$, $NPT$, $PMT$. Действительно, выполнение условий 1 и~2 в~этом случае проверяется непосредственно. С~другой стороны, если треугольники $ABC_1$, $BCA_1$, $CAB_1$ удовлетворяют условиям теоремы, а~треугольники $MNT$ и~$NPT$ подобны соответственно треугольникам $ABC_1$ и~$BCA_1$, то и~треугольник $PMT$ подобен $CAB_1$.

Конструкция из обобщённой теоремы Наполеона интересна, в~ней можно обнаружить ещё несколько красивых фактов, например, окружности $ABC_1$, $BCA_1$, $CAB_1$ и~$A_1B_1C_1$ имеют общую точку (отсюда можно понять, что на самом деле в~этой конструкции треугольники $ABC_1$, $BCA_1$, $CAB_1$ и~$A_1B_1C_1$ равноправны).

\newpage
\setcounter{page}{284}
\refstepcounter{section}
\refstepcounter{subsection}
\refstepcounter{subsection}
\refstepcounter{subsection}
\refstepcounter{subsection}
\refstepcounter{subsection}

\subsection{Теоремы Птолемея и~Кези (3$^*$). \emph{А.~Д.~Блинков}, \emph{А.~А.~Заславский}}\label{s10.6}
\sectionmark{Окружность}

%%\Opensolutionfile{_hintAholder}
%\Opensolutionfile{_hintBholder}

\subsubsection{Теорема Птолемея}

\begin{pr}\label{pr2-6-1}
(a) Докажите, что для любых четырёх различных точек $A$, $B$, $C$, $D$ выполнено \emph{неравенство Птолемея}
$$
AB\cdot CD+AD\cdot BC\geq AC\cdot BD.
$$
\index{Неравенство!Птолемея}

(b) \textbf{Теорема Птолемея.} Это неравенство обращается в~равенство тогда и~только тогда, когда $ABCD$ "--- вписанный четырёхугольник. \index{Теорема!Птолемея}
\end{pr}

\begin{pr}\label{pr2-6-2}
В остроугольном треугольнике $ABC$ обозначим $|BC|=a$, $|AC|=b$.
Найдите $|AB|$, если радиус окружности, описанной около $\triangle
ABC$, равен $R$.
\end{pr}

\begin{pr}\label{pr2-6-3}
Биссектриса угла $A$ треугольника $ABC$ пересекает описанную около
него окружность в~точке $W$.

(a) Выразите отношение $AW/IW$, где $I$ "--- центр окружности,
вписанной в~треугольник $ABC$, через длины сторон треугольника.

(b) Докажите, что $AW>\frac{AB+AC}2$.
\end{pr}

\begin{pr}\label{pr2-6-4}
На гипотенузе $AB$ прямоугольного треугольника $ABC$ во внешнюю
сторону построен квадрат, $O$ "--- его центр. Найдите $|OC|$, если
$a$ и~$b$ "--- катеты треугольника.
\end{pr}

\begin{pr}\label{pr2-6-5}
Дан правильный треугольник $ABC$ и~точка $P$.

(a) Докажите, что если точка $P$ лежит на описанной около треугольника
окружности, то расстояние от неё до одной из вершин треугольника
равно сумме расстояний до двух других вершин.

(b) \textbf{Теорема Помпейю}. Для любой точки $P$, не лежащей на описанной окружности, из отрезков $PA$, $PB$, $PC$ можно составить треугольник. \index{Теорема!Помпейю}
\end{pr}

\begin{pr}\label{pr2-6-6}
Сумма расстояний от точки $X$, выбранной вне квадрата, до двух его ближайших соседних вершин равна $m$. Найдите наибольшее значение суммы расстояний от $X$ до двух других вершин квадрата.
\end{pr}

\begin{pr}\label{pr2-6-7}
Точки $M$ и~$N$ "--- середины диагоналей $AC$ и~$BD$ вписанного четырёхугольника $ABCD$. Известно, что $\angle ABD=\angle MBC$. Докажите, что $\angle BCA=\angle NCD$. (\emph{Кубок Колмогорова\textup{, 1999}~г.})
\end{pr}

\begin{pr}\label{pr2-6-8}
(a) Точки $A$, $B$, $C$ и~$D$ "--- четыре последовательные вершины
правильного семиугольника. Докажите, что
$\frac1{AB}=\frac1{AC}+\frac1{AD}$.

(b) Докажите, что
$\frac1{\sin(\pi/7)}=\frac1{\sin(2\pi/7)}+\frac1{\sin(3\pi/7)}$.
\end{pr}

\begin{pr}\label{pr2-6-9}
В выпуклом шестиугольнике $ABCDEF$ известно, что $AB=BC=a$, $CD=DE=b$,
$EF=FA=c$. Докажите, что
$\frac{a}{BE}+\frac{b}{AD}+\frac{c}{CF}\geq\frac32$.
\end{pr}

\begin{pr}\label{pr2-6-10}
Стороны вписанного четырёхугольника равны $a$, $b$, $c$, $d$.
Найдите его диагонали.
\end{pr}

\begin{pr}\label{pr2-6-11}
Выведите из теоремы Птолемея формулу Карно (см.~п.\;\ref{s9.4} <<Формула
Карно>>). \index{Формула!Карно}
\end{pr}

\subsubsection{Теорема Кези}

\begin{pr}\label{pr2-7-1}
\textbf{Обобщённая теорема Птолемея, или теорема Кези.} (a)~Даны четыре непересекающихся круга, ограниченных окружностями
$\alpha$, $\beta$, $\gamma$, $\delta$. Докажите, что окружность,
касающаяся их внешним образом, или прямая, касающаяся их всех так,
что круги лежат относительно неё в~одной полуплоскости, существует
тогда и~только тогда, когда
$$
l_{\alpha\beta}l_{\gamma\delta}+l_{\alpha\delta}l_{\beta\gamma}=l_{\alpha\gamma}l_{\beta\delta},
$$
где $l_{\alpha\beta}$ "--- длина общей внешней касательной к~окружностям $\alpha$, $\beta$ и~т.\,д.

(b) Сформулируйте теорему Кези для случая, когда искомая окружность
касается некоторых из данных окружностей внутренним образом.
\index{Теорема!Кези}
\end{pr}

\begin{pr}\label{pr2-7-2}
Сформулируйте утверждение, аналогичное теореме Кези, для случая,
когда

(a) одна;

(b) две из данных окружностей вырождаются в~прямые;

(c) какие-то из данных окружностей вырождаются в~точки.
\end{pr}

\begin{pr}\label{pr2-7-3}
Пусть на сторонах $AC$ и~$BC$ треугольника $ABC$ взяты такие точки $X$,
$Y$, что $XY\parallel AB$. Докажите, что существует
окружность, проходящая через $X$, $Y$ и~касающаяся одинаковым
образом вневписанных окружностей треугольника, вписанных в~углы $A$
и~$B$.
\end{pr}

\begin{pr}\label{pr2-7-4}
Докажите \textbf{теорему Фейербаха}: окружность, проходящая через
середины сторон треугольника, касается его вписанной и~вневписанных
окружностей. \index{Теорема!Фейербаха}
\end{pr}

\begin{pr}\label{pr2-7-5}
Докажите, что три окружности, каждая из которых касается внутренним
образом одной из вневписанных окружностей треугольника и~внешним
образом двух других, пересекаются в~одной точке.
\end{pr}

\begin{pr}\label{pr2-7-6}
Даны две окружности, лежащие одна вне другой. Произвольная
окружность, касающаяся их одинаковым образом, пересекает одну из их
общих внутренних касательных в~точках $A$ и~$A'$, а~другую "--- в~точках $B$ и~$B'$. Докажите, что среди прямых $AB$, $AB'$, $A'B$,
$A'B'$ найдутся две, параллельные общим внешним касательным к~данным
окружностям.
\end{pr}

\begin{pr}\label{pr2-7-7} Даны две концентрические окружности $a_1$ и~$a_2$. Каждая из окружностей $b_1$ и~$b_2$ касается внешним образом окружности $a_1$ и~внутренним "--- $a_2$, а~каждая из окружностей $c_1$ и~$c_2$ касается внутренним образом обеих окружностей $a_1$ и~$a_2$. Оказалось, что окружности $b_1$, $b_2$ пересекают $c_1$, $c_2$ в~восьми точках. Докажите, что эти точки лежат на двух окружностях или прямых, отличных от $b_1$, $b_2$, $c_1$,~$c_2$. (\emph{В.~Протасов\textup{, III} Олимпиада им.~И.~Ф.~Шарыгина}.)
\end{pr}

%%\Closesolutionfile{_hintAholder}
%\Closesolutionfile{_hintBholder}

%\subsubsection*{Подсказки}
%%\Readsolutionfile{_hintAholder}

\sseccol{Указания, ответы и~решения}
%\Readsolutionfile{_hintBholder}

%%\ssection{Указания, ответы, решения}
%%
%%\ssection{1. Теорема Птолемея}
%%\ssection{2. Теорема Кези}

\paragraph*{\ref{pr2-6-1}.}
Сделайте инверсию с~центром $A$ и~воспользуйтесь утверждениями задач 11.10.2, 11.10.5%\ref{pr2-16-2}, \ref{pr2-16-5}
. Заметим, что неравенство Птолемея верно даже для точек, не лежащих в~одной плоскости.\index{Инверсия}

\paragraph*{\ref{pr2-6-2}.}
\emph{Ответ}: $\frac{a\sqrt{4R^2-b^2}+b\sqrt{4R^2-a^2}}{2R}$.

Проведите диаметр $CD$ и~примените к~полученному четырёхугольнику
теорему Птолемея.

\newpage
\setcounter{page}{292}
\section{Геометрические преобразования}\label{preobr}

В данном параграфе задачи расположены так, чтобы сначала новые
понятия (геометрических преобразований) использовались для решения
интересных задач, формулируемых без этих понятий, и~только потом эти
новые (но уже мотивированные) понятия изучались сами по себе.

%%%!!! перенесено в \S1

Подробнее о~геометрических преобразованиях см., например,
\cite{zasl} (теорема Шаля "--- \S\,1.2, подобие и~гомотетия "---
\S\,1.3, аффинные преобразования "--- гл.\,2, проективные
преобразования "--- гл.\,3, инверсия "--- гл.\,4, комплексная
интерпретация движений и~подобий "--- \S\,6.1, комплексная
интерпретация инверсии "--- \S\,6.2), \cite{prasolov} и~\cite{Yaglom}.

\subsection{Применения движений. (1) \emph{А.~Д.~Блинков}}\label{s11.1}

\emph{Поворотом} вокруг точки $O$ на угол $\varphi$ называется преобразование плоскости, оставляющее точку $O$ на месте и~переводящее любую отличную от $O$ точку $X$ в~такую точку $X'$, что $|OX|=OX'$ и~ориентированный угол между векторами $\overrightarrow{OX}$ и~$\overrightarrow{OX'}$ равен $\varphi$. Поворот на $180^{\circ}$ называется \emph{центральной симметрией}.
\index{Поворот|textbf}
\index{Симметрия!центральная|textbf}

\emph{Параллельным переносом} на вектор $\overrightarrow{m}$ называется преобразование плоскости, переводящее любую точку $X$ в~такую точку $X'$, что $\overrightarrow{XX'}=\overrightarrow{m}$.
\index{Параллельный перенос|textbf}

\emph{Осевой симметрией} относительно прямой $l$ называется преобразование плоскости, переводящее любую точку $X$ в~такую точку $X'$, что $XX'\perp l$ и~точки $X$, $X'$ лежат по разные стороны от прямой $l$ и~равноудалены от неё.
\index{Симметрия!осевая|textbf}

\begin{pr}\circpr\label{pr2-8-1}
Параллелограмм имеет ровно четыре оси симметрии. Какое из следующих
утверждений верно?

%\begin{arab}
1) %\item
это прямоугольник, отличный от квадрата;

2) %\item
это ромб, отличный от квадрата;

3) %\item
это квадрат;

4) %\item
такого параллелограмма не существует.
%\end{arab}
\end{pr}

\comment
\begin{pr}\circpr\label{pr2-8-2}
Треугольник имеет центр симметрии.
Какое из следующих утверждений верно?

%\begin{arab}
%\item
1) он разносторонний;

%\item
2) он равносторонний;

%\item
3) он равнобедренный;

%\item
4) такого треугольника не существует.
%\end{arab}
\end{pr}

\endcomment

\newpage
\setcounter{page}{312}
\refstepcounter{subsection}
\refstepcounter{subsection}
\refstepcounter{subsection}
\subsection{Поворотная гомотетия (2). \emph{П.~А.~Кожевников}}\label{s11.5}

\subsubsection{Вводные задачи: немного о~велосипедистах}

\begin{pr}\label{pr2-11-2}
По двум окружностям, пересекающимся в~точках $P$ и~$Q$, одновременно
начали движение из точки $P$ по часовой стрелке с~равными угловыми
скоростями два велосипедиста $A$ и~$B$.

(a) Докажите, что прямая $AB$ всё время проходит через $Q$.

(b) Докажите, что треугольники $PAB$ всё время подобны друг другу и~треугольнику $PO_1O_2$, где $O_1$ и~$O_2$ "--- центры окружностей.

(c) Найдите ГМТ (траекторию движения) середин отрезков $AB$; центров
вписанных окружностей треугольников $PAB$; любых соответственных
точек подобных треугольников $PAB$.

(d) \textbf{Задача о~велосипедистах.} Докажите, что $A$ и~$B$ всё
время равноудалены от фиксированной точки. %(См. \cite{Kvant}.)
\end{pr}

\begin{pr}\label{pr2-11-3}
По трём окружностям, имеющим общую точку $O$ и~попарно различные
точки пересечения $P$, $Q$ и~$R$, одновременно начали движение из
точки $O$ по часовой стрелке с~равными угловыми скоростями три
велосипедиста $A$, $B$ и~$C$.

(a) Докажите, что все треугольники $ABC$ подобны между собой и~треугольнику $O_1O_2O_3$, где $O_1$, $O_2$ и~$O_3$ "--- центры окружностей.

(b) Какова траектория движения центра масс треугольника $ABC$?
\end{pr}

\begin{pr}\label{pr2-11-4}
Два велосипедиста $P$ и~$Q$ едут равномерно по двум прямым, пересекающимся в~точке~$O$.

(a) Найдите траекторию середины отрезка $PQ$.

(b) Докажите, что если скорости велосипедистов равны, то середина дуги (одной
из дуг) $PQ$ окружности $OPQ$ неподвижна.

(c) Докажите, что если велосипедисты проходят $O$ не одновременно, то
окружности $OPQ$ имеют вторую общую точку, отличную от~$O$.
\end{pr}

\begin{pr}\label{pr2-11-5}
Дан фиксированный треугольник $ABC$. По прямым $BC$, $CA$, $AB$ едут
соответственно велосипедисты $P$, $Q$, $R$ так, что углы между $RP$ и~$PQ$, $PQ$ и~$QR$, $QR$ и~$RP$ фиксированные.

(a) Докажите, что точка пересечения окружностей $RAQ$, $RBP$, $PCQ$ неподвижна.

(b) Найдите ГМТ центров вписанных окружностей треугольников $PQR$.
\end{pr}

\subsubsection{Основные задачи}

\begin{pr}\label{pr2-11-6}
Три велосипедиста $P$, $Q$ и~$R$ едут равномерно по трём прямым.
Известно, что в~некоторые два момента времени треугольник $PQR$ был
подобен с~сохранением ориентации фиксированному треугольнику $XYZ$.
Докажите, что это условие будет выполняться в~любой момент времени.
\end{pr}

\begin{pr}\label{pr2-11-7}
В треугольник $ABC$ вписан подобный ему треугольник $PQR$ ($P\in
BC$, $Q\in CA$, $R\in AB$, $\angle P=\angle A$, $\angle Q=\angle B$,
$\angle R=\angle C$).

(a) Докажите, что центр описанной окружности треугольника $ABC$
совпадает с~ортоцентром треугольника $PQR$.

(b) Найдите максимальное значение выражения
$\frac{S_{ABC}}{S_{PQR}}$.

(c) Докажите, что центр описанной окружности треугольника $PQR$
равноудалён от центра описанной окружности и~ортоцентра треугольника
$ABC$.
\end{pr}

\begin{pr}\label{pr2-11-8}
Через вершины треугольника $ABC$ проводятся три произвольные
параллельные прямые $d_a$, $d_b$, $d_c$. Прямые $d'_a$, $d'_b$,
$d'_c$, симметричные $d_a$, $d_b$, $d_c$ относительно $BC$, $CA$,
$AB$ соответственно, образуют треугольник $XYZ$. Найдите
геометрическое место центров вписанных окружностей таких
треугольников.
%(А.~Заславский, олимпиада им. И.~Ф.~Шарыгина 2009 г.)
\end{pr}

\begin{pr}\label{pr2-11-9}
Дан выпуклый четырёхугольник $ABCD$, стороны $BC$ и~$AD$ которого равны, но не параллельны. Пусть $E$ и~$F$ "--- внутренние точки отрезков $BC$ и~$AD$ соответственно, удовлетворяющие условию $BE=DF$. Прямые $AC$ и~$BD$ пересекаются в~точке~$P$, прямые $BD$ и~$EF$ пересекаются в~точке~$Q$, прямые $EF$ и~$AC$ пересекаются в~точке~$R$. Докажите, что для всевозможных способов выбора точек $E$, $F$ окружности $PQR$ имеют общую точку, отличную от~$P$. %(См. \cite{IMO}, 2005~г.)
\end{pr}

\begin{pr}\label{pr2-11-10}
Пусть $O$ и~$I$ "--- центры описанной и~вписанной окружностей
треугольника $ABC$ соответственно. Точки $D$, $E$ и~$F$ выбраны на
сторонах $BC$, $CA$ и~$AB$ соответственно так, что $BD+BF=CA$ и~$CD+CE=AB$. Описанные окружности треугольников $BDF$ и~$CDE$
пересекаются в~точках $D$ и~$P$. Докажите, что $OP=OI$. %(См. \cite{IMO}, 2012~г.)
\end{pr}

\subsubsection{Дополнительные задачи}

\begin{pr}\label{pr2-11-11}
Впишите в~данный остроугольный треугольник равносторонний треугольник c~минимальной стороной.
\end{pr}

\begin{pr}\label{pr2-11-12}
На пол положили правильный треугольник $ABC$, выпиленный из фанеры.
В пол вбили три гвоздя (по одному вплотную к~каждой стороне треугольника) так, что треугольник невозможно повернуть, не отрывая от пола. Первый гвоздь делит сторону $AB$ в~отношении $1:3$, считая от вершины~$A$, а~второй делит сторону $BC$ в~отношении $2:1$, считая от вершины~$B$. В~каком отношении делит сторону $AC$ третий гвоздь? %(\emph{Московская математическая олимпиада \textup{1998}~г.})
\end{pr}

\begin{pr}\label{pr2-11-13}
Выпуклый многоугольник $M$ можно поместить в~треугольник~$T$. Докажите, что это можно сделать так, чтобы одна из сторон многоугольника $M$ лежала на стороне треугольника~$T$.
\end{pr}

\emph{Поворотной гомотетией} называют преобразование %$H_O^{k,\varphi}:=
$H_O^k\circ R_O^\varphi$. \index{Гомотетия!поворотная}

\begin{pr}\label{pr2-11-1}
(a) Окружности $\alpha$ и~$\beta$ пересекаются в~точках $A$ и~$B$.
Пусть $H$ "--- поворотная гомотетия с~центром в~точке~$A$, переводящая
$\alpha$ в~$\beta$. Докажите, что для любой точки $X\in\alpha$ точка
$H(X)$ получена пересечением прямой $BX$ с~окружностью $\beta$.
(См. \cite[19.27]{prasolov}.)

(b) Окружности $S_1,\ldots,S_n$ проходят через точку $O$. Кузнечик из точки $X_i\in S_i$ прыгает в~точку $X_{i+1}\in S_{i+1}$ так, что прямая $X_iX_{i+1}$ проходит через вторую точку пересечения окружностей $S_i$ и~$S_{i+1}$. Докажите, что после $n$ прыжков (с~$S_1$ на $S_2,\ldots,$ с~$S_n$ на $S_1$) кузнечик вернётся в~исходную точку. (См. \cite[19.28]{prasolov}.)

\comment
(c) Пусть концы отрезков $AB$ и~$CD$ попарно различны, а~$P$ "---
точка пересечения прямых $AB$ и~$CD$. Центром поворотной гомотетии,
переводящей $AB$ в~$CD$, является (отличная от $P$) точка
пересечения описанных окружностей треугольников $ACP$ и~$BDP$.
(См. \cite[19.41\,(б)]{prasolov}.)

\endcomment
\end{pr}

%\subsubsection*{Подсказки}

\newpage
\setcounter{page}{343}
\refstepcounter{section}
\refstepcounter{subsection}
\refstepcounter{subsection}
\subsection{Полярное соответствие (2). \emph{А.~А.~Гаврилюк}, \emph{П.~А.~Кожевников}}\label{s12.3}
\sectionmark{Аффинная и проективная геометрия}

%%\Opensolutionfile{_hintAholder}
%\Opensolutionfile{_hintBholder}

Традиционно при изучении полярного соответствия существенно используются свойства проективных преобразований. Мы же делаем попытку познакомиться с~полярным соответствием и~применением его свойств без привлечения проективной геометрии.

Введём нужные нам определения и~обозначения. Пусть на плоскости фиксированы точка $O$ и~окружность $\omega$ радиуса~$R$ с~центром в~$O$.

%\begin{deff}
Для каждой точки $X\neq O$ на луче $OX$ строим такую точку~$X'$, что
$OX \cdot OX' = R^2$. (Говорят, что $X'$ и~$X$ \emph{инверсны} относительно окружности~$\omega$.) Через точку $X'$ проведём прямую $x$, перпендикулярную~$OX'$. Прямая $x$ называется \emph{полярой} точки~$X$, а~точка $X$ называется \emph{полюсом} прямой~$x$.
Соответствие $X \leftrightarrow x$ является взаимно однозначным соответствием между точками, отличными от~$O$, и~прямыми, не проходящими через~$O$. Это соответствие и~называется \emph{полярным соответствием}.

Ниже мы обозначаем точки, отличные от $O$ (полюсы), большими латинскими
буквами, а~их поляры "--- соответствующими маленькими буквами: $A
\leftrightarrow a$, $B\leftrightarrow b$, $C \leftrightarrow c$,
$\ldots $
%\end{deff}
\index{Полярное!соответствие} \index{Полюс} \index{Поляра}

\subsection*{Основные свойства и~вводные задачи}

Установите два основных \emph{свойства полярного соответствия.}

\renewcommand{\thepr}{П\arabic{pr}}

\begin{pr}\label{pr2-p1}
\textbf{Двойственность.} Включение $A\in b$ выполняется тогда и~только тогда
$B\in a$, т.\,е. поляра любой точки является геометрическим местом полюсов проходящих через неё прямых. \index{Двойственность}
\end{pr}

\begin{pr}\zvezda\label{pr2-p2}
Пусть две прямые $m$ и~$l$, проходящие через произвольную точку $A\notin \omega$, пересекают $\omega $ в~точках $M_1$, $M_2$ и~$L_1$, $L_2$. Тогда $M_1L_1\cap M_2L_2 \in a$ или $M_1L_1\parallel M_2L_2 \parallel a$.
\end{pr}

Докажите следующие факты.

\setcounter{pr}{0}
\renewcommand{\thepr}{B\arabic{pr}}

\begin{pr}\label{pr2-b1}
Если $A\in \omega$, то $a$ "--- это касательная к~$\omega $, проведённая через $A$.
\end{pr}

\begin{pr}\label{pr2-b2}
Если точка $A$ расположена вне окружности $\omega$, то
$a$ проходит через точки касания с~$\omega$ касательных, проведённых
через $A$.
\end{pr}

\begin{pr}\label{pr2-b3}
Если $O, A, B$ не лежат на одной прямой, то $a\cap b \leftrightarrow AB$.
\end{pr}

\begin{pr}\label{pr2-b4}
Точки $A, B, C$ лежат на одной прямой тогда и~только тогда, когда $a, b, c$ проходят через одну точку или параллельны.
\end{pr}

\subsection*{Основные задачи}

\setcounter{pr}{0}
\renewcommand{\thepr}{\arabic{section}.\arabic{subsection}.\arabic{pr}}

%Следующие задачи предназначены для закрепления материала

\begin{pr}\circpr\label{pr2-19-15}
Даны окружность и~её хорда $AB$. Где лежит точка пересечения поляр
точек $A$ и~$B$?

1) внутри окружности;\quad 2) вне окружности;\quad 3) на окружности.
\end{pr}

\begin{pr}\circpr\label{pr2-19-16}
Пусть $C$ "--- середина хорды $AB$. Тогда поляра точки $C$

%\begin{arab}
%\item
1) параллельна $AB$;

%\item
2) перпендикулярна $AB$;

%\item
3) касается окружности.
%\end{arab}
\end{pr}

\begin{pr}\circpr\label{pr2-19-17}
При полярном соответствии относительно вписанной окружности треугольник переходит

%\begin{arab}
%\item
1) в~серединный треугольник;

%\item
2) в~ортотреугольник;

%\item
3) в~треугольник, образованный точками касания сторон с~вписанной окружностью.
%\end{arab}
\end{pr}

\begin{pr}\label{pr2-19-1}
Даны окружность $\omega$ и~прямая $l$, не имеющие общих точек. Из
точки $X$, которая движется по прямой $l$, проводятся касательные
$XA$, $XB$ к~$\omega$. Докажите, что все хорды $AB$ имеют общую
точку.
\end{pr}

\begin{pr}\label{pr2-19-2}
\textbf{Симметричная бабочка.} (a) Дана точка $A$ на диаметре $BC$
полуокружности $\omega $. Точки $X, Y$ на $\omega $ таковы, что
$\angle XAB = \angle YAC$. Докажите, что прямые $XY$ проходят через
одну точку или параллельны.

(b) Точки $A$ и~$A'$ инверсны относительно окружности~$\omega $, причём точка $A'$ расположена внутри $\omega $. Через $A'$ проводятся хорды~$XY$. Докажите, что центры вписанной и~одной из вневписанных окружностей
треугольника $AXY$ фиксированны. %(\emph{С.~Маркелов, см.}~\cite{shar}.)
\end{pr}

\begin{pr}\label{pr2-19-3}
\textbf{Основное свойство симедианы.} Касательные к~описанной
окружности треугольника $ABC$, проведённые через точки $B$ и~$C$,
пересекаются в~точке~$P$. Докажите, что $AP$ "--- симедиана (т.\,е.
прямая, симметричная медиане $AM$ относительно биссектрисы угла
$A$). \index{Симедиана}
\end{pr}

\begin{pr}\label{pr2-19-4}
\textbf{Гармонический четырёхугольник.} Пусть четырёхугольник $ABCD$ вписан в~окружность $\omega$. Известно,
что касательные к~$\omega$, проведённые в~точках $A$ и~$C$,
пересекаются на прямой $BD$ или параллельны $BD$. Докажите, что
касательные к~$\omega$, проведённые в~точках $B$ и~$D$, пересекаются
на прямой $AC$ или параллельны $AC$.
\end{pr}

В следующих трёх задачах дан четырёхугольник $ABCD$, у~которого диагонали пересекаются в~точке~$P$, продолжения сторон $AB$ и~$CD$  "--- в~точке~$R$, продолжения сторон $BC$ и~$DA$ "--- в~точке~$Q$.

\begin{pr}\label{pr2-19-5}
\textbf{Вписанный четырёхугольник.} Пусть четырёхугольник $ABCD$ вписан в~окружность с~центром $O$. Докажите, что четвёрка точек $O$, $P$, $Q$, $R$ ортоцентрическая (т.\,е. каждая точка является ортоцентром треугольника с~вершинами в~оставшихся трёх точках).
\end{pr}

\begin{pr}\label{pr2-19-6}
\textbf{Описанный четырёхугольник.} Пусть четырёхугольник $ABCD$ описан около окружности; $K$, $L$, $M$, $N$ "--- точки касания с~окружностью сторон~$AB$, $BC$, $CD$, $DA$ соответственно; прямые $KL$ и~$MN$ пересекаются в~точке~$S$, а~прямые $LM$ и~$NK$ "--- в~точке~$T$.

(a) Докажите, что точки $Q$, $R$, $S$, $T$ лежат на одной прямой.

(b) Докажите, что $KM$ и~$LN$ пересекаются в~точке~$P$.
\end{pr}

\begin{pr}\label{pr2-19-7}
\textbf{Вписанно-описанный четырёхугольник.} Четырёхугольник $ABCD$ описан около окружности $\omega $ с~центром $I$ и~вписан в~окружность~$\Omega$ с~центром~$O$.

(a) Докажите, что $O$, $I$, $P$ лежат на одной прямой.

(b) Зафиксируем $\omega $ и~$\Omega$ и~рассмотрим всевозможные четырёхугольники $ABCD$, описанные около окружности $\omega $ и~вписанные в~окружность~$\Omega $. Докажите, что для всех таких четырёхугольников точки $P$ совпадают, а~также что прямые $QR$ совпадают. \index{Четырёхугольник!вписанно-опи\-сан\-ный}

\emph{Комментарий}. Согласно теореме Понселе если существует хотя бы один четырёхугольник, описанный около окружности $\omega$ и~вписанный в~окружность~$\Omega$, то существует бесконечно много таких четырёхугольников. \index{Теорема!Понселе}
\end{pr}

\sectionmark{Комплексные числа и~геометрия}
\comment
\section{Комплексные числа и~геометрия (3). \emph{А.~А.~Заславский}}\label{compgeom}

%\Author{}

Этот параграф посвящён применению комплексных чисел для решения
геометрических задач. Для его изучения достаточно знать определение
комплексных чисел и~действий над ними (см.~п.\;\ref{s:compl}). Желательно также знание
определений и~основных свойств движений, подобий, аффинных
преобразований для первого пункта и~инверсии для второго.
\endcomment

\newpage
\setcounter{page}{352}
\refstepcounter{section}
\subsection{Комплексные числа и~элементарная геометрия.}\label{compgeom-1}

%\Name{А.~А.~Заславский}]{Комплексные числа и~элементарная геометрия}

%%\Opensolutionfile{_hintAholder}
%\Opensolutionfile{_hintBholder}

Пусть на плоскости задана система координат. Тогда комплексному числу $z=x+yi$ соответствует точка плоскости $Z$ с~координатами $(x,y)$. При этом модуль числа $z$ равен расстоянию от $Z$ до начала координат~$O$, а~аргумент равен ориентированному углу между положительным направлением оси $Ox$ и~вектором $\overrightarrow{OZ}$, т.\,е. углу, на который надо повернуть против часовой стрелки ось~$Ox$, чтобы совместить её положительное направление с~направлением вектора $\overrightarrow{OZ}$. Оси $Ox$ и~$Oy$ называют \emph{действительной} и~\emph{мнимой} осями.

\begin{pr}\label{pr2-20-0}
(Загадка.) Выясните геометрический смысл сложения комплексных чисел.
\end{pr}

\begin{pr}\label{pr2-20-1}
(a) Каким геометрическим преобразованием комплексной плоскости
получается число~$iz$ из~числа~$z$?

(b) (Загадка.) Обозначим $e^{i\varphi}:=\cos\varphi+i\sin\varphi$. Каков
геометрический смысл умножения на~$e^{i\varphi}$? А~на~$re^{i\varphi}$, где $r$ "--- вещественное число (см.~определение тригонометрической формы комплексного числа в~п.\;4.5%\ref{s:compl}
)?

(c) Выразите число~$w$, полученное из~числа~$z$ поворотом на~угол~$\varphi$ против часовой стрелки относительно центра~$z_0$, через $z$, $z_0$ и~$\varphi$.

(d) Докажите, что композиция поворотов плоскости (с~различными центрами) "--- поворот или параллельный перенос.

(e) Докажите, что точки $z_1$, $z_2$, $z_3$ лежат на одной прямой тогда и~только тогда, когда отношение $(z_3-z_1)/(z_2-z_1)$ вещественно.
\end{pr}

\emph{Комментарий.} Задача \ref{pr2-20-1}\,(b) легко решается с~помощью
тригонометрических формул сложения. Однако можно поступить наоборот: решить эту задачу геометрически, доказать, что при умножении комплексных чисел их модули перемножаются, а~аргументы складываются, а~затем, используя этот результат, получить доказательство формул сложения, не требующее перебора различных случаев.

\begin{pr}\circpr\label{pr2-20-8}
Какое преобразование плоскости задаётся формулой $z\mapsto 2z+2$?

\comment
%\begin{arab}
%\item
1) параллельный перенос на вектор длины 2;

%\item
2) гомотетия с~коэффициентом 2;

%\item
3) сжатие к~прямой с~коэффициентом 2.
%\end{arab}
\endcomment
\end{pr}

\newpage
\setcounter{page}{355}
\subsection{Комплексные числа и~круговые преобразования.}\label{s13.2}

%%\Opensolutionfile{_hintAholder}
%\Opensolutionfile{_hintBholder}

Преобразование круговой плоскости, сохраняющее обобщённые окружности,
называется \emph{круговым}. Произвольное отличное от подобия
круговое преобразование может быть представлено как композиция
инверсии и~движения. \index{Круговое!преобразование}

\begin{pr}\circpr\label{pr2-20-9}
Четвёрка комплексных чисел $z_1, z_2, z_3, z_4$ удовлетворяет
равенству $\frac{(z_1-z_3)(z_2-z_4)}{(z_1-z_4)(z_2-z_3)}=2$. Что
можно сказать о~четвёрке точек плоскости, соответствующих числам
$z_1, z_2, z_3, z_4$?

%\begin{arab}
%\item
1) Они являются вершинами параллелограмма.

%\item
2) Они лежат на одной прямой или на одной окружности.

%\item
3) Площадь треугольника $0z_1z_2$ равна площади треугольника
$0z_3z_4$ (точка $0$ "--- начало координат).
%\end{arab}
\end{pr}

\begin{pr}\label{pr2-21-1}
Докажите, что преобразование комплексной плоскости является круговым
тогда и~только тогда, когда оно задаётся дробно-линейной функцией
вида $f(z)=(az+b)/(cz+d)$ или $f(z)=(a\bar z+b)/(c\bar z+d)$, где
$ad-bc\not=0$.
\end{pr}

\begin{pr}\label{pr2-21-2}
Докажите, что для любых шести различных точек $A$, $B$, $C$, $A'$,
$B'$, $C'$ существует ровно два круговых преобразования, переводящих
$A$ в~$A'$, $B$ в~$B'$, $C$ в~$C'$.
\end{pr}

\emph{Двойным отношением} четырёх комплексных чисел $a$, $b$, $c$,
$d$, где $a\not= d$, $b\not= c$, называется комплексное число
$(a,b,c,d)= \frac{(a-c)(b-d)}{(a-d)(b-c)}$.
\index{Отношение!двойное}

\begin{pr}\label{pr2-21-3}
Докажите, что для данных восьми различных точек $A$, $B$, $C$, $D$;
$A'$, $B'$, $C'$, $D'$ круговое преобразование, переводящее $A$ в~$A'$, $B$ в~$B'$, $C$ в~$C'$, $D$ в~$D'$, существует тогда и~только тогда, когда для соответствующих комплексных чисел выполняется равенство
$(a,b,c,d)=(a',b',c',d')$ или $\overline{(a,b,c,d)}=(a',b',c',d')$.
\end{pr}

\begin{pr}\label{pr2-21-4}
Даны два треугольника $ABC$ и~$A'B'C'$. Докажите, что существует инверсия, переводящая треугольник $ABC$ в~треугольник, равный $A'B'C'$.
\end{pr}

\begin{pr}\label{pr2-21-5}
Дан четырёхугольник $ABCD$. Докажите, что существует инверсия, переводящая его вершины в~вершины параллелограмма, причём все параллелограммы, полученные в~результате таких инверсий, подобны.
\end{pr}

См.~также задачу 11.10.4%\ref{pr2-16-4}
\,(c) п.~<<Инверсия>>.

\subsection*{Дополнительные задачи}

\begin{pr}\label{pr2-21-6}
(a) Пусть $a$, $b$, $c$ "--- комплексные числа, соответствующие не
лежащим на одной прямой точкам $A$, $B$, $C$; $f(z)=(z-a)(z-b)(z-c)$. Докажите, что две точки, соответствующие корням производной $f'(z)$, изогонально сопряжены относительно треугольника~$ABC$.

(b)$^*$ \emph{Эллипсом Штейнера} треугольника $ABC$ называется
эллипс наибольшей площади, лежащий внутри треугольника. Докажите,
что фокусы эллипса Штейнера соответствуют корням производной~$f'(z)$.
\index{Эллипс!Штейнера}
\end{pr}

\begin{pr}\label{pr2-21-7}
Пусть $a$, $b$, $c$ "--- комплексные числа, соответствующие точкам $A$, $B$, $C$, причём $|a|=|b|=|c|=1$. Докажите, что точки $Z_1$, $Z_2$ изогонально сопряжены относительно треугольника $ABC$ тогда и~только тогда, когда соответствующие комплексные числа удовлетворяют соотношению
$$
z_1+z_2+abc\bar{z_1}\bar{z_2}=a+b+c.
$$
\end{pr}

\begin{pr}\label{pr2-21-8}
Как известно, расстояние между центрами $O$ и~$I$ описанной и~вписанной окружностей треугольника выражается через радиусы $R$, $r$ этих окружностей с~помощью \emph{формулы Эйлера}: $OI^2=R^2-2Rr$ (см.~задачу \ref{pr1-4-10} п.~<<Формула Карно>>). Докажите обобщение формулы Эйлера: если в~треугольник вписан эллипс с~фокусами $F_1$, $F_2$ и~малой осью $l$, то
$$
R^2l^2=(R^2-OF_1^2)(R^2-OF_2^2).
$$
\index{Формула!Эйлера}
\end{pr}

\begin{pr}\label{pr2-21-9}(А.~Акопян).
Даны треугольник $ABC$ и~точка $P$. Найдите ГМТ, изогонально сопряжённых точке~$P$ относительно всех треугольников, имеющих те же описанную и~вписанную окружности, что и~$ABC$.
\end{pr}

\newpage
\setcounter{page}{373}
\refstepcounter{section}
\refstepcounter{subsection}
\refstepcounter{subsection}
\subsection{Построения. Ящик инструментов (2). \emph{А.~А.~Гаврилюк}}\label{s14-3}

\sectionmark{Построения и геометрические места точек}
%%\Opensolutionfile{_hintAholder}
%\Opensolutionfile{_hintBholder}

При изучении материала этого раздела желательно знакомство
с~\S\;10%\ref{okruzhn}
 <<Окружность>> и~рекомендованной в~нем литературой.

\begin{pr}\label{pr3-3-1}
Даны два отрезка с~длинами $x$, $y$. С~помощью циркуля и~линейки постройте отрезок длины $\sqrt{3xy+y\sqrt[4]{xy^3}}$.
\end{pr}

\begin{pr}\label{pr3-3-2}
(a) Даны две параллельные прямые, на одной из которых дан отрезок.
С~помощью одной линейки разделите его пополам.

(b) Даны две параллельные прямые, на одной из которых дан отрезок.
С~помощью одной линейки удвойте его.

(c) Даны две параллельные прямые, на одной из которых дан отрезок.
С~помощью одной линейки разделите его на $n$ равных частей.

Ср. с~задачей 11.9.5%\ref{pr2-15-5}
 п.~<<Центральная проекция и~проективные
преобразования".
\end{pr}

\begin{pr}\label{pr3-3-3}
Даны окружность $\omega$, её диаметр $AB$ и~точка $X$. С~помощью
одной линейки постройте перпендикуляр из точки~$X$ на~$AB$, если точка~$X$
лежит

(a) не на окружности; \quad (b) на окружности.
\end{pr}

\begin{pr}\label{pr3-3-4}
Даны окружность $\omega$ и~точка $X$. С~помощью одной линейки
постройте (все возможные) касательные, проведённые из точки $X$ к~окружности, если точка~$X$ лежит

(a) вне окружности; \quad (b) на окружности.
\end{pr}

\begin{pr}\label{pr3-3-5}
При помощи только циркуля постройте образ данной точки $X$ при инверсии относительно данной окружности~$\omega$.
\end{pr}

\begin{pr}\label{pr3-3-6}
Дана окружность на плоскости. С~помощью двусторонней линейки
постройте её центр. (С помощью двусторонней линейки можно проводить
прямую через две точки, проводить прямую, параллельную проведённой
ранее прямой и~отстоящую от неё на расстояние, равное ширине
линейки, а~также проводить через две точки, расстояние между
которыми не меньше ширины линейки, две параллельные прямые,
расстояние между которыми равно ширине линейки.)
\end{pr}

\begin{pr}\label{pr3-3-7}
Даны прямая $l$ и~отрезок $OA$, ей параллельный. С~помощью
двусторонней линейки постройте точки пересечения прямой $l$
с~окружностью радиуса $OA$ и~с центром в~точке $O$.
\end{pr}

\begin{pr}\label{pr3-3-8}
При помощи только циркуля постройте окружность, проходящую через три данные точки.
\end{pr}

\begin{pr}\label{pr3-3-9}
\textbf{Задача Аполлония.} Постройте окружность, касающуюся трёх данных, при помощи циркуля и~линейки. \index{Задача!Аполлония}
\end{pr}

В последующих задачах этого пункта \emph{построением} будем
называть некоторую последовательность следующих элементарных операций:

"--* с~помощью линейки провести прямую через две данные или
ранее построенные точки;

"--* с~помощью циркуля построить окружность с~центром $A$ и~радиусом $BC$, где $A$, $B$, $C$ "--- данные или ранее построенные
точки;

"--* найти точки пересечения двух данных или ранее
построенных линий (прямых или окружностей).

В последующих теоремах никакие другие операции не разрешаются (в~отличие от предыдущих задач, где разрешена, например, операция <<взять произвольную точку уже построенного множества>>). В~частности, если изначально не даны хотя бы две точки, ничего построить нельзя.

\begin{pr}\zvezda\label{pr3-3-10}
\textbf{Теорема.} С~помощью циркуля и~линейки можно осуществлять те и~только те построения, которые <<сводятся>> к~арифметическим операциям и~операции извлечения квадратного корня, то есть если на плоскости фиксирована система координат, то координаты всех построимых точек выражаются через координаты исходных точек с~помощью указанных операций.

\emph{Комментарий.} Из этой теоремы, в~частности, следует, что если дан отрезок длины~1, то для любых отрезков с~длинами $a$, $b$ можно построить отрезки с~длинами $a+b$, $a-b$, $ab$, $a/b$, $\sqrt{a}$ и~длина любого отрезка, который можно построить, выражается через $a$ и~$b$ с~помощью указанных операций (ср. с~основной теоремой из п.~5.2.3%\ref{motcon}
).
\end{pr}

\begin{pr}\zvezda\label{pr3-3-11}
\textbf{Теорема} (Мор"--~Маскерони). Любое построение, осуществимое циркулем
и~линейкой, можно осуществить одним циркулем (прямая считается построенной, если построены две различные лежащие на ней точки, см.~\cite{FuksD}).
\index{Теорема!Мора"--~Маскерони}
\end{pr}

\begin{pr}\zvezda\label{pr3-3-12}
\textbf{Теорема} (Штейнер). Любое построение, осуществимое циркулем и~линейкой, можно осуществить одной линейкой, если начерчена одна окружность и~отмечен её центр (окружность считается построенной, если построены её центр и~лежащая на ней точка, см.~\cite{Smog}).\index{Теорема!Штейнера}
\end{pr}

Следующая задача предназначена для закрепления материала.

\begin{pr}\circpr\label{pr3-3-14}
Пользуясь теоремами Мора"--~Маскерони и~Штейнера, определите, какие инструменты необходимы для построения центра данной окружности.

1) циркуль и~линейка;\quad 2) только линейка;\quad 3) только циркуль.
\end{pr}

%%\Closesolutionfile{_hintAholder}
%\Closesolutionfile{_hintBholder}

%\subsubsection*{Подсказки}
%%\Readsolutionfile{_hintAholder}

\sseccol{Указания, ответы и~решения}
%\Readsolutionfile{_hintBholder}

\paragraph*{\ref{pr3-3-1}.}
Высота прямоугольного треугольника является средним геометрическим отрезков, на которые она делит гипотенузу. Поэтому если даны отрезки с~длинами $a$, $b$, то, построив полуокружность с~диаметром $a+b$ и~найдя её пересечение с~прямой, перпендикулярной диаметру и~делящей его на отрезки длины $a$ и~$b$, получим отрезок длины~$\sqrt{ab}$. Для решения данной задачи достаточно последовательно построить отрезки с~длинами $z_1=\sqrt{xy}$, $z_2=\sqrt{yz_1}$, $z_3=3x+z_2$, $z=\sqrt{yz_3}$.

\newpage
\setcounter{page}{385}
\section{Стереометрия}\label{stereo}

\begin{flushright}

{\small

Чужбина так же сродственна отчизне,\\
Как тупику соседствует пространство.\\[1mm]
\emph{И.~Бродский.}

}

\end{flushright}

\subsection{Рисование (2). \emph{А.~Б.~Скопенков}}

%%\Opensolutionfile{_hintAholder}
%\Opensolutionfile{_hintBholder}

\begin{pr}\label{pr4-1-1}
Куб с~ребром 3 разбит на 27 единичных кубиков. На\-ри\-суйте

(a) ежа (объединение центрального кубика и~имеющих с~ним общую грань);

(b) то, что получается при выкидывании ежа из куба;

(c)$^*$ то, что получается при выкидывании угловых кубиков из куба.
\end{pr}

\begin{pr}\label{pr4-1-2}
Можно ли пространство заполнить непересекающимися ежами?
\end{pr}

\begin{pr}\label{pr4-1-3}
Правильные многоугольники c каким числом сторон могут получиться при
пересечении куба плоскостью?
\end{pr}

\begin{pr}\label{pr4-1-4}
(a) Нарисуйте объединение куба $A\ldots D_1$ с~кубом, полученным из
него поворотом на $\pi/3$ относительно большой диагонали.

(b) Нарисуйте объединение тетраэдра $ABCD$ с~тетраэдром, полученным
из него поворотом на $\pi/2$ относительно бимедианы, т.\,е. прямой,
соединяющей середины противоположных рёбер.
\end{pr}

\begin{pr}\label{pr4-1-5}
На плоскости стоят куб и~каркас треугольной пирамиды, высота которой
больше высоты куба. Нарисуйте тень от каркаса пирамиды на кубе, если
пучок света параллелен прямой, соединяющей вершину пирамиды
с~центром верхней грани куба.
\end{pr}

\begin{pr}\label{pr4-1-6}
(a) Нарисуйте тело, проекции которого на три взаимно
ортогональные плоскости являются треугольником, квадратом и~кругом
соответственно.

(b) Проекции пространственной фигуры на две пересекающиеся плоскости
являются прямыми линиями. Обязательно ли эта фигура "--- прямая
линия?
\end{pr}

\newpage
\setcounter{page}{392}
\refstepcounter{subsection}
\subsection{Многомерье (4$^*$). \emph{А.~Я.~Канель-Белов}}

\subsubsection{Простейшие многогранники в~многомерном пространстве. \emph{Ю.~М.~Бурман}, \emph{А.~Я.~Канель-Белов}}\label{s15.3}

Хорошо известно, что точке плоскости можно сопоставить пару чисел "--- её декартовых координат (для этого нужно предварительно выбрать систему координат, то есть начало координат и~оси). Тем самым плоскость можно понимать просто как множество всевозможных пар $(x_1,x_2)$ действительных чисел. Аналогично трёхмерное пространство можно считать просто множеством всевозможных троек $(x_1,x_2,x_3)$. Накладывая на числа различные ограничения, мы получим описание разнообразных подмножеств плоскости и~пространства (плоских фигур и~трёхмерных тел).

\begin{pr}\circpr\label{pr4-4-1}
Даны три набора условий на числа $x_1, x_2, x_3$:

%\begin{arab}
%\item
1) $x_1 = x_2 = 2x_3$;

%\item
2) $x_1 + 2x_2 + 3x_3 = 0$, $3x_1 + 2x_2 + x_1 = 1$;

%\item
3) $x_1^2 + x_3^2 - 2x_3 = -1$.
%\end{arab}

Какие из них задают прямую в~трёхмерном пространстве?
\end{pr}

Когда измерений больше, чем три, координатный подход становится ведущим: удобно \emph{определить}, скажем, четырёхмерное пространство как множество всевозможных наборов $(x_1, x_2, x_3, x_4)$ из четырёх действительных чисел.

В этом пункте \emph{отрезком} мы будем называть множество $[-1,1] =
\{x\colon |x|\le 1\}$ чисел, по модулю не превосходящих $1$; \emph{квадратом} "--- множество $[-1,1]^2 = \{(x_1,x_2)\colon |x_1|, |x_2|\le 1\}$ пар чисел, каждое из которых по модулю не превосходит $1$; \emph{кубом} "--- множество $[-1,1]^3 = \{(x_1,x_2,x_3)\colon |x_1|, |x_2|, |x_3|\le 1\}$ троек таких чисел; \emph{четырёхмерным кубом} "--- четвёрок и~т.\,д. См.~рисунок~\ref{Fg:Cubes}.
\index{Куб!многомерный}

\begin{figure}[ht]\centering
\includegraphics[scale=.78]{burman.1}
\caption{Кубы размерностей 1, 2, 3, 4}
\label{Fg:Cubes}
\end{figure}

\begin{pr}\label{pr4-4-2}
%\quest\label{Qu:CoordVert}
(a) Расставьте на рисунке координаты вершин.
%\quest\label{Qu:NumVert}

(b) Сколько вершин у~$n$"~мерного куба?
%\quest\label{Qu:Joined}

(c) Какие вершины $n$"~мерного куба соединены между собой ребром, а~какие нет?

\comment
Выпуклый многогранник с~$n+1$ вершиной, попарные расстояния между которыми равны, называется $n$"~\emph{мерным правильным тетраэдром} или \emph{симплексом}, а~выпуклый многогранник с~$2n$ вершинами в~точках $(x_1=\ldots=x_{i-1}=x_{i+1}=\ldots=x_n=0, x_i=\pm a), i=1,\ldots,n$, "--- \emph{правильным $n$"~мерным октаэдром}.
\index{Симплекс!правильный} \index{Октаэдр!многомерный}
\end{pr}
\endcomment

%%%!!!
\newpage
\setcounter{page}{397}
\subsubsection{Многомерные объёмы}

Объём $n$"~мерного многогранника определяется аналогично площади фигуры на плоскости (см.~п.\,25.5%\ref{ss:dirgeo2} 
 <<Принцип Дирихле и~его применения в~геометрии>>).

\emph{Объёмом} $n$"~мерных многогранников называется заданная на множестве многогранников неотрицательная функция~$V$, удовлетворяющая следующим условиям:

"--* если многогранник $M_1$ можно движением перевести в~многогранник $M_2$, то $V(M_1)=V(M_2)$;

"--* $V(M_1\cup M_2)=V(M_1)+V(M_2)-V(M_1\cap M_2)$;

"--* объём любого подмножества $(n-1)$"~мерной гиперплоскости
равен нулю;

"--* объём куба с~ребром $a$ равен $a^n$. \index{Объём}

Используя эти свойства и~при необходимости верхние и~нижние оценки, можно найти объём любого многогранника. Например, объём $n$"~мерной пирамиды задаётся формулой $V=\frac1{n}Sh$, где $S$ "--- $(n-1)$"~мерный объём основания пирамиды, а~$h$ "--- её высота. Можно также находить объёмы некоторых $n$"~мерных тел (т.\,е. ограниченных подмножеств $n$"~мерного пространства), не являющихся многогранниками.

\begin{pr}\label{pr4-5-1}
У 100-мерного арбуза (шара) радиус равен 1 метру, а~толщина корки  "--- 1~см.~Какой процент его объёма занимает мякоть?
\end{pr}

\begin{pr}\label{pr4-5-2}
Докажите, что в~единичный куб достаточно большой размерности можно поместить здание МГУ, т.\,е. существует трёхмерная плоскость, в~пересечение которой с~кубом можно поместить это здание.
\end{pr}

\begin{pr}\label{pr4-5-3}
Укажите какое-нибудь такое $n$, что в~$n$"~мерный единичный куб можно поместить круг радиуса~$R$.
\end{pr}

\begin{pr}\label{pr4-5-4}
Укажите какое-нибудь такое $n$, что в~$n$"~мерный единичный куб можно поместить шар радиуса~$R$.
\end{pr}

\begin{pr}\label{pr4-5-5}
Укажите какое-нибудь такое $n$, что в~$n$"~мерный единичный куб можно поместить $n$"~мерный шар радиуса~$R$.
\end{pr}

\begin{pr}\label{pr4-5-6}
К чему стремится объём $n$"~мерного шара радиуса 2015 при $n\to\infty$?

Известно, что объём $n$"~мерного шара радиуса $R$ равен
 $$
 B_n=\frac{\pi^{n/2}}{\Gamma(n/2+1)},
 $$
где $\Gamma(z)=\int\limits_0^\infty y^ze^{-y}dy$, $z>0$ "--- знаменитая \emph{гамма-функция} Эйлера. Она доопределяет факториал на комплексную плоскость: $\Gamma(k)=(k+1)!$ при целом $k$ и~$\Gamma(z)=\Gamma(z-1)z$. Последнее равенство позволяет доопределить $\Gamma(z)$ также и~при $\Ree(z)<0$. Известно, что $\Gamma(x)\Gamma(1-x)= \pi/\sin(\pi z)$, в~частности $\Gamma(1/2)=\sqrt{\pi}/2$.
\end{pr}

\begin{pr}\label{pr4-5-7}
Найдите площадь поверхности $n$"~мерного шара единичного объёма.
\end{pr}

\begin{pr}\label{pr4-5-8}
Найдите объём $n$"~мерного симплекса с~единичным ребром. Найдите ребро $n$"~мерного симплекса с~единичным объёмом. (Определения $n$"~мерных симплекса и~октаэдра приведены в~п.\,\ref{s15.3} <<Комбинаторная геометрия в~многомерном пространстве>>.)
\end{pr}

\emph{Диаметром} ограниченного подмножества $M$ $n$"~мерного
пространства называется $\sup\{|XY|, X,Y\in M\}$, где $\dist(X, Y)$ "--- расстояние между точками $X$ и~$Y$.
\index{Диаметр!множества}

\begin{pr}\label{pr4-5-9}
Найдите объём $n$"~мерного октаэдра с~единичным ребром. Найдите диаметр $n$"~мерного симплекса с~единичным объёмом.
\end{pr}

\endgroup

\begingroup

\chapter{Комбинаторика}

%убрать \label{*} в~\begin{hintB}\label{*}

%\Closesolutionfile{_hintAholder}
%\Closesolutionfile{_hintBholder}

\setcounter{page}{441}
\refstepcounter{section}
\section{Подсчеты в~комбинаторике}\label{s17}
\index{Подсчет} \index{Комбинаторика}

Этот параграф посвящен в~основном вопросу <<Сколько существует
объектов с~данными свойствами?>>. В~нем собраны материалы для самого
первого знакомства с~подсчетами в~комбинаторике. Продолжить их
изучение мы рекомендуем по главе~1 книги \cite{GDI2}.

\subsection{Подсчеты числа способов (1). \emph{А.~А.~Гаврилюк}, \emph{Д.~А.~Пермяков}}\label{s:count-sim}

%\Opensolutionfile{_hintAholder}
%\Opensolutionfile{_hintBholder}

Этот пункт не требует никаких знаний и~подходит для первого знакомства с~комбинаторикой.

\begin{pr}\label{pro1-1-1}
(a) Назовем натуральное число \emph{симпатичным}, если в~его записи встречаются только четные цифры. Выпишите все двузначные симпатичные числа и~подсчитайте их количество.

(b) Сколько существует пятизначных симпатичных чисел?

(c) Сколько существует шестизначных чисел, в~записи которых есть
хотя бы одна четная цифра?

(d) Каких семизначных чисел больше: тех, в~записи которых есть
единица, или остальных?
\end{pr}

\begin{pr}\label{pro1-1-2}
Из двух математиков и~десяти экономистов надо составить комиссию из
восьми человек. Сколькими способами можно составить комиссию, если в~нее должен входить хотя бы один математик?
\end{pr}

\begin{pr}\label{pro1-1-3}
(a) Найдите сумму всех семизначных чисел, которые можно получить
всевозможными перестановками цифр $1,\ldots,7$.

(b) Из цифр $1,2,3,\ldots, 9$ составлены все четырехзначные числа, не
содержащие повторяющихся цифр. Найдите сумму этих чисел.

(c) Найдите сумму всех четырехзначных чисел, не содержащих
повторяющихся цифр.
\end{pr}

\begin{pr}\label{pro1-1-4}
(a) На двух клетках шахматной доски стоят черный и~белый короли. За
один ход можно пойти любым королем (короли дружат, так что могут
стоять в~соседних клетках, но не в~одной и~той же). Могут ли в~результате их передвижений встретиться все возможные варианты
расположения этих королей, причем ровно по одному разу?

(b) Тот же вопрос, если короли разучились ходить по диагонали.
\end{pr}

\begin{pr}\label{pro1-1-5}
(a) Найдите сумму всех $6$"~значных чисел, получаемых при всех
перестановках цифр $4, 5, 5, 6, 6, 6$.

(b) Найдите сумму всех 10-значных чисел, получаемых при всех
перестановках цифр $4, 5, 5, 6, 6, 6, 7, 7, 7, 7$.
\end{pr}

\begin{pr}\label{pro1-1-6}
(a) Тому Сойеру поручили покрасить забор из $8$~досок в~белый цвет. В~силу  своей лени он покрасит не более $3$~досок. Сколько у~него способов это сделать?

(b) А~сколько способов покрасить не более $5$~досок?

(c) А~сколько способов покрасить любое количество досок?
\end{pr}

%\Closesolutionfile{_hintAholder}
%\Closesolutionfile{_hintBholder}

\sseccol{Указания, ответы и~решения}

%\Readsolutionfile{_hintBholder}

\paragraph*{\ref{pro1-1-1}.}
\emph{Ответы}: (b) 2500; \quad (c) $884\,375$; \quad (d)
в которых есть единица.

(b) \emph{Решение} (написано А. Колоченковым). Первой цифрой симпатичного числа может быть 2, 4, 6, или 8 "--- всего 4 варианта. Для каждой цифры со второй по пятую есть 5 вариантов: 0, 2, 4, 6, 8. Значит, всего симпатичных чисел
$4\cdot5\cdot5\cdot5\cdot5=2500$.

Это рассуждение в комбинаторике называется \emph{правилом произведения} и~подробно обсуждается в статье \cite{Vi71}.

\newpage
\setcounter{page}{448}
\refstepcounter{subsection}
\subsection{Формула включений и~исключений (2). \emph{Д.~А.~Пермяков}}\label{0inex}

%\Opensolutionfile{_hintAholder}
%\Opensolutionfile{_hintBholder}

Этот пункт посвящен доказательству и~использованию формулы включений и~исключений. Она позволяет отвечать на вопрос <<Сколько существует объектов с~данными свойствами?>> во многих непростых случаях. Потребуются базовые навыки решения задач по комбинаторике. В~частности, нужно уметь приводить строгие доказательства с~использованием взаимно однозначных соответствий, правил суммы и~произведения. Например, полезно прорешать п.\;\ref{s:count-sim} <<Подсчеты числа способов>> или задачи из статьи \cite{Vi71}.

\begin{pr}\label{inex-perm}
Сколькими способами можно переставить числа от $1$~до $n$, чтобы

(a) и~$1$, и~$2$~не оказались на своем месте;

(b) ровно одно из чисел $1$, $2$ и~$3$ оказалось на своем месте;

(c) каждое из чисел $1$, $2$ и~$3$ оказалось не на своем месте;

(d) каждое из чисел $1$, $2$, $3$ и~$4$ оказалось не на своем месте?
\end{pr}

Обозначим через $\varphi(n)$ функцию Эйлера, т.\,е. количество чисел от $1$ до $n$, взаимно простых с~числом $n$. \index{Функция!Эйлера}

\begin{pr}\label{vkl-iskl-4isla}
(a) Найдите количество целых чисел от 1 до $1001$, не делящихся ни на одно из чисел $7$, $11$, $13$.

(b) Найдите $\varphi(1)$, $\varphi(p)$, $\varphi(p^2)$, $\varphi(p^{\alpha})$, где $p$ "--- простое число, $\alpha>2$.

%{vkl-iskl-euler_pt}

%\label{vkl-iskl-euler}

(c) Докажите, что $\varphi(n) = n \Big(1 - \frac{1}{p_1}\Big) \ldots \Big(1-\frac{1}{p_s}\Big)$, где $n=p_1^{\alpha_1}\cdot \ldots \cdot p_s^{\alpha_s}$ "--- каноническое разложение числа $n$.
\end{pr}

\begin{pr}\label{prolov-coat}
(a) На полу комнаты площадью $24\;\text{м}^2$ расположены три ковра (произвольной формы) площадью $12\;\text{м}^2$ каждый. Тогда площадь
пересечения некоторых двух ковров не меньше $4\;\text{м}^2$.

(b) На кафтане расположено пять заплат (произвольной формы). Площадь
каждой из них больше трех пятых площади кафтана. Тогда площадь общей
части некоторых двух заплат больше одной пятой площади кафтана.

(c)$^*$ То же, что в~п.~(b), если площадь каждой заплаты больше  \emph{половины} площади кафтана.
\end{pr}

В этом пункте предлагаются задачи следующего типа: даны конечное множество $U$ и~набор свойств (подмножеств) $A_k\subset U$, $k=1,\ldots,n$. Требуется найти количество элементов, для которых выполнено хотя бы одно из свойств $A_k$ (т.\,е. $|A_1\cup\ldots\cup A_n|$), либо количество элементов, для которых не выполнено ни одно из свойств $A_k$ (т.\,е. $|U - (A_1\cup\ldots\cup A_n)|$). Для этого используются два варианта формулы включений и~исключений
(см.~задачу \ref{vkl-iskl-m}\,(b)). При этом если во всех пересечениях множеств набора число элементов зависит только от количества пересекаемых множеств, то формулу можно упростить (см.~задачу \ref{vkl-iskl-m}\,(a)).

\begin{pr}\label{inex-34}\label{pro1-3-5}
Рассмотрим подмножества $A_1,A_2,A_3, A_4$ конечного множества~$U$. Докажите равенства

(a) $A_1\cup A_2=(A_1\backslash A_2)\sqcup(A_1\cap
A_2)\sqcup(A_2\backslash A_1)$;

(b) $|A_1\cup A_2|=|A_1|+|A_2|-|A_1\cap A_2|$;

(c) $|A_1\cup A_2\cup A_3|=|A_1|+|A_2|+|A_3|-|A_1\cap A_2|-|A_2\cap
A_3|-|A_1\cap A_3|+|A_1\cap A_2\cap A_3|$.

(d) Количество элементов в~$U$, не принадлежащих ни одному из подмножеств $A_1$, $A_2$, $A_3$, равно
$$
|U|-|A_1|-|A_2|-|A_3|+|A_1\cap A_2|+|A_2\cap A_3|+|A_1\cap
A_3|-|A_1\cap A_2\cap A_3|.
$$

(e) Для $k=1,2,3,4$ обозначим
$$
M_k:=\sum\limits_{1\le i_1<\ldots<i_k\le 4}|A_{i_1}\cap A_{i_2}\cap
\ldots\cap A_{i_k}|.
$$
Докажите, что количество элементов в~$A$, не принадлежащих ни одному
из $A_i$, равно $|U|-M_1+M_2-M_3+M_4$.

(f) В~условиях п.\,(e) количество элементов, принадлежащих ровно одному из множеств $A_i$, равно $M_1-2M_2+3M_3-4M_4$.
\end{pr}

\begin{pr}\label{vkl-iskl-m}
\textbf{Формула включений и~исключений.} \index{Формула!включений и~исключений} Рассмотрим подмножества $A_1, \ldots, A_n$ конечного множества $U$. Положим по определению $\Big|\bigcap\limits_{j\in \emptyset}A_j\Big|:=U$.

(a) Пусть число $\alpha_{|S|}:=\Big|\bigcap\limits_{j\in S}A_j\Big|$ зависит только от размера $|S|$ набора $S \subset \{1,\ldots,n\}$ индексов, а~не от самого набора. Тогда
\begin{align*}
|A_1\cup\ldots\cup A_n| &= \sum_{k=1}^n (-1)^{k+1} \binom{n}{k} \alpha_k ,\\
|U - (A_1\cup\ldots\cup A_n)| &= \sum_{k=0}^n (-1)^{k} \binom{n}{k} \alpha_k.
\end{align*}

(b) Обозначим $M_k:=\sum\limits_{S \in \binom{n}{k}}\Big|\bigcap\limits_{j\in S}A_j\Big|$, где суммирование производится по всем $k$-элементным подмножествам множества $\{1,\ldots,n\}$. В~частности, $M_0:=|U|$. Тогда
\begin{align*}
|A_1\cup\ldots\cup A_n| &= M_1 - M_2 + M_3 - \ldots +(-1)^{n+1} M_n,\\
|U - (A_1\cup\ldots\cup A_n)| &= M_0 - M_1 + M_2 + \ldots
+(-1)^{n} M_n.
\end{align*}

(c) \textbf{Неравенства Бонферрони.} \index{Неравенство!Бонферрони}
Для любого $0\le s<n/2$ справедливы неравенства
\begin{multline*}
M_1-M_2+M_3-\ldots-M_{2s}\le|A_1\cup\ldots\cup A_n|\le M_1-M_2+M_3-\ldots+M_{2s+1},\\
 \shoveleft{M_0-M_1+M_2-\ldots+M_{2s}\ge|U - (A_1\cup\ldots\cup A_n)|\ge}\\
\ge M_0-M_1+M_2-\ldots-M_{2s+1}.
\end{multline*}

(d) Число элементов, принадлежащих ровно $r$ из подмножеств $A_1, \ldots, A_n$, равно $\sum\limits_{k=r}^n(-1)^{k-r}\binom{k}{r}M_k$.
\end{pr}

\begin{pr}\label{vkl-iskl-knigi}
На полке стоят $10$ различных книг.

(a) Сколькими способами их можно переставить так, чтобы ни одна
книга не осталась на своем месте?

(b) Количество таких перестановок книг, при которых на месте
остаётся ровно $4$~книги, больше $50\,000$.
\end{pr}

В следующей задаче в~ответе можно использовать суммы (аналогично формуле включений и~исключений).

\begin{pr}\label{vkl-iskl-turisti}\label{pro1-3-7}
(a) Сколькими способами можно расселить $20$~туристов по $5$~различным домикам, чтобы ни один домик не оказался пустым?

(b) Сколько существует различных сюръекций $f\colon \Z_{k}\to \Z_n$?
\end{pr}

\begin{pr}\label{inex-rou}\label{pro1-3-8}
По кругу стоят числа $1,2,\ldots,n$. Найдите число способов выбрать $k$ из них, чтобы никакие два выбранных числа не стояли рядом.

(b) Найдите число способов рассадить $n$ пар враждующих рыцарей за круглый стол с~нумерованными местами, чтобы никакие два враждующих рыцаря не сидели рядом.
\end{pr}

\begin{pr}\label{inex-cube}\label{pro1-3-9}
Куб с~ребром длины $20$~разбит на $8000$~единичных кубиков, и~в~каждом кубике записано число. Известно, что в~каждом столбике из $20$~кубиков, параллельном ребру куба, сумма чисел равна~$1$ (рассматриваются столбики всех трех направлений). В~некотором кубике записано число 10. Через этот кубик проходят три слоя $1\times20\times20$, параллельные граням куба. Найдите сумму всех
чисел вне этих слоев.
\end{pr}

\begin{pr}\zvezda\label{inex-tick}
Сколько существует шестизначных трамвайных билетов, в~которых нет
двух семерок рядом и~всего

(a) не более трех семерок; \quad

(b) не более четырех семерок; \quad

(c) сколько угодно семерок?
\end{pr}

\begin{pr}\zvezda\label{vkl-iskl-mnogo4len}
Докажите следующую формулу:
\vspace{-0.3cm}
\begin{multline*}
n! \cdot x_1 x_2\ldots x_n = (x_1 + x_2+ \ldots + x_n)^n -{}\\
{}- \sum_{1 \le i_1<i_2<\ldots < i_{n-1} \le n} (x_{i_1} + x_{i_2}+ \ldots + x_{i_{n-1}})^n +{}\\
+ \sum_{1\le i_1<i_2<\ldots < i_{n-2} \le n} (x_{i_1} +
x_{i_2}+ \ldots + x_{i_{n-2}})^n - \ldots + (-1)^{n-1} \sum_{i=1}^n
x_i^n.
\end{multline*}
%\begin{hintB}
%В качестве множеств $A_i$ рассмотрим упорядоченные мономы степени $n$, т.\,е. выражения
%вида $x_{k_1}x_{k_2}\ldots x_{k_n}$, которые не содержат переменную $x_i$.
%Тогда первое слагаемое в~правой части равенства "--- это множество всех мономов такого вида,
%второе "--- объединение по $A_i$, третье "--- объединение по пересечениям, и~т.~д.

%В итоге, мы получим своеобразную формулу включения-исключения (где вместо количества элементов выступают сами элементы множеств),
%а в~левой части мы получим все упорядоченные мономы, где все элементы различны, т.\,е. они по сути равны $x_1 x_2\ldots x_n$,
%а количество равно количеству способов упорядочить числа от $1$ до $n$, т.\,е. $n!$.
%\end{hintB}
\end{pr}

\newpage
\setcounter{page}{494}
\refstepcounter{section}
\refstepcounter{section}
\section{Конструкции и~инварианты}\label{s20}

Эта тема доступна и~для учеников 6"--~7 классов, но тогда нужно
пользоваться не этим параграфом, а~статьёй~\cite{Lvovsky-Toom-89} и~соответствующим разделом книги~\cite{Genkin-etal-94}.

\subsection{Конструкции\protect\footnotemark\ (1). \emph{А.~В.~Шаповалов\protect\footnotemark}}\index{Конструкция}

%\Opensolutionfile{_hintAholder}
%\Opensolutionfile{_hintBholder}

\addtocounter{footnote}{-1}
\footnotetext{Эта подборка задач составлена по книгам
\cite{Shapovalov-14} и~\cite{Shapovalov-15}.}

\addtocounter{footnote}{1}
\footnotetext{\url{http://www.ashap.info}.}

Если на вопрос <<Может ли?>> вы подозреваете ответ <<Может>>, то
стоит спросить себя: <<Как такое может быть?>>. Уточните вопрос:
<<Какими свойствами эта конструкция должна обладать?>>. Дополнительное знание поможет сильно сузить круг поисков. Задавайте
себе вопросы на протяжении всего построения. Вы с~удивлением увидите, как много конструкций окажутся логичными и~единственно возможными.

Часто примеров много, а~нужен только один. Избыток свободы может сбивать с~толку: неясно, с~чего начинать. Примените \emph{здравый
смысл}, \emph{естественные соображения}. Они ограничивают поле для
поиска примера, но зато поиск убыстряется и~облегчается. Вообще, ваш
опыт гораздо больше, чем вы думаете. Ответом может оказаться \emph{хорошо знакомый объект}, просто надо посмотреть на него под
нужным углом.

\begin{pr}\label{pri1-1-1}
У двух треугольников равны по две стороны, а~также равны высоты,
проведённые к~третьей стороне. Обязательно ли эти треугольники равны?
\end{pr}

\begin{pr}\label{pri1-1-2}
Верно ли, что в~вершинах любого треугольника можно поставить по
положительному числу так, чтобы длина каждой стороны была равна
сумме чисел в~её концах?
\end{pr}

\begin{pr}\label{pri1-1-3}
В~кружке у~каждого участника ровно по 6 друзей. Может ли у~каждой
пары участников быть ровно по два общих друга?
\end{pr}

Конструкцию с~большим числом деталей проще строить из одинаковых
<<кирпичей>>. Даже если все они одинаковыми быть не могут,
попробуйте взять одинаковых побольше. Можно ещё выбрать два вида
деталей и~посчитать, сколько нужно тех и~других.

Ну, а~если детали <<для сборки>> заданы и~они разные? Тогда стоит
попытаться объединить эти части в~\emph{одинаковые блоки},
и~строить из блоков.

\begin{pr}\label{pri1-1-4}
Назовём неотрицательное целое число \emph{зеброй}, если в~его
записи строго чередуются чётные и~нечётные цифры и~среди цифр есть
не менее трёх различных. Может ли разность двух $100$"~значных зебр
быть $100$"~значной зеброй?
\end{pr}

\begin{pr}\label{pri1-1-5}
Грани параллелепипеда со сторонами 3, 4 и~5 разбиты на единичные
клетки. В~каждую клетку вписали по натуральному числу. Рассмотрим
всевозможные кольца шириной в~одну клетку, параллельные
какой-нибудь грани. Может ли сумма чисел в~каждом таком кольце быть
одной и~той же?
\end{pr}

В задачах, где требуются равные части, приходится выбирать форму
частей. Тут может помочь такое соображение: части заведомо равны,
если они получаются друг из друга симметрией, сдвигом или поворотами. Так, для квадрата популярны разрезания, переходящие в~себя при повороте на~$90^{\circ}$, а~для правильного треугольника "--- при
повороте на~$120^{\circ}$. Для симметричных объектов поиск примера
начинают с~симметричных или <<почти симметричных>> конструкций.
Симметрия и~идея <<расположить объекты по кругу>> применима и~в~негеометрических задачах. \index{Симметрия}

\begin{pr}\label{pri1-1-6}
Можно ли рёбра куба занумеровать числами $-6$, $-5$, $-4$, $-3$,
$-2$, $-1$, $1$, $2$, $3$, $4$, $5$, $6$ так, чтобы для каждой
тройки рёбер, выходящих из одной вершины, сумма была одинакова?
\end{pr}

\begin{pr}\label{pri1-1-7}
Круг разрезали на несколько равных частей. Обязательно ли граница
каждой части проходит через центр круга?
\end{pr}

Если к~конструкции предъявляются противоречивые требования, присмотритесь внимательнее. Часто эти противоречия мнимые. Так, \emph{большой} периметр не противоречит \emph{малой} площади. Вообще, словам <<много>>, <<мало>>, <<сильно>> нужно уметь придать в~решении точный математический смысл с~помощью уравнений и~неравенств.

\begin{pr}\zvezda\label{pri1-1-8}
В~море плавает айсберг в~форме выпуклого многогранника. Может ли
случиться, что $90\,\%$ его объёма находится ниже уровня воды и~при
этом больше половины его поверхности находится выше уровня воды?
\end{pr}

\begin{pr}\label{pri1-1-9}
Есть три игральных кубика с~нестандартными наборами чисел на гранях.
Скажем, что кубик А~\emph{выигрывает} у~кубика B, если при их
одновременном бросании число на A будет больше числа на B
с~вероятностью \emph{больше} $0{,}5$. Может ли первый кубик
выигрывать у~второго, второй "--- у~третьего, а~третий "---
у~первого? \index{Вероятность}

(Приведём равносильную формулировку этой же задачи, не использующую
понятие вероятности: для пары кубиков A и~B составим 36
упорядоченных пар вида (грань A, грань B). Заменим в~каждой паре
грань на число, стоящее на грани. Кубик А~\emph{выигрывает} у~B,
если более чем в~половине пар первое число больше второго.)
\end{pr}

Помешать решить задачу могут невидимые барьеры в~голове решателя.
Если очевидного решения не видно, надо расширять список вариантов,
по возможности до полного. \emph{Инерция мышления} проявляется
в~том, что ключевой вариант пропускают либо не подозревают, что
вариантов более одного. Примените <<метод Шерлока Холмса>>:
отбросьте все невозможные случаи, тогда \emph{последний вариант}
окажется возможным, каким бы невероятным он ни казался.

\begin{pr}\label{pri1-1-10}
\begin{figure}[ht]\centering
\includegraphics%[width=2.5cm]
{construct-01.mps}
\caption{}
\label{ris-construct-01}
\end{figure}
На столе лежат 9 яблок, образуя 10 рядов по 3 яблока в~каждом (см.~рис.~\ref{ris-construct-01}). Известно, что у~девяти рядов веса
одинаковы, а~вес десятого ряда отличается. Есть электронные весы, на
которых за рубль можно узнать вес любой группы яблок. Какое
наименьшее число рублей надо заплатить, чтобы узнать, вес какого
именно ряда отличается?
\end{pr}

\begin{pr}\label{pri1-1-11}
Может ли прямая разбить какой-нибудь шестиугольник на 4 равных
треугольника?
\end{pr}

\index{Редукция}
\emph{Редукция} "--- это свед\'ение сложной задачи к~более простой. Так, если сложную конструкцию не удаётся сразу построить целиком, постройте её \emph{необходимую часть}. Даже если эту часть не удастся потом достроить до целого, решение упрощённой задачи может послужить разминкой, после чего вы вернётесь к~сложной задаче уже с~накопленным опытом.

\begin{pr}\label{pri1-1-12}
Барон Мюнхгаузен говорит, что у~него есть многозначное
число-\emph{палиндром} (т.\,е. оно читается одинаково слева направо и~справа
налево). Написав его на бумажной ленте, барон сделал несколько
разрезов между цифрами. Лента распалась на $N$ кусков. Переложив
куски в~другом порядке, барон увидел, что на кусках по разу записаны
числа $1$, $2$, $\ldots$, $N$. Могут ли слова барона быть правдой?
\end{pr}

При построении конструкции может мешать неоднозначность выбора.
В~\emph{узком месте} всё однозначно или неопределённость
минимальна, что сокращает перебор.\index{Узкое место} Начав
с~узкого места, мы либо быстро придём к~противоречию, либо построим
большой кусок конструкции. Как искать узкие места? Присмотритесь:
они служат препятствиями к~построению конструкции или кажутся
таковыми.

\begin{pr}\label{pri1-1-13}
Записав числа $1$, $\frac{1}{2}$, $\frac{1}{3}$, $\ldots$,
$\frac{1}{10}$ в~некотором порядке, соедините их знаками четырёх
арифметических действий так, чтобы полученное выражение равнялось 0.
(Скобки использовать нельзя.)
\end{pr}

\begin{pr}\label{pri1-1-14}
Существуют ли три равных семиугольника, все вершины которых
совпадают, но никакие стороны не совпадают?
\end{pr}

\begin{pr}\label{pri1-1-15}
Можно ли разрезать какой-нибудь треугольник на четыре выпуклые
фигуры: треугольник, четырёхугольник, пятиугольник и~шестиугольник?
\end{pr}

При \emph{постепенном конструировании }к примеру идут через
цепочку вспомогательных конструкций-\emph{заготовок}. На каждом
шаге очередная конструкция \emph{улучшается} до следующей.
В~заготовке требования к~окончательной конструкции выполнены лишь
частично. Оставляем \emph{принципиальные} условия, временно
забываем или ослабляем \emph{технические}.
\index{Конструкция!постепенная}

\begin{pr}\label{pri1-1-16}
Могут ли в~остроугольном треугольнике все стороны и~высоты
измеряться целым числом сантиметров?
\end{pr}

\begin{pr}\label{pri1-1-17}
Докажите, что существует палиндром, делящийся на $6^{100}$. (Напомним, что \emph{палиндром} "--- это число, которое не меняется при записи его цифр в~обратном порядке.)
\index{Палиндром|textbf}
\end{pr}

Наконец, при \emph{конструкции по индукции} результат получается
постепенно, но уже за бесконечное число шагов. Таким конструкциям
посвящён п.\;18.5%\ref{ss:ks}
 <<Конечное и~счётное>>.

Продолжить знакомство с~конструкциями можно по
статье~\cite{Genkin-etal-90} и~книгам \cite{Shapovalov-14,
Shapovalov-15, Shapovalov-08}.

%\Closesolutionfile{_hintAholder}
%\Closesolutionfile{_hintBholder}

\sseccol{Указания, ответы и~решения}

Решения задач и~\emph{пути к~решению} тщательно разделены.
Решение "--- это то, что решающий задачу в~идеале должен написать. Путь
к~решению должен остаться в~голове, здесь он поясняет, как это
решение можно было придумать. В~задачах на конструкцию решение
и~путь к~решению обычно имеют мало общего.

Решение в~задаче на конструкцию состоит из двух частей: \emph{примера}, то есть описания конструкции, и~\emph{доказательства} того, что она удовлетворяет условию задачи. Для наших задач вторая часть не представляет труда и~обычно опускается. Но иногда из многих возможных примеров нужно ещё выбрать тот, для которого доказательство проще.

%\Readsolutionfile{_hintBholder}

\paragraph*{\ref{pri1-1-1}.}
\emph{Ответ}: не обязательно.

\emph{Решение}. Рассмотрим равнобедренный треугольник $ACD$ и~точку $B$ на продолжении основания $DC$. У~треугольников $ABC$ и~$ABD$ сторона $AB$ и~высота $AH$ общие, стороны $AC$ и~$AD$ равны. Однако эти треугольники не равны: один "--- часть другого.

\emph{Путь к~решению}. Попробуем \emph{построить} треугольник по
двум сторонам $b$, $c$ и~высоте $h$, проведённой к~третьей стороне.
Для этого проведём прямую $l$ (на ней будет лежать третья сторона)
и построим вершину $A$ на расстоянии $h$ от $l$. Две другие вершины
треугольника должны лежать на этой прямой на расстояниях $b$ и~$c$
от точки~$A$. Проведя окружности указанных радиусов с~центром в~точке~$A$, получим (при $b>h$ и~$c>h$) по две точки пересечения каждой из окружностей с~$l$. Видим, что с~точностью до симметрии есть два принципиально разных треугольника: когда вершины выбираются по одну сторону от ближайшей к~$A$ точки прямой и~по разные стороны от неё.

\newpage
\setcounter{page}{524}
\refstepcounter{subsection}
\refstepcounter{subsection}
\refstepcounter{subsection}
\subsection{Полуинварианты\protect\footnotemark (1). \emph{А.~В.~Шаповалов}}%
\label{ss:polu}

%\Opensolutionfile{_hintAholder}
%\Opensolutionfile{_hintBholder}

\footnotetext{Идейными предшественниками подборок <<Инварианты>> и~<<Полуинварианты>> были, среди прочего, соответствующие параграфы
книги~\cite{Kanel-etal-08}.}

Если слово <<инвариант>> означает <<неизменный>>, то <<полуинвариант>> "--- неизменный наполовину.

Бывает так, что мы меняем конструкцию, а~какая-то связанная с~этой конструкцией величина может меняться только в~одну сторону, то есть либо только увеличиваться, либо только уменьшаться. Ещё возможно, что мы делаем ходы и~в~одну сторону меняется величина, связанная с~позицией. Например, при игре в~крестики-нолики число заполненных клеток с~каждым ходом увеличивается. На ограниченной доске из этого следует, что рано или поздно игра закончится. При игре на бесконечной доске игра может не закончиться никогда, но зато мы можем гарантировать, что позиция не повторится, "--- ведь число заполненных клеток каждый раз новое!

Чуть более формально: пусть мы меняем конструкции (или позиции) с~помощью \emph{разрешённых операций} (или \emph{ходов}) и~нам
удалось связать с~каждой конструкцией/позицией \emph{величину},
значение которой при любом разрешённом преобразовании либо не
меняется, либо меняется всегда в~одну и~ту же сторону. Тогда эта
величина называется \emph{полуинвариантом}\footnote{Эта фраза не
является формальным определеним полуинварианта. Но для решения
задач формальное определение этого понятия не нужно.}. Если
полуинвариант меняется при каждой операции/ходе, он называет
\emph{строгим}, иначе "--- \emph{нестрогим}.
\index{Полуинвариант!строгий} \index{Полуинвариант!нестрогий}

В типовых задачах <<на полуинвариант>> доказывают невозможность а)
повторения позиций; б) бесконечного числа ходов; в) построения
конструкций. Для последнего находят полуинвариант и~проверяют, что
для получения искомой конструкции из исходной полуинвариант должен
был бы \emph{меняться не в~ту сторону}.

Но как найти полуинвариант? Начните с~проверки типовых величин:
сумм, произведений, площадей, периметров и~их комбинаций. Если
конструкция зависит от целых чисел, то полуинвариантом может быть
НОД или НОК.

В следующих двух задачах важно, что полуинвариант целочисленный и~не
может быть больше определённого числа.

\begin{pr}\label{pri1-7-1}
На шахматной доске $100\times 100$ королю разрешено ходить вправо,
вверх или вправо-вверх по диагонали. Какое наибольшее число ходов он
может сделать?
\end{pr}

\begin{pr}\label{pri1-7-2}
В~клетках таблицы $99\times 99$ расставлены целые числа. Если
в~каком-то ряду (строке или столбце) сумма отрицательна, разрешается
в~этом ряду поменять знаки всех чисел на противоположные. Докажите,
что в~итоге можно сделать лишь конечное число таких операций.
\end{pr}

{Если полуинвариант не целочисленный, то его ограниченность ещё не
гарантирует окончания процесса (например, убывающий положительный
полуинвариант мог бы бесконечно долго принимать значения $1$, $1/2$,
$1/3$, $1/4$, $\ldots,$ $1/n$, $\ldots$). В~этих случаях прекращение
ходов гарантируется конечным числом позиций.}

\begin{pr}\label{pri1-7-3}
Дано 10 чисел. За одну операцию можно два неравных числа заменить на
два равных с~той же суммой. Может ли этот процесс для какого-то
исходного набора чисел

(a) продолжаться бесконечно долго;

(b) зациклиться (то есть может ли один и~тот же набор чисел возникнуть дважды)?
\end{pr}

\begin{pr}\label{pri1-7-4}
По кругу выписано несколько чисел. Если для некоторых четырёх идущих
подряд чисел $a,b,c,d $ оказывается, что $(a-d)(b-c)<0$, то числа
$b$ и~$c$ можно поменять местами. Докажите, что такую операцию можно
проделать лишь конечное число раз.
\end{pr}

Очень часто положение, в~котором нет разрешённых операций,
и~является искомым.

\begin{pr}\label{pri1-7-5}
В~клетки прямоугольной таблицы вписаны числа. Разрешается
одновременно менять знак у~всех чисел некоторого столбца или
некоторой строки. Докажите, что многократным повторением этой
операции можно превратить данную таблицу в~такую, у~которой суммы
чисел в~любой строке или любом столбце неотрицательны.
\end{pr}

В комбинаторных задачах полуинвариантом часто служит число комбинаций, например пар, троек, подмножеств или перестановок какого-то вида.

\begin{pr}\label{pri1-7-6}
В тридевятом царстве все города подняли над ратушами флаги "--- голубые либо оранжевые. Каждый день жители узнают цвета флагов у~соседей в~радиусе 100~км. Один из городов, где у~большинства соседей флаги другого цвета, меняет свой флаг на этот другой цвет. Докажите, что со временем смены цвета флагов прекратятся.
\end{pr}

{Некоторые конструкции создаются <<методом последовательного
улучшения>>. Мы берём несовершенную конструкцию и~начинаем её
преобразовывать. Полуинвариант гарантирует завершение процесса
и~достижение нужного эффекта в~конце.}

\begin{pr}\label{pri1-7-7}
В~парламенте каждый депутат имеет не более трёх врагов. Докажите,
что парламент можно так разбить на две палаты, что у~каждого
депутата в~его палате будет не более одного врага.
\end{pr}

\begin{pr}\label{pri1-7-8}
На плоскости дано $100$~красных и~$100$~синих точек, никакие три из
которых не лежат на одной прямой. Докажите, что можно провести
$100$~непересекающихся отрезков с~концами разных цветов.
\end{pr}

{Полуинвариант может быть и~\emph{нестрогим}, т.\,е. не меняться
при некоторых ходах. Тогда полезно найти ещё один полуинвариант,
который строго меняется как раз тогда, когда первый остаётся
неизменным.}\index{Полуинвариант!нестрогий}

\begin{pr}\label{pri1-7-9}
На шахматной доске $100\times 100$ королю разрешено ходить вправо,
вверх, вправо-вверх или вправо-вниз по диагонали. Докажите, что он
может сделать лишь конечное число ходов.
\end{pr}

{Если и~второй полуинвариант оказывается нестрогим, то приходится
рассматривать и~третий, и~четвёртый и~т.\,д. В~этом случае естественно рассматривать наборы значений полуинвариантов как строки, упорядоченные \emph{лексикографически} (как слова в~словаре:
сравниваются первые элементы, при равенстве "--- вторые и~т.\,д. и~так до первого несовпадения).} \index{Порядок!лексикографический}

\comment
\begin{pr}\label{pri1-7-10}
В~колоде часть карт лежит рубашкой вниз. Время от времени Петя
вынимает из колоды пачку из нескольких подряд идущих карт, в~которой
верхняя и~нижняя карты лежат рубашкой вниз (в~частности, может
вынуть просто одну карту рубашкой вниз), переворачивает эту пачку
как одно целое и~вставляет в~то же место колоды. Докажите, что
независимо от того, как Петя выбирает пачки, в~конце концов все
карты лягут рубашкой вверх.
\end{pr}
\endcomment

\newpage
\setcounter{page}{536}
\section{Алгоритмы}\label{s:alg} \index{Алгоритм}

%\setcounter{secnumdepth}{1}
%\tableofcontents

\subsection{Игры (1)\protect\footnotemark. \emph{Д.~А.~Пермяков}, \emph{М.~Б.~Скопенков}, \emph{А.~В.~Шаповалов}}

%%%!!!

%\Opensolutionfile{_hintAholder}
%\Opensolutionfile{_hintBholder}

\footnotetext{Подпункты <<Симметричная стратегия>>, <<Выращивание
дерева позиций>>, <<Передача хода>> написаны Д.~А.~Пермяковым и~М.~Б.~Скопенковым, <<Игры-шутки>>, <<Игра на опережение>>, <<Накопление преимущества>> "--- А.~В.~Шаповаловым, <<Смесь>> "--- всеми тремя авторами.}

\index{Игра} \index{Симметрия}

На конкретных примерах мы познакомимся с~некоторыми красивыми идеями теории игр.
Общие методические указания по теме <<Игры>> можно найти в~соответствующем разделе книги~\cite{Genkin-etal-94}.

\subsubsection*{Симметричная стратегия}

Самая распространённая стратегия в~играх "--- \emph{симметричная}
(а~также её обобщение "--- \emph{дополняющая}). Для решения
последующих задач полезно знакомство с~п.\;20.2%\ref{ss:inv1} 
 <<Инварианты I>>, поскольку многие стратегии в~играх основаны на инвариантах (пример инварианта "--- симметричность позиции). \index{Инвариант}
\index{Стратегия!cимметричная}\index{Стратегия!дополняющая}

\begin{pr}\label{pa1-1-1}
(a) Двое по очереди выкладывают доминошки на шахматную доску. Каждая
доминошка покрывает ровно две клетки доски, каждая клетка может быть
покрыта не более чем одной доминошкой. Проигрывает тот игрок,
который не может положить очередную доминошку. Кто выигрывает при
правильной игре? Как он должен для этого играть?

(b) То же для доски $8\times 9$.
\end{pr}

Вот что означают вопросы этой задачи. В~ответе на первый вопрос
нужно назвать игрока, который выигрывает при \emph{любой} игре
своего противника. В~ответе на второй вопрос нужно привести
\emph{алгоритм} действий этого игрока, который гарантирует выигрыш
(\emph{выигрышную стратегию}). Важно чётко отделять \emph{сам}
алгоритм от \emph{доказательства} того, что алгоритм приводит к~желаемому результату. \index{Стратегия!выигрышная} \index{Алгоритм}

Второй пункт этой задачи показывает, что не всегда симметричность
позиции гарантирует, что симметричная стратегия работает.

\begin{pr}\circpr\label{pa1-1-2} (Загадка.)
К~какому результату приведёт попытка чёрных зеркально-симметрично
копировать ходы противника в~обычных шахматах при правильной игре
белых? Выберите верный вариант ответа:

1) к~ничьей; \quad 2) к~выигрышу белых; \quad 3) к~выигрышу чёрных.
\end{pr}

Ключевой идеей является не столько симметрия, сколько разбиение всех
возможных позиций на пары. \emph{Дополняющая стратегия} состоит в~том, чтобы на ход противника отвечать ходом во вторую позицию
соответствующей пары.

\begin{pr}\label{pa1-1-3}
На шахматной доске стоит король. Двое по очереди ходят им.
Проигрывает игрок, после хода которого король оказывается в~клетке,
в~которой побывал ранее. Кто выигрывает при правильной игре и~как он
должен для этого играть?
\end{pr}

\index{Игра!-шутка}\index{Игра!на опережение}\index{Накопление
преимущества}

\subsubsection*{Игра на опережение}

\emph{Игра на опережение} "--- распространённый приём
в~нематематических играх. Но и~в~математических играх бывает, что
выигрыш достаётся тому, кто первый сумеет занять ключевое положение.
После этого, как правило, работает дополняющая стратегия.

\begin{pr}\label{pa1-1-4}
Есть $9$ запечатанных прозрачных коробок соответственно с~$1,2,3,\ldots,9$
фишками. Двое играющих по очереди берут по одной фишке из любой
коробки, распечатывая, если необходимо, коробку. Проигрывает тот,
кто последним распечатает коробку. Кто из них может всегда выиграть
независимо от игры противника? \index{Четность}
\end{pr}

\begin{pr}\label{pa1-1-5}
В~одном из углов шахматной доски лежит плоский картонный квадрат
$2\times 2$, а~в~противоположном "--- квадрат $1\times1$. Двое
играющих по очереди перекатывают каждый свой квадрат через сторону:
Боря "--- большой квадрат, а~Миша "--- маленький. Боря выигрывает,
если не позднее 100"~го хода Мишин квадрат окажется на клетке,
накрытой Бориным квадратом. Может ли Боря выиграть независимо от
игры Миши, если

(a) первым ходит Боря;

(b) первым ходит Миша?
\end{pr}

\subsubsection*{Накопление преимущества}

\emph{Накопление преимущества} "--- тоже весьма распространённый
приём в~нематематических играх. В~математических играх накопление
обычно связано с~каким-нибудь полуинвариантом. Поэтому для изучения
таких игр полезно знакомство с~п.\;\ref{ss:polu} <<Полуинварианты>>. При этом
надо придумать алгоритм, ведущий к~накоплению независимо от
сопротивления соперника. \index{Полуинвариант}

\begin{pr}\label{pa1-1-7}
Миша стоит в~центре круглой лужайки радиуса $100$~метров. Каждую
минуту он делает шаг длиной $1$~метр. Перед каждым шагом он
объявляет направление, в~котором хочет шагнуть. Катя имеет право
заставить его сменить направление на противоположное. Может ли Миша
действовать так, чтобы в~какой-то момент обязательно выйти
с~лужайки, или Катя всегда сможет ему помешать?
\end{pr}

\begin{pr}\label{pa1-1-8}
На клетчатой доске $1\times100\,000$ (вначале пустой) двое ходят
по очереди. Первый может за ход выставить два крестика в~любые
два свободных поля доски. Второй может стереть любое количество
крестиков, идущих подряд "--- без пустых клеток между ними. Если
после хода первого образуется~$13$ или более крестиков подряд, он
выиграл. Может ли первый игрок выиграть при правильной игре обеих сторон?
\end{pr}

\begin{pr}\label{pa1-1-9}
Двое играющих по очереди ломают палку: первый на две части, затем
второй ломает любой из кусков на две части, затем первый "---
любой из кусков на две части и~т.\,д. Один из игроков выигрывает,
если сможет после какого-то из своих ходов сложить из $6$ кусков
два равных треугольника. Может ли другой ему помешать?
\end{pr}

\subsubsection*{Игры-шутки}

В~\emph{играх-шутках} побеждает всегда одна из сторон независимо
от её желания.

\begin{pr}\label{pa1-1-10}
(a) На столе лежат $2015$~кучек по одному ореху. За один ход разрешается объединить две кучки в~одну. Двое играющих делают ходы по очереди, кто не сможет сделать ход, тот проигрывает. Кто выиграет?

(b) То же, но разрешается объединять кучки только с~одинаковым числом орехов.
\end{pr}

\begin{pr}\label{pa1-1-11}
Дана клетчатая полоса $1\times N$. Двое играют в~следующую игру. На
очередном ходу первый игрок ставит в~одну из свободных клеток крестик, а~второй "--- нолик. Не разрешается ставить в~соседние клетки два крестика или два нолика. Проигрывает тот, кто не может сделать ход. Кто из игроков выигрывает при правильной игре? Как он
должен для этого играть?
\end{pr}

Кроме игр-шуток бывают и~\emph{почти шутки}, где выигрышная
стратегия такова: если есть выигрыш в~один ход, его надо сделать,
иначе можно делать любой ход. Или, наоборот: делать любой ход, кроме
тех, которые проигрывают в~один ход. В~таких играх важно догадаться,
кому обязательно представится возможность сделать выигрышный ход
или кто будет вынужден сделать проигрывающий ход, "--- и~доказать
это. Кроме того, выигрышная стратегия может состоять в~достижении
позиции, после которой игра превращается в~игру-шутку с~нужным
исходом.

\begin{pr}\label{pa1-1-12}
В десяти корзинах лежат яблоки: 1, 3, 5, \ldots, 19 яблок. Сначала
берёт одно яблоко из любой корзины Вася, потом "--- Гена, потом
Лёва, потом опять Вася и~т.\,д. по кругу. Проигрывает тот, после чьего хода в~каких-то корзинах станет яблок поровну. Кто из них не может избежать проигрыша?
\end{pr}

\begin{pr}\label{pa1-1-13}
Из спичек сложен клетчатый квадрат $9\times9$, сторона каждой клетки
"--- одна спичка. Петя и~Вася по очереди убирают по спичке, начинает
Петя. Выиграет тот, после чьего хода не останется целых квадратиков
$1\times1$. Кто может действовать так, чтобы обеспечить себе победу,
как бы ни играл его соперник?
\end{pr}

\subsubsection*{Выращивание дерева позиций}
\index{Дерево!позиций}\index{Передача хода}

Один из универсальных способов анализа игры "--- \emph{выращивание
дерева позиций}.

\begin{pr}\label{pa1-1-14}
\emph{Ферзя "--- в~угол, или <<цзяньшицзы>>.} Ферзь стоит на d1.
Двое по очереди ходят им по направлению вверх, вправо или
вправо-вверх. Выигрывает тот, кто поставит его на h8. Кто выигрывает
при правильной игре и~как он должен для этого играть?
\end{pr}

Если не получается, подумайте сначала над следующим вопросом.

\begin{pr}\circpr\label{pa1-1-15}
Кто выигрывает в~игре из предыдущей задачи, если в~начальный момент
ферзь стоит на клетке f4? Выберите верный вариант ответа:

1) первый игрок;\ \  2) второй игрок.
\end{pr}

Выращивание дерева позиций означает полный анализ игры. Перейдём
теперь к~более сложной идее \emph{передачи хода}, которая помогает
даже тогда, когда для полного анализа нет никакой возможности.

\subsubsection*{Передача хода}

\begin{pr}\label{pa1-1-16}
В~\emph{двухходовых} шахматах фигуры ходят по обычным правилам, только за каждый ход разрешается сделать ровно два хода одной фигурой. Цель игры "--- съесть короля соперника. Правила троекратного повторения позиции и~$50$~ходов не действуют\footnote{Если не знаете, что это за правила, игнорируйте это предложение.}. Докажите, что белые в~двухходовых шахматах могут играть так, что заведомо не проиграют (т.\,е. либо выиграют, либо сыграют вничью).
\end{pr}

\begin{pr}\circpr\label{pa1-1-17}
Правила \emph{шахмат без цугцванга}\footnote{\emph{Цугцвангом} в~шахматах называется такая позиция для игрока, в~которой любой его
ход эту позицию ухудшает.} отличаются от правил обычных шахмат
только добавлением возможности пропустить свой ход для каждого из
игроков. Могут ли чёрные выиграть при правильной игре белых? Выберите верный вариант ответа:

1) могут;\ \ 2) не могут.

\end{pr}

%\chapter{Миникурс по теории вероятности. (3)}

\newpage
\setcounter{page}{579}
\section{Вероятность\protect\footnotemark. \emph{А.~А.~Заславский}}\index{Вероятность}\label{s22}

\footnotetext{Автор благодарен Ю.~Н.~Тюрину за полезное обсуждение.}

Данный параграф посвящён простейшим понятиям и~применениям теории
вероятности. Для его изучения необходимо знакомство с~основами
комбинаторики, например с~п.\;\ref{s:count-sim} <<Подсчёт числа способов>> и~\ref{0inex} <<Формула включений и~исключений>> данной книги. Кроме того, знакомство с~теорией вероятностей полезно начинать на <<физическом>> уровне строгости, как в~книгах [Shen] и [KzhP]%\cite{Shen}, \cite{Kolm}
. Здесь же мы сразу даём <<математические>> определения. Однако мы приводим многие задачи на <<практическом>> языке и~показываем на примерах, как их формализовать. Формализацию остальных задач оставляем читателю. Такая формализация является первым шагом решения, от которого может зависеть ответ. См., например, задачи \ref{probas-pha}\,(b, c).

\subsection{Классическое определение вероятности (1).}

%%\Opensolutionfile{_hintAholder}
%%\Opensolutionfile{_hintBholder}

Рассмотрим эксперимент, имеющий $m$ равновозможных исходов, например бросание игральной кости, вытаскивание карты из колоды
и~т.\,д. Если интересующее нас событие (например, выпадение
шестёрки, вытаскивание туза и~т.\,д.) происходит в~$a$ из этих
исходов, то \emph{вероятность} события считают равной $p=a/m$.

Это пояснение полезно для начинающего, но не является математическим
определением. Вот математическое определение.

\emph{Вероятностью} подмножества $A$ конечного множества $M$
называется число\index{Вероятность}
$$
P(A)=P_M(A):=|A|/|M|.
$$
Далее, если не оговорено противное, множество $M$ фиксировано и~пропускается из обозначений. Тогда вероятность определена для всех его подмножеств. Их часто называют \emph{событиями}.

\begin{pr}\label{probas52}
Из колоды в~52 карты вытаскивается одна карта. Найдите вероятность того, что она окажется

(a) чёрной масти; \quad (b) тузом; \quad (c) картинкой;

(d) дамой пик; \quad (e) королём или бубной.
\end{pr}

Например, в~задаче \ref{probas52}\,(с) множество $M$ (<<всех возможных исходов>>) совпадает с~множеством карт в~колоде, а~множество $A$ (<<исходов, в~которых происходит рассматриваемое событие>>) "--- с~множеством картинок. Так эта и~многие другие вероятностные задачи могут быть строго сформулированы на комбинаторном языке.

\begin{pr}
Монета бросается 3 раза. Найдите вероятность выпадения

(a) трёх орлов; \quad (b) двух орлов и~решки.
\end{pr}

\begin{pr}
Найдите вероятность того, что при бросании двух игральных костей

(a) на первой выпадет больше очков, чем на второй;

(b) сумма выпавших очков составит $2, 3,\ldots, 12$.
\end{pr}

\begin{pr}\label{probas-chisla}
Найдите вероятность того, что случайное целое число от 1 до 105

(a) делится на 5; \quad (b) делится на 7; \quad (c) делится на 35.

(a$'$, b$'$, c$'$) То же для случайного целого числа от 1 до 100.
\end{pr}

\begin{pr} \label{probas-exam}
Федя знает ответы на 10 вопросов из 30. Билет состоит из двух вопросов. С~какой вероятностью Федя ответит на оба вопроса?
\end{pr}

Для решения некоторых из вышеприведённых задач полезны следующие.

\begin{pr}\label{probas-rule}
(a) \textbf{Правило сложения.} Пусть $A\cap B=\emptyset$. Выразите
$P(A \cup B)$ через $P(A)$ и~$P(B)$.
\index{Правило!сложения вероятностей}

(b) Выразите вероятность $P(A \cup B)$ через $P(A)$, $P(B)$ и~$P(A\cap B)$.

(c) \textbf{Правило умножения.} Выразите вероятность $P_{M\times N}(A\times B)$ через $P_M(A)$ и~$P_N(B)$.
\index{Правило!умножения вероятностей}

\emph{Комментарий}: $P_M(A)=P_{M\times N}(A\times N)$ и~$P_N(B)=P_{M\times N}(M\times B)$.

\end{pr}

\begin{pr}\label{noski}
(a) В~ящике лежат красные и~чёрные носки. Какое минимальное количество носков может быть в~ящике, если вероятность того, что два случайно вытянутых носка красные, равна 1/2?

(b) То же, если дополнительно известно, что число чёрных носков
чётно.
\end{pr}

\comment

\begin{pr}\zvezda\label{procon-geom}
(a) С~какой вероятностью треугольник, образованный тремя случайными
вершинами правильного $2n$"~угольника, будет прямоугольным; остроугольным; тупоугольным?

(Если эта задача не получается, то см.~следующий пункт.)

(b) Найдите пределы полученных вероятностей при $n\to\infty$. (Подумайте о~смысле полученных результатов. Ср. с~задачей \ref{probas-pha}\,(с).)
\end{pr}

%%\Closesolutionfile{_hintAholder}
%%\Closesolutionfile{_hintBholder}

\sseccol{Указания, ответы и~решения}
%%\Readsolutionfile{_hintAholder}

\paragraph*{\ref{probas-chisla}.}
\emph{Ответы}: (a) $0{,}2$;\ (b) $\frac{1}{7}$;\ (c) $\frac{1}{35}$;\ (a$'$) $0{,}2$;\ \ (b$'$) $0{,}14$;\
(c$'$) $0{,}05$.

(a) \emph{Решение} (написано Е.~Павловым). Пусть $M=\{1,2, \dots, 105\}$ "--- множество всех возможных исходов, $A=\{5,10,\dots, 105\}=\{x\in M\colon 5\mid x\}$ "--- множество благоприятных исходов. Тогда по определению вероятность множества $A$ равна $P(A)=\frac{|A|}{|M|}=\frac{\lfloor\frac{105}{5}\rfloor}{105}=0{,}2$.

\paragraph*{\ref{probas-exam}.}
\emph{Ответ}: $\frac{3}{29}.$
\emph{Решение} (написано П.~Белопашенцевой). Обозначим через $M$ множество всех неупорядоченных пар различных чисел от 1 до 30. Это множество соответствует множеству всех возможных билетов. Количество элементов в $M$ равно $|M|={30\choose 2}$.

Обозначим через $A$ множество всех неупорядоченных пар различных чисел от 1 до 10. Множество $A$ соответствует множеству выученных Федей билетов. Количество элементов в $A$ есть $|A|={10\choose 2}$.
Вероятность подмножества $A$ в множестве $M$ есть по определению отношение $ P_M(A)= |A|/|M|={10\choose 2}/{30\choose 2}=\frac{10\cdot 9/2}{30\cdot 29/2}=\frac{3}{29}.$
\endcomment

\newpage
\setcounter{page}{583}
\subsection{Более общее определение вероятности (1)}
\label{ss:gendefprob}

%%\Opensolutionfile{_hintAholder}
%%\Opensolutionfile{_hintBholder}

\begin{pr}\label{strelki} (a) Один стрелок попадает в~цель с~вероятностью $0{,}8$, другой "--- $0{,}7$. Найдите вероятность поражения цели, если оба стреляют одновременно.

(В этой и некоторых других задачах этого пункта формализация приводится после условий.)

(b) Рабочий обслуживает три станка. Вероятности их остановки равны
соответственно $0{,}1$;\ $0{,}2$;\ $0{,}15$. Найдите вероятность
безотказной работы всех станков.
\end{pr}

Для формализации вышеприведённых задач необходимо следующее более
общее определение. Пусть задано множество $M$ и каждому $m\in M$ поставлено в~соответствие неотрицательное число $P(m)$, причём сумма всех этих чисел равна~1. Тогда \emph{вероятностью} события $A$ называется сумма чисел $P(m)$ по всем $m\in A$.

Например, в~вышеприведённой задаче разумно считать, что множество
$M$ состоит из четырёх элементов: оба стрелка попали, первый попал и~второй промахнулся, первый промахнулся и~второй попал, оба
промахнулись.

\begin{pr}\label{probas-gen}
Сформулируйте и~докажите аналоги правил суммы и~произведения для
вышеприведённого обобщения.
\end{pr}

Приведённое определение можно обобщить на случай бесконечного множества~$M$. (В этом случае для всех $m\in M$, кроме счётного числа, $P(m)=0$.) Ещё более интересно следующее обобщение.

\begin{pr}\label{probas-geo}
Найдите вероятность того, что случайная точка правильного
треугольника лежит

(a) в~треугольнике, образованном средними линиями;

(b) во вписанном круге.
\end{pr}

Пусть $A\subset M$ "--- подмножества прямой (или плоскости, или пространства), имеющие длину. Не все подмножества имеют длину (или
площадь или объём), см.~замечание в~п.\,25.5%\ref{ss:dirgeo2}
<<принцип Дирихле и его применения в~геометрии>>. Тогда \emph{вероятностью} подмножества $A$ в~$M$
называется число
$$
P(A)=P_M(A):=L(A)/L(M),
$$
где $L(A),L(M)$ "--- длины подмножеств.

Пусть $A\subset M$ "--- подмножества плоскости (или пространства),
имеющие площадь. Тогда \emph{вероятностью} подмножества $A$ в~$M$
называется
$$
P(A)=P_M(A):=S(A)/S(M),
$$
где $S(A),S(M)$ "--- площади подмножеств. Аналогично определяется вероятность для подмножеств $A\subset M$ пространства, имеющих объёмы.

Как и~в~дискретном случае, когда множество $M$ фиксировано, его
подмножества, имеющие длину (площадь, объём), часто называются
\emph{событиями}.
\index{Вероятность!геометрическая}
\index{Вероятность|textbf}

\begin{pr}\zvezda\label{probas-geng}
Сформулируйте и~докажите аналоги правил суммы и~произведения для
вышеопределенных <<геометрических>> вероятностей.
\end{pr}

\begin{pr}\zvezda\label{probas-pha}
(a) Дуэли в~городе Осторожности редко кончаются печальным исходом.
Дело в~том, что каждый дуэлянт прибывает на место встречи в~случайный момент времени между 5 и~6 часами утра и, прождав соперника 5 минут, удаляется. В~случае же прибытия последнего в~эти 5 минут дуэль состоится. Какая часть дуэлей действительно заканчивается поединком?

(b) Стержень случайным образом ломают на три части. С~какой
вероятностью из этих частей можно составить треугольник?

(c) Найдите вероятность того, что случайный треугольник является
остроугольным.
\end{pr}

В задаче \ref{probas-pha}\,(b) за $M$ можно принять равносторонний
треугольник с~высотой, равной длине стержня. Так как для каждой
точки внутри треугольника сумма расстояний от неё до сторон равна
высоте, эти расстояния можно считать равными длинам получившихся
при разломе частей стержня. Парадоксально, что у~этой задачи (и~у~других, например, \ref{probas-pha}\,(с)) имеются другие естественные формализации, дающие другой ответ!

\sseccol{Указания, ответы и~решения}

\paragraph*{\ref{strelki}.}
(a) \emph{Ответ}: $1-(1-0{,}7)(1-0{,}8)=0{,}94$.

\subsection{Независимость и~условная вероятность (1)}

%%\Opensolutionfile{_hintAholder}
%%\Opensolutionfile{_hintBholder}

Следующее определение обобщает ситуацию правила умножения
\ref{probas-rule}\,(c). Подмножества (т.\,е. события) $A$ и~$B\ne\emptyset$ конечного множества $M$ \emph{независимы}, если
доля (т.\,е. вероятность) множества $A\cap B$ в~$B$ равна доле
(т.\,е. вероятности) множества $A$ в~$M$. Приведём симметричную
переформулировку, которая работает и~для $B=\emptyset$.
Подмножества $A$ и~$B$ конечного множества $M$ называются
\emph{независимыми}, если\index{Независимость!}
$$
|A\cap B|\cdot|M| =|A|\cdot|B|.
$$
Основной пример независимых подмножеств "--- в~множестве всех клеток шахматной доски подмножество клеток в~первых трёх её строках и~подмножество клеток в~последних четырёх её столбцах, или, более строго, $A\times N$ и~$M\times B$ в~$M\times N$.

\begin{pr}\label{prolov-indep}
Зависимы ли следующие подмножества? (Мы называем \emph{зависимыми}
подмножества, не являющиеся независимыми.)

(a) Подмножества $\{1,2\}\subset\{1,2,3,4\}$ и~$\{1,3\}\subset\{1,2,3,4\}$.

(b) Подмножества $\{1,2\}\subset\{1,2,3,4,5,6\}$ и~$\{1,3\}\subset\{1,2,3,4,5,6\}$.
\end{pr}

\begin{pr}\label{prolov-indep105}
Зависимы ли следующие подмножества множества целых чисел от 1 до
105?

(a) Подмножество чисел, делящихся на 5, и~подмножество чисел,
делящихся на 7.

(b) Подмножество чисел, делящихся на 15, и~подмножество чисел,
делящихся на 21.

(c) Подмножество чисел, делящихся на 15, и~подмножество чисел,
делящихся на 5.

(d) Подмножество чисел, делящихся на 10, и~подмножество чисел,
делящихся на 7.
\end{pr}

\begin{pr}\label{prolov-com}
Подмножества $A$ и~$B$ конечного множества независимы тогда и~только
тогда, когда $B$ и~$A$ независимы.
\end{pr}

\begin{pr}\label{procon-div}
Два дворянина из свиты короля в~ожидании выхода его Величества
решили сыграть в~кости. Они сделали одинаковые ставки и~договорились, что тот, кто первым выиграет 10 партий, получает все
деньги. При счёте 9:8 появился король и~игру пришлось закончить. Как
следует поделить деньги?
\end{pr}

Это одна из задач, положивших начало теории вероятностей. (Решить её
вам будет проще после задачи \ref{proind-bay}.) В~XVII в. её предложил великому французскому математику Блезу Паскалю его знакомый "--- один из тех дворян, о~которых говорится в~задаче. Паскаль понял, что следует поделить деньги пропорционально шансам, которые имели игроки на окончательную победу в~момент остановки игры. Он нашёл способ вычисления этих шансов (для любого счёта). Другой метод решения задачи, приводящий к~тому же результату, нашёл другой великий математик XVII~в. Пьер Ферма. Их методы основаны на следующем понятии.

\emph{Условной вероятностью} подмножества $A$ при условии подмножества $B$, для которого $P(B)\ne0$, называется отношение\index{Вероятность!условная|textbf}
$$
P(A|B)=P(A \cap B)/P(B).
$$
Ясно, что независимость подмножеств $A$ и~$B$ равносильна тому, что $P(A|B)=P(A)$.

\begin{pr}\label{deti}
(a) Известно, что при броске игральной кости выпало чётное число.
Найдите вероятность того, что оно меньше 5.

(b) В~семье два ребёнка. Известно, что один из них мальчик. Найдите
вероятность того, что второй ребёнок тоже мальчик. (Мы предполагаем,
что вероятности рождения мальчика и~девочки равны половине и~что
пол второго ребёнка не зависит от пола первого.)
\end{pr}

\begin{pr}\label{proind-lamp}
Лампочки выпускаются двумя заводами, причём первый из них производит
$70\,\%$ всей продукции. Лампочки, произведённые первым заводом,
горят с~вероятностью $0{,}98$, вторым "--- $0{,}95$. Найдите
вероятность того, что купленная лампочка горит.
\end{pr}

Решение этой задачи обобщает следующий факт.

\begin{pr}\label{proind-full}
\textbf{Формула полной вероятности.} Если $M=B_1\sqcup\ldots\sqcup
B_n$ и~$P(B_j)\ne0$ (говорят, что $B_1,\ldots,B_n$ "--- \emph{полная система событий}), то \index{Формула!полной вероятности}
$$
P(A)=P(A|B_1)P(B_1)+ \ldots +P(A|B_n)P(B_n).
$$
\end{pr}

\begin{pr}\label{probas}
Победитель в~поединке двух боксёров определяется большинством голосов трёх судей. Двое судей выносят верное решение с~вероятностью~$p$, а~третий голосует, бросая монету. Найдите вероятность принятия
судьями верного решения.
\end{pr}

\begin{pr}
Отец, мать и~сын увлекаются шахматами. Отец обещает сыну приз, если
он выиграет две партии подряд из трёх, сыгранных поочерёдно с~отцом
и матерью. Сын знает, что отец играет лучше матери. С~кем ему
выгоднее играть первую партию?
\end{pr}

\begin{pr}\zvezda\label{procon-sum}
Правила распространённой в~ряде стран игры следующие: игрок
бросает две кости. Он выигрывает, если сумма выпавших очков равна 7
или 11, и~проигрывает, если она равна 2, 3 или 12. Во всех остальных
случаях он бросает кости до тех пор, пока не выиграет, выбросив
первоначальную сумму, или не проиграет, выбросив 7. Найти
вероятность выигрыша.
\end{pr}

\begin{pr}\label{proind-lamp2}
Лампочки выпускаются двумя заводами, причём первый из них производит
$70\,\%$ всей продукции. Лампочки, произведённые первым заводом,
горят с~вероятностью $0{,}98$, вторым "--- $0{,}95$. Купленная лампочка оказалась бракованной. Найдите вероятность того, что она
выпущена первым заводом.
\end{pr}

Решение этой задачи обобщает следующий факт.

\begin{pr}\label{proind-bay} \textbf{Формула Байеса.} Справедливо равенство $$P(B|A)=P(A|B)P(B)/P(A).$$
\index{Формула!Байеса}
\end{pr}

\vspace{-0.5cm}
Часто применяется следствие формул \ref{proind-full} и~\ref{proind-bay}:
$$
P(X|A)=\frac{P(A|X)P(X)}{P(A|B_1)P(B_1)+ \ldots +P(A|B_n)P(B_n)}.
$$

\begin{pr}\label{proind-test}
Вероятность того, что изделие бракованное, равна $0{,}04$. Если
изделие бракованное, то оно пройдёт тест с~вероятностью $0{,}05$, а~иначе "--- с~вероятностью $0{,}98$. Найдите (с~точностью до
$0{,}0001$) вероятность того, что изделие, дважды выдержавшее тест, бракованное.
\end{pr}

\newpage
\setcounter{page}{591}
\subsection{Случайные величины (3)}\label{ss:vel}

%\Opensolutionfile{_hintAholder}
%\Opensolutionfile{_hintBholder}

Пусть дано конечное или счётное множество $M$ и~для каждого
элемента $m\in M$ задано число (вероятность) $P(m)\geq 0$,
$\sum\limits_{m\in M}P(m)=1$. Числовая функция $X$, заданная на $M$,
называется \emph{случайной величиной}. Множество пар $(x_i,p_i),
i=1,2,\ldots$, где $\{x_1,x_2,\ldots\}$ "--- множество возможных
значений случайной величины $X$, а~$p_i=P(\{m\in M\colon X(m)=x_i\})$, $i=1,2,\ldots,$ "--- соответствующие им вероятности, называется \emph{распределением} случайной величины $X$.
\index{Случайная!величина|textbf}\index{Распределение|textbf}

\emph{Комментарий}. Как правило, при изучении случайной величины
$X$ не требуется знать, на каком множестве она определена. Достаточно знать только её распределение.

Событие $\{m\in M\colon X(m)=x_i\}$ в~дальнейшем сокращённо обозначается $X=x_i$.

\begin{pr}\label{prt1-4-1}
Монета подбрасывается 5 раз. Найдите распределение числа выпавших орлов.

\smallskip
\emph{Математическим ожиданием} или \emph{средним значением}
случайной величины $X$ называется сумма
$$
E(X)=\sum x_ip_i =x_1P(X=x_1)+x_2P(X=x_2)+\ldots.
$$\index{Математическое!ожидание|textbf}\index{Среднее!значение|textbf}

\emph{Комментарий}. Если множество значений случайной величины
бесконечно, то это определение нуждается в~уточнении. Сумма ряда в~правой части называется \emph{математическим ожиданием}, только когда этот ряд сходится абсолютно. В~противном случае говорят, что у~величины $X$ \emph{не существует математического ожидания}. Например, пусть случайная величина $X$ принимает значение $n\in \N$ с~вероятностью $p_n=\frac1{n(n+1)}$. Тогда ряд $\sum
np_n=\sum\frac1{n+1}$ расходится, т.\,е. $E(X)$ не существует. В~дальнейшем мы предполагаем, что для всех рассматриваемых случайных
величин математические ожидания существуют, т.\,е. ряд $\sum
x_iP(X=x_i)$ сходится абсолютно.
\end{pr}

\begin{pr}\label{prt1-4-2}
a) Докажите, что математическое ожидание случайной величины $X$,
заданной на множестве $M$, равно $\sum\limits_{m\in M} X(m)P(m)$.

b) Докажите, что если $ E(X) \le x $, то существует $m\in M\colon
X(m) \le x $.

c) Пусть случайная величина $X$ при всех $m\in M$ принимает
одно и~то же значение $\mu$: $X(m)=\mu$. Найдите $E(X)$.

d) Докажите \emph{линейность} математического ожидания:
$E(aX+bY)=aE(X)+bE(Y)$, где $a$, $b$ "--- вещественные числа, а~$X$, $Y$ "--- случайные величины.
\end{pr}

Случайные величины $X$ и~$Y$ называются \emph{независимыми}, если
события $X=x_i$ и~$Y=y_j$ независимы при любых $x_i$, $y_j$, т.\,е.
$$P(\{m\in M\colon X(m)=x_i\text{\ и\ }Y(m)=y_j\})=P(X=x_i)P(Y=y_i).$$
\index{Случайные величины!независимые|textbf}

Неформально независимость означает, что значения одной из случайных
величин не влияют на распределение другой.

\begin{pr}\label{prt1-4-3}
Докажите, что если случайные величины $X$ и~$Y$ независимы, то
математическое ожидание их произведения равно произведению их математических ожиданий: $E(XY)=E(X)E(Y)$.
\end{pr}

\emph{Дисперсией} случайной величины $X$ называется число
$D(X)=E\big((X-E(X))^2\big)$. \index{Дисперсия|textbf}

\emph{Комментарий}. Если множество значений случайной величины
бесконечно, то дисперсия может не существовать. В~дальнейшем
предполагается, что для всех рассматриваемых случайных величин
дисперсия существует.

\begin{pr}\label{prt1-4-4}
Докажите, что $D(X)=E(X^2)-E(X)^2$.
\end{pr}

\begin{pr}\label{prt1-4-5}
Докажите, что если $X$ и~$Y$ независимы, то $D(X+Y)=D(X)+D(Y)$.
\end{pr}

\begin{pr}\label{prt1-4-6}
\textbf{Неравенство Чебышёва.} Докажите, что для любой случайной
величины $X$ и~любого $\varepsilon>0$ выполняется неравенство
$$
P(|X-E(X)| \geq \varepsilon) \leq D(X)/\varepsilon^2.
$$\index{Неравенство!Чебышёва}
\end{pr}

\begin{pr}\label{prt1-4-7}
Федя знает ответы на 20 из 30 вопросов. В~билет входят 3 вопроса.
Найдите распределение числа вопросов, на которые Федя сможет ответить.
\end{pr}

\begin{pr}\label{prt1-4-8}
Две одинаковые колоды карт перетасовываются, и~карты последовательно
парами выкладываются на стол. Найдите среднее значение числа пар,
карты в~которых совпадают.
\end{pr}

\begin{pr}\label{prt1-4-9}
В городе $N$ предприниматели обязаны предоставлять всем рабочим
выходной, если хотя бы у~одного из них день рождения. Остальные дни
являются рабочими. Сколько человек следует принять на работу, чтобы
среднее значение числа рабочих человеко-дней было максимальным?
\end{pr}

\begin{pr}\label{prt1-4-10}
В задаче \ref{prt1-4-7} найдите среднее значение Фединой оценки (если Федя ответит на 3 вопроса, он получит 5, на 2 "--- 4 и~т.\,д.).
\end{pr}

\begin{pr}\label{prt1-4-11}
В распространённой азартной игре игрок может делать ставку на один
из номеров от 1 до 6. Бросаются 3 кости, и~если выбранный номер
выпал хотя бы на одной, то игрок получает свою ставку плюс столько
же за каждое появление выбранного номера. Выгодна ли игра для
игрока?
\end{pr}

\begin{pr}\label{prt1-4-12}
(Загадка.) Площадка имеет форму квадрата со стороной $350$~м. При
измерении стороны вероятность ошибки $\pm10$~м равна $0{,}16$,
$\pm20$~м "--- $0{,}08$, $\pm30$~м "--- $0{,}05$. Найдите среднее
значение измеренной площади.

\emph{Комментарий.} На самом деле ответ на этот вопрос зависит от
того, как формализовано понятие измерения площади. Если независимо
измерить каждую из сторон квадрата и~перемножить полученные
значения, то по задаче~\ref{prt1-4-3} среднее значение будет равно $350^2$~м$^2$. Если же измерить только одну сторону и~возвести результат в~квадрат, то ответ будет другим.
\end{pr}

\begin{pr}\label{prt1-4-13}
В ряд в~случайном порядке выписаны $m$ единиц и~$n$ нулей. Найдите
среднее число серий из $k$ одинаковых цифр подряд.
\end{pr}

\begin{pr}\label{prt1-4-14}
Из колоды в~52 карты вынимаются карты до первого туза. Сколько карт
в среднем будет вынуто?
\end{pr}

\begin{pr}\label{prt1-4-15}
По узкой дороге в~одном направлении едут $n$ машин. Вначале скорости
всех машин различны. Каждая машина едет с~постоянной скоростью, пока
не догонит едущую впереди, после чего едет со скоростью передней
машины. В~результате через достаточно большое время машины
разбиваются на несколько групп. Найдите среднее значение числа групп.
\end{pr}

%\Closesolutionfile{_hintAholder}
%\Closesolutionfile{_hintBholder}

\sseccol{Указания, ответы и~решения}

%\Readsolutionfile{_hintBholder}

\paragraph*{\ref{prt1-4-8}.}
\emph{Ответ}: 1.

\paragraph*{\ref{prt1-4-9}.}
\emph{Ответ}: 364 или 365.

Вероятность того, что в~данный день ни у~одного из $n$ рабочих не
будет дня рождения, равна $(364/365)^n$. Следовательно, среднее число рабочих человеко-дней равно $365n(364/365)^n$. Это выражение
достигает максимального значения при $n$, равном 364 или 365
%(см.~задачу 34 из книги~\cite{sb0}).

\newpage
\setcounter{page}{595}
\subsection{Испытания Бернулли (3)}

%\Opensolutionfile{_hintAholder}
%\Opensolutionfile{_hintBholder}

\index{Испытания Бернулли|textbf} \emph{Испытаниями Бернулли}
называется последовательность $n$ независимых случайных величин,
каждая из которых принимает два значения: $1$ с~вероятностью $p$ и~$0$ с~вероятностью $q=1-p$. Обычно появление 1 называют
\emph{успехом}, а~$0$ "--- \emph{неудачей}. \index{Успех|textbf}
\index{Неудача|textbf}

Приведём другое определение испытаний Бернулли. Пусть $M$ "---
множество $n$"~мерных векторов, все координаты которых равны~$0$
или~$1$, и~для каждого $x\in M$ задана вероятность $P(x)=\Pi_{i=1}^n p_i$, где $p_i=p$, если $x_i=1$, и~$p_i=q=1-p$, если $x_i=0$. Элементы множества $M$ тоже назовём \emph{испытаниями Бернулли}.

\emph{Комментарий}. Оба определения являются эквивалентными в~следующем смысле. Очевидно, что определённые на множестве $M$
случайные величины $x_i$ независимы и~каждая из них принимает
значение 1 с~вероятностью $p$ и~0 с~вероятностью $q$. Поэтому каждый
вектор $x\in M$ можно рассматривать как набор значений $n$
независимых случайных величин $x_i$.

Случайная величина $X=\sum x_i$ называется \emph{числом успехов}.
\index{Число!успехов}

\begin{pr}\label{prt1-5-1}
В $n$ испытаниях Бернулли с~вероятностью успеха $p$ найдите

а) вероятность ровно $k$ успехов;

б) среднее значение числа успехов;

в) дисперсию числа успехов;

г) наиболее вероятное значение числа успехов.
\end{pr}

\begin{pr}\label{prt1-5-2}
\textbf{Закон больших чисел.} Пусть $X$ "--- число успехов в~$n$
испытаниях Бернулли с~вероятностью успеха $p$. Пусть $t>0$.
Докажите, что
$$
P\bigg(\Big|\frac{X}{n}-p\Big|\geq t \sqrt{\frac{pq}{n}}\bigg) \leq
\frac1{t^2}
$$\index{Закон!больших чисел}

\emph{Указание.} Примените неравенство Чебышёва.
\end{pr}

Закон больших чисел означает, что при большом числе испытаний
вероятность того, что частота успеха сильно отличается от его
вероятности, мала. На самом деле этот закон справедлив не только для
испытаний Бернулли: если наблюдать много независимых реализаций
произвольной случайной величины, то их среднее с~большой
вероятностью будет мало отличаться от её математического ожидания.
Этот закон позволяет, например, проводить социологические
исследования, в~которых на основе опроса некоторого количества
случайно выбранных людей (достаточно большого, но составляющего
малую часть всего населения) делаются выводы о~распространённости в~обществе тех или иных мнений и~предпочтений.

\begin{pr}\label{prt1-5-3}
Пассажиру купейного вагона удобно, если все его попутчики одного с~ним пола. Какая часть пассажиров испытывает удобства?
\end{pr}

\begin{pr}\label{prt1-5-4}
Вероятность рождения мальчика равна $0{,}515$. Найдите вероятность того, что среди 6 детей не более 2 девочек.
\end{pr}

\begin{pr}\label{prt1-5-5}
Кооператив отгружает железные балки. Средняя длина балки $3$~м, дисперсия $0{,}09$~м$^2$. Сколько балок надо заказать, чтобы с~вероятностью, не меньшей чем $0{,}999$, хотя бы 1000 из них имели
длину не менее $2$~м?
\end{pr}

\begin{pr}\label{prt1-5-6}
Найдите среднее число испытаний до первого успеха, если вероятность
успеха равна~$p$.
\end{pr}

\begin{pr}\label{prt1-5-7}
Проводятся независимые испытания с~вероятностью успеха $0{,}8$.
Испытания проводятся до первого успеха, но не более четырёх раз. Найдите среднее число испытаний.
\end{pr}

\begin{pr}\label{prt1-5-8}
(Загадка.) Старик ловил неводом рыбу ровно тридцать лет и~три года.
Каждый день он ловил ровно 7 рыб, которых как раз хватало на ужин.
Живущий у~старухи кот-долгожитель ест только макрель, которая
ловится вдвое реже остальных рыб. В~результате он 700 раз оставался
голодным. Плавает ли макрель в~море косяками или поодиночке?

\emph{Комментарий.} Конечно, точно ответить на поставленный вопрос
невозможно. Однако можно оценить, какая из двух гипотез лучше
согласуется с~данными.
\end{pr}

%\Closesolutionfile{_hintAholder}
%\Closesolutionfile{_hintBholder}

\sseccol{Указания, ответы и~решения}

%\Readsolutionfile{_hintBholder}

\paragraph*{\ref{prt1-5-1}.}
\emph{Ответы}: а) $\smash[b]{\binom{n}{k}}p^kq^{n-k}$;

\comment

б) $np$;

в) $npq$;

г) $\lfloor np\rfloor$, если $\{np\}\leq q$, и~$\lfloor np\rfloor+1$, если $\{np\}\geq q$
(при равенстве соответствующие вероятности совпадают).
\emph{Указание}. Найдите отношение вероятностей ровно $k$ успехов
и ровно $k+1$ успехов.

\endcomment

%kommentarii vernuty, ne ubirat!

\newpage
\setcounter{page}{620}
\section{Перестановки. \emph{А.~Б.~Скопенков}}\label{0gro}

%\Opensolutionfile{_hintAholder}
%%\Opensolutionfile{_hintBholder}

Задачи этого параграфа не требуют для решения предварительных знаний.
%AS peremestil frazu
Они естественным образом подводят читателя к~понятию \emph{группы преобразований}, которое явно вводится в~\S\,\ref{0gro-cyc}. Миникурс <<Рождение понятия группы>> можно составить из этого параграфа, \S\,\ref{s:feresi} <<Умножение по простому модулю>>, \S\,\ref{preobr} <<Геометрические преобразования>>, \S\,\ref{s:rad} <<Разрешимость в~радикалах>>, \S\,\ref{0gro-cyc} <<Группы>> и~статьи \cite{Sk15}.

% Например, лемма Бернсайда (п. \ref{0grobur}) и теорема о порядках циклических групп (п. \ref{0gro-cyc}).
%<<Когда любая группа из $N$ элементов циклическая?>>
%(там приведено простое доказательство {\it основного} результата книги [Al]),

\subsection{Порядок, тип, сопряжённость (1)}\label{0groper}

\begin{pr}\label{dobryj}
Пятнадцать школьников сидят на пятнадцати пронумерованных стульях.
Каждую минуту добрый преподаватель пересаживает их по следующей схеме:
$$
{\setcounter{MaxMatrixCols}{15}
\begin{pmatrix}
1 & 2 & 3 & 4 & 5 & 6 & 7 & 8 & 9 & 10 & 11 & 12 & 13 & 14 & 15\\
3 & 5 & 10 & 8 & 11 & 14 & 15 & 6 & 13 & 1 & 4 & 9 & 7 & 2 & 12
\end{pmatrix}}.
$$
Через сколько минут все школьники впервые окажутся на своих первоначальных местах?
\end{pr}

\emph{Перестановка} множества "--- запись элементов этого множества в~некотором порядке. Если говорить более строго, \emph{перестановкой} множества называется взаимно однозначное отображение этого множества на себя (т.\,е. биекция).\index{Перестановка|textbf}
(Перестановку $f$ удобно изображать в~виде {\it ориентированного графа}, вершины которого "--- элементы множества, а~рёбра идут из вершины $a_k$ в~вершину~$f(a_k)$).
Перестановка множества $\{a_1,a_2,\ldots,a_n\}$, переводящая $a_k$ в~$f(a_k)$, записывается в~виде
$$
\begin{pmatrix}
a_1 & a_2 & \ldots &a_n\\
f(a_1) & f(a_2) & \ldots & f(a_n)
\end{pmatrix};
$$
обычно $a_k=k$ для всех $k=1,\ldots,n$.

\emph{Обратной} к~$f$ перестановкой называется перестановка $f^{-1}$, определённая формулой $f(f^{-1}(x))=x$.
Она записывается в~виде
$$
\begin{pmatrix}
f(a_1) & f(a_2) & \ldots & f(a_n)\\
a_1 & a_2 & \ldots &a_n
\end{pmatrix}.
$$\index{Перестановка!обратная|textbf}

\emph{Композицией} перестановок $f$ и~$g$ называется перестановка $f\circ g$, определённая формулой $(f\circ g)(x):=f(g(x))$.
\index{Композиция!перестановок|textbf}

\begin{pr}\label{g-comp0}
Найдите композиции

(a) $\begin{pmatrix}
 1 & 2 & 3\\ 2 & 1 & 3
 \end{pmatrix}
 \circ
 \begin{pmatrix}
 1 & 2 & 3\\ 3 & 1 & 2
 \end{pmatrix}$;\qquad
(b) $\begin{pmatrix}
 1 & 2 & 3\\ 2 & 3 & 1
 \end{pmatrix}
 \circ
 \begin{pmatrix}
 1 & 2 & 3\\ 3 & 1 & 2
 \end{pmatrix}$.
\end{pr}

\emph{Циклом} $(a_1,a_2,\ldots,a_n)$ называется перестановка
$$
\begin{pmatrix}
a_1 & a_2 & \ldots & a_{n-1} & a_n\\
a_2 & a_3 & \ldots & a_n &a_1
\end{pmatrix}
$$
множества, содержащего элементы $a_1,a_2,\ldots,a_n$, которая переводит $a_n$ в~$a_1$ и~$a_i$ в~$a_{i+1}$ для любого $i<n$, а~каждый из остальных элементов переводит в~себя.
\index{Цикл|textbf}

На этом языке результаты задачи \ref{g-comp0} можно коротко выразить так: $(12)\circ (13)=(132)$ и~$(123)\circ (132)=(1)$.

\begin{pr}\label{g-comp}
Найдите композиции (перестановок на множестве цифр)

(a) $(12)\circ (23)$; \quad (b) $(23)\circ (12)$; \quad
(c) $(12)\circ (13)\circ (12)$; \quad

(d) $(12345)\circ (12)$; \quad (e) $(12345)\circ(56789)$.

Ответ дайте в~виде композиции непересекающихся циклов.
Например, $(123)\circ(234)=(12)\circ(34)$.
\end{pr}

Далее знак композиции опускается.

\begin{pr}\label{orderexi}
Для любой перестановки $f$ существует $n>0$, для которого $f^n=\id$
(т.\,е. после $n$"~кратного применения перестановки $f$ каждый элемент перейдёт в~себя).
\end{pr}

\emph{Порядком} $\ord f$ перестановки $f$ называется наименьшее целое положительное число~$n$, для которого $f^n=\id$.
\index{Порядок!перестановки|textbf}
% (если такое $n$ существует).

\begin{pr}\label{orderex}
Существуют ли перестановки 9-элементного множества порядков 7; 10; 12; 11?
\end{pr}

\begin{pr}\label{order}
Чему равен порядок композиции непересекающихся циклов из $n_1,\ldots,n_k$ элементов соответственно?
\end{pr}

\begin{figure}[ht]\centering
\includegraphics{skopenkov-005.mps}
\caption{Перестановка типа $\langle1,2,3,4\rangle$}
%\label{fig3-2}
\end{figure}

Перестановки $(n_1+\ldots+n_k)$"~элементного множества из задачи \ref{order} называются перестановками \emph{типа} $\langle n_1,\ldots,n_k\rangle$.
Например, перестановки $(14)(253)$, $(15)(432)$ типа $\langle2,3\rangle$, а~перестановка $(1)(3)(245)$ "--- другого типа $\langle1,1,3\rangle$.
\index{Тип перестановки|textbf}

\begin{pr}\label{numper}
 Найдите число перестановок типа

(a) $\langle2,3\rangle$; \quad (b) $\langle3,3\rangle$; \quad (c) $\langle1,2,3,4\rangle$.
\end{pr}

Перестановки $a$ и~$b$ называются \emph{сопряжёнными}, если $a=xbx^{-1}$ для некоторой перестановки $x$.
\index{Сопряжённые перестановки|textbf}

\begin{pr}\label{conju}
(a) Перестановки $a$ и~$b$ сопряжены тогда и~только тогда, когда их типы одинаковы.

(b) Пусть $a$ и~$x$ "--- произвольные перестановки $n$"~элементного множества. Тогда
$$
 xax^{-1}=\begin{pmatrix}
 x(1) & x(2) & \ldots & x(n)\\
 x(a(1))& x(a(2)) & \ldots & x(a(n))
\end{pmatrix}.
$$
Иными словами, циклическое разложение перестановки $xax^{-1}$ получается из циклического разложения перестановки $a$ заменой каждого элемента на его $x$-образ: если
$a=\prod\limits_{j=1}^q(i_{j,1},i_{j,2},\ldots,i_{j,s_j})$, то
 $$
 xax^{-1}=\prod\limits_{j=1}^q(x(i_{j,1}),x(i_{j,2}),\ldots,x(i_{j,s_j})).
 $$

(c) Найдите $gf^{-1}g^{-1}f$ для $f:=(1,2,\ldots,N)$ и~$g := (N,N+1,\ldots,L).$

(d) Вращения куба вокруг больших диагоналей порождают сопряжённые перестановки множества его вершин.
\end{pr}

\begin{pr}\label{decom}
Любая перестановка представляется в~виде композиции

(a) непересекающихся циклов;

(b) \emph{транспозиций}, т.\,е. перестановок, каждая из которых меняет местами
некоторые два элемента, а~остальные оставляет на месте (иными словами, циклов длины 2);
\index{Транспозиции|textbf}

(c) транспозиций $(1i)$, $i=2,3,\ldots,n$.
\end{pr}

\begin{pr}\label{twocom}
Найдите \emph{две} перестановки, композициями которых можно получить любую перестановку $n$"~элементного множества.
\end{pr}

%{\bf Зач\"етные задачи для 8 класса:} \ref{comp}.def, \ref{orderex}.12, \ref{order}, \ref{numper}.ab,
%\ref{decom}.ab, \ref{twocom}.a, \ref{conju}.a.
%Письменно: \ref{numper}.b или \ref{decom}.a по выбору школьника.

%{\bf Зач\"етные задачи для 9 класса:} \ref{comp}.efg, \ref{orderex}.11, \ref{order}, \ref{numper}.ac,
%\ref{decom}.bc, \ref{twocom}.b, \ref{conju}.b.
%Письменно: \ref{numper}.c или \ref{decom}.a по выбору школьника.

%{\it Для школьников, готовящих доклады, письменные задачи не обязательны.}

\newpage
\setcounter{page}{626}
\refstepcounter{subsection}
\subsection{Комбинаторика классов эквивалентности (2)}\label{0grobur}

%\Opensolutionfile{_hintAholder}
%%\Opensolutionfile{_hintBholder}

Этот пункт посвящён подсчёту числа классов эквивалентости (т.\,е. раскрасок и~т.\,д.). Такой подсчёт подводит читателя к~важному понятию \emph{группы преобразований} и~к~элементарной формулировке \emph{леммы Бёрнсайда}. Формулировка и~доказательство этого и~других результатов на языке абстрактной теории групп делает их менее доступными. Ср. \S\,28.%\ref{s:motiv1}.

%Третий раздел не зависит от первых двух.
%Первая подборка посвящена простейшим свойствам перестановок.

%Изложение в этом разделе улучшено по сравнению с
%[S4, Z, глава 10, раздел <<Комбинаторика классов эквивалентности>>].

Не требуется, чтобы в раскраске присутствовали все данные цвета. Раскраски, совмещающиеся вращением пространства (т.\,е. движением пространства, сохраняющим ориентацию и~имеющим неподвижную точку), считаются одинаковыми (кроме задачи \ref{cleq1}\,(с)).

Следующие определения используются только в~задачах \ref{cleq1}.(b), \ref{per-raskr}.(e), 23.3.11%\ref{bernsgraph}
 (и потому могут быть пропущены при решения остальных задач).

\emph{Изоморфизм} между графами"--- такая биекция между
множествами их вершин, что для любых двух вершин эти вершины соединены ребром тогда и~только тогда, когда их образы при биекции соединены ребром.
\emph{Автоморфизм} графа"--- его изоморфизм на себя.

%%%!!!перестановка

\begin{pr}\label{cleq1}
Сколько существует

(a) раскрасок граней куба в~красный и~серый цвета;

(b) различных (т.\,е. неизоморфных) неориен\-ти\-ро\-ванных графов с~4 вершинами;

(c) раскрасок в~$r$ цветов вершин правильного тетраэдра?

%Вернул, как было, чтобы не было вложенных скобок

\noindent Здесь раскраски, совмещающиеся движением пространства (не обязательно сохраняющим ориентацию), считаются одинаковыми.

%%%!!!
%(d) раскрасок вершин полного графа с~4 вершинами в~$r$ цветов\blue{?}
%\noindent \blue{Здесь} раскраски, совмещающиеся перестановкой вершин этого графа, считаются одинаковыми.

\end{pr}

% Найдите количество закрасок в $r$ цветов
% (a) правильного треугольника, разбитого средними линиями на 4 равных треугольника.
% (b) то же для 9 треугольников.
% (c) квадрата, разбитого средними линиями на 4 равных квадрата.
% (d) то же для 9 равных квадратов.

\begin{pr}\label{cleq2}
Для простого $p$ найдите количество замкнутых ориентированных связных $p$"~звенных ломаных (возможно, самопересекающихся), проходящих через все вершины данного правильного $p$"~угольника.

Здесь ломаные, совмещающиеся поворотом, неотличимы.
\end{pr}

Задачи \ref{cleq1} и~\ref{cleq2} простые, их можно решить без идей, приводящих к~лемме Бёрнсайда.

\begin{pr}\label{cleq3}
Найдите количество раскрасок карусели из $n$ незанумерованных вагончиков в~$r$ цветов (т.\,е. количество раскрасок вершин правильного $n$"~угольника в~$r$ цветов, если раскраски, совмещающиеся поворотом, неотличимы) для

(a) $n=5$; \quad (b) $n=4$; \quad (c) $n=6$.
\end{pr}

Задачу \ref{cleq3} для произвольного $n$ можно решить способом, аналогичным придуманному вами для малых~$n$. Однако решение будет громоздким. Приведём более простой (для <<очень непростых>> $n$) способ на примере решения задачи \ref{cleq3}\,(с).

Назовём (\emph{раскрашенным}) \emph{поездом} раскраску карусели из \emph{занумерованных} вагончиков в~$r$ цветов. Тогда всего имеется $r^6$ поездов из 6 вагончиков.

Распределим поезда по вокзалам так, чтобы на каждом вокзале находились все поезда, полученные из некоторой одной раскраски карусели всевозможными разрубаниями, т.\,е. искомое количество $Z$ раскрасок равно количеству вокзалов.

%На каждый вокзал пригоним все поезда, получающиеся разрубанием в некотором месте из той раскраски карусели, которая соответствует этому вокзалу.

Назовем \emph{периодом} $T(\alpha)$ поезда $\alpha$ наименьшую положительную величину циклического сдвига, переводящего поезд $\alpha$ в~себя.

\begin{pr}\label{burn-d}
Количество поездов на вокзале равно периоду каждого из поездов, стоящих на этом вокзале. В~частности, периоды поездов, стоящих на одном вокзале, равны.
\end{pr}

На каждом вокзале выберем один поезд. Посадим в~него 6 пассажиров и~выдадим им билеты с~числами 0, 1, 2, 3, 4, 5. Тогда нужно найти общее число $6Z$ пассажиров.

По команде каждый пассажир переходит в (раскрашенный) поезд, полученный из выбранного поезда циклическим сдвигом на число, указанное в~билете пассажира. Ясно, что каждый пассажир остается на прежнем вокзале.

%Число пассажиров, оставшихся в выбранном поезде, равно количеству циклическим сдвигов, переводящих этот поезд в себя.

\begin{pr}\label{burn-d}
(a) В выбранном поезде $\alpha$ останется $6/T(\alpha)$ пассажиров.
Более формально, количество тех $s\in\{0,1,2,3,4,5\}$, для которых циклический сдвиг на $s$ переводит поезд $\alpha$ в себя,
равно $6/T(\alpha)$.

(b) В каждом поезде $\alpha$ окажется $6/T(\alpha)$ пассажиров.
\end{pr}

Значит, общее число $6Z$ пассажиров равно количеству всех пар $(\alpha,s)$, в~которых $s\in\{0,1,2,3,4,5\}$ и $\alpha$ --- поезд,
переходящий в~себя при циклическом сдвиге на $s$ вагончиков. Циклический сдвиг на $s$ переводит в себя ровно $r^{\gcd(s,6)}$ поездов. Поэтому
 $$
 6Z= r^6+r+r^2+r^3+r^2+r.
 $$
Приведенный план решения можно представить в виде формулы
 $$
 6Z=\sum\limits_x T(x)\cdot\frac6{T(x)} =
 \sum\limits_\alpha \frac6{T(\alpha)}=r^6+r+r^2+r^3+r^2+r.
 $$
Здесь первое суммирование происходит по всем по всем раскраскам $x$ каруселей, а второе "--- по всем поездам~$\alpha$.

\begin{figure}[ht]\centering
\includegraphics{skopenkov-002.mps}
\caption{Граф $K_{3,3}$}
\label{f:k33}
\end{figure}

\begin{pr}\label{per-raskr}
Найдите количество

(a) раскрасок карусели из $n$ вагончиков в~$r$ цветов (см.~формализацию и~другое решение в~\cite[\S\,1.5]{GDI});

(b) $r$"~цветных ожерелий из $n=2k+1$ бусин (ожерелья считаются одинаковыми, если они совмещаются либо поворотом вокруг центра ожерелья, либо осевой симметрией ожерелья);

(c) раскрасок незанумерованных граней куба в~$r$ цветов;

(d) раскрасок незанумерованных вершин куба в~$r$ цветов;

(e) раскрасок незанумерованных вершин графа $K_{3,3}$ (рис.~\ref{f:k33}) в~$r$ цветов (раскраски считаются одинаковыми, если они совмещаются автоморфизмом этого графа).

%%%В этом графе 6 вершин, поделенных на 2 группы по 3 вершины.
%%%Ребро между двумя вершинами проведено тогда и только тогда, когда эти вершины из разных групп.

\end{pr}

\begin{pr}
Перечислите все вращения куба (т.\,е. вращения пространства, переводящие куб в~себя). (Эта задача разбита на шаги в~п.~15.2.2%\ref{s15.2.2}
 <<Самосовмещения>>.)
\end{pr}

Приведем план решения задачи \ref{per-raskr}\,(c). (Пункты (b)--(e) решаются аналогично. Пункт~(b) решается и~без этого указания.)

Назовём {\it (раскрашенной) коробкой} (или {\it замороженной раскраской}) раскраску занумерованых граней куба в~$r$ цветов. Тогда всего имеется $r^6$ коробок.

Распределим коробки по комнатам так, чтобы в~каждой комнате находились все коробки, полученные из некоторой одной коробки всевозможными вращениями, т.\,е. искомое количество $Z$ раскрасок равно количеству комнат.

В каждой комнате выберем одну коробку. Посадим в~нее 24 таракана, соответствующих вращениям куба. Тогда нужно найти общее число тараканов~$24Z$.

По команде каждый таракан переползает в коробку, полученную из выбранной тем вращением, которое соответствует этому таракану. Ясно, что каждый таракан остается в прежней комнате. Число тараканов, оставшихся в выбранной коробке, равно количеству вращений куба, переводящих эту коробку в себя. Обозначим через $\st\alpha$ количество вращений куба, переводящих (раскрашенную) коробку (т.\,е. замороженную раскраску) $\alpha$ в~себя.

\begin{pr}\label{burn-cube} (a) Число тараканов, оказавшихся в~коробке~$\alpha$, равно~$\st\alpha$. Более формально, если существует вращение, переводящее замороженную раскраску $\alpha$ в~замороженную раскраску~$\alpha'$, то количество таких вращений равно~$\st\alpha$.

(b) В любой другой коробке из выбранной комнаты окажется столько же тараканов, сколько в~выбранной коробке. Более формально, для любых двух замороженных раскрасок $\alpha$ и~$\alpha'$, переходящих друг в~друга при некотором вращении, выполняется равенство $\st\alpha=\st\alpha'$. (Эти равные числа обозначаются~$\st x$, где $x$ "--- соответствующая раскраска незанумерованных граней куба.)
\end{pr}

%(Пункты очень важны для рождения понятия группы.)

Поэтому общее число тараканов равно количеству $P$ всех пар $(\alpha,s)$, в~которых $s$ "--- вращение куба и~$\alpha$ "---
коробка, переходящая в себя при вращении~$s$. Поэтому осталось решить следующую задачу.

\begin{pr}\label{burn-fix}
Для каждого вращения куба $s$ найдите количество $\fix s$ коробок (т.\,е. замороженных раскрасок), переходящих в~себя при вращении~$s$.
\end{pr}

Обозначим через $N_x$ количество замороженных раскрасок, отвечающих раскраске~$x$. Тогда для любой раскраски $x$ число $\st x\cdot N_x$ равно количеству вращений куба, т.\,е.~24. Поэтому приведенный план решения можно представить в виде формулы
$$
 24Z=\sum\limits_x\st x\cdot N_x=\sum\limits_\alpha\st\alpha=\sum\limits_s\fix s.
 $$
Здесь первое суммирование происходит по всем раскраскам $x$ незанумерованных граней, второе "--- по всем замороженным раскраскам~$\alpha$, а~третье "--- по всем вращениям куба~$s$.

Как сформулировать общий результат, который можно было применять вместо повторения намеченных решений задач \ref{per-raskr}\,(a),~(c)?

\begin{pr}
{\bf Лемма Бёрнсайда}. Пусть заданы конечное множество $M$ и~семейство $\{g_1,g_2,\ldots,g_n\}$ преобразований этого множества, замкнутое относительно взятия композиции и~взятия обратного элемента. Назовём элементы множества $M$ \emph{эквивалентными}, если один из них можно перевести в~другой одним из данных преобразований. Тогда количество классов эквивалентности равно $\frac 1n\sum\limits_{k=1}^n \fix(g_k)$, где $\fix(g_k)$ "--- количество элементов множества~$M$, которые преобразование $g_k$ переводит в~себя.
\index{Лемма!Бёрнсайда}
\end{pr}

\newpage
\setcounter{page}{632}
\section{Группы. \emph{В.~Брагин}, \emph{А.~Клячко}, \emph{А.~Скопенков}} \label{0gro-cyc}

%\section{%Когда любая группа из $N$ элементов циклическая?
%\chapter

%\footnotetext{
%AS: Изменил для согласованности и ибо неясно, с чего сноска

Обновляемую версию см.~на \url{http://arxiv.org/abs/1108.5406}.
Благодарим Д.~Баранова, М.~Вялого, П.~Кожевникова, К.~Кохася, А.~Сгибнева, Б.~Френкина и~А.~Шеня и~за полезные замечания. Этот текст основан на цикле задач [BKK]%\cite{BKK}.

%}

% ЗАКОММЕНТИРОВАЛ МС СДЕЛАВ СНОСКОЙ к НАЗВАНИЮ
%\footnotetext{}

\subsection{Зачем, для кого и~как устроен этот параграф}\label{s:grocycwhy}

Мы хотели бы привлечь внимание к~теории групп широкого круга людей, интересующихся математикой и~программированием: учителей, руководителей кружков, студентов и~старшеклассников. В~этой теории есть доступные и~интересные им результаты-жемчужины. Формулировки таких результатов кратки и~используют лишь простейшие определения;
доказательства красивы и~похожи на решения сложных олимпиадных задач. Именно с~таких жемчужин полезно начинать изучение теории, на примере их доказательства показывая, как появляются её основные понятия.
К~сожалению, в~б\'ольшей части существующей литературы эти жемчужины погребены под огромным количеством немотивированного материала, что делает их неинтересными и~недоступными.

\begin{th*}{Основной вопрос} Дано семейство $G$ из $n$ перестановок некоторого множества, замкнутое относительно композиции и~взятия обратной перестановки. Для каких $n$ обязательно найдётся такая перестановка $g\in G$, что $G= \{g,g^2,\ldots,g^n\}$\textup{?}
\end{th*}

Этот параграф предназначен для тех, кому понятна и~интересна формулировка этого вопроса, ср. конец п.\;\ref{s:grostat}. На примере исследования этого просто формулируемого вопроса мы покажем, как \emph{появляются} некоторые основные понятия теории групп. Ср. \S\,\ref{phil-met} <<Начинать с~языка или содержания?>>. Мы дадим простое доказательство теоремы, отвечающей на этот вопрос. Оно не претендует на новизну, хотя мы не видели такого доказательства в~литературе. (Ввиду элементарности вопроса выяснить новизну не представляется возможным.)

Этот параграф может быть интересен читателю, не знакомому с~основами  абстрактной теории групп, но изучавшему перестановки и~основы теории чисел, например, по \S\,\ref{s:numb},\ref{s:feresi},\ref{0gro}. В~частности, этот параграф не должен быть единственным и~даже первым шагом в~теорию групп. Он может быть интересен и~читателю, знакомому с~этими основами, ибо ответ на сформулированный вопрос нетривиален. Такому читателю может быть достаточно прочитать~п.\;24.3.%\ref{s:grocycpro}.

%Пункты \ref{s:grohow}  и \ref{s:grocycpro} формально независимы друг от друга.

%Вернул, как было

В п.\;\ref{s:grohow} проиллюстрировано в задачах, как придумать ответ и~доказательство. Хотя придумать их непросто, {\it изложить} их можно коротко. В~п.\;24.3%\ref{s:grocycpro}
 приведено доказательство, формально независимое от п.\;\ref{s:grohow}. Освобождение доказательства от деталей, возникших при его придумывании, но не нужных для него самого, "--- важная часть его проверки.

%%%!!!
%Задачи, приведённые в~п.\;\ref{s:grohow}, иллюстрируют, как придумать ответ и~доказательство.
%Хотя это и~непросто, \emph{изложить} полученные результаты??? можно коротко.

Для понимания доказательства необходим опыт работы с~перестановками и~числами, включая теорему Ферма"--~Эйлера (задачи \ref{numbfer-fer} и~\ref{numbfer-eul}). Знаний по теории групп и~опыта работы с~определением абстрактной группы не требуется\footnote{В частности, наше доказательство не привлекает явно понятия факторгруппы, в~отличие от более традиционных доказательств, см., например, [Br]%\cite{Br}
. Конструкция факторгруппы "--- одна из простейших конструкций, которая всё-таки уже настолько сложна, что для неё удобнее общее понятие группы вместо группы преобразований. Мы также не используем теорем Силова, хотя наш разбор второго случая похож на их доказательство.}. Небольшое количество необходимых понятий вводятся (и~могут быть освоены читателем) в~процессе доказательства. Конечно, читателю, не знакомому с~основами теории групп, нужно будет самостоятельно доказывать некоторые факты. Хотя эти факты просты, они могут касаться новых для читателя объектов, и~тогда ему нужно будет потрудиться. Такие упражнения "--- важная часть изучения этих понятий. Выполнять их интереснее ради красивых результатов, формулировки которых ясны и~доступны неспециалисту (в~частности, не используют этих понятий), но в~доказательствах которых эти понятия возникают. Это предпочтительнее долгого немотивированного изучения теории. Этот параграф будет особенно интересен читателю, предпочитающему изучить доказательство красивого результата на несколько страниц, самостоятельно разбираясь в~деталях, чем прочитать сотню страниц более лёгкого материала, не мотивированных таким результатом. Подробнее см.~п.\;\ref{s:intpro} <<Изучение путём решения и~обсуждения задач>> и~\S\,28%\ref{s:motiv1}
 <<О необходимости мотивировок>>.

%, чем ради сдачи зачета или изучения немотивированной теории.
%Т.е. в качестве части {\it пути к} пониманию теории, а не в качестве {\it награды за}

Опыт работы с~абстрактными группами как раз появится при изучении данного параграфа, хотя в~нём формально не используется это понятие. Читатель увидит, что помимо всех рассматриваемых объектов есть ещё одно множество (на котором действуют перестановки). Странным образом оно так никогда и~не выходит из тени. В результате естественно возникает общее понятие \emph{группы}. Итак, \emph{этот параграф посвящён мотивировке важного общего понятия группы} (здесь оно не используется, но его и~основы соответствующей теории можно найти в~\cite{Al, KS85}). Хороший опыт в~работе с~основными понятиями теории групп получит и~тот, кто не дойдёт до полного доказательства основного результата.

\vspace{-0.2cm}

\subsection{Как придумать}\label{s:grohow}

 %\footnote{Это пункт может быть пропущен читателем, знакомым с основами теории групп.}

\subsubsection{Постановка задачи (2)}\label{s:grostat}

\begin{pr}\label{open}
Дано семейство $G$ из 11 перестановок некоторого множества, замкнутое относительно композиции и~взятия обратной перестановки (т.\,е. если $f,g\in G$, то $f\circ g\in G$ и~$f^{-1}\in G$). Тогда найдётся перестановка $g\in G$, для которой $G= \{g,g^2,\ldots,g^{11}\}$.
\end{pr}

\vspace{-0.3cm}

\begin{pr}\label{pr-cycle}
Верен ли аналог предыдущего утверждения для аналогичного семейства $G$ из $n$ перестановок при $n=2$; 3; 4; 5; 6; 7; 8; 9; 10; 12; 15; 21; 1001? (Ответ может быть разным для разных~$n$.)
\end{pr}

\vspace{-0.2cm}

\index{Группа|textbf}
\emph{Группой преобразований} называется непустое семейство $G$ преобразований (т.\,е. перестановок) некоторого множества, замкнутое относительно композиции и~взятия обратного преобразования (т.\,е. если $f,g\in G$, то $f\circ g\in G$ и~$f^{-1}\in G$). Мы будем опускать слово <<преобразований>> (поскольку это определение <<равносильно>> обычному определению \emph{группы} ввиду \emph{теоремы Кэли}; ср. с~цитатой из книги В.~И.~Арнольда в~п.\;28.1%\ref{s:motiv1gen}
).

Если в~конечной группе $G$ найдётся перестановка $g$, из всех возможных степеней которой состоит $G$ (т.\,е. $G=\{g,g^2,\ldots,g^n,\ldots\}$), то эта группа называется \emph{циклической\/}.\index{Группа!циклическая|textbf} Примеры циклических и~нециклических групп вы привели при решении задачи~\ref{pr-cycle}.

На этом языке основной вопрос из п.\;\ref{s:grocycwhy} формулируется так: \emph{для каких $n$ любая группа из $n$ элементов циклическая?}

\comment
\begin{pr}\label{pr-anyn} Для любого $n$ имеется циклическая группа из $n$ элементов.
\end{pr}
\endcomment

%\newlength{\wrfigwidth}
\begingroup

\newpage

%\section{Комбинаторная геометрия}

\sectionmark{Комбинаторная геометрия}

\setcounter{page}{694}
\refstepcounter{section}
\refstepcounter{subsection}
\refstepcounter{subsection}
\refstepcounter{subsection}
\refstepcounter{subsection}
\refstepcounter{subsection}
\refstepcounter{subsection}
\refstepcounter{subsection}
\subsection{Собери квадрат (3*). \emph{М.~Б.~Скопенков}, \emph{О.~А.~Малиновская}, \emph{С.~А.~Дориченко}, \emph{Ф.~А.~Шаров}}
\label{ss:soberi}

%\Opensolutionfile{_hintAholder}
%\Opensolutionfile{_hintBholder}

\index{Разрезание} \index{Подобие} \index{Электрическая цепь}
Этот пункт посвящён решению такой задачи (для некоторых частных
случаев).

\paragraph*{Задача.}
Когда из прямоугольников, подобных данному, можно составить квадрат?

\smallskip

В процессе решения мы познакомимся с~красивыми применениями алгебры
в комбинаторной геометрии, а~именно "--- систем линейных уравнений и~многочленов с~целыми коэффициентами. Для решения задач необходимо
первоначальное знакомство с~этими темами. Желательно также
первоначальное знакомство с~задачами на разрезание, см., например,
\cite{Savin-87}.

Наш подход к~решению развивает идеи книги \cite{Ya}.

Другой подход к~решению "--- это физическая интерпретация,
использующая электрические цепи (хотя без неё решать проще).
Познакомиться с~этой физической интерпретацией и~её применением к~решению поставленной задачи можно в~статьях \cite{PS-12, SMD-15}.
Увлекательный рассказ об истории её возникновения можно прочитать в~книге~\cite{G}.

\subsubsection*{Наводящие вопросы}

\begin{flushright}

{\small

\parbox{8cm}{"--* У~меня есть мысль! "--- сказал удав,
открывая глаза. "--- Мысль. И~я~е\" е думаю.\par
"--* Какая мысль? "--- спросила мартышка.\par
"--* Так сразу не скажешь...\par
"--* Ух ты! "--- подпрыгнула мартышка. "--- Ох, какая
хорошая мысль. А~можно я~е\" е тоже немножко подумаю?\par
\hfill \emph{Г.~Остёр}. Бабушка удава}

}

\end{flushright}

\begin{pr}\circpr\label{p1}
Верно ли, что при любых натуральных $m$ и~$n$ из нескольких
прямоугольников $m\times n$ можно сложить квадрат? Выберите верный вариант ответа:

1) верно; \quad 2) неверно.

\end{pr}

\begin{pr}\label{p-9}
Дизайнеру заказали рамы для квадратного окна. На проектах (рис.~\ref{ris25.8.1}\,A,~B) показано, как должны примыкать стёкла друг к~другу и~как они
должны быть ориентированы (короткой или длинной стороной вверх).
Можно ли сделать все стёкла в~каждой раме подобными прямоугольниками?

\begin{figure}[ht]\centering
\begin{minipage}[b]{0.4\textwidth}\centering
\includegraphics[scale=.85]{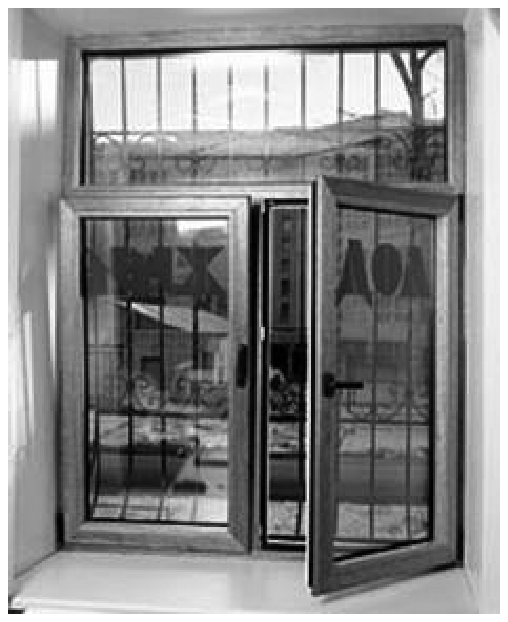}\\[1mm]
{\small A}
\end{minipage}\hfil
\begin{minipage}[b]{0.58\textwidth}\centering
\includegraphics[scale=.85]{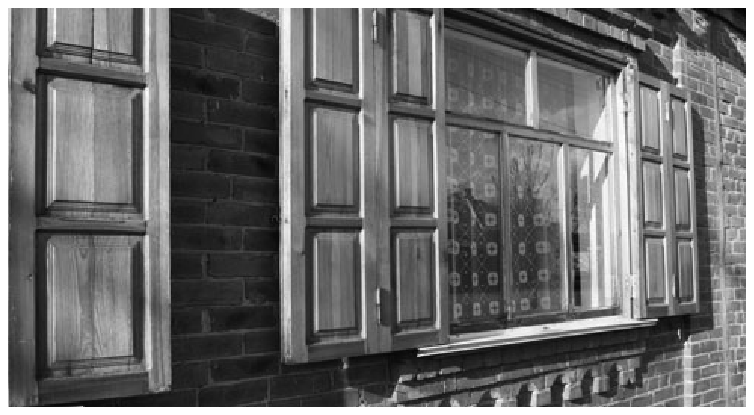}\\[1mm]
{\small B}
\end{minipage}
\caption{Проекты оконных рам; см. задачу \ref{p-9}}
\label{ris25.8.1}
\end{figure}
\end{pr}

\begin{pr}\label{pr3}
Можно ли разрезать квадрат на три подобных, но неравных
прямоугольника?

\end{pr}

\begin{pr}\label{p4}
Можно ли разрезать квадрат на $5$ квадратов?
\end{pr}

\begin{pr} \label{p-2}
Все полки у~шкафа на рис.~\ref{ris25.8.2}\,C, как и~все лоскутки, из которых сшито одеяло на рис.~\ref{ris25.8.2}\,D "--- квадратные. Являются ли квадратными сами шкаф и~одеяло?

\begin{figure}[ht]\centering
\begin{minipage}[b]{0.44\textwidth}\centering
\includegraphics[scale=.9]{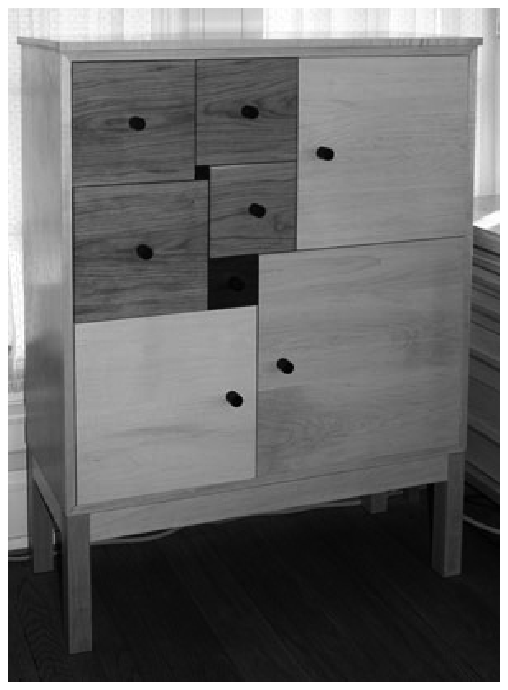}\\[1mm]
{\small C}
\end{minipage}\hfil
\begin{minipage}[b]{0.52\textwidth}\centering
\includegraphics[scale=.9]{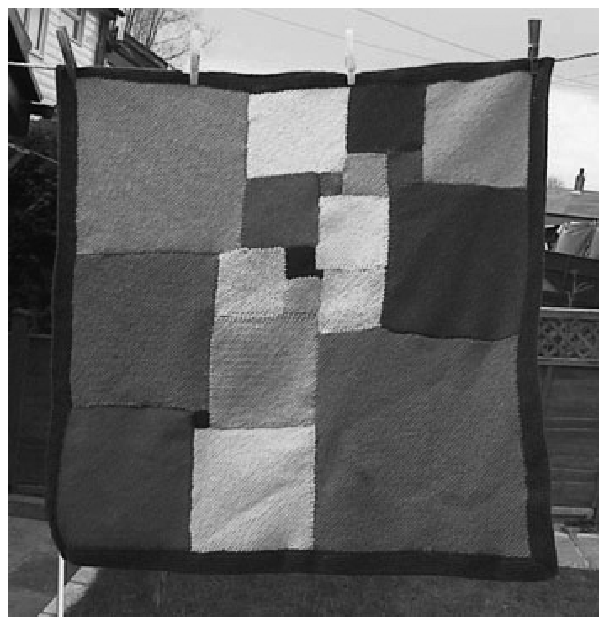}\\[1mm]
{\small D}
\end{minipage}
\caption{Шкаф и одеяло; см. задачу \ref{p-2}}
\label{ris25.8.2}
\end{figure}
\end{pr}

\begin{pr}\label{p-3}
Можно ли замостить всю плоскость попарно различными квадратами, длины сторон которых "--- %рациональные
целые числа?

\end{pr}

\begin{pr}\label{pr9}
Можно ли разрезать квадрат на прямоугольники с~отношением сторон
$2+\sqrt{2}$? То же для $2-\sqrt{2}$, для $3+2\sqrt{2}$ и~для $3-2\sqrt{2}$.

\end{pr}

\begin{pr}\label{px}
Является ли $1+\sqrt{2}$ суммой квадратов чисел вида $a+b\sqrt{2}$, где $a$ и~$b$ рациональны?

\end{pr}

\begin{deff}{} Пусть на прямоугольном листе бумаги нарисовано разбиение на прямоугольники. Разрешается разрезать лист вдоль любого  отрезка на два прямоугольника, потом произвести такие операции по отдельности с~каждой из получившихся частей и~так далее. Если таким образом можно реализовать исходное разбиение, то назовём его \emph{тривиальным}. Например, разбиения на рис.~\ref{ris25.8.1} тривиальные, а~на рис.~\ref{ris25.8.2} нетривиальные.\index{Разрезание!тривиальное}

Следующие $4$ задачи предлагается сначала решить для тривиальных
разбиений, а~уже потом подумать над произвольными разбиениями. В~последующих подпунктах будут даны подсказки к~решению этих трудных
задач.
\end{deff}

\begin{pr}\label{p9}
Какие прямоугольники можно (тривиально) разрезать на прямоугольники
со стороной~$1$?

\end{pr}

\begin{pr}\label{p13}
Какие прямоугольники можно (тривиально) разрезать на квадраты?

\end{pr}

\begin{pr}\label{p-12}
Можно ли квадрат (тривиально) разрезать на прямоугольники с~отношением сторон~$\sqrt2$?
То же для $1+\sqrt2$.

\end{pr}

Все числа, которые можно представить в~виде $x=a+b\sqrt{2}$ с~рациональными $a$ и~$b$, назовём \emph{хорошими}.

\begin{pr}\label{p12}(Основная задача.)
При каких хороших $x$ квадрат можно (тривиально) разрезать на
прямоугольники с~отношением сторон $x$?

\end{pr}

\subsubsection*{Прямоугольник из квадратов.}

%%%!!!\subsubsection*{От разрезаний к~системам линейных уравнений.}

\begin{flushright}

{\small

Ты, дорога, иду по тебе и~гляжу, но мне думается,\\
Мне думается, в~тебе много такого, чего не увидишь глазами.\\[2pt]
\emph{Уолт Уитмен}. Песня большой дороги

}

\end{flushright}

В этом подпункте мы наметим новый вариант элементарного решения задач~\ref{p13} и~\ref{p12}. В~этом подпункте латинские буквы $a$, $b$, $c$, $d$ и~эти же буквы с~индексами обозначают \textit{рациональные} числа.

\begin{pr}\label{p23}
Можно ли прямоугольник $1\times\sqrt{2}$ разрезать на квадраты с~рациональными сторонами? А~со сторонами, которые либо рациональны, либо имеют вид $b\sqrt{2}$? А~со сторонами, которые являются произвольными хорошими числами? Те же вопросы для прямоугольников $1\times(1+\sqrt{2}\,)$ и~$1\times(2+\sqrt{2}\,)$.

\end{pr}

Для доказательства невозможности разрезаний естественно использовать площадь и её \textit{аддитивность}: площадь целого равна сумме площадей частей. Вряд ли получится ответить на вопросы задачи~\ref{p23} для прямоугольника $1\times(2+\sqrt{2}\,)$ без следующего обобщения понятия площади (мы обобщаем понятие площади так, чтобы площадь этого прямоугольника стала отрицательной, а площади квадратов оставались неотрицательными).

\begin{deff}{}
Пусть $x$ "--- действительное число. Назовём \emph{$x$-площадью} (или \textit{площадью Гамеля}) прямоугольника $(a+b\sqrt{2}\,)\times(c+d\sqrt{2}\,)$ число $(a+bx)(c+dx)$. Число $\bar s:=a-b\sqrt{2}$ назовём \textit{сопряжённым} к~числу $s=a+b\sqrt{2}$.
\end{deff}

\begin{pr} Обычная площадь прямоугольника $(a+b\sqrt{2}\,)\times(c+d\sqrt{2}\,)$ и сопряжённое к ней число "--- это одни из его $x$-площадей. Чему равно $x$ в каждом из случаев?
\end{pr}

\begin{pr} Найдите все прямоугольники вида $(a+b\sqrt{2}\,)\times(c+d\sqrt{2}\,)$, $x$-площади которых неотрицательны при всех~$x$.
\end{pr}

\begin{pr} \textbf{Аддитивность $x$-площади.} \label{pr-additive}
Если прямоугольник разрезан на конечное число прямоугольников, стороны которых "--- хорошие числа, то для любого $x\in\mathbb{R}$ $x$-площадь разрезаемого прямоугольника равна сумме $x$-площадей прямоугольников, на которые он разрезан.
\end{pr}

\textit{Указание.} Начните со случая разрезания на $2$ прямоугольника.

\begin{pr}\label{std}
Решите задачи~\ref{p13} и~\ref{p12} для частного случая, когда стороны всех квадратов и всех прямоугольников, участвующих в разрезании, "--- хорошие числа (разрезание не обязательно тривиально).
\end{pr}

В следующих трёх задачах мы считаем, что
прямоугольник $s_0\times t_0$ разрезан на прямоугольники ${s_1\times t_1}$, ${s_2\times t_2}$, \ldots, ${s_N\times t_N}$, причем $s_0$ и $t_0$ несоизмеримы.

\begin{pr} \label{10-basis}
Обозначим
$$
P=\{s_0,t_0,s_1,t_1,\ldots,s_N,t_N\}.
$$
Тогда можно выбрать такие числа ${e_1,e_2,\ldots,e_n\in
P}$, чтобы любое число ${p\in P}$ единственным образом представлялось
в виде
$$
p=as_0+bt_0+a_1e_1+a_2e_2+\ldots+a_ne_n.
$$
\end{pr}

\textit{Указание.} Начните с примера, изображённого на рис.~\ref{fig-basis}.

\begin{figure}[htbp]\centering
\unitlength=1mm
\begin{picture}(53,20)\put(2,-60){
\begin{picture}(45,83)
\thicklines
%\put(0,80){\line(1,0){35}}
%\put(0,60){\line(1,0){35}}
\put(0,60){\line(0,1){15}}
\put(25,60){\line(0,1){15}}
\put(51,60){\line(0,1){15}}
\put(0,75){\line(1,0){51}}
\put(0,60){\line(1,0){51}}
\put(0,70){\line(1,0){25}}
%\put(35,60){\line(0,1){20}}
%\put(14,60){\line(0,1){20}}
%\put(5.5,68){$S_1$}
%\put(22.5,68){$S_2$}
\put(-2,66){1}
\put(20,75.5){$2+\sqrt{2}$}
\put(7.5,69.8){$^{1/3\times\sqrt{3}}$}
\put(7.5,62.3){$^{2/3\times\sqrt{3}}$}
\put(27.8,64.8){$^{1\times\left(2+\sqrt{2}-\sqrt{3}\right)}$}
%\put(0.5,80){$^{c_1+d_1\sqrt{2}}$}
%\put(18.5,80){$^{c_2+d_2\sqrt{2}}$}
%\put(48,80){$^{\gamma+\delta\sqrt{p}}$}

%\put(42,80){\line(1,0){21}}
%\put(42,60){\line(1,0){21}}
%\put(42,60){\line(0,1){20}}
%\put(63,60){\line(0,1){20}}
%\put(42,72){\line(1,0){21}}
%\put(50.5,74.5){$S_3$}
%\put(50.5,64.5){$S_4$}
%\put(63.5,74){$^{\alpha_1+\beta_1\sqrt{p}}$}
%\put(63.5,64){$^{\alpha_2+\beta_2\sqrt{p}}$}

%\put(17.5,55){\ris\label{ri3}}
\end{picture}}\end{picture}
\unitlength=1pt
\caption{К построению базиса}
\label{fig-basis}
\end{figure}

\index{Базис}
Зафиксируем набор чисел $s_0$, $t_0$, $e_1$, $e_2$, \ldots, $e_n$ из задачи~\ref{10-basis}. Он называется \emph{базисом}.

\begin{deff}{}
Пусть $y$ "--- действительное число. Назовём \textit{$y$-площадью} прямоугольника со сторонами
 $$
 as_0+bt_0+a_1e_1+a_2e_2+\ldots+a_ne_n\,\,\,\text{и}\,\,\,cs_0+dt_0+c_1e_1+c_2e_2+\ldots+c_ne_n
 $$
число $(a+by)(c+dy)$.
\end{deff}

Обратите внимание на то, что при $y=x$ и хороших несоизмеримых  $s_0$, $t_0$ это определение не всегда эквивалентно определению $x$-площади выше!

\begin{pr} Вычислите $y$-площадь разрезаемого прямоугольника $s_0\times t_0$. Является ли она неотрицательной при всех~$y$?
\end{pr}

\begin{pr} Докажите, что для любого $y$ $y$-площадь разрезаемого прямоугольника $s_0\times t_0$ равна сумме $y$-площадей прямоугольников, на которые он разрезан.
\end{pr}

\begin{pr}\label{p26}
\textbf{Теорема Дена}. Если прямоугольник разрезан на квадраты (не
обязательно равные), то отношение его сторон рационально.
\end{pr}

\comment
\begin{pr} Если квадрат $1\times 1$ разрезан на прямоугольники, отношение сторон каждого из которых "--- хорошее число, то и сами стороны всех прямоугольников "--- хорошие числа.
\end{pr}
\endcomment

\endgroup

\endgroup

\newpage
%%%!!!

\chapter[О преподавании. \emph{А.~Б.~Скопенков}]{О преподавании}
\label{s:tea}
\setcounter{page}{724}

\vspace*{-3.5\baselineskip}

\begin{center}
{\large\itshape\bfseries А.~Б.~Скопенков}
\end{center}

\vspace*{3\baselineskip}

Редакторы считают важным обсуждение вопросов и~идей, затронутых в~следующих статьях.
При этом мнение авторов статей может не совпадать с~мнением редакторов.

%\newpage
\section[Олимпиады и~математика]{Олимпиады и~математика}\label{oim}

\begin{flushright}

{\small

To him a thinking man's job was not to deny one reality\\
at the expense of the other, but to include and to connect\\[2pt]
\emph{U.~K.~Le Guin}. The Dispossessed\footnotemark

\footnotetext{Для него работой мыслителя было не отрицание одной реальности за счёт другой, а~взаимовключение и~взаимосвязь. \emph{У.~К.~Ле Гуин}, <<Обделённые>> (пер. автора).}

}

\end{flushright}

%\medskip

%\epigraph{at the expense of the other, but to include and to connect.}%
%{To him a thinking man's job was not to deny one reality\\
%at the expense of the other, but to include and to connect.}%
%{\emph{U.~K.~Le Guin}. The Dispossessed\footnotemark}
%
%\footnotetext{Для него работой мыслителя было не отрицание одной реальности за счёт другой, а~взаимовключение и~взаимосвязь. \emph{У.~К.~Ле Гуин}, <<Обделённые>> (пер. автора).}

%\smallskip

Перед школьниками, их учителями и~руководителями кружков встаёт вопрос: готовиться к~олимпиадам или к~<<серьёзной>> математике? Некоторые думают, что для первого надо прорешивать задачи последних олимпиад, для второго надо читать вузовские учебники, и~что ввиду принципиальной разницы первого и~второго бессмысленно пытаться достичь и~того, и~другого. Я~придерживаюсь распространённого мнения о~том, что эти подходы недостаточно эффективны и~приводят к~вредным
<<побочным эффектам>>: школьники либо чрезмерно увлекаются \emph{спортивным} элементом в~решении задач, либо изучают \emph{язык} математики вместо её содержания\footnote{Имеется обширная литература, в~которой в~первую очередь излагается содержание, а~язык появляется по ходу дела. Однако часто такая \emph{популярная} литература недооценивается ввиду её <<недостаточной серьёзности>> по сравнению с~учебниками \emph{для университета}. Подробнее см. \S\,\ref{phil-met} <<Начинать с~языка или содержания?>>.

Кроме того, даже чтение хороших книг без решения задач, как правило, неэффективно.}.

По моему мнению, основу математического образования должно составлять \emph{решение и~обсуждение интересных ученику задач, в~процессе которых он знакомится с~важными математическими идеями и~теориями}.
Это одновременно подготовит школьника и~к~математической науке, и~к~олимпиадам и~не нанесёт вред его развитию в~целом. Это будет более эффективно и~для достижения успеха только в~олимпиадах или только в~науке (если не учитывать большого количества других факторов, кроме разумной организации занятий).

Как и~при естественном развитии самой математики, каждая следующая задача должна быть мотивирована либо практикой, либо уже решёнными задачами (см.~подробнее \S\,\ref{phil-met} <<Начинать с~языка или содержания?>> и~28%\ref{s:motiv1}
 <<О~необходимости мотивировок>>). Ученик, занимающийся <<мотивированной для него>> математикой (обычно более элементарной, но содержательной и~потому сложной) вместо <<немотивированной для него>> математики (обычно менее элементарной, но языковой и~потому тривиальной), имеет преимущество в~дальнейшей учёбе и~научной работе. А.~Н.~Колмогоров говорил, что до тридцати лет математику разумнее всего заниматься решением конкретно поставленных задач. А~значит, умение решать сложные задачи является одним из важнейших для молодого математика.

Олимпиадных задач очень много; большинство из них интересны школьнику, и~среди них много математически содержательных. Такие задачи могут составить основу изучаемого материала. Однако решение олимпиадных задач без изучения математических идей и~теорий недостаточно эффективно для подготовки к~олимпиадам (на долгих "--- год и~более "--- промежутках времени, как и~вообще решение сиюминутных задач без фундаментального развития). А~решение олимпиадных задач \emph{вместе} с~изучением стоящих за ними математических идей и~теорий более эффективно. Это также позволит по-настоящему разобраться в~идеях и~теориях.

Кроме того, большинству людей легче достичь успеха на олимпиадах в~том случае, когда они не считают успех главной целью. Задачу легче решить, если спокойно думать о~самой задаче, а~не о~награде, которая последует за её решением. Поэтому школьник, мотивированный более высокой целью, чем успех на олимпиаде, имеет на этой олимпиаде психологическое преимущество.

См.~также п.\;\ref{s:intpro} <<Изучение путём решения и~обсуждения задач>>.

%\hfill\emph{Круг мог, нацелясь в~стаю самых }
%\hfill\emph{признанных и~возвышенных человеческих мыслей, }
%\hfill\emph{вмиг ссадить ворону в~павлиньих перьях.}
%Но убить смерть он не мог.
%\hfill\emph{В.~Набоков, Под знаком незаконнорождённых.}

\setcounter{page}{726}
\section[Начинать с~языка или содержания?]{Начинать с~языка или содержания?}\label{phil-met}

По моему мнению, именно с~\emph{новых идей\/}, изложенных на уже имеющемся языке, а~не с~\emph{введения нового языка\/}, полезно \emph{начинать} изучение любой теории. Удачно представлять основные идеи на <<олимпиадных>> примерах: на простейших частных случаях,
свободных от технических деталей. Как правило, такие идеи наиболее ярко выражаются доказательствами, подобными приведённым в~\S\,\ref{s:rad}, \ref{0gro-cyc} и~других частях этой книги. Имеется много других ярких примеров, упомянем только фейнмановские лекции по физике (там приводятся физические рассуждения, а~не доказательства).

<<\emph{Мы стараемся свести к~минимуму число понятий, откладывая определения до момента, когда они напрашиваются сами собой, и~избегая задач на понимание и~применение формальных определений \textup{(}типа\textup{ \glqq}\,является ли множество целых чисел группой по сложению}?\grqq)>>[Shen].%\cite{Sh}.

<<\emph{При изложении материала нужно ориентироваться на объекты, которые основательнее всего укореняются в~человеческой памяти. Это "--- отнюдь не системы аксиом и~не логические приёмы в~доказательстве теорем. Изящное решение красивой задачи, формулировка которой ясна и~доступна, имеет больше шансов удержаться в~памяти студента, нежели абстрактная теория. Скажем больше, именно по такому решению, при наличии некоторой математической культуры, студент впоследствии сможет восстановить теоретический материал. Обратное же, как показывает опыт, практически невозможно}>> [Kol,предисловие].%\cite[предисловие]{Ko}.

Такой стиль изложения не только делает материал более доступным, но позволяет сильным ученикам (для которых доступно даже абстрактное изложение) приобрести математический вкус и~стиль. Это означает разумный выбор проблем для исследования и~их мотивировки. Например, математик, понимающий, что теория Галуа мотивируется более важными и~более сложными проблемами, чем построимость правильных многоугольников и~разрешимость алгебраических уравнений в~радикалах, вряд ли станет мотивировать созданную им теорию приложениями, которые можно получить и~без его теории. Это означает также ясное изложение собственных открытий, не скрывающее ошибку или известность полученного результата за чрезмерным формализмом. К~сожалению, такое "--- обычно непреднамеренное "--- сокрытие ошибки часто происходит с~математиками, воспитанными на чрезмерно формальных курсах. Происходило это и~с~автором этих строк; к~счастью, все мои серьёзные ошибки исправлялись \emph{перед} публикациями.

Мода на искусственно формализованное изложение привела к~следующему парадоксу. По данному \emph{известному понятию} высшей математики зачастую непросто восстановить \emph{конкретный красивый результат}, для которого это понятие действительно необходимо (и~при получении которого это понятие возникло).

Доказательство c использованием некоторого нового термина имеют свои преимущества: оно подготавливает читателя к~доказательству тех теорем, которые уже трудно или невозможно доказать без этого термина. Однако такие доказательства, как правило, не должны быть \emph{первыми} доказательствами данного результата (легко себе представить результат \emph{первого} знакомства с~теоремой Пифагора на основе понятий векторного пространства и~скалярного умножения). Кроме того, при приведении <<терминологического>> доказательства полезно оговорить его мотивированность не доказываемым результатом, а~обучением полезному новому методу. Ну и, конечно, важно соблюсти баланс между доказываемым результатом и~уровнем предлагаемой абстракции. 

\printindex

\end{document}